\documentclass[journal]{IEEEtran}
\usepackage{PerlazaTIT2023}
\begin{document}

\title{Empirical Risk Minimization with \\ Relative Entropy Regularization}

\author{Samir~M.~Perlaza, Gaetan~Bisson, I\~{n}aki~Esnaola, Alain~Jean-Marie, and Stefano~Rini.
\thanks{Samir M. Perlaza is with INRIA, Centre Inria d'Université Côte d'Azur, Sophia Antipolis 06902, France; also with the ECE Dept. at Princeton University, Princeton N.J. 08544, USA; and also with the GAATI Laboratory at the Université de la Polynésie Française, Faaa 98702, French Polynesia. }
\thanks{Gaetan Bisson is with the GAATI Laboratory at the Université de la Polynésie Française, Faaa 98702, French Polynesia.}
\thanks{Iñaki Esnaola is with the ACSE Dept. at The University of Sheffield, Sheffield S1 3JD, UK; and also with the ECE Dept. at Princeton University, Princeton N.J. 08544, USA.}
\thanks{Alain Jean-Marie is with INRIA, Centre Inria d'Université Côte d'Azur, Sophia Antipolis 06902, France.}
\thanks{Stefano Rini is with the ECE Dept. at the National Yang Ming Chiao Tung University (NYCU), Hsinchu, Taiwan 30010, ROC.}
\thanks{This work was presented in part at the IEEE International Symposium on Information Theory (ISIT), Espoo, Finland, 2022 in \cite{Perlaza-ISIT-2022}; and appears as an INRIA technical report in \cite{InriaRR9454}. }
%
}



\maketitle 

\begin{abstract} 
The empirical risk minimization (ERM) problem with relative entropy regularization (ERM-RER) is investigated under the assumption that the reference measure is a~$\sigma$-finite measure, and not necessarily a probability measure. 
Under this assumption, which leads to a generalization of the ERM-RER problem allowing a larger degree of flexibility for incorporating prior knowledge, numerous relevant properties are stated. 
Among these properties, the solution to this problem, if it exists, is shown to be a unique probability measure, mutually absolutely continuous with the reference measure. Such a solution exhibits a probably-approximately-correct guarantee for the ERM problem independently of whether the latter possesses a solution.
For a fixed dataset and under a specific condition, the empirical risk is shown to be a sub-Gaussian random variable when the models are sampled from the solution to the ERM-RER problem. 
The generalization capabilities of the solution to the ERM-RER problem (the Gibbs algorithm) are studied via the sensitivity of the expected empirical risk to deviations from such a solution towards alternative probability measures. 
Finally, an interesting  connection between sensitivity, generalization error, and lautum information is established.  
\end{abstract}


\section{Introduction}

\IEEEPARstart{I}{n} statistical machine learning, the problem of empirical risk minimization (ERM) with relative entropy regularization (ERM-RER) has been the workhorse for building probability measures on the set of models, without any additional assumption on the statistical description of the datasets. See for instance \cite{catoni2004statistical, zdeborova2016statistical, alquier2016properties} and \cite{Young-Book-2004}. 
Instead of additional statistical assumptions on the datasets, which are typical in Bayesian methods \cite{robert2007bayesian}, relative entropy  regularization requires a reference probability measure on the set of models, which is external to the ERM problem.  
Often, such a reference measure represents prior knowledge or side information and is chosen for guiding the search of models towards those inducing low empirical risks with high probability over seen and unseen datasets. 
From this perspective, the reference measure can be seen as an additional degree of freedom to improve the generalization capabilities of machine learning algorithms based on ERM-RER, e.g, Gibbs algorithms \cite{aminian2021information, aminian2021exact, jiang2008gibbs, alquier2016properties, Perlaza-ISIT2023b, InriaRR9515, bu2022characterizing, he2023does,InriaRR9539, hellstrom2023generalization,zou2024Generalization} and \cite{aminian2021jensen}. 
This new degree of freedom is one of the main motivations for regularizing the ERM problem using relative entropy, or more generally, any $f$-divergence regularization, as discussed in \cite{masiha2023f, alquier2021non} and~\cite{InriaRR9521}.
Beyond probability measures, as shown in this paper, the reference measure can be any $\sigma$-finite measure with arbitrary support. 
The flexibility introduced by this generalization becomes particularly relevant for the case in which priors are available in the form of probability distributions that can be evaluated up to some normalizing factor, cf.   \cite{robert2004monte}, or cannot be represented by probability distributions, e.g., equal preferences among elements of infinite countable sets.
For some specific choices of $\sigma$-finite reference measures, the ERM-RER boils down to particular cases of special interest: $(i)$ the information-risk minimization problem presented in \cite{zhang2006information}; $(ii)$ the ERM with differential entropy regularization (ERM-DiffER); and $(iii)$ the ERM with discrete entropy regularization (ERM-DisER). See for instance  \cite{mazuelas2022generalized, kapur1989bookMaxEntropy} and references therein. From this perspective, the proposed ERM-RER formulation yields a unified mathematical framework that comprises a large class of problems. 

When the reference measure is a probability measure, the solution to the ERM-RER problem is known to be unique and correspond to a Gibbs probability measure. Such a Gibbs probability measure has been studied using measure theoretic and information theoretic notions in \cite{de2021quantitative, AthreyaGibbs, hwang1980laplace, hasenpflug2022wasserstein, aminian2021exact, raginsky2016information, russo2019much, zhang2006information, asadi2020chaining, xu2017information}; statistical physics in \cite{catoni2004statistical};  PAC (Probably Approximatively Correct)-Bayesian learning theory in \cite{shawe1997pac,mcallester2003pac, haddouche2020pacbayes, guedj2019free};   and proved to be of particular interest in classification problems in \cite{alquier2021non, alquier2016properties, mazuelas2022generalized, jaakkola1999maximum, zhu2009maximum, InriaRR9515} and~\cite{lecun2006tutorial}.
In the general case in which the reference is a $\sigma$-finite measure, a solution to the ERM-RER problem does not always exist. Nonetheless, if it exists,  it is shown to be a unique Gibbs probability  despite the fact that its partition function is defined with respect to a $\sigma$-finite measure.
The condition for the existence is mild and is always satisfied when the reference measure is a probability measure, as highlighted above. 
Interestingly, such a solution is mutually absolutely continuous with the reference measure and most of the properties known for the classical ERM-RER problem are shown to hold in the most general case.  
For instance, under certain conditions, the empirical risk observed when models are sampled from the ERM-RER-optimal probability measure is a sub-Gaussian random variable that exhibits a PAC guarantee for the ERM problem without regularization. 

When the solution to the ERM-RER problem is used to sample models to label unseen patterns, the process is known as the Gibbs algorithm. One of the traditional performance metrics to evaluate the generalization capabilities of machine learning algorithms is the generalization error, for which closed-form expressions in terms of information measures are presented in \cite{InriaRR9539}. 
When the reference measure is a probability measure, a closed-form expression for the generalization error of the Gibbs algorithm is presented in \cite{aminian2021exact}, while upper bounds have been derived  in \cite{e24091178, Jiang2020Fantastic, zhang2006ep, zhang2006information, jiao2017dependence, xu2017information, wang2019information, issa2019strengthened, russo2019much, bu2020tightening, asadi2018chaining, lopez2018generalization,  asadi2020chaining, hafez2020conditioning, haghifam2020sharpened, rodriguez2021tighter,  esposito2021generalization, aminian2021jensen, aminian2022information, aminian2022tighter, shawe1997pac, mcallester2003pac, haddouche2020pacbayes, guedj2019free}, and references therein.
In this work, a new performance metric coined \emph{sensitivity}, which quantifies the variations of the expected empirical risk due to deviations from the solution of the ERM-RER problem is introduced. The sensitivity is defined as the difference between two quantities:~$(a)$ The expectation of the empirical risk with respect to the solution to the ERM-RER problem; and~$(b)$ the expectation of the empirical risk with respect to an alternative measure.  
The absolute value of the sensitivity is shown to be upper bounded by a term that is proportional to the squared-root of the relative entropy of the alternative measure with respect to the ERM-RER-optimal measure. 
Such bound allows providing lower and upper bounds on the expected empirical risk after a deviation from the ERM-RER-optimal measure towards an alternative probability measure. 
More interestingly, the expectation (with respect to the probability distribution of the datasets) of the sensitivity to deviations to a specific measure is shown to be equal to the generalization error of the Gibbs algorithm. Using this result, the closed-form expression for the generalization error  of the Gibbs algorithm  presented in~\cite{aminian2021exact} is shown to hold even in the case in which the reference measure is a $\sigma$-finite measure.
Moreover, the generalization error is shown to be upper bounded by a term that is proportional to the squared-root of the lautum information between the models and the datasets, cf.~\cite{palomar2008lautum}.
This bound is reminiscent of the result in \cite[Theorem~$1$]{xu2017information} in which a similar bound is presented using the mutual information instead of the lautum information. While \cite[Theorem~$1$]{xu2017information} follows immediately from the variational representation of relative entropy, c.f., \cite[Lemma $4.18$ (Transportation Lemma)]{boucheron2013book}, the new result follows from the fact that the empirical risk  when models are sampled from the ERM-RER-optimal probability measure is a sub-Gaussian random variable. 
 Interestingly, the new upper-bound does not require any of the conditions in~\cite[Theorem~$1$]{xu2017information}.

The remainder of this work is organized as follows. Section~\ref{SecERM} introduces two optimization problems: the ERM and the ERM-RER. The asymmetry of the relative entropy is analyzed in the context of the ERM-RER and two variants,  coined \mbox{Type-I} and \mbox{Type-II}, are distinguished. The former considers the case in which the regularization is the relative entropy of the optimization measure with respect to the reference measure. The latter considers a regularization by the relative entropy of the reference measure with respect to the optimization measure.
Section~\ref{SecProperties} presents the solution to the ERM-RER problem  in the general case and introduces its main properties. 
Section~\ref{SecRefMes} introduces two new classes of reference measures and the solution of the ERM-RER problem is shown to exhibit different properties for each class. This section ends by studying the ERM-RER problem in the special case in which the reference measure is a Gibbs probability measure.  This special case exhibits a solution that is identical to the solution to an ERM-RER problem whose reference measure is the same used to build the above mentioned Gibbs measure.
Section~\ref{SecLogPartition} studies the properties of the log-partition function of the ERM-RER-optimal probability measure. The first, second, and third cumulants of the empirical risk when the models are sampled from the ERM-RER-optimal measure and the reference measure are respectively characterized.  
Section~\ref{SecExpectation} and Section~\ref{SecVariance} study the properties of the expectation and variance of the empirical risk when the models are sampled from the ERM-RER-optimal probability measure. This mean and variance are compared with the mean and variance of the empirical risk when models are sampled from the reference measure. 
Section~\ref{SecSubGaussianity} introduces several  explicit expressions for the cumulant generating function of the empirical risk when the models are sampled from the ERM-RER-optimal measure. Using these equivalent expressions, under a specific condition, it is shown that the empirical risk is a sub-Gaussian random variable when models are sampled from the ERM-RER-optimal measure. 
Section~\ref{SecMonotonicConcentration} describes the monotonic concentration of the ERM-RER-optimal probability measure when the regularization factor tends to zero. 
Section~\ref{SecDeltaEpsilonOptimality}  show that the empirical risk when the models are sampled from the ERM-RER-optimal probability measure exhibits a PAC-type guarantee with respect to the ERM problem without regularization. 
Finally, Section~\ref{SecSubSensitivity} studies the sensitivity of the expected empirical risk  with respect to deviations from the ERM-RER-optimal measure to alternative measures and shows connections with the generalization error and the lautum information.
Section~\ref{SecDiscussion} ends this work with conclusions and a discussion on the results.
  
\section{Empirical Risk Minimization (ERM)}\label{SecERM}

Let~$\set{M}$,~$\set{X}$ and~$\set{Y}$, with~$\set{M} \subseteq \reals^{d}$ and~$d \in \ints$, be sets of \emph{models}, \emph{patterns}, and \emph{labels}, respectively.  
A pair~$(x,y) \in \mathcal{X} \times \mathcal{Y}$ is referred to as a \emph{labeled pattern} or as a \emph{data point}.
Given~$n$ data points, with~$n \in \ints$,  denoted by~$\left(x_1, y_1 \right)$, $\left( x_2, y_2\right)$, $\ldots$, $\left( x_n, y_n \right)$, the corresponding dataset is represented by the tuple
\begin{equation}\label{EqTheDataSet}
\vect{z} = \big(\left(x_1, y_1 \right), \left(x_2, y_2 \right), \ldots, \left(x_n, y_n \right)\big)  \in \left( \set{X} \times \set{Y} \right)^n.
\end{equation}  

Let the function~$f: \set{M} \times \mathcal{X} \rightarrow \mathcal{Y}$ be such that the label assigned to the pattern $x$ according to the model $\vect{\theta} \in \set{M}$ is $f(\vect{\theta}, x)$.
%
Let also the function 
\begin{equation}\label{EqEll}
\ell: \set{Y} \times \set{Y} \rightarrow [0, +\infty]
\end{equation} 
be such that given a data point~$(x, y) \in \set{X} \times \set{Y}$, the  risk induced by a model~$\vect{\theta} \in \set{M}$ is ~$\ell\left( f(\vect{\theta}, x), y \right)$.  
In the following, the risk function~$\ell$ is assumed to be nonnegative and  for all~$y \in \set{Y}$, ~$\ell\left( y , y\right) = 0$.

The \emph{empirical risk} induced by the model~$\vect{\theta}$, with respect to the dataset $\vect{z}$ in~\eqref{EqTheDataSet} is determined by the  function~$\mathsf{L}_{\vect{z}}: \set{M} \rightarrow [0, +\infty ]$, which satisfies  
\begin{IEEEeqnarray}{rCl}
\label{EqLxy}
\mathsf{L}_{\vect{z}} \left(\vect{\theta} \right)  & = & 
\frac{1}{n}\sum_{i=1}^{n}  \ell\left( f(\vect{\theta}, x_i), y_i\right).
\end{IEEEeqnarray}
Using this notation, the ERM consists of the following optimization problem:
\begin{equation}\label{EqOriginalOP}
\min_{\vect{\theta} \in \set{M}} \mathsf{L}_{\vect{z}} \left(\vect{\theta} \right).
\end{equation}
Let the set of solutions to the ERM problem in~\eqref{EqOriginalOP} be denoted by
\begin{equation}\label{EqHatTheta}
\set{T}\left( \vect{z} \right) \triangleq \arg\min_{\vect{\theta} \in \set{M}}    \mathsf{L}_{\vect{z}} \left(\vect{\theta} \right).
\end{equation}
Note that if the set $\set{M}$ is finite, the ERM problem in~\eqref{EqOriginalOP} always possesses a solution, and thus, $\abs{\set{T}\left( \vect{z} \right)} > 0$. Nonetheless, in general, the ERM problem might not necessarily possess a solution, i.e., $\abs{\set{T}\left( \vect{z} \right) } = 0$.

\subsection{Notation and Main Assumptions} 
In the following, given a measurable space~$\left( \Omega , \mathscr{F} \right)$, the notation~$\triangle\left( \Omega , \mathscr{F} \right)$ is used to represent the set of~$\sigma$-finite measures  that can be defined over $\left( \Omega , \mathscr{F} \right)$. Given a measure~$Q \in \triangle\left( \Omega , \mathscr{F} \right)$, the subset~$\triangle_{Q}\left( \Omega , \mathscr{F} \right)$ of~$\triangle\left( \Omega , \mathscr{F} \right)$ contains all $\sigma$-finite measures that are absolutely continuous with respect to the measure~$Q$. Alternatively, the subset~$\bigtriangledown_{Q}\left( \Omega , \mathscr{F} \right)$ of~$\triangle\left( \Omega , \mathscr{F} \right)$ contains all probability measures $P$ such that $Q$ is absolutely continuous with respect to~$P$.
Given a set~$\set{A} \subset \reals^d$, the Borel~$\sigma$-field over~$\set{A}$ is denoted by~$\BorSigma{\set{A}}$.

The main assumption adopted in this work is that the function $\mathsf{L}_{\vect{z}}$ in~\eqref{EqLxy} is measurable with respect to the Borel measurable spaces $\left( \set{M}, \BorSigma{\set{M}} \right)$ and $\left([0, +\infty], \BorSigma{[0, +\infty]} \right)$.

\subsection{Relative Entropy Extended to $\sigma$-Finite Measures} 

In this work, the \emph{relative entropy}, which is usually defined for probability measures, is extended to $\sigma$-finite measures.
\begin{definition}[Generalized Relative Entropy]\label{DefRelEntropy}
Given two~$\sigma$-finite measures~$P$ and~$Q$ on the same measurable space, such that~$P$ is absolutely continuous with respect to~$Q$,  the relative entropy of~$P$ with respect to~$Q$ is
\begin{equation}
\label{EqGKL}
\KL{P}{Q} = \int \frac{\mathrm{d}P}{\mathrm{d}Q}(x)  \log\left( \frac{\mathrm{d}P}{\mathrm{d}Q}(x)\right)  \mathrm{d}Q(x),
\end{equation}
where the function~$\frac{\mathrm{d}P}{\mathrm{d}Q}$ is the Radon-Nikodym derivative of~$P$ with respect to~$Q$.
\end{definition}
The relative entropy exhibits a property often referred to as the \emph{information inequality} \cite[Theorem~$2.6.3$]{cover2006book} in the case of probability measures on $\left( \Omega , \mathscr{F} \right)$, with $\Omega$ a countable set. The following theorem explores this property in a more general scenario.
\begin{theorem}\label{LemmaREA}
If $P$ and $Q$ are both probability measures on a general measurable space $\left( \Omega , \mathscr{F} \right)$, with $P$ absolutely continuous with respect to $Q$, then, 
\begin{IEEEeqnarray}{rcl}
\KL{P}{Q} & \geqslant & 0,
\end{IEEEeqnarray}
with equality if and only if $P$ and $Q$ are identical.
\end{theorem}
\begin{IEEEproof}
Consider the function $f:[0,\infty) \to \reals$ such that for all $x \in (0, +\infty)$, $f(x) = x \log(x)$ and $f(0) = 0$.  Note that $f$ is strictly convex.  
If $P$ and $Q$ are both probability measures on the measurable space $\left( \Omega , \mathscr{F} \right)$, the following holds:
\begin{IEEEeqnarray}{rcl}
\label{EqInTheAirplane101}
\KL{P}{Q} & = & \int  \frac{\mathrm{d}P}{\mathrm{d}Q}(x) \log\left(  \frac{\mathrm{d}P}{\mathrm{d}Q}(x) \right) \mathrm{d}Q(x)\\
\label{EqInTheAirplane102}
& = & \int f\left(  \frac{\mathrm{d}P}{\mathrm{d}Q}(x) \right) \mathrm{d}Q(x)\\
\label{EqInTheAirplane102DC}
& \geqslant &  f\left( \int \frac{\mathrm{d}P}{\mathrm{d}Q}(x)  \mathrm{d}Q(x) \right)\\
\label{EqInTheAirplane103}
& = &  f\left( 1 \right)\\
\label{EqInTheAirplane104}
& = &  0,
\end{IEEEeqnarray}
where the inequality~\eqref{EqInTheAirplane103} follows from Jensen's inequality \cite[Section $6.3.5$]{ash2000probability}. Equality in~\eqref{EqInTheAirplane103} holds if and only if for all $x \in \supp Q$, $\frac{\mathrm{d}P}{\mathrm{d}Q}(x) = 1$, which implies that both $P$ and $Q$ are identical.
This completes proof.
\end{IEEEproof}
If $Q$ is not a probability measure, then it might be observed that $\KL{P}{Q} < 0$. Consider for instance the case in which $P$ is a zero-mean Gaussian probability measure with variance $\sigma^2$ and $Q$ is the Lebesgue measure on $\Bormeaspace{\reals}$. 
Hence, the Radon-Nikodym derivative $\frac{\mathrm{d}P}{\mathrm{d}Q}$ is the Gaussian probability density function  such that for all $x \in~\reals$,
\begin{IEEEeqnarray}{rcl}
\frac{\mathrm{d}P}{\mathrm{d}Q}(x) & = & \frac{1}{\sqrt{2 \pi \sigma^2}}\exp\left(- \frac{x^2}{2 \sigma^2}\right).
\end{IEEEeqnarray}
Under this assumption, the relative entropy of $P$ with respect to $Q$ is the negative of the differential entropy of $P$. That is,
\begin{IEEEeqnarray}{rcl}
\label{EqKeyExampleZero}
\KL{P}{Q} & = & -\frac{1}{2} \log\left(2 \pi  \epsilon \sigma^2 \right),
\end{IEEEeqnarray}
with $\epsilon$ being Néper's constant. See for instance  \cite[Example~$8.1.2$]{cover2006book}. 
Hence, $\KL{P}{Q}$ is negative for all $\sigma^2 \in \left(\frac{1}{2 \pi  \epsilon} , + \infty \right)$ and nonnegative for all $\sigma^2 \in \left( 0,  \frac{1}{2 \pi  \epsilon} \right]$. 
Finally, note also that 
\begin{IEEEeqnarray}{rcl}
\label{EqNoisyMiami}
\lim_{\sigma^2 \to 0}\KL{P}{Q} & = & +\infty, \mbox{ and} \\
\label{EqRafaelIsBorn}
\lim_{\sigma^2 \to +\infty}\KL{P}{Q} & = & -\infty.
\end{IEEEeqnarray}

A central observation from~\eqref{EqKeyExampleZero} is that the equality $\KL{P}{Q}  = 0$ does not necessarily imply that $P$ and $Q$ are identical measures. For instance,  when $\sigma^2 = \frac{1}{2 \pi  \epsilon}$ in~\eqref{EqNoisyMiami}, it holds that $\KL{P}{Q} = 0$, while $P$ is a Gaussian probability measure and $Q$ is the Lebesgue measure. 

The following  property, known for the case of probability measures as the \emph{joint-convexity of the relative entropy}, is extended by the following theorem. 

\begin{theorem}\label{TheoremCasaBlanca}
Let $P_1$ and $P_2$ be two probability measures and $Q_1$ and $Q_2$ be  two $\sigma$-finite measures, all on the same measurable space. For all $i \in \lbrace 1,2 \rbrace$, let $P_{i}$ be   absolutely continuous with  respect to $Q_{i}$. Then, for all $\lambda \in [0,1]$,
\begin{IEEEeqnarray}{rcl}
\nonumber
&  & \KL{\lambda P_1 + (1-\lambda) P_2}{\lambda Q_1 + (1-\lambda) Q_2} \\
\label{EqDontWorry}
& \leqslant &  \lambda \KL{ P_1 }{Q_1} + (1 - \lambda) \KL{ P_2}{Q_2}.
\end{IEEEeqnarray}
Equality in \eqref{EqDontWorry} holds if and only if $P_{1} = P_{2}$ and $Q_{1} = Q_{2}$.
\end{theorem}
\begin{IEEEproof}
The proof is presented in Appendix~\ref{AppProofOfTheoremCasaBlanca}.
\end{IEEEproof}

 \subsection{ERM with Relative Entropy Regularization}

Given a dataset, the \emph{expected empirical risk} induced by a measure $P \in \Delta\Bormeaspace{\set{M}}$ is defined as follows.
\begin{definition}[Expected Empirical Risk]\label{DefEmpiricalRisk}
Let $P$ be a probability measure in $\Delta\Bormeaspace{\set{M}}$. The expected empirical risk with respect to the dataset~$\vect{z}$ in~\eqref{EqTheDataSet} induced by the measure $P$ is 
\begin{equation}
\label{EqRxy}
\mathsf{R}_{\vect{z}}\left( P  \right) = \int \mathsf{L}_{ \vect{z} } \left(\vect{\theta} \right)  \mathrm{d} P(\vect{\theta}),
\end{equation}
where the function~$\mathsf{L}_{\vect{z}}$ is defined in~\eqref{EqLxy}. 
\end{definition}

The ERM-RER problem is parametrized by a~$\sigma$-finite measure in~$\triangle\Bormeaspace{\set{M}}$ and a positive real, which are referred to as the \emph{reference measure} and the \emph{regularization factor}, respectively.
Let~$Q \in \triangle\Bormeaspace{\set{M}}$ be a~$\sigma$-finite measure  and let~$\lambda$ be a positive real. The ERM-RER problem, with parameters~$Q$ and~$\lambda$, consists of the following optimization problem:
\begin{subequations}\label{EqERMRER}
\begin{IEEEeqnarray}{CCl}
\label{EqERMRERa1}
    \min_{P \in \triangle_{Q}\Bormeaspace{\set{M}}} & &  \mathsf{R}_{\vect{z}} \left( P \right)  + \lambda \KL{P}{Q},\\
\label{EqERMRERb1}
    \mathrm{s.~t.} & & \int \d P \left( \vect{ \theta} \right) = 1,
\end{IEEEeqnarray}
\end{subequations}
where the dataset~$\vect{z}$ is in~\eqref{EqTheDataSet}, and the functional~$\mathsf{R}_{\vect{z}}$ is defined in~\eqref{EqRxy}. 

\subsection{Type-I and Type-II Relative Entropy Regularization}

The optimization problem in~\eqref{EqERMRER} is coined \mbox{Type-I} ERM-RER in~\cite{Perlaza-ISIT2023a} in the aim of distinguishing it from the optimization problem
\begin{subequations}\label{EqERMRERbis}
\begin{IEEEeqnarray}{CCl}
\label{EqERMRERbisa}
    \min_{P \in \bigtriangledown_{Q}\Bormeaspace{\set{M}}} & &  \mathsf{R}_{\vect{z}} \left( P \right)  + \lambda \KL{Q}{P},\\
\label{EqERMRERbisb}
    \mathrm{s.~t.} & & \int \d P \left( \vect{ \theta} \right) = 1,
\end{IEEEeqnarray}
\end{subequations}
which is coined \mbox{Type-II} ERM-RER.

The \mbox{Type-II} ERM-RER problem in~\eqref{EqERMRERbis}, when $Q$ is a probability measure, exhibits a solution that is identical to the solution to the following \mbox{Type-I} ERM-RER problem~\cite[Theorem~$1$]{Perlaza-ISIT2023a}:
\begin{subequations}\label{EqERMRERTypeI}
\begin{IEEEeqnarray}{CCl}
\label{EqERMRERTypeIa}
    \min_{P \in \triangle_{Q}\Bormeaspace{\set{M}}} & &  \int \log\left( \beta + \mathsf{L}_{\vect{z}} \left(\vect{\nu} \right) \right) \mathrm{d}P \left( \vect{\nu}  \right)  + \KL{P}{Q},
\middlesqueezeequ
\IEEEeqnarraynumspace\\
\label{EqERMRERTypeIb}
    \mathrm{s.~t.} & & \int \d P \left( \vect{ \theta} \right) = 1,
\end{IEEEeqnarray}
where $\beta$ is a constant chosen to satisfy
\begin{IEEEeqnarray}{rcl}
\int \frac{\lambda}{\beta + \mathsf{L}_{\vect{z}} \left(\vect{\nu} \right)} \mathrm{d}Q(\vect{\nu}) &=& 1.
\end{IEEEeqnarray}
\end{subequations}
Essentially, by appropriately transforming the objective function, an equivalence can be established between \mbox{Type-I} and \mbox{Type-II} ERM-RER problems. 
Hence, without loss of generality, the remainder of this work focuses exclusively on \mbox{Type-I} ERM-RER, which is simply referred to as ERM-RER.

\section{The Solution to the ERM-RER Problem}\label{SecProperties}

The solution to the ERM-RER problem in~\eqref{EqERMRER} is presented in terms of two objects. First, the function~$K_{Q, \vect{z}}: \reals \rightarrow \reals \cup \lbrace +\infty\rbrace$  such that for all~$t \in \reals$,
\begin{IEEEeqnarray}{rcl}
\label{EqK}
K_{Q, \vect{z}}\left(t \right) & = &  \log\left( \int \exp\left( t \; \mathsf{L}_{\vect{z}}\left(\vect{\theta}\right)  \right) \mathrm{d}Q(\vect{\theta}) \right),
\end{IEEEeqnarray} 
with~$\mathsf{L}_{\vect{z}}$ in~\eqref{EqLxy}. Second, the set~$\set{K}_{Q, \vect{z}} \subset (0, +\infty)$, which is defined by
\begin{IEEEeqnarray}{rcl}
\label{EqSetKxy}
\set{K}_{Q, \vect{z}} & \triangleq &\left\lbrace s \in (0, +\infty): \; K_{Q, \vect{z}}\left(-\frac{1}{s} \right)  < +\infty \right\rbrace.
\end{IEEEeqnarray}
The notation for the function~$K_{Q, \vect{z}}$ and the set~$\set{K}_{Q, \vect{z}}$ are chosen such that their parametrization by (or dependence on) the  dataset~$\vect{z}$ in~\eqref{EqTheDataSet} and the~$\sigma$-finite measure~$Q$  in~\eqref{EqERMRER} are highlighted.

The following lemma describes the set~$\set{K}_{Q, \vect{z}}$.
\begin{lemma}\label{LemmaSubSetK}
The set~$\set{K}_{Q, \vect{z}}$ in~\eqref{EqSetKxy} is a convex subset of $\reals$. If the measure~$Q$ in~\eqref{EqERMRER} is a probability measure, then, the set~$\set{K}_{Q, \vect{z}}$ in~\eqref{EqSetKxy}  satisfies
\begin{equation}
\set{K}_{Q, \vect{z}} = (0, +\infty).
\end{equation}

\end{lemma}
\begin{IEEEproof}
The proof is presented in Appendix~\ref{AppProofLemmaSubSetK}.
\end{IEEEproof}
Using this notation, the solution to the ERM-RER problem in~\eqref{EqERMRER} is presented by the following theorem.
\begin{theorem}\label{TheoremOptimalModel}
If $\lambda \in \set{K}_{Q, \vect{z}}$, with $\set{K}_{Q, \vect{z}}$ in~\eqref{EqSetKxy}, the solution to the optimization problem in~\eqref{EqERMRER} is a unique probability measure, denoted by $P^{\left(Q, \lambda\right)}_{\vect{\Theta} | \vect{Z} = \vect{z}}$, which satisfies for all~$\vect{\theta} \in \supp Q$,
\begin{IEEEeqnarray}{rcl}\label{EqGenpdf}
\frac{\mathrm{d}P^{\left(Q, \lambda\right)}_{\vect{\Theta} | \vect{Z} = \vect{z}}}{\mathrm{d}Q} \left( \vect{\theta} \right) 
  & =& \exp\left( - K_{Q, \vect{z}}\left(- \frac{1}{\lambda} \right) - \frac{1}{\lambda} \mathsf{L}_{\vect{z}}\left( \vect{\theta}\right)\right),
\end{IEEEeqnarray}
where the function~$\mathsf{L}_{\vect{z}}$ is defined in~\eqref{EqLxy} and the function~$K_{Q, \vect{z}}$ is defined in~\eqref{EqK}.
\end{theorem}
\begin{IEEEproof}
The proof is presented in Appendix~\ref{AppProofTheoremOptimalModel}.
\end{IEEEproof}

Contrary to the ERM problem in~\eqref{EqOriginalOP}, which does not necessarily possess a solution, the ERM-RER problem in~\eqref{EqERMRER} always possess a solution when $Q$ is a probability measure. This is essentially because the set $\set{K}_{Q, \vect{z}}$ is the set of all positive reals (Lemma~\ref{LemmaSubSetK}), and thus, the condition $\lambda \in \set{K}_{Q, \vect{z}}$ is always verified. On the contrary, when $Q$ is a $\sigma$-finite measure, the solution to the ERM-RER problem in~\eqref{EqERMRER} depends on whether $\lambda \in \set{K}_{Q, \vect{z}}$. If the solution exists, it is~$P^{\left(Q, \lambda\right)}_{\vect{\Theta} | \vect{Z} = \vect{z}}$ in~\eqref{EqGenpdf}, which is unique and corresponds to a Gibbs probability measure~\cite{georgii2011gibbs}.  
The function~$K_{Q, \vect{z}}$ is often referred to as the  \emph{log-partition function}, see for instance, \cite[Section~$7.3.1$]{dembo2009large}.

The following lemma shows that the Radon-Nikodym derivative in~\eqref{EqGenpdf} is both nonnegative and finite. 
\begin{lemma}\label{CorPositive}
The Radon-Nikodym derivative~$\frac{\mathrm{d}P^{\left(Q, \lambda\right)}_{\vect{\Theta} | \vect{Z} = \vect{z}}}{\mathrm{d}Q}$ in~\eqref{EqGenpdf} satisfies for all~$\vect{\theta} \in \supp Q$  that
    \begin{equation}\label{EqpdfAT}
\frac{\mathrm{d}P^{\left(Q, \lambda\right)}_{\vect{\Theta} | \vect{Z} = \vect{z}}}{\mathrm{d}Q}\left( \vect{\theta}\right) < +\infty.
    \end{equation}  
Moreover, it is strictly positive almost surely with respect to~$Q$.
\end{lemma}
\begin{IEEEproof}
The proof is presented in Appendix~\ref{AppProofCorPositive}.
\end{IEEEproof}

Theorem~\ref{TheoremOptimalModel} has shown that the probability measure~$P^{\left(Q, \lambda\right)}_{\vect{\Theta} | \vect{Z} = \vect{z}}$ is absolutely continuous with respect to the measure~$Q$.
The following lemma shows that the converse is also true.
 \begin{lemma}\label{LemmaMutualAC}
The $\sigma$-finite measure~$Q$ and the probability measure~$P^{\left(Q, \lambda\right)}_{\vect{\Theta} | \vect{Z} = \vect{z}}$ in~\eqref{EqGenpdf} are mutually absolutely continuous.
\end{lemma}
\begin{IEEEproof}
The proof  is presented in Appendix~\ref{AppProofLemmaMutualAC}.
\end{IEEEproof}
The relevance of Lemma~\ref{LemmaMutualAC} is that it shows that if~$\lambda \in \set{K}_{Q, \vect{z}}$, the corresponding collections of negligible sets with respect to the measures~$P^{\left(Q, \lambda\right)}_{\vect{\Theta} | \vect{Z} = \vect{z}}$ and~$Q$ are identical.  
The following lemma shows that the negligible sets with respect to the measure~$P^{\left(Q, \lambda \right)}_{\vect{\Theta}| \vect{Z} = \vect{z}}$ in~\eqref{EqGenpdf} are invariant with respect to~$\lambda$.
\begin{lemma}\label{LemmaMutualAlphaBeta}
For all~$\left(\alpha, \beta\right) \in \set{K}_{Q, \vect{z}}\times\set{K}_{Q, \vect{z}}$, with~$\set{K}_{Q, \vect{z}}$ in~\eqref{EqSetKxy}, assume that the probability measures~$P^{\left(Q, \alpha\right)}_{\vect{\Theta} | \vect{Z} = \vect{z}}~$ and~$P^{(Q, \beta)}_{\vect{\Theta}| \vect{Z} = \vect{z}}$ satisfy~\eqref{EqGenpdf} with~$\lambda = \alpha$ and~$\lambda = \beta$, respectively. Then,~$P^{\left(Q, \alpha\right)}_{\vect{\Theta} | \vect{Z} = \vect{z}}~$ and~$P^{(Q, \beta)}_{\vect{\Theta}| \vect{Z} = \vect{z}}$  are mutually   absolutely continuous.
\end{lemma}
\begin{IEEEproof}
The proof  is presented in Appendix~\ref{AppProofLemmaMutualAlphaBeta}.
\end{IEEEproof} 

Particular assumptions on the set~$\set{M}$ and the reference measure~$Q$ lead to well-known instances of the ERM-RER problem in~\eqref{EqERMRER}, as discussed hereunder.

\subsection{Examples}
Three examples are of particular interest: 
$(a)$ The set~$\set{M} \subset  \reals^{d}$ is countable and the measure~$Q$ is the counting measure  in $\Bormeaspace{\set{M}}$, which leads to the ERM-DisER problem;
$(b)$ The set~$\set{M}$ is an uncountable subset of~$\reals^{d}$, and~$Q$ is the Lebesgue measure on~$\Bormeaspace{\set{M}}$, which leads to the ERM-DiffER problem; and
$(c)$ The set~$\set{M}$ and the measure~$Q$ form a Borel probability measure space~$\left(\set{M}, \BorSigma{\set{M}}, Q \right)$, which leads to the information-risk minimization problem.

\subsubsection{ERM with Discrete Entropy Regularization}

When the set~$\set{M} \subset \reals^d$  is countable and the~$\sigma$-finite measure~$Q$ in~\eqref{EqERMRER} is the counting measure in $\Bormeaspace{\set{M}}$, given a probability measure $P \in \triangle\Bormeaspace{\set{M}}$,  the Radon-Nikodym derivative~$\frac{\mathrm{d}P}{\mathrm{d}Q}$ is a probability mass function, denoted by $p$. Thus, the relative entropy $\KL{P}{Q}$ is equivalent to the negative of the discrete entropy induced by $p$~\cite[Chapter~$2$]{cover2006book}, denoted by $H(p)$.
In this case, the ERM-RER in~\eqref{EqERMRER} can be re-written as the following ERM-DisER problem:
\begin{IEEEeqnarray}{CCCl}
\label{EqKillerEntropyLearningAPingPong}
&  \min_{p} & &  \sum_{\vect{\theta} \in \set{M}} \mathsf{L}_{\vect{z}} \left(\vect{\theta} \right) p\left(\vect{\theta} \right)    - \lambda H\left( p \right),
\end{IEEEeqnarray}
where the optimization domain in~\eqref{EqKillerEntropyLearningAPingPong} is the set of probability mass functions that can be defined over the measure space~$\triangle\Bormeaspace{\set{M}}$.  
In this special case, the probability measure $P^{\left(Q, \lambda\right)}_{\vect{\Theta} | \vect{Z} = \vect{z}}$ in~\eqref{EqGenpdf} whose probability mass function is the solution to the ERM-DisER problem in~\eqref{EqKillerEntropyLearningAPingPong} satisfies
\begin{IEEEeqnarray}{rCl}
\label{EqpmfGibbsSS}
\frac{\mathrm{d}P^{\left(Q, \lambda\right)}_{\vect{\Theta} | \vect{Z} = \vect{z}}}{\mathrm{d}Q} \left( \vect{\theta} \right) & = & \frac{\exp\left( -\frac{\mathsf{L}_{\vect{z}}\left( \vect{\theta}\right)}{\lambda}\right)}{\displaystyle\sum_{\vect{\nu} \in \set{M}} \exp\left( -\frac{\mathsf{L}_{\vect{z}}\left( \vect{\nu} \right)}{\lambda} \right)},
\end{IEEEeqnarray}
which describes the discrete Gibbs probability measure on~$\triangle\Bormeaspace{\set{M}}$, with temperature parameter~$\lambda$, and energy function~$\mathsf{L}_{\vect{z}}$ in~\eqref{EqLxy}.

\subsubsection{ERM with Differential Entropy Regularization}
When~$\set{M} \subseteq \reals^d$ is uncountable and the~$\sigma$-finite measure~$Q$ in~\eqref{EqERMRER} is the Lebesgue measure in $\Bormeaspace{\set{M}}$, for all probability measures $P\in\Delta_Q\Bormeaspace{\set{M}}$, the Radon-Nikodym derivative~$\frac{\mathrm{d}P}{\mathrm{d}Q}$ is a probability density function, denoted by $g$. 
Thus, the relative entropy $\KL{P}{Q}$ is equivalent to the negative of the differential entropy induced by $g$~\cite[Chapter~$8$]{cover2006book}, denoted by $h(g)$.
In this special case, the ERM-RER in~\eqref{EqERMRER} can be re-written as the following ERM-DiffER problem:
\begin{IEEEeqnarray}{cCCl}
\label{EqKillerEntropyLearningAgadir}
& \min_{g} & & \int_{\set{M}} \mathsf{L}_{\vect{z}} \left(\vect{\theta} \right) g\left(\vect{\theta} \right) \mathrm{d}\vect{\theta}  - \lambda h\left( g \right),
\end{IEEEeqnarray}
where the optimization domain in~\eqref{EqKillerEntropyLearningAgadir} is the set of probability density functions  that can be defined over the measure space~$\Bormeaspace{\set{M}}$.  
The probability measure $P^{\left(Q, \lambda\right)}_{\vect{\Theta} | \vect{Z} = \vect{z}}$ in \eqref{EqGenpdf} whose probability density function  is the solution to the ERM-RER problem in~\eqref{EqKillerEntropyLearningAgadir} satisfies
\begin{IEEEeqnarray}{rCl}
\label{EqpdfGibbsSS}
\frac{\mathrm{d}P^{\left(Q, \lambda\right)}_{\vect{\Theta} | \vect{Z} = \vect{z}}}{\mathrm{d}Q} \left( \vect{\theta} \right) & = & \frac{\exp\left( -\frac{\mathsf{L}_{\vect{z}}\left( \vect{\theta}\right)}{\lambda}\right)}{\displaystyle\int_{\set{M}}\exp\left( -\frac{\mathsf{L}_{\vect{z}}\left( \vect{\nu} \right)}{\lambda} \right)  \mathrm{d} \vect{\nu} },
\end{IEEEeqnarray}
which describes the absolutely continuous Gibbs probability measure with temperature parameter~$\lambda$ and energy function~$\mathsf{L}_{\vect{z}}$ in~\eqref{EqLxy}.

Both, the ERM-DiffER and ERM-DisER problems are closely related to those typically arising while using Jayne's maximum entropy principle~\cite{jaynes1957maxEntropy1, jaynes1957maxEntropy2} for classification problems such as those in~\cite{jaakkola1999maximum, zhu2009maximum, mazuelas2022generalized}, and  ~\cite{geman1984stochastic}. 

\subsubsection{Information-Risk Minimization}
When~$Q$ is a probability measure, the ERM-RER in~\eqref{EqERMRER} is equivalent to the  \emph{information-risk minimization} (IRM) problem in~\cite{zhang2006information}.
The IRM problem in~\eqref{EqERMRER} is known to possess a unique solution equal to the Gibbs probability measure in~\eqref{EqGenpdf}, as independently shown in \cite{zhang2006information, georgii2011gibbs, gibbs1902bookElementary, xu2017information, catoni2007pac} and \cite{guedj2019primer}. 

\subsection{Bounds on the Radon-Nikodym Derivative}\label{SecFelicidadeA}
The Radon-Nikodym derivative~$\frac{\mathrm{d}P^{\left(Q, \lambda\right)}_{\vect{\Theta} | \vect{Z} = \vect{z}}}{\mathrm{d}Q}$ in~\eqref{EqGenpdf} is greater for models inducing smaller empirical risks, as shown by the following corollary of Theorem~\ref{TheoremOptimalModel}. 
\begin{corollary}\label{CorMode}
The Radon-Nikodym derivative~$\frac{\mathrm{d}P^{\left(Q, \lambda\right)}_{\vect{\Theta} | \vect{Z} = \vect{z}}}{\mathrm{d}Q}$ in~\eqref{EqGenpdf} satisfies  for all~$\left( \vect{\theta}_1, \vect{\theta}_2 \right) \in \supp Q \times \supp Q$,  with $\mathsf{L}_{\vect{z}} \left(\vect{\theta}_2 \right) \leqslant \mathsf{L}_{\vect{z}} \left(\vect{\theta}_1 \right) $, that
    \begin{equation}\label{EqpdfA}
        \frac{\mathrm{d}P^{\left(Q, \lambda\right)}_{\vect{\Theta} | \vect{Z} = \vect{z}}}{\mathrm{d}Q}\left( \vect{\theta}_1 \right)\leqslant   \frac{\mathrm{d}P^{\left(Q, \lambda\right)}_{\vect{\Theta} | \vect{Z} = \vect{z}}}{\mathrm{d}Q}\left( \vect{\theta}_2\right),
    \end{equation}
with equality if and only if $\mathsf{L}_{\vect{z}} \left(\vect{\theta}_1 \right) = \mathsf{L}_{\vect{z}} \left(\vect{\theta}_2 \right)$.
\end{corollary}
The intuition that follows from corollary~\ref{CorMode} is that under the assumption that the ERM problem in~\eqref{EqOriginalOP} possesses a solution in the support of the reference measure, i.e., $\set{T}\left( \vect{z} \right) \cap \supp Q$ is not empty, with $\set{T}\left( \vect{z} \right)$  in~\eqref{EqHatTheta}, the maximum of the  function~$\frac{\mathrm{d}P^{\left(Q, \lambda\right)}_{\vect{\Theta} | \vect{Z} = \vect{z}}}{\mathrm{d}Q}$ in~\eqref{EqGenpdf} is achieved by the models in $\set{T}\left( \vect{z} \right) \cap \supp Q$.  
When the Radon-Nikodym derivative~$\frac{\mathrm{d}P^{\left(Q, \lambda\right)}_{\vect{\Theta} | \vect{Z} = \vect{z}}}{\mathrm{d}Q}$ in~\eqref{EqGenpdf}  is either  the probability mass function in~\eqref{EqpmfGibbsSS} or  the probability density function in~\eqref{EqpdfGibbsSS},  Corollary~\ref{CorMode} shows that the elements of the set~$\set{T}\left( \vect{z} \right) \cap \supp Q$  are the \emph{modes} of the corresponding probability density function or probability mass function.

\subsection{Asymptotes of the Radon-Nikodym Derivative}\label{SecFelicidadeB}

The following lemma describes the asymptotic behavior of the Radon-Nikodym derivative~$\frac{\mathrm{d}P^{\left(Q, \lambda\right)}_{\vect{\Theta} | \vect{Z} = \vect{z}}}{\mathrm{d}Q}$ in~\eqref{EqGenpdf} when the regulariation factor increases, i.e.,~$\lambda \rightarrow +\infty$ and the reference measure $Q$ is a probability measure.

\begin{lemma}\label{CorAssymptoticProperties}
Let the measure~$Q$ in~\eqref{EqERMRER} be a probability measure. Then, for all~$\vect{\theta} \in  \supp Q$, the Radon-Nikodym derivative~$\frac{\mathrm{d}P^{\left(Q, \lambda\right)}_{\vect{\Theta} | \vect{Z} = \vect{z}}}{\mathrm{d}Q}$ in~\eqref{EqGenpdf} satisfies
\begin{IEEEeqnarray}{rcl}
\label{EqAssymptoticInftyA1p23458}
\lim_{\lambda \rightarrow +\infty} \frac{\mathrm{d}P^{\left(Q, \lambda\right)}_{\vect{\Theta} | \vect{Z} = \vect{z}}}{\mathrm{d}Q} \left( \vect{\theta} \right)  & =& 1.
\end{IEEEeqnarray}
\end{lemma}
\begin{IEEEproof}
From Theorem~\ref{TheoremOptimalModel}, it follows that for all~$\vect{\theta} \in \supp Q$,
\begin{IEEEeqnarray}{rcl}
\lim_{\lambda \rightarrow +\infty} \frac{\mathrm{d}P^{\left(Q, \lambda\right)}_{\vect{\Theta} | \vect{Z} = \vect{z}}}{\mathrm{d}Q} \left( \vect{\theta} \right)  
& = & \lim_{\lambda \rightarrow +\infty} \frac{\exp\left( -\frac{\mathsf{L}_{\vect{z}}\left( \vect{\theta}\right)}{\lambda}\right)}{\displaystyle\int\exp\left( -\frac{\mathsf{L}_{\vect{z}}\left( \vect{\nu} \right)}{\lambda} \right)  \mathrm{d} Q\left(\vect{\nu} \right)} \IEEEeqnarraynumspace \\
& = &
  \frac{1}{\displaystyle\int \mathrm{d} Q\left(\vect{\nu} \right)} \\
& = & 
 1, 
\end{IEEEeqnarray} 
where the function~$\mathsf{L}_{\vect{z}}$ is defined in~\eqref{EqLxy}. This completes the proof.
\end{IEEEproof}

Lemma~\ref{CorAssymptoticProperties} unveils the fact that, when~$Q$ is a probability measure, in the limit when~$\lambda \rightarrow +\infty$, both probability measures $P^{\left(Q, \lambda\right)}_{\vect{\Theta} | \vect{Z} = \vect{z}}$ and $Q$ are identical.
This is consistent with the fact that when~$\lambda$ tends to infinity, the optimization problem in~\eqref{EqERMRER} boils down to exclusively minimizing the relative entropy. Such minimum is zero and is observed  when both probability measures~$P^{\left(Q, \lambda\right)}_{\vect{\Theta} | \vect{Z} = \vect{z}}$ and~$Q$ are identical (Theorem~\ref{LemmaREA}).
Such intuition breaks when the reference measure is a $\sigma$-finite measure, but not a probability measure. In such a case, the relative entropy term in~\eqref{EqERMRER} might be negative and a minimum might not exist. See for instance, the case of the relative entropy between a Gaussian measure and the Lebesgue measure in \eqref{EqKeyExampleZero}, which satisfies \eqref{EqRafaelIsBorn}.

The limit  of the Radon-Nikodym derivative~$\frac{\mathrm{d}P^{\left(Q, \lambda\right)}_{\vect{\Theta} | \vect{Z} = \vect{z}}}{\mathrm{d}Q}$ in~\eqref{EqGenpdf}, when~$\lambda$ tends to zero from the right,  can be studied using the following set 
\begin{IEEEeqnarray}{rcl}
\label{EqSetL}
\set{L}_{\vect{z}} \left( \delta \right)  
& \triangleq &   \left\lbrace  \vect{\theta}  \in \set{M}:   \mathsf{L}_{\vect{z}}\left( \vect{\theta} \right)  \leqslant  \delta  \right\rbrace,
\end{IEEEeqnarray}
where the function~$\mathsf{L}_{\vect{z}}$ is defined in~\eqref{EqLxy} and~$\delta \in [0, +\infty)$.  
In particular consider the nonnegative real 
\begin{equation}\label{EqDeltaStar}
\delta_{Q,\vect{z}}^{\star} \triangleq \inf\left\lbrace \delta \in [0, +\infty) : Q\left( \set{L}_{\vect{z}} \left( \delta \right) \right) > 0 \right\rbrace.
\end{equation}
Let also~$\set{L}^{\star}_{Q,\vect{z}}$ be the following level set of the empirical risk function~$\mathsf{L}_{\vect{z}}$ in~\eqref{EqLxy}:
\begin{IEEEeqnarray}{rcl}
\label{EqSetLStar}
\set{L}^{\star}_{Q,\vect{z}}  
&  \triangleq & \left\lbrace  \vect{\theta}  \in \supp Q:   \mathsf{L}_{\vect{z}}\left( \vect{\theta} \right)  =  \delta_{Q,\vect{z}}^{\star}\right\rbrace.
\end{IEEEeqnarray}
Using this notation, the limit of the Radon-Nikodym derivative~$\frac{\mathrm{d}P^{\left(Q, \lambda\right)}_{\vect{\Theta} | \vect{Z} = \vect{z}}}{\mathrm{d}Q}$ in~\eqref{EqGenpdf}, when~$\lambda$ tends to zero from the right, is described by the following lemma.  

\begin{lemma}\label{CorAssymptoticZero}
If~$Q\left( \set{L}^{\star}_{Q,\vect{z}} \right) > 0$, with the set~$\set{L}^{\star}_{Q,\vect{z}}$ in~\eqref{EqSetLStar}  and~$Q$ the~$\sigma$-finite measure in~\eqref{EqERMRER}, then for all~$\vect{\theta} \in  \supp Q$, the Radon-Nikodym derivative~$\frac{\mathrm{d}P^{\left(Q, \lambda\right)}_{\vect{\Theta} | \vect{Z} = \vect{z}}}{\mathrm{d}Q}$ in~\eqref{EqGenpdf} satisfies
\begin{IEEEeqnarray}{rcl}
\label{EqAssymptotic0aUp}
\lim_{\lambda \rightarrow 0^{+}} \frac{\mathrm{d}P^{\left(Q, \lambda\right)}_{\vect{\Theta} | \vect{Z} = \vect{z}}}{\mathrm{d}Q} \left( \vect{\theta} \right)  & = & \frac{1}{Q\left( \set{L}^{\star}_{Q,\vect{z}}\right)} \ind{\vect{\theta} \in \set{L}^{\star}_{Q,\vect{z}} }.
\end{IEEEeqnarray}
Alternatively, if~$Q\left( \set{L}^{\star}_{Q,\vect{z}} \right) = 0$. Then, for all~$\vect{\theta} \in  \supp Q$,  
\begin{IEEEeqnarray}{rcl}
\label{EqAssymptotic0aDown}
\lim_{\lambda \rightarrow 0^{+}}  \frac{\mathrm{d}P^{\left(Q, \lambda\right)}_{\vect{\Theta} | \vect{Z} = \vect{z}}}{\mathrm{d}Q} \left( \vect{\theta} \right) & = & 
\left\lbrace
\begin{array}{cl}
+\infty & \text{ if } \vect{\theta} \in  \set{L}^{\star}_{Q,\vect{z}}\\
0 & \text{ otherwise}.
\end{array}
\right.
\IEEEeqnarraynumspace
\end{IEEEeqnarray}
\end{lemma}
\begin{IEEEproof}
The proof  is presented in Appendix~\ref{AppProofCorAssymptoticZero}.
\end{IEEEproof}

  Consider that~$Q\left( \set{L}^{\star}_{Q,\vect{z}}  \right) >  0$, with~$\set{L}^{\star}_{Q,\vect{z}}$ in~\eqref{EqSetLStar}. Under this assumption, from  Lemma~\ref{CorAssymptoticZero},  it holds that  the probability measure~$P^{\left(Q, \lambda\right)}_{\vect{\Theta} | \vect{Z} = \vect{z}}$ asymptotically concentrates on the set~$\set{L}^{\star}_{Q,\vect{z}}$ when~$\lambda$ tends to zero from the right.
More specifically, note that for all measurable sets~$\set{A} \subseteq \set{L}^{\star}_{Q,\vect{z}} \cap \supp Q$, it holds that
\begin{IEEEeqnarray}{rcl}
\label{EqSummer2022a}
\lim_{\lambda \rightarrow 0^{+}} P^{\left(Q, \lambda\right)}_{\vect{\Theta} | \vect{Z} = \vect{z}} \left( \set{A} \right) & = & \lim_{\lambda \rightarrow 0^{+}}  \int_{\set{A}} \frac{\mathrm{d}P^{\left(Q, \lambda\right)}_{\vect{\Theta} | \vect{Z} = \vect{z}}}{\mathrm{d}Q} \left( \vect{\theta} \right) \d Q \left( \vect{\theta} \right)
\IEEEeqnarraynumspace \\
 \label{EqSummer2022c}
 & = & \int  \lim_{\lambda \rightarrow 0^{+}}  \frac{\mathrm{d}P^{\left(Q, \lambda\right)}_{\vect{\Theta} | \vect{Z} = \vect{z}}}{\mathrm{d}Q} \left( \vect{\theta} \right) \ind{\vect{\theta} \in \set{A}} \d Q \left( \vect{\theta} \right)
\squeezeequ  \IEEEeqnarraynumspace \\
 \label{EqSummer2022d}
 & = & \int \frac{1}{Q\left( \set{L}^{\star}_{Q,\vect{z}} \right)} \ind{\vect{\theta} \in \set{L}^{\star}_{Q,\vect{z}} } \ind{\vect{\theta} \in \set{A}} \d Q \left( \vect{\theta} \right) 
 \Tsupersqueezeequ
 \IEEEeqnarraynumspace \\
  \label{EqSummer2022e}
  & = &  \frac{1}{Q\left(  \set{L}^{\star}_{Q,\vect{z}} \right)} \int \ind{\vect{\theta} \in \set{A}} \d Q \left( \vect{\theta} \right)\\
  \label{EqSummer2022f}
  & = &  \frac{Q\left( \set{A} \right)}{Q\left( \set{L}^{\star}_{Q,\vect{z}} \right)},
\end{IEEEeqnarray}
where the equality in~\eqref{EqSummer2022c} follows from Lemma~\ref{CorPositive} and  the dominated convergence theorem~\cite[Theorem~$2.6.9$]{ash2000probability}. The equality in~\eqref{EqSummer2022d} follows from Lemma~\ref{CorAssymptoticZero}.
In the particular case in which~$\set{A} =  \set{L}^{\star}_{Q,\vect{z}}$ in~\eqref{EqSummer2022f}, it holds that $\displaystyle\lim_{\lambda \rightarrow 0^{+}} P^{\left(Q, \lambda\right)}_{\vect{\Theta} | \vect{Z} = \vect{z}} \left(   \set{L}^{\star}_{Q,\vect{z}}  \right)$ $=$ $1$, which verifies the asymptotic concentration of the probability measure~$P^{\left(Q, \lambda\right)}_{\vect{\Theta} | \vect{Z} = \vect{z}}$ on the set~$  \set{L}^{\star}_{Q,\vect{z}}$.

Another interesting observation is that the Radon-Nikodym derivative~$\frac{\mathrm{d}P^{\left(Q, \lambda\right)}_{\vect{\Theta} | \vect{Z} = \vect{z}}}{\mathrm{d}Q}$ in~\eqref{EqGenpdf} is a constant among the elements of the set~$\set{L}^{\star}_{Q,\vect{z}}$. 
This can be assimilated to a uniform distribution of the probability among the elements of the set~$\set{L}^{\star}_{Q,\vect{z}}$ in the limit when~$\lambda$ tends to zero from the right, as previously highlighted in \cite{de2021quantitative, AthreyaGibbs, hwang1980laplace} and \cite{hasenpflug2022wasserstein}. This becomes more  evident in the case in which the set~$\set{M}$ is finite and~$Q$ is the counting measure. In such a case, the asymptotic probability  of each of the elements in~$\set{L}^{\star}_{Q,\vect{z}}$ when~$\lambda$ tends to zero from the right  is~$\frac{1}{\abs{ \set{L}^{\star}_{Q,\vect{z}}}}$.  

Consider now that~$Q\left( \set{L}^{\star}_{Q,\vect{z}}  \right) =  0$, with~$\set{L}^{\star}_{Q,\vect{z}}$ in~\eqref{EqSetLStar}. Under this assumption, in the asymptotic regime when $\lambda \to 0$, the measure $P^{\left(Q, \lambda\right)}_{\vect{\Theta} | \vect{Z} = \vect{z}}$ is not a probability measure but either the trivial measure or the infinite measure. 
This is typically the case in which~$\set{M} = \reals^d$, the measure~$Q$ is absolutely continuous with respect to the Lebesgue measure,  and the solution to the ERM problem in~\eqref{EqOriginalOP} has a unique solution on the support of $Q$, i.e.,~$\set{L}^{\star}_{Q,\vect{z}} = \set{T}\left( \vect{z} \right)$ and $\abs{ \set{T}\left( \vect{z} \right)} = 1$, which implies~$Q(\set{L}^{\star}_{Q,\vect{z}})~=~0$. 

An interesting question, which is left out of the scope of this paper, is the rate at which $P^{\left(Q, \lambda\right)}_{\vect{\Theta} | \vect{Z} = \vect{z}}$ converges to such limiting measure. The interested reader is referred to \cite{de2021quantitative, hasenpflug2022wasserstein}, and references therein. 
 
The following lemma shows that independently of whether the set~$\set{L}^{\star}_{Q,\vect{z}}$ is negligible with respect to the measure~$Q$, the limit when~$\lambda$ tends to zero from the right of~$P^{\left(Q, \lambda\right)}_{\vect{\Theta} | \vect{Z} = \vect{z}} \left( \set{L}^{\star}_{Q,\vect{z}} \right)$ is equal to one.  
\begin{lemma}\label{TheoSetLStarConcentration}
The measure~$P^{\left(Q, \lambda\right)}_{\vect{\Theta} | \vect{Z} = \vect{z}}$ in~\eqref{EqGenpdf} and the set~$\set{L}^{\star}_{Q,\vect{z}}$ in~\eqref{EqSetLStar} satisfy,  
\begin{IEEEeqnarray}{rcl}
 \lim_{\lambda \rightarrow 0^{+}} P^{\left(Q, \lambda\right)}_{\vect{\Theta} | \vect{Z} = \vect{z}}\left( \set{L}^{\star}_{Q,\vect{z}}  \right) & = & 1.
\end{IEEEeqnarray}
\end{lemma}
\begin{IEEEproof}
The proof is presented in Appendix~\ref{AppProofTheoSetLStarConcentration}.
\end{IEEEproof} 
Note that if the ERM problem in~\eqref{EqOriginalOP} possesses at least one solution and such solution is within the support of the measure~$Q$, i.e.,  $ \set{T}\left( \vect{z} \right) \cap \supp Q \neq \emptyset$, then, when~$\lambda$ tends to zero from the right, the probability measure~$P^{\left(Q, \lambda\right)}_{\vect{\Theta} | \vect{Z} = \vect{z}}$ asymptotically concentrates on the solution (or the set of solutions within the support of $Q$) to the ERM problem in~$\eqref{EqOriginalOP}$. 
Alternatively, in the case in which~$\set{L}^{\star}_{Q,\vect{z}} \cap \set{T}\left( \vect{z} \right) = \emptyset$, when~$\lambda$ tends to zero from the right, the probability measure~$P^{\left(Q, \lambda\right)}_{\vect{\Theta} | \vect{Z} = \vect{z}}$ asymptotically concentrates on a set that does not contain the set  of solutions to the ERM problem in~$\eqref{EqOriginalOP}$. 
This observation leads to the introduction to two new classes of reference measures, namely, \emph{coherent} and \emph{consistent} measures, in the following section.

\section{Reference Measures}\label{SecRefMes}

This section introduces two classes of reference measures, namely \emph{coherent} and \emph{consistent} measures, and discusses the special case of Gibbs reference measures. 

\subsection{Coherent and Consistent Reference Measures} 

A class of reference measures of particular importance to establish connections between the set of solutions to the ERM problem in~\eqref{EqOriginalOP} and the solution to the ERM-RER problem in~\eqref{EqERMRER} is that of \emph{coherent} measures. Let $\rho^{\star} \geqslant 0$ be the infimum of the empirical risk~$\mathsf{L}_{\vect{z}}$ in~\eqref{EqLxy}. That is,
\begin{IEEEeqnarray}{rcl}
\label{EqRhoStar}
\rho^{\star} \triangleq \inf \lbrace \mathsf{L}_{\vect{z}}\left( \vect{\theta} \right): \vect{\theta} \in \set{M} \rbrace.
\end{IEEEeqnarray}
Using this notation, coherent measures are defined as follows.
\begin{definition}[Coherent Measures]\label{DefCoherentMess}
The~$\sigma$-finite measure~$Q$ in~\eqref{EqERMRER} is said to be coherent if, for all~$\delta \in \left( \rho^{\star}, +\infty \right)$, with $\rho^{\star}$ in~\eqref{EqRhoStar}, it holds that
\begin{equation}
Q \left( \set{L}_{\vect{z}} \left( \delta \right)  \right) > 0,
\end{equation}
where the set~$\set{L}_{\vect{z}} \left( \delta \right)$ is defined in~\eqref{EqSetL}.
\end{definition}
When the reference measure $Q$ in the EMR-RER problem in~\eqref{EqERMRER} is a coherent measure, it holds that for all $\delta > \rho^{\star}$, the set $\set{L}_{\vect{z}} \left( \delta \right)$ in~\eqref{EqSetL} exhibits positive probability with respect to the probability measure~$P^{\left(Q, \lambda\right)}_{\vect{\Theta} | \vect{Z} = \vect{z}}$ in~\eqref{EqGenpdf}. 
The following lemma highlights this observation.
\begin{lemma}\label{LemmaHotelAmourA}
The probability measure~$P^{\left(Q, \lambda\right)}_{\vect{\Theta} | \vect{Z} = \vect{z}}$ in~\eqref{EqGenpdf}  satisfies for all~$\delta \in \left( \rho^{\star}, +\infty \right)$, with $\rho^{\star}$ in~\eqref{EqRhoStar},  that
\begin{IEEEeqnarray}{ccl}
\label{EqCrazyEight}
P^{\left(Q, \lambda\right)}_{\vect{\Theta} | \vect{Z} = \vect{z}}\left( \set{L}_{\vect{z}} \left( \delta \right)  \right) & > & 0, 
\end{IEEEeqnarray}
with~$\set{L}_{\vect{z}} \left( \delta \right)$ in~\eqref{EqSetL},
if and only if the~$\sigma$-finite measure~$Q$ in~\eqref{EqERMRER} is coherent.
\end{lemma}
\begin{IEEEproof}
The proof is presented in Appendix~\ref{AppProofLemmaHotelAmourA}.
\end{IEEEproof}

Under the assumption that the ERM problem in~\eqref{EqOriginalOP} possesses a solution, it holds that 
\begin{IEEEeqnarray}{rcl}\label{EqInfMinHappy}
\min_{\theta \in \set{M}} \mathsf{L}_{\vect{z}}\left( \vect{\theta} \right) & = &\inf \lbrace \mathsf{L}_{\vect{z}}\left( \vect{\theta} \right): \vect{\theta} \in \set{M} \rbrace.
\end{IEEEeqnarray}
Hence, when the~$\sigma$-finite measure~$Q$ in~\eqref{EqERMRER} is coherent, then~
\begin{IEEEeqnarray}{rCl}
\delta_{Q,\vect{z}}^{\star} = \rho^{\star},
\end{IEEEeqnarray} 
with~$\delta_{Q,\vect{z}}^{\star}$ in~\eqref{EqDeltaStar} and~$\rho^{\star}$ in~\eqref{EqRhoStar}, which implies that
\begin{IEEEeqnarray}{rCl}
\label{EqInclusionThor}
\set{L}^{\star}_{Q,\vect{z}} & \subseteq &\set{T}\left( \vect{z} \right) ,
\end{IEEEeqnarray}
with~$\set{T}\left( \vect{z} \right)$ in~\eqref{EqHatTheta} and~$\set{L}^{\star}_{Q,\vect{z}}$ in~\eqref{EqSetLStar}. 
This observation, together with Lemma~\ref{TheoSetLStarConcentration}, leads to the following result.

\begin{lemma}\label{LemmaSetLStarConcentrationC}
Assume that the ERM problem in~\eqref{EqOriginalOP} possesses a solution. Then, the probability measure~$P^{\left(Q, \lambda\right)}_{\vect{\Theta} | \vect{Z} = \vect{z}}$ in~\eqref{EqGenpdf} and the sets~$\set{T}\left( \vect{z} \right)$ in~\eqref{EqHatTheta} and~$\set{L}^{\star}_{Q,\vect{z}}$  in~\eqref{EqSetLStar} satisfy
\begin{IEEEeqnarray}{ccl}
\label{EqHotSummer}
 \lim_{\lambda \rightarrow 0^{+}} P^{\left(Q, \lambda\right)}_{\vect{\Theta} | \vect{Z} = \vect{z}}\left(\set{L}^{\star}_{Q,\vect{z}} \cap \set{T}\left( \vect{z}\right) \right) & = & 1, 
\end{IEEEeqnarray}
if and only if the~$\sigma$-finite measure~$Q$ in~\eqref{EqERMRER} is coherent.
\end{lemma}
\begin{IEEEproof}
The proof follows by observing that if $Q$ is a coherent measure and the ERM problem in~\eqref{EqOriginalOP} possesses a solution, the inclusion in~\eqref{EqInclusionThor} holds. 
Thus, from Lemma~\ref{TheoSetLStarConcentration}, the equality in~\eqref{EqHotSummer} holds.
Alternatively, when the measure~$Q$ in~\eqref{EqERMRER} is noncoherent,  then~$\delta_{Q,\vect{z}}^{\star} > \rho^{\star}$, which implies that~$\set{L}^{\star}_{Q,\vect{z}} \cap \set{T}\left( \vect{z} \right) = \emptyset$. Hence, from Lemma~\ref{TheoSetLStarConcentration}, it follows that 
\begin{IEEEeqnarray}{ccl}
\label{EqUberHotSummer}
 \lim_{\lambda \rightarrow 0^{+}} P^{\left(Q, \lambda\right)}_{\vect{\Theta} | \vect{Z} = \vect{z}}\left( \set{L}^{\star}_{Q,\vect{z}} \cap \set{T}\left( \vect{z}\right)  \right) & = & 0, 
\end{IEEEeqnarray}
and completes the proof.
\end{IEEEproof}

The relevance of coherent measures in ERM-RER problems is well highlighted by Lemma~\ref{LemmaSetLStarConcentrationC}. Essentially, when the ERM problem in~\eqref{EqOriginalOP} possesses at least one solution,  the concentration of the probability measure $P^{\left(Q, \lambda\right)}_{\vect{\Theta} | \vect{Z} = \vect{z}}$ in~\eqref{EqGenpdf} on the set (or a subset) of solutions to the ERM problem in~\eqref{EqOriginalOP} occurs asymptotically when $\lambda$ tends to zero from the right,  if only if the reference measure $Q$ in~\eqref{EqERMRER} is coherent.
Nonetheless, such asymptotic concentration is not a guarantee that for strictly positive values of $\lambda$ in~\eqref{EqERMRER}, the set~$\set{T}\left( \vect{z} \right)$ in~\eqref{EqHatTheta} and the measure $P^{\left(Q, \lambda\right)}_{\vect{\Theta} | \vect{Z} = \vect{z}}$ in~\eqref{EqGenpdf} satisfy $P^{\left(Q, \lambda\right)}_{\vect{\Theta} | \vect{Z} = \vect{z}}\left(\set{T}\left( \vect{z} \right)  \right) > 0$. In order to ensure this, another class of reference measures, known as \emph{consistent measures},  is introduced.

\begin{definition}[Consistent Measure]\label{DefConsistentMess}
The~$\sigma$-finite measure~$Q$ in~\eqref{EqERMRER} is said to be consistent if~$Q \left( \set{L}^{\star}_{Q,\vect{z}} \right) > 0$, with~$\set{L}^{\star}_{Q,\vect{z}}$ in~\eqref{EqSetLStar}.
\end{definition}
 
Note that every consistent measure is not necessarily coherent. For instance, if $Q$ is consistent but $\delta_{Q,\vect{z}}^{\star} > \rho^{\star}$, with $\rho^{\star}$ in~\eqref{EqRhoStar} and $\delta_{Q,\vect{z}}^{\star}$ in~\eqref{EqDeltaStar}, then, for all $\delta \in (\rho^{\star}, \delta_{Q,\vect{z}}^{\star})$, it follows that $Q \left(  \set{L}_{\vect{z}} \left( \delta \right) \right)$ = 0, and thus, $Q$ is not coherent. 
Alternatively, every coherent measure is not necessarily consistent. For instance, if $\abs{\set{L}^{\star}_{Q,\vect{z}}} = 1$ and $Q$ is coherent and absolutely continuous with respect to the Lebesgue measure, it follows that $Q \left( \set{L}^{\star}_{Q,\vect{z}} \right) = 0$, and thus, $Q$ is not consistent. 
    
The relevance of consistent measures is highlighted by the following lemma.

\begin{lemma}\label{LemmaSetLStarConcentrationCaos}
The  probability measure~$P^{\left(Q, \lambda\right)}_{\vect{\Theta} | \vect{Z} = \vect{z}}$ in~\eqref{EqGenpdf} and the set~$\set{L}^{\star}_{Q,\vect{z}}$ in~\eqref{EqSetLStar} satisfy
\begin{IEEEeqnarray}{ccl}
P^{\left(Q, \lambda\right)}_{\vect{\Theta} | \vect{Z} = \vect{z}}\left( \set{L}^{\star}_{Q,\vect{z}} \right) & > & 0, 
\end{IEEEeqnarray}
if and only if the~$\sigma$-finite measure~$Q$ in~\eqref{EqERMRER} is consistent.
\end{lemma}
\begin{IEEEproof}
When~$Q$ is nonconsistent, it holds that~$Q\left( \set{L}^{\star}_{Q,\vect{z}} \right) =0$ and thus, from the fact that the measure~$P^{\left(Q, \lambda\right)}_{\vect{\Theta} | \vect{Z} = \vect{z}}$ in~\eqref{EqGenpdf} is absolutely continuous with respect to~$Q$, it holds that  $P^{\left(Q, \lambda\right)}_{\vect{\Theta} | \vect{Z} = \vect{z}}\left( \set{L}^{\star}_{Q,\vect{z}} \right) =0$.
When~$Q$ is consistent, it holds that~$Q\left( \set{L}^{\star}_{Q,\vect{z}} \right) > 0$. Moreover, for all~$\vect{\theta} \in  \set{L}^{\star}_{Q,\vect{z}}$, it holds that~$\mathsf{L}_{\vect{z}}\left(\vect{\theta}\right) < +\infty$ and thus, from Lemma~\ref{CorPositive}, it follows that $\frac{\mathrm{d}P^{\left(Q, \lambda\right)}_{\vect{\Theta} | \vect{Z} = \vect{z}}}{\mathrm{d}Q}\left( \vect{\theta}\right) > 0$. Hence, 
\begin{IEEEeqnarray}{CCL}
\label{EqCrazyFuzz}
P^{\left(Q, \lambda\right)}_{\vect{\Theta} | \vect{Z} = \vect{z}}\left( \set{L}^{\star}_{Q,\vect{z}}  \right)  
& = & \int_{\set{L}^{\star}_{Q,\vect{z}} } \d P^{\left(Q, \lambda\right)}_{\vect{\Theta} | \vect{Z} = \vect{z}}\left( \vect{\theta} \right)\\
& = & \int_{\set{L}^{\star}_{Q,\vect{z}}} \frac{\mathrm{d}P^{\left(Q, \lambda\right)}_{\vect{\Theta} | \vect{Z} = \vect{z}}}{\mathrm{d}Q}\left( \vect{\theta}\right)  \d Q\left( \vect{\theta} \right) > 0, 
\end{IEEEeqnarray}
which completes the proof.
\end{IEEEproof}

The following lemma highlights a central property of consistent measures when the ERM problem in~\eqref{EqOriginalOP} possesses a solution.

\begin{lemma}\label{LemmaSetLStarConcentration230}
Assume that the ERM problem in~\eqref{EqOriginalOP} possesses a solution in the support of $Q$.
The probability measure~$P^{\left(Q, \lambda\right)}_{\vect{\Theta} | \vect{Z} = \vect{z}}$ in~\eqref{EqGenpdf} and the sets~$\set{T}\left( \vect{z} \right)$ in~\eqref{EqHatTheta} and~$\set{L}^{\star}_{Q,\vect{z}}$  in~\eqref{EqSetLStar} satisfy
\begin{IEEEeqnarray}{ccl}
P^{\left(Q, \lambda\right)}_{\vect{\Theta} | \vect{Z} = \vect{z}}\left( \set{L}^{\star}_{Q,\vect{z}} \cap \set{T}\left( \vect{z}\right) \right) & > & 0, 
\end{IEEEeqnarray}
if and only if the~$\sigma$-finite measure~$Q$ in~\eqref{EqERMRER} is consistent.
\end{lemma}
 \begin{IEEEproof}
The proof follows from Lemma~\ref{LemmaSetLStarConcentrationCaos} by noticing that when the ERM problem in~\eqref{EqOriginalOP} possesses a solution in the support of $Q$, the inclusion in~\eqref{EqInclusionThor} holds. 
 \end{IEEEproof}

The distinction between coherent and consistent measures becomes more evident under certain conditions. 
Consider the case in which  $\set{M}$ is finite. In this case, if the solution to the ERM problem in~\eqref{EqOriginalOP} is in the support of the $\sigma$-finite measure~$Q$, then~$Q$ is both coherent and consistent. This is essentially because all measurable singletons (models)  in $\supp Q$ exhibit positive measure with respect to $Q$. Alternatively, if the solution to the ERM problem in~\eqref{EqOriginalOP} is not in the support of $Q$, then $Q$ is consistent but not coherent. 
Consider the case in which $\set{M}$ is the set $\reals^{d}$; the loss function $\ell$ in~\eqref{EqEll} is continuous; and the ERM problem in~\eqref{EqOriginalOP} admits a unique solution. In this case, any probability measure $Q$ absolutely continuous with respect to the Lebesgue measure is a coherent measure, but it is not a consistent measure. Alternatively, if the set of solutions to the ERM problem in~\eqref{EqOriginalOP} exhibits positive Lebesgue measure, then, the measure $Q$ is both coherent and consistent. 

\subsection{Gibbs Reference Measures}\label{SecGibbsRefMes}

In model selection, a natural idea is to proceed by successive approximations in the seek of lower computation complexity. From this perspective, one might wonder whether the solution to a
	current instance of an ERM-RER problem might serve as reference measure for
	the next instance. In this section, it is shown that this yields no
	benefit. Composing two successive ERM-REM problems boils down
	to a unique ERM-RER problem with the initial reference measure and a particular regularization factor.
Under the assumption that $\lambda \in \set{K}_{Q, \vect{z}}$, with~$\set{K}_{Q, \vect{z}}$ in~\eqref{EqSetKxy}, the problem of interest is:
\begin{subequations}\label{EqERMRERb}
\begin{IEEEeqnarray}{CCl}
\label{EqERMRERbz}
    \min_{P \in \triangle_{Q}\Bormeaspace{\set{M}}} & &  \mathsf{R}_{\vect{z}} \left( P \right)  + \alpha \KL{P}{P^{\left(Q, \lambda\right)}_{\vect{\Theta} | \vect{Z} = \vect{z}}},\\
\label{EqERMRERbv}
    \mathrm{s.~t.} & & \int \d P \left( \vect{ \theta} \right) = 1,
\end{IEEEeqnarray}
\end{subequations}
where $\alpha > 0$;  the reference measure $P^{\left(Q, \lambda\right)}_{\vect{\Theta} | \vect{Z} = \vect{z}}$, which satisfies~\eqref{EqGenpdf},  is the solution of the ERM-RER problem in~\eqref{EqERMRER}; and the functional~$\mathsf{R}_{\vect{z}}$ is defined in~\eqref{EqRxy}. 
From Theorem~\ref{TheoremOptimalModel}, the solution to the ERM-RER problem in~\eqref{EqERMRERb}, which is denoted by $P^{\left(P^{\left(Q, \lambda\right)}_{\vect{\Theta} | \vect{Z} = \vect{z}}, \alpha \right)}_{\vect{\Theta} | \vect{Z} = \vect{z}}$, satisfies for all $\theta \in \supp Q$ that
\begin{IEEEeqnarray}{rcl}
 \label{EqGenpdfb}
 \frac{\mathrm{d}P^{\left(P^{\left(Q, \lambda\right)}_{\vect{\Theta} | \vect{Z} = \vect{z}}, \alpha \right)}_{\vect{\Theta} | \vect{Z} = \vect{z}}}{\mathrm{d}P^{\left(Q, \lambda\right)}_{\vect{\Theta} | \vect{Z} = \vect{z}}} \left( \vect{\theta} \right) \middlesqueezeequ 
  & =& \exp\left( - K_{P^{\left(Q, \lambda\right)}_{\vect{\Theta} | \vect{Z} = \vect{z}}, \vect{z}}\left(- \frac{1}{\alpha} \right) - \frac{1}{\alpha} \mathsf{L}_{\vect{z}}\left( \vect{\theta}\right)\right). \squeezeequ\IEEEeqnarraynumspace
\end{IEEEeqnarray}

The log-partition functions~$K_{Q,\vect{z}}$ in~\eqref{EqK}  and~$K_{ P^{\left(Q, \lambda\right)}_{\vect{\Theta}| \vect{Z} = \vect{z}},\vect{z}}$ in~\eqref{EqGenpdfb} are strongly related, as shown by the following lemma.
\begin{lemma}\label{LemmaSophia101} 
The functions~$K_{Q,\vect{z}}$ in~\eqref{EqK}  and~$K_{ P^{\left(Q, \lambda\right)}_{\vect{\Theta}| \vect{Z} = \vect{z}},\vect{z}}$ in~\eqref{EqGenpdfb} satisfy for all $t \in \reals$,
\begin{IEEEeqnarray}{rcl}
\label{EqITookBusA}
K_{ P^{\left(Q, \lambda\right)}_{\vect{\Theta}| \vect{Z} = \vect{z}},\vect{z}}\left(t\right)  & = &  K_{Q,\vect{z}}\left(t - \frac{1}{\lambda} \right) - K_{Q,\vect{z}}\left(- \frac{1}{\lambda} \right). \IEEEeqnarraynumspace 
\end{IEEEeqnarray}
Moreover, for all $t \leqslant 0$, 
\begin{IEEEeqnarray}{rcl}
\label{EqITookBuseta}
K_{ P^{\left(Q, \lambda\right)}_{\vect{\Theta}| \vect{Z} = \vect{z}},\vect{z}}\left(t\right) & \leqslant & 0.
\end{IEEEeqnarray}
\end{lemma}
\begin{IEEEproof}
The proof of~\eqref{EqITookBusA} relies on the fact that for all $t \in \left\lbrace \nu \in \reals: K_{ P^{\left(Q, \lambda\right)}_{\vect{\Theta}| \vect{Z} = \vect{z}},\vect{z}}(\nu) < \infty \right\rbrace$, the function $K_{ P^{\left(Q, \lambda\right)}_{\vect{\Theta}| \vect{Z} = \vect{z}},\vect{z}}$ in~\eqref{EqGenpdfb}  satisfies 
\begin{IEEEeqnarray}{rcl}
\nonumber
& & K_{ P^{\left(Q, \lambda\right)}_{\vect{\Theta}| \vect{Z} = \vect{z}},\vect{z}}\left( t \right) \\
\label{EqITookBusANightmare101}
 & = & \log\left( \int \exp\left(t \; \mathsf{L}_{\vect{z}}\left(\vect{\theta}\right)  \right) \mathrm{d} P^{\left(Q, \lambda\right)}_{\vect{\Theta}| \vect{Z} = \vect{z}}(\vect{\theta}) \right)
 \\
\label{EqITookBusANightmare102}
 & = & \log\left( \int \exp\left( t\; \mathsf{L}_{\vect{z}}\left(\vect{\theta}\right)  \right) \frac{\mathrm{d} P^{\left(Q, \lambda\right)}_{\vect{\Theta}| \vect{Z} = \vect{z}}}{\mathrm{d} Q}(\vect{\theta})  \; \mathrm{d} Q (\vect{\theta}) \right) 
 \\
 \label{EqITookBusANightmare103}
 & = & \log\left( \int \exp\left( \left( t -  \frac{1}{\lambda}  \right) \; \mathsf{L}_{\vect{z}}\left(\vect{\theta}\right)  - K_{Q,\vect{z}}\left(- \frac{1}{\lambda} \right)  \right) \; \mathrm{d} Q (\vect{\theta}) \right)\middlesqueezeequ \IEEEeqnarraynumspace
 \\
 \label{EqITookBusANightmare104}
 & = & \log\left( \int \exp\left( \left( t -  \frac{1}{\lambda}  \right) \; \mathsf{L}_{\vect{z}}\left(\vect{\theta}\right)    \right) \; \mathrm{d} Q (\vect{\theta}) \right) - K_{Q,\vect{z}}\left(- \frac{1}{\lambda} \right)\middlesqueezeequ
 \\
 \label{EqITookBusANightmare105}
 & = &  K_{Q,\vect{z}}\left(t - \frac{1}{\lambda} \right) - K_{Q,\vect{z}}\left(- \frac{1}{\lambda} \right),
\end{IEEEeqnarray}
where the equality in~\eqref{EqITookBusANightmare103} follows from~\eqref{EqGenpdf}. 
Moreover, from Lemma~\ref{LemmaContinuousK}, it follows that the function $K_{ P^{\left(Q, \lambda\right)}_{\vect{\Theta}| \vect{Z} = \vect{z}},\vect{z}}$ is continuous and nondecresing.  
Let $s^{\star} \in \reals \cup \lbrace +\infty \rbrace$ be defined by
\begin{IEEEeqnarray}{rcl}
s^{\star} & \triangleq & \sup \left\lbrace \nu \in \reals: K_{ P^{\left(Q, \lambda\right)}_{\vect{\Theta}| \vect{Z} = \vect{z}},\vect{z}}(\nu) < \infty \right\rbrace.
\end{IEEEeqnarray}
If $s^{\star} = + \infty$,  then for all $t \in \reals$, $K_{ P^{\left(Q, \lambda\right)}_{\vect{\Theta}| \vect{Z} = \vect{z}},\vect{z}}\left( t \right) < + \infty$, and the proof of~\eqref{EqITookBusA} is completed.

Alternatively, if $s^{\star}  < +\infty$, it follows that for all $t > s^{\star}$, $K_{ P^{\left(Q, \lambda\right)}_{\vect{\Theta}| \vect{Z} = \vect{z}},\vect{z}}(t) = + \infty$, which implies that $K_{Q,\vect{z}}\left(t - \frac{1}{\lambda} \right) = +\infty$, as the function $K_{Q,\vect{z}}$ is also continuous (Lemma~\ref{LemmaContinuousK}) and $ K_{Q,\vect{z}}\left(- \frac{1}{\lambda} \right) < \infty$ (due to the choice of $\lambda$). Hence, in this case, the equality in \eqref{EqITookBusA} is of the form $+\infty = +\infty$.
This completes the proof of~\eqref{EqITookBusA}.

The proof of~\eqref{EqITookBuseta} follows by noticing that for all $t\leqslant 0$ and for all $\vect{\theta} \in \supp Q$, it holds that $\exp\left( t \; \mathsf{L}_{\vect{z}}\left(\vect{\theta}\right)\right)\leqslant 1$. Hence, 
\begin{IEEEeqnarray}{rcl}
 K_{ P^{\left(Q, \lambda\right)}_{\vect{\Theta}| \vect{Z} = \vect{z}},\vect{z}}\left( t \right)  
 & = & \log\left( \int \exp\left(t \; \mathsf{L}_{\vect{z}}\left(\vect{\theta}\right)  \right) \mathrm{d} P^{\left(Q, \lambda\right)}_{\vect{\Theta}| \vect{Z} = \vect{z}}(\vect{\theta}) \right)\IEEEeqnarraynumspace\\
 & \leqslant & \log\left( \int  \mathrm{d} P^{\left(Q, \lambda\right)}_{\vect{\Theta}| \vect{Z} = \vect{z}}(\vect{\theta}) \right)\\
 & = & 0,
\end{IEEEeqnarray}
 which completes the proof.
\end{IEEEproof}

The following lemma establishes that the solution to the ERM-RER problem in~\eqref{EqERMRERb} is identical to the solution to another ERM-RER problem of the form
\begin{subequations}\label{EqERMRERc}
\begin{IEEEeqnarray}{CCl}
\label{EqERMRERca}
    \min_{P \in \triangle_{Q}\Bormeaspace{\set{M}}} & &  \mathsf{R}_{\vect{z}} \left( P \right)  + \left( \frac{1}{\frac{1}{\alpha} + \frac{1}{\lambda}} \right) \KL{P}{Q},\IEEEeqnarraynumspace\\
\label{EqERMRERcb}
    \mathrm{s.~t.} & & \int \d P \left( \vect{ \theta} \right) = 1,
\end{IEEEeqnarray}
\end{subequations}
with $\lambda \in \set{K}_{Q, \vect{z}}$, with~$\set{K}_{Q, \vect{z}}$ in~\eqref{EqSetKxy}, and whose solution, denoted by $P^{\left(Q, \frac{1}{\frac{1}{\lambda} + \frac{1}{\alpha}} \right)}_{\vect{\Theta}| \vect{Z} = \vect{z}}$, satisfies for all $\vect{\theta} \in \supp  Q$, 
\begin{IEEEeqnarray}{rcl}
\nonumber
&  &\frac{\mathrm{d} P^{\left(Q, \frac{1}{\frac{1}{\lambda} + \frac{1}{\alpha}} \right)}_{\vect{\Theta}| \vect{Z} = \vect{z}} }{\mathrm{d}Q}(\vect{\theta})\\
 \label{EqGenpdfc}
  & =& \exp\left( - K_{Q, \vect{z}}\left(- \frac{1}{\lambda} - \frac{1}{\alpha} \right) - \left(\frac{1}{\lambda} + \frac{1}{\alpha} \right) \mathsf{L}_{\vect{z}}\left( \vect{\theta}\right)\right).
\end{IEEEeqnarray}
The formal statement is as follows.

\begin{lemma}\label{LemmaSophia101A} 
Let $\alpha \in \left(0,+\infty \right)$ and $\lambda \in \set{K}_{Q, \vect{z}}$, with~$\set{K}_{Q, \vect{z}}$ in~\eqref{EqSetKxy}. Then, the probability measures~$P^{\left(P^{\left(Q, \lambda\right)}_{\vect{\Theta} | \vect{Z} = \vect{z}}, \alpha \right)}_{\vect{\Theta} | \vect{Z} = \vect{z}}$ in~\eqref{EqGenpdfb} and~$P^{\left(Q, \frac{1}{\frac{1}{\lambda} + \frac{1}{\alpha}} \right)}_{\vect{\Theta}| \vect{Z} = \vect{z}}$ in~\eqref{EqGenpdfc}  are identical.
\end{lemma}
\begin{IEEEproof}
For all $\vect{\theta} \in \supp  Q$, 
\begin{IEEEeqnarray}{rcl}
\nonumber
&   &\frac{\mathrm{d} P^{\left(P^{\left(Q, \lambda\right)}_{\vect{\Theta}| \vect{Z} = \vect{z}}, \alpha \right)}_{\vect{\Theta}| \vect{Z} = \vect{z}}}{\mathrm{d}Q}(\vect{\theta}) \\
\label{EqSilencioCervezaYTequila9877}
& = &  \frac{\mathrm{d} P^{\left(P^{\left(Q, \lambda\right)}_{\vect{\Theta}| \vect{Z} = \vect{z}}, \alpha \right)}_{\vect{\Theta}| \vect{Z} = \vect{z}}}{\mathrm{d}P^{\left(Q, \lambda\right)}_{\vect{\Theta}| \vect{Z} = \vect{z}}}(\vect{\theta})     \frac{\mathrm{d}P^{\left(Q, \lambda\right)}_{\vect{\Theta}| \vect{Z} = \vect{z}}}{\mathrm{d}Q} \left( \vect{\theta} \right) \\
\nonumber
& = & \exp \Bigg( - K_{P^{\left(Q, \lambda\right)}_{\vect{\Theta}| \vect{Z} = \vect{z}},\vect{z}}\left(- \frac{1}{\alpha} \right) - K_{Q,\vect{z}}\left(- \frac{1}{\lambda} \right) \\
\label{EqSilencioCervezaYTequila9878}
& &- \left( \frac{1}{\alpha} + \frac{1}{\lambda} \right)\mathsf{L}_{\vect{z}}\left( \vect{\theta}\right)\Bigg)\middlesqueezeequ \IEEEeqnarraynumspace\\
\label{EqSilencioCervezaYTequila9879}
& = & \exp\left( - K_{Q,\vect{z}}\left(- \frac{1}{\alpha} - \frac{1}{\lambda}\right) - \left( \frac{1}{\alpha} + \frac{1}{\lambda} \right)\mathsf{L}_{\vect{z}}\left( \vect{\theta}\right)\right) \\
\label{EqSilencioCervezaYTequila98710}
& = &\frac{\mathrm{d} P^{\left(Q, \frac{1}{\frac{1}{\lambda} + \frac{1}{\alpha}} \right)}_{\vect{\Theta}| \vect{Z} = \vect{z}} }{\mathrm{d}Q}(\vect{\theta}),
\end{IEEEeqnarray}
where 
the equality in~\eqref{EqSilencioCervezaYTequila9877} follows from the fact that the measure $P^{\left(P^{\left(Q, \lambda\right)}_{\vect{\Theta}| \vect{Z} = \vect{z}}, \alpha \right)}_{\vect{\Theta}| \vect{Z} = \vect{z}}$ is absolutely continuous with respect to $P^{\left(Q, \lambda\right)}_{\vect{\Theta}| \vect{Z} = \vect{z}}$ and $P^{\left(Q, \lambda\right)}_{\vect{\Theta}| \vect{Z} = \vect{z}}$ is absolutely continuous with respect to the measure $Q$;
the equality in~\eqref{EqSilencioCervezaYTequila9878} follows from Lemma~\ref{LemmaSophia101}; and 
the equality in~\eqref{EqSilencioCervezaYTequila98710} follows from Theorem~\ref{TheoremOptimalModel}.

For all measurable subsets $\set{A}$ of $\set{M}$, the following holds:
\begin{IEEEeqnarray}{rcl}
P^{\left(P^{\left(Q, \lambda\right)}_{\vect{\Theta}| \vect{Z} = \vect{z}}, \alpha \right)}_{\vect{\Theta}| \vect{Z} = \vect{z}} \left( \set{A} \right) & = &\int_{\set{A}} \frac{\mathrm{d} P^{\left(P^{\left(Q, \lambda\right)}_{\vect{\Theta}| \vect{Z} = \vect{z}}, \alpha \right)}_{\vect{\Theta}| \vect{Z} = \vect{z}}}{\mathrm{d}Q}(\vect{\theta}) \mathrm{d}Q (\vect{\theta}) \\
\label{EqWithThorInAgadir}
& = &\int_{\set{A}} \frac{\mathrm{d} P^{\left(Q, \frac{1}{\frac{1}{\lambda} + \frac{1}{\alpha}} \right)}_{\vect{\Theta}| \vect{Z} = \vect{z}} }{\mathrm{d}Q} \mathrm{d}Q (\vect{\theta})  \\
& = &\int_{\set{A}}  \mathrm{d}P^{\left(Q, \frac{1}{\frac{1}{\lambda} + \frac{1}{\alpha}} \right)}_{\vect{\Theta}| \vect{Z} = \vect{z}}(\vect{\theta}) \\
& = & P^{\left(Q, \frac{1}{\frac{1}{\lambda} + \frac{1}{\alpha}} \right)}_{\vect{\Theta}| \vect{Z} = \vect{z}} (\set{A}),
\end{IEEEeqnarray}
where the equality in~\eqref{EqWithThorInAgadir} follows from~\eqref{EqSilencioCervezaYTequila98710}. This completes the proof. 
\end{IEEEproof}

The following theorem establishes a relation between the solutions to the following optimization problems
\begin{subequations}\label{EqImprovementp235a8}
\begin{IEEEeqnarray}{ccl}
\min_{P \in \triangle_{Q}\Bormeaspace{\set{M}}} & &  \mathsf{R}_{\vect{z}} \left( P \right) ,\\
\text{s. t.} & \quad & \KL{P}{P^{\left(Q, \lambda\right)}_{\vect{\Theta}| \vect{Z} = \vect{z}}} \leqslant c, \mbox{ and } \IEEEeqnarraynumspace\\
& \quad & \int\mathrm{d} P (\vect{\theta}) = 1,
\end{IEEEeqnarray}
\end{subequations}
and 
\begin{subequations}\label{EqERMRERzeta}
\begin{IEEEeqnarray}{CCl}
    \min_{P \in \triangle_{Q}\Bormeaspace{\set{M}}} & &  \mathsf{R}_{\vect{z}} \left( P \right)  + \omega \KL{P}{Q},\\
    \mathrm{s.~t.} & & \int \d P \left( \vect{ \theta} \right) = 1,
\end{IEEEeqnarray}
\end{subequations}
 with 
$c > 0$ and $\omega \in \set{K}_{Q, \vect{z}}$, with~$\set{K}_{Q, \vect{z}}$ in~\eqref{EqSetKxy}, two constants;
$P^{\left(Q, \lambda\right)}_{\vect{\Theta}| \vect{Z}=\vect{z}}$  the probability measure  in~\eqref{EqGenpdf}; and  
$\mathsf{R}_{\vect{z}}$ the functional in~\eqref{EqRxy}.

From Theorem~\ref{TheoremOptimalModel}, the solution to the ERM-RER problem in~\eqref{EqERMRERzeta}, which is denoted by $P^{\left(Q, \omega \right)}_{\vect{\Theta} | \vect{Z} = \vect{z}}$, satisfies for all $\theta \in \supp Q$ that
\begin{IEEEeqnarray}{rcl}\label{EqGenSolB}
\frac{\mathrm{d}P^{\left(Q, \omega\right)}_{\vect{\Theta}| \vect{Z} = \vect{z}}}{\mathrm{d}Q} \left( \vect{\theta} \right) 
  & =& \exp\left( - K_{Q,\vect{z}}\left(- \frac{1}{\omega} \right) - \frac{1}{\omega} \mathsf{L}_{\vect{z}}\left( \vect{\theta}\right)\right),
\end{IEEEeqnarray}
where the function~$K_{Q,\vect{z}}$ is in~\eqref{EqK}.

The following theorem formalizes the relation between both optimization problems.
\begin{theorem}\label{TheoremSensitivityB}
Assume that $c$ and $\omega$ in~\eqref{EqImprovementp235a8} and~\eqref{EqERMRERzeta} satisfy
\begin{IEEEeqnarray}{rcl}
\label{EqBombasticElastic}
\KL{P^{\left(Q, \omega\right)}_{\vect{\Theta}| \vect{Z} = \vect{z}}}{P^{\left(Q, \lambda\right)}_{\vect{\Theta}| \vect{Z} = \vect{z}}} = c,
\end{IEEEeqnarray}
with $P^{\left(Q, \lambda\right)}_{\vect{\Theta} | \vect{Z} = \vect{z}}$ and~$P^{\left(Q, \omega\right)}_{\vect{\Theta}| \vect{Z} = \vect{z}}$ being the probability measures in~\eqref{EqGenpdf} and ~\eqref{EqGenSolB}, respectively. 
Then, the solution to the optimization problem in~\eqref{EqImprovementp235a8} is the probability measure $P^{\left(Q, \omega\right)}_{\vect{\Theta}| \vect{Z} = \vect{z}}$.
\end{theorem}
\begin{IEEEproof}
The proof is presented in Appendix~\ref{AppProofTheoremSensitivityB}.
\end{IEEEproof}

\section{The Log-Partition Function}\label{SecLogPartition}

This section introduces some properties of the log-partition function~$K_{Q, \vect{z}}$ in~\eqref{EqK}  using the notion of \emph{separable} empirical risk functions. 

\subsection{Separable Empirical Risk Functions}

Separable empirical risk functions are defined with respect to a measure~$P \in \triangle\left( \set{M} \right)$.

\begin{definition}[Separable Empirical Risk Function]\label{DefSeparableLxy}
The empirical risk function~$\mathsf{L}_{\vect{z}}$ in~\eqref{EqLxy} is said to be separable with respect to a~$\sigma$-finite measure~$P\in \triangle\left( \set{M} \right)$, if there exist a positive real~$c > 0$  and two subsets~$\set{A}$ and~$\set{B}$ of~$\set{M}$  that are nonnegligible with respect to~$P$, and for all~$(\vect{\theta}_1,\vect{\theta}_2) \in \set{A} \times \set{B}$, 
\begin{IEEEeqnarray}{rcl}
\label{EqTwoNonnegligibleSets}
 \mathsf{L}_{\vect{z}} \left( \vect{\theta}_1 \right) &< c <& \mathsf{L}_{\vect{z}}\left(\vect{\theta}_2\right) < +\infty.
\end{IEEEeqnarray}
\end{definition}
In a nutshell, a nonseparable empirical risk function  with respect to the measure~$Q$ is a constant almost surely. More specifically, there exists a real~$a \geqslant 0$, such that
\begin{equation}
Q\left( \left\lbrace \vect{\theta} \in \set{M}: \mathsf{L}_{\vect{z}}\left(\vect{\theta}\right)  = a \right\rbrace\right) = 1.
\end{equation}
From this perspective, nonseparable empirical risk functions exhibit little practical interest for model selection. 
 
The definition of separability in Definition~\ref{DefSeparableLxy} and Lemma~\ref{LemmaMutualAC} lead to the following lemma.
 \begin{lemma}\label{CorAC}
The empirical risk function~$\mathsf{L}_{\vect{z}}$ in~\eqref{EqLxy} is separable with respect to the~$\sigma$-finite measure~$Q$  in~\eqref{EqERMRER} if and only if it is separable with respect to the probability measure~$P^{\left(Q, \lambda\right)}_{\vect{\Theta} | \vect{Z} = \vect{z}}$ in~\eqref{EqGenpdf}.
\end{lemma}
\begin{IEEEproof}
Consider first that the function~$\mathsf{L}_{\vect{z}}$ is separable with respect to the~$\sigma$-finite measure~$Q$. Hence, there exist a positive real~$c > 0$  and two subsets~$\set{A}$ and~$\set{B}$ of~$\set{M}$  that are nonnegligible with respect to~$Q$, such that for all~$(\vect{\theta}_1,\vect{\theta}_2) \in \set{A} \times \set{B}$ the inequality in~\eqref{EqTwoNonnegligibleSets} holds. Hence, from~\eqref{EqTwoNonnegligibleSets} the following inequalities hold: 
\begin{IEEEeqnarray}{rcl}
-\frac{1}{\lambda} \mathsf{L}_{\vect{z}} \left( \vect{\theta}_1 \right) &>&  -\frac{c}{\lambda}  >  -\frac{1}{\lambda}  \mathsf{L}_{\vect{z}}\left(\vect{\theta}_2\right)  >  -\infty, \mbox{ and }\IEEEeqnarraynumspace\\
\exp\left(-\frac{1}{\lambda} \mathsf{L}_{\vect{z}} \left( \vect{\theta}_1 \right) \right) &> & \exp\left(-\frac{c}{\lambda}\right)   >   \exp\left(-\frac{1}{\lambda}  \mathsf{L}_{\vect{z}}\left(\vect{\theta}_2\right) \right)   >  0. \squeezeequ\IEEEeqnarraynumspace
\end{IEEEeqnarray}
This implies that
\begin{IEEEeqnarray}{rCl}
\label{EqCalorIntenso}
\frac{\mathrm{d}P^{\left(Q, \lambda\right)}_{\vect{\Theta} | \vect{Z} = \vect{z}}}{\mathrm{d}Q} \left( \vect{\theta}_1 \right)  &> & \exp\left(- K_{Q, \vect{z}}\left(- \frac{1}{\lambda} \right)  -\frac{c}{\lambda}\right)  \\
& > &  \frac{\mathrm{d}P^{\left(Q, \lambda\right)}_{\vect{\Theta} | \vect{Z} = \vect{z}}}{\mathrm{d}Q} \left( \vect{\theta}_2 \right)  \\
& > & 0. \squeezeequ
\end{IEEEeqnarray}
Using the inequality in~\eqref{EqCalorIntenso} and the facts that~$Q \left( \set{A} \right) >0$ and~$Q \left( \set{B} \right) >0$, the following holds
\begin{IEEEeqnarray}{rcl}
\label{EqShouldBeOnTheBeach1}
P^{\left(Q, \lambda\right)}_{\vect{\Theta} | \vect{Z} = \vect{z}}\left( \set{A} \right) & = & \int_{\set{A}} \frac{\mathrm{d}P^{\left(Q, \lambda\right)}_{\vect{\Theta} | \vect{Z} = \vect{z}}}{\mathrm{d}Q} \left( \vect{\theta} \right) \d Q \left( \vect{\theta} \right) 
> 0,
\end{IEEEeqnarray}
and
\begin{IEEEeqnarray}{rcl}
\label{EqShouldBeOnTheBeach2}
P^{\left(Q, \lambda\right)}_{\vect{\Theta} | \vect{Z} = \vect{z}}\left( \set{B} \right) & = & \int_{\set{B}} \frac{\mathrm{d}P^{\left(Q, \lambda\right)}_{\vect{\Theta} | \vect{Z} = \vect{z}}}{\mathrm{d}Q} \left( \vect{\theta} \right) \d Q \left( \vect{\theta} \right) > 0.
\end{IEEEeqnarray}
which implies that the function~$\mathsf{L}_{\vect{z}}$ is separable with respect to the probability measure~$P^{\left(Q, \lambda\right)}_{\vect{\Theta} | \vect{Z} = \vect{z}}$.

Consider now that the function~$\mathsf{L}_{\vect{z}}$ is separable with respect to the probability measure~$P^{\left(Q, \lambda\right)}_{\vect{\Theta} | \vect{Z} = \vect{z}}$. Hence, there exist a positive real~$c > 0$  and two subsets~$\set{A}$ and~$\set{B}$ of~$\set{M}$  that are nonnegligible with respect to~$P^{\left(Q, \lambda\right)}_{\vect{\Theta} | \vect{Z} = \vect{z}}$, such that for all~$(\vect{\theta}_1,\vect{\theta}_2) \in \set{A} \times \set{B}$ the inequality in~\eqref{EqTwoNonnegligibleSets} holds. More specifically, $P^{\left(Q, \lambda\right)}_{\vect{\Theta} | \vect{Z} = \vect{z}}\left( \set{A} \right) > 0$ and $P^{\left(Q, \lambda\right)}_{\vect{\Theta} | \vect{Z} = \vect{z}}\left( \set{B} \right) > 0$. 
From Lemma~\ref{CorPositive} and the inequality in~\eqref{EqTwoNonnegligibleSets}, it follows that for all pairs~$(\vect{\theta}_1,\vect{\theta}_2) \in \set{A} \times \set{B}$,~$\frac{\mathrm{d}P^{\left(Q, \lambda\right)}_{\vect{\Theta} | \vect{Z} = \vect{z}}}{\mathrm{d}Q} \left( \vect{\theta}_1 \right) > 0$ and~$\frac{\mathrm{d}P^{\left(Q, \lambda\right)}_{\vect{\Theta} | \vect{Z} = \vect{z}}}{\mathrm{d}Q} \left( \vect{\theta}_2 \right) > 0$.
Hence, from the fact that $P^{\left(Q, \lambda\right)}_{\vect{\Theta} | \vect{Z} = \vect{z}}\left( \set{A} \right) > 0$ and $P^{\left(Q, \lambda\right)}_{\vect{\Theta} | \vect{Z} = \vect{z}}\left( \set{B} \right) > 0$, it follows that~$Q \left( \set{A} \right) >0$ and~$Q \left( \set{B} \right) >0$, which implies that the function~$\mathsf{L}_{\vect{z}}$ is separable with respect to the~$\sigma$-finite measure~$Q$.
This completes the proof.
\end{IEEEproof}
Lemma~\ref{CorAC} shows that separable empirical risk functions, and only these functions, lead to ERM-RER-optimal probability measures from which models are sampled with different probabilities. For the case of nonseparable empirical risk functions, all models are sampled from the ERM-RER-optimal probability measure with the same probability.

\subsection{Properties of the Log-Partition Function}

The log-partition function~$K_{Q, \vect{z}}$ in~\eqref{EqK} is a nondecreasing continuous convex function as shown by the following lemmas.
\begin{lemma}\label{LemmaContinuousK}
The function~$K_{Q, \vect{z}}$ in~\eqref{EqK} is nondecreasing and differentiable infinitely many times in the interior of $\left\lbrace t \in \reals : K_{Q, \vect{z}}(t) < +\infty \right\rbrace$. 
\end{lemma}
\begin{IEEEproof}
The proof is presented in Appendix~\ref{AppProofLemmaContinuousK}.
\end{IEEEproof}

\begin{lemma}\label{LemmaConvexK}
The function~$K_{Q, \vect{z}}$ in~\eqref{EqK} is convex in $\left\lbrace t \in \reals : K_{Q, \vect{z}}(t) < +\infty \right\rbrace$. Moreover, it is strictly convex if and only if the empirical risk function~$\mathsf{L}_{\vect{z}}$ in~\eqref{EqLxy} is separable with respect to the~$\sigma$-finite measure~$Q$  in~\eqref{EqERMRER}. 
\end{lemma}
\begin{IEEEproof}
The proof is presented in Appendix~\ref{AppProofLemmaConvexK}.
\end{IEEEproof}

In Lemma~\ref{LemmaContinuousK}, it has been established that the log-partition function~$K_{Q, \vect{z}}$ in~\eqref{EqK} is differentiable infinitely many times in the interval $\left\lbrace t \in \reals : K_{Q, \vect{z}}(t) < +\infty \right\rbrace$. 
Let the~$m$-th derivative of the function~$K_{Q, \vect{z}}$ in~\eqref{EqK} be denoted by~$K^{(m)}_{Q, \vect{z}}: \reals \rightarrow \reals$, with~$m \in \ints$. Hence, for all~$s \in  \set{K}_{Q, \vect{z}}$,
\begin{IEEEeqnarray}{rcl}
\label{EqK1}\label{EqK123}
K^{(m)}_{Q, \vect{z}} \left(- \frac{1}{s} \right) & \triangleq &  \frac{\mathrm{d}^{m}}{\mathrm{d}t^{m}} K_{Q, \vect{z}}\left(t \right)\Bigr|_{t = - \frac{1}{s}}.
\end{IEEEeqnarray}

The following lemma provides explicit expressions for the first, second and third derivatives of the function~$K_{Q, \vect{z}}$ in~\eqref{EqK}.
\begin{lemma}\label{CorDerivatives}
The first, second and third derivatives of the function~$K_{Q, \vect{z}}$ in~\eqref{EqK}, denoted respectively by~$K^{(1)}_{Q, \vect{z}}$,~$K^{(2)}_{Q, \vect{z}}$, and~$K^{(3)}_{Q, \vect{z}}$, satisfy for all~$\lambda \in \mathrm{int} \set{K}_{Q, \vect{z}}$, with~$\set{K}_{Q, \vect{z}}$ in~\eqref{EqSetKxy}, 
\begin{IEEEeqnarray}{rcl}
\label{EqNiceK1}
K^{(1)}_{Q, \vect{z}} \left(-\frac{1}{\lambda} \right) & = &   \int \mathsf{L}_{\vect{z}}\left( \vect{\theta} \right)  \mathrm{d}P^{\left(Q, \lambda\right)}_{\vect{\Theta} | \vect{Z} = \vect{z}}  (\vect{\theta}), \\
\label{EqNiceK2}
 K^{(2)}_{Q, \vect{z}}\left(-\frac{1}{\lambda} \right)  & = &   \int \left( \mathsf{L}_{\vect{z}}\left( \vect{\theta} \right)  - K^{(1)}_{Q, \vect{z}}\left(-\frac{1}{\lambda} \right)  \right)^2 \mathrm{d}P^{\left(Q, \lambda\right)}_{\vect{\Theta} | \vect{Z} = \vect{z}} (\vect{\theta}), \squeezeequ \IEEEeqnarraynumspace   \\
\label{EqNiceK3}
 K^{(3)}_{Q, \vect{z}}\left(-\frac{1}{\lambda} \right)  & = &   \int \left( \mathsf{L}_{\vect{z}}\left( \vect{\theta} \right)  - K^{(1)}_{Q, \vect{z}}\left(-\frac{1}{\lambda} \right)  \right)^3 \mathrm{d}P^{\left(Q, \lambda\right)}_{\vect{\Theta} | \vect{Z} = \vect{z}} (\vect{\theta}), \middlesqueezeequ\IEEEeqnarraynumspace
\end{IEEEeqnarray}
where the function~$\mathsf{L}_{\vect{z}}$ is defined in~\eqref{EqLxy} and  the measure~$P^{\left(Q, \lambda\right)}_{\vect{\Theta} | \vect{Z} = \vect{z}}$ satisfies~\eqref{EqGenpdf}.
\end{lemma}
\begin{IEEEproof}
The proof is presented in Appendix~\ref{ProofCorDerivatives}.
\end{IEEEproof}
From Lemma~\ref{CorDerivatives}, it follows that if~$\vect{\Theta} \sim P^{\left(Q, \lambda\right)}_{\vect{\Theta} | \vect{Z} = \vect{z}}$, with $P^{\left(Q, \lambda\right)}_{\vect{\Theta} | \vect{Z} = \vect{z}}$ in~\eqref{EqGenpdf}, the random variable 
\begin{equation}\label{EqW}
W \triangleq \mathsf{L}_{\vect{z}}\left( \vect{\Theta} \right),
\end{equation}
with the function~$\mathsf{L}_{\vect{z}}$ in~\eqref{EqLxy}, possesses a mean,  variance, and third cumulant that are equivalent to~$K^{(1)}_{Q, \vect{z}}\left(-\frac{1}{\lambda}\right)$ in~\eqref{EqNiceK1}, $K^{(2)}_{Q, \vect{z}}\left( -\frac{1}{\lambda} \right)$ in~\eqref{EqNiceK2}, and~$K^{(3)}_{Q, \vect{z}}\left( -\frac{1}{\lambda} \right)$ in~\eqref{EqNiceK3}, respectively.

Note that if there exists a $\delta > 0$ such that the log-partition function $K_{Q, \vect{z}}$ is differentiable within the open interval $(-\delta, \delta)$ and $Q$ in \eqref{EqERMRER} is a probability measure, the function $K_{Q, \vect{z}}$ is the cumulant generating function of the random variable
\begin{IEEEeqnarray}{rcl}
\label{EqV}
V \triangleq \mathsf{L}_{\vect{z}}\left( \vect{\Theta} \right), \mbox{ with $\vect{\Theta} \sim Q$}.
\end{IEEEeqnarray}
The following lemma leverages this observation.
\begin{lemma}\label{CorDerivativesQ}
Assume that $Q$ in \eqref{EqERMRER} is a probability measure and that  there exists real $\delta >0$ such that the log-partition function~$K_{Q, \vect{z}}$ in~\eqref{EqK} is differentiable within  $(-\delta, \delta)$. Then, the first, second and third derivatives of $K_{Q, \vect{z}}$,  denoted respectively by~$K^{(1)}_{Q, \vect{z}}$,~$K^{(2)}_{Q, \vect{z}}$, and~$K^{(3)}_{Q, \vect{z}}$, satisfy 
\begin{IEEEeqnarray}{rcl}
\label{EqNiceKQ1}
K^{(1)}_{Q, \vect{z}} \left( 0 \right) & = &   \int \mathsf{L}_{\vect{z}}\left( \vect{\theta} \right)  \mathrm{d}Q (\vect{\theta}), \\
\label{EqNiceKQ2}
 K^{(2)}_{Q, \vect{z}}\left( 0 \right)  & = &   \int \left( \mathsf{L}_{\vect{z}}\left( \vect{\theta} \right)  - K^{(1)}_{Q, \vect{z}}\left( 0 \right)  \right)^2 \mathrm{d} Q  (\vect{\theta}), \middlesqueezeequ \IEEEeqnarraynumspace   \\
\label{EqNiceKQ3}
 K^{(3)}_{Q, \vect{z}}\left( 0 \right)  & = &   \int \left( \mathsf{L}_{\vect{z}}\left( \vect{\theta} \right)  - K^{(1)}_{Q, \vect{z}}\left( 0 \right)  \right)^3 \mathrm{d}Q (\vect{\theta}), \middlesqueezeequ\IEEEeqnarraynumspace
\end{IEEEeqnarray}
where the function~$\mathsf{L}_{\vect{z}}$ is defined in~\eqref{EqLxy}.
\end{lemma}
\begin{IEEEproof}
The proof follows along the same arguments of the proof of Lemma~\ref{CorDerivatives}.
\end{IEEEproof}

The mean,  variance, and third cumulant of the random variable $V$ in~\eqref{EqV} are~$K^{(1)}_{Q, \vect{z}}\left( 0 \right)$ in~\eqref{EqNiceKQ1}, $K^{(2)}_{Q, \vect{z}}\left( 0 \right)$ in~\eqref{EqNiceKQ2}, and~$K^{(3)}_{Q, \vect{z}}\left( 0 \right)$ in~\eqref{EqNiceKQ3}, respectively.

\section{Expectation of the Empirical Risk}\label{SecExpectation}

The mean of the random variable~$W$ in~\eqref{EqW} is equivalent to  the expectation of the  empirical risk function~$\mathsf{L}_{\vect{z}}$ with respect to the probability measure~$P^{\left(Q, \lambda\right)}_{\vect{\Theta} | \vect{Z} = \vect{z}}$ in~\eqref{EqGenpdf}, which is equal to~$\mathsf{R}_{\vect{z}}\left(P^{\left(Q, \lambda\right)}_{\vect{\Theta} | \vect{Z} = \vect{z}} \right)$, with the functional~$\mathsf{R}_{\vect{z}}$ in~\eqref{EqRxy}. Often,~$\mathsf{R}_{\vect{z}}\left(P^{\left(Q, \lambda\right)}_{\vect{\Theta} | \vect{Z} = \vect{z}} \right)$ is referred to as the ERM-RER-optimal expected empirical risk to emphasize that this is the expected value of the empirical risk when models are sampled from the solution of the ERM-RER problem in~\eqref{EqERMRER}.
The following corollary of Lemma~\ref{CorDerivatives} formalizes this observation.

\begin{corollary}\label{CorRK}
The probability measure~$P^{\left(Q, \lambda\right)}_{\vect{\Theta} | \vect{Z} = \vect{z}}$ in~\eqref{EqGenpdf} verifies that  
\begin{equation}\label{EqRK}
\mathsf{R}_{\vect{z}}\left(P^{\left(Q, \lambda\right)}_{\vect{\Theta} | \vect{Z} = \vect{z}} \right) = K^{(1)}_{Q, \vect{z}}\left(-\frac{1}{\lambda}\right),
\end{equation}
where the functional~$\mathsf{R}_{\vect{z}}$ and the function~$K^{(1)}_{Q, \vect{z}}$ are defined in~\eqref{EqRxy} and~\eqref{EqNiceK1}, respectively.
\end{corollary}

The expected empirical risk~$\mathsf{R}_{\vect{z}}\left(P^{\left(Q, \lambda\right)}_{\vect{\Theta} | \vect{Z} = \vect{z}} \right)$ in~\eqref{EqRK} exhibits the following property.
\begin{theorem}\label{CorDecreasingAverage}
The expected empirical risk~$\mathsf{R}_{\vect{z}}\left(P^{\left(Q, \lambda\right)}_{\vect{\Theta} | \vect{Z} = \vect{z}} \right)$ in~\eqref{EqRK} is nondecreasing with~$\lambda \in \set{K}_{Q, \vect{z}}$, with~$\set{K}_{Q, \vect{z}}$ in~\eqref{EqSetKxy}.
Moreover,~$\mathsf{R}_{\vect{z}}\left(P^{\left(Q, \lambda\right)}_{\vect{\Theta} | \vect{Z} = \vect{z}} \right)$ is strictly increasing  with~$\lambda \in \set{K}_{Q, \vect{z}}$  if and only if the function~$\mathsf{L}_{\vect{z}}$ in~\eqref{EqLxy} is separable with respect to the measure~$Q$.
\end{theorem}
\begin{IEEEproof}
The proof is presented in Appendix~\ref{AppProofCorDecreasingAverage}.
\end{IEEEproof}
The expected empirical risk~$\mathsf{R}_{\vect{z}}\left(P^{\left(Q, \lambda\right)}_{\vect{\Theta} | \vect{Z} = \vect{z}} \right)$ in~\eqref{EqRK}  has been shown to be nondecreasing with $\lambda$ in \cite[Appendix E.4]{aminian2021exact} for the special case in which $Q$ is a probability measure.

A question that arises from Theorem~\ref{CorDecreasingAverage} is whether  the value~$\mathsf{R}_{\vect{z}}\left(P^{\left(Q, \lambda\right)}_{\vect{\Theta} | \vect{Z} = \vect{z}} \right)$ in~\eqref{EqRK} can be made arbitrarily close to~$\delta_{Q,\vect{z}}^{\star}$, with~$\delta_{Q,\vect{z}}^{\star}$ in~\eqref{EqDeltaStar}, by making~$\lambda$  arbitrarily small. 
The following lemma shows that  the  value~$\mathsf{R}_{\vect{z}}\left(P^{\left(Q, \lambda\right)}_{\vect{\Theta} | \vect{Z} = \vect{z}} \right)$ is often bounded away from  $\delta_{Q,\vect{z}}^{\star}$, even for arbitrarily small values of~$\lambda$. 
\begin{lemma}\label{LemmaRzLowerBoundB}
The expected empirical risk~$\mathsf{R}_{\vect{z}}\left(P^{\left(Q, \lambda\right)}_{\vect{\Theta} | \vect{Z} = \vect{z}} \right)$ in~\eqref{EqRK} satisfies,
\begin{equation}\label{EqLimitKB}
\mathsf{R}_{\vect{z}}\left(P^{\left(Q, \lambda\right)}_{\vect{\Theta} | \vect{Z} = \vect{z}} \right) \geqslant \delta_{Q,\vect{z}}^{\star},
\end{equation} 
where~$\delta_{Q,\vect{z}}^{\star}$ is defined in~\eqref{EqDeltaStar}.
Moreover, the inequality in~\eqref{EqLimitKB} is strict if and only if the function~$\mathsf{L}_{\vect{z}}$ in~\eqref{EqLxy} is separable with respect to the measure~$Q$ in~\eqref{EqERMRER}. 
\end{lemma}
\begin{IEEEproof}
The proof is presented in Appendix~\ref{AppProofLemmaRzLowerBoundB}.
\end{IEEEproof}

In the asymptotic regime when~$\lambda$ tends to zero,  the expected empirical risk $\mathsf{R}_{\vect{z}}\left(P^{\left(Q, \lambda\right)}_{\vect{\Theta} | \vect{Z} = \vect{z}} \right)$ in~\eqref{EqRK} is equal to~$\delta_{Q,\vect{z}}^{\star}$, as shown by the following lemma.

\begin{theorem}\label{CorAssymptoticMean}
The expected empirical risk~$\mathsf{R}_{\vect{z}}\left(P^{\left(Q, \lambda\right)}_{\vect{\Theta} | \vect{Z} = \vect{z}} \right)$ in~\eqref{EqRK} satisfies,  
\begin{IEEEeqnarray}{rcl}
\label{EqAssymptotic0brrqqr}
\lim_{\lambda \rightarrow 0^{+}} \mathsf{R}_{\vect{z}}\left(P^{\left(Q, \lambda\right)}_{\vect{\Theta} | \vect{Z} = \vect{z}} \right) 
& = & \delta_{Q,\vect{z}}^{\star},
\end{IEEEeqnarray}
where~$\delta_{Q,\vect{z}}^{\star}$ is defined in~\eqref{EqDeltaStar}.
\end{theorem}
\begin{IEEEproof}
The proof is presented in Appendix~\ref{AppProofCorAssymptoticMean}.
\end{IEEEproof} 
 
The following lemma determines the value of the objective function of the ERM-RER problem in~\eqref{EqERMRER} when it is evaluated at its solution. This result appeared first in \cite[Lemma~$3$]{Perlaza-ISIT2023b}.  

\begin{lemma}[Lemma~$3$ in~\cite{Perlaza-ISIT2023b}]\label{LemmaAgadir} 
The probability measure~$P^{\left(Q, \lambda\right)}_{\vect{\Theta}| \vect{Z} = \vect{z}}$ in~\eqref{EqGenpdf} and the $\sigma$-finite measure $Q$ in~\eqref{EqERMRER} satisfy
\begin{IEEEeqnarray}{rcl}
\label{EqAgadirEarlyMorningA}
- \lambda K_{Q,\vect{z}}\left(- \frac{1}{\lambda} \right) & = & 
\mathsf{R}_{\vect{z}}\left( P^{\left(Q, \lambda\right)}_{\vect{\Theta}| \vect{Z} = \vect{z}} \right)  + \lambda \KL{P^{\left(Q, \lambda\right)}_{\vect{\Theta}| \vect{Z} = \vect{z}}}{Q} \middlesqueezeequ  \IEEEeqnarraynumspace\\
\label{EqAgadirEarlyMorningB}
& = & \mathsf{R}_{\vect{z}}\left( Q  \right) - \lambda \KL{Q}{P^{\left(Q, \lambda\right)}_{\vect{\Theta}| \vect{Z} = \vect{z}}} 
\end{IEEEeqnarray}
where
the functional~$\mathsf{R}_{\vect{z}}$ is defined in~\eqref{EqRxy};  and 
the function~$K_{Q,\vect{z}}$ is defined in~\eqref{EqK}.
\end{lemma}
\begin{IEEEproof}
From Theorem~\ref{TheoremOptimalModel}, it follows that for all~$\vect{\theta} \in \supp Q$, 
\begin{IEEEeqnarray}{rcl}
\label{EqDeniedVisaUSA}
\log\left(\frac{\mathrm{d}P^{\left(Q, \lambda\right)}_{\vect{\Theta}| \vect{Z} = \vect{z}}}{\mathrm{d}Q} \left( \vect{\theta} \right) \right)
& = &  - K_{Q,\vect{z}}\left(- \frac{1}{\lambda} \right) - \frac{1}{\lambda} \mathsf{L}_{\vect{z}}\left( \vect{\theta}\right),
\end{IEEEeqnarray}
where the function~$\mathsf{L}_{\vect{z}}$ is defined in~\eqref{EqLxy}. 
Thus, 
\begin{IEEEeqnarray}{rcl}
\KL{P^{\left(Q, \lambda\right)}_{\vect{\Theta}| \vect{Z} = \vect{z}}}{Q} \squeezeequ& = & \int  \log\left(\frac{\mathrm{d}P^{\left(Q, \lambda\right)}_{\vect{\Theta}| \vect{Z} = \vect{z}}}{\mathrm{d}Q} \left( \vect{\theta} \right)\right) \mathrm{d}P^{\left(Q, \lambda\right)}_{\vect{\Theta}| \vect{Z} = \vect{z}} \left( \vect{\theta} \right) \\
& = &  - K_{Q,\vect{z}}\left(- \frac{1}{\lambda} \right) \hspace{-0.5ex} - \hspace{-0.5ex} \frac{1}{\lambda}\int \hspace{-0.5ex} \mathsf{L}_{\vect{z}}\left( \vect{\theta}\right) \mathrm{d}P^{\left(Q, \lambda\right)}_{\vect{\Theta}| \vect{Z} = \vect{z}} \left( \vect{\theta} \right) \squeezeequ \IEEEeqnarraynumspace\\
\label{EqTheProofAgadiri1}
& = &  - K_{Q,\vect{z}}\left(- \frac{1}{\lambda} \right) - \frac{1}{\lambda}\mathsf{R}_{\vect{z}}\left( P^{\left(Q, \lambda\right)}_{\vect{\Theta}| \vect{Z} = \vect{z}} \right),
\end{IEEEeqnarray}
where the functional~$\mathsf{R}_{\vect{z}}$ is defined in~\eqref{EqRxy}. This completes the proof of~\eqref{EqAgadirEarlyMorningA}.

From Lemma~\ref{LemmaMutualAC} and~\eqref{EqDeniedVisaUSA}, it follows that
\begin{IEEEeqnarray}{rcl}
\KL{Q}{P^{\left(Q, \lambda\right)}_{\vect{\Theta}| \vect{Z} = \vect{z}}}& = & - \int  \log\left(\frac{\mathrm{d}P^{\left(Q, \lambda\right)}_{\vect{\Theta}| \vect{Z} = \vect{z}}}{\mathrm{d}Q} \left( \vect{\theta} \right)\right) \mathrm{d}Q\left( \vect{\theta} \right)\IEEEeqnarraynumspace \\
& = &   K_{Q,\vect{z}}\left(- \frac{1}{\lambda} \right) + \frac{1}{\lambda}\int \mathsf{L}_{\vect{z}}\left( \vect{\theta}\right) \mathrm{d}Q \left( \vect{\theta} \right) \\
\label{EqXYcaprice1}
& = &  K_{Q,\vect{z}}\left(- \frac{1}{\lambda} \right) + \frac{1}{\lambda}\mathsf{R}_{\vect{z}}\left( Q \right),
\end{IEEEeqnarray}
which completes the proof of ~\eqref{EqAgadirEarlyMorningB}.
\end{IEEEproof} 
The following corollary of Lemma~\ref{LemmaAgadir} characterizes the difference between the expected values of the random variables $W$ and $V$ in \eqref{EqW} and \eqref{EqV}, respectively. 
\begin{corollary}\label{EqLastMinuteMean}
If  measures $Q$ and~$P^{\left(Q, \lambda\right)}_{\vect{\Theta}| \vect{Z} = \vect{z}}$ in~\eqref{EqGenpdf} are both probability measures, then, 
\begin{IEEEeqnarray}{rcl}
\nonumber
& & \mathsf{R}_{\vect{z}}\left( Q  \right) - \mathsf{R}_{\vect{z}}\left( P^{\left(Q, \lambda\right)}_{\vect{\Theta}| \vect{Z} = \vect{z}} \right)  \\
\label{EqNiceEarlyMorning5432}
& = & \lambda \left( \KL{Q}{P^{\left(Q, \lambda\right)}_{\vect{\Theta}| \vect{Z} = \vect{z}}} + \KL{P^{\left(Q, \lambda\right)}_{\vect{\Theta}| \vect{Z} = \vect{z}}}{Q}  \right).
\end{IEEEeqnarray}
\end{corollary}
The right-hand side of \eqref{EqNiceEarlyMorning5432} is a symmetrized Kullback-Liebler divergence, also known as Jeffrey's divergence~\cite{jeffreys1946invariant}, between the measures $Q$ and $P^{\left(Q, \lambda\right)}_{\vect{\Theta}| \vect{Z} = \vect{z}}$.
More importantly, when $Q$ is a probability measure, it follows that  $\KL{P^{\left(Q, \lambda\right)}_{\vect{\Theta}| \vect{Z} = \vect{z}}}{Q}\geqslant 0$ and $\KL{Q}{P^{\left(Q, \lambda\right)}_{\vect{\Theta}| \vect{Z} = \vect{z}}}\geqslant0$, which leads to the following corollary from Lemma~\ref{LemmaAgadir}.

\begin{corollary}\label{CorSensitivityIneq}
If the $\sigma$-finite measure $Q$ in~\eqref{EqERMRER} is a probability measure, then, the probability measure~$P^{\left(Q, \lambda\right)}_{\vect{\Theta}| \vect{Z} = \vect{z}}$ in~\eqref{EqGenpdf} satisfies
\begin{IEEEeqnarray}{l} 
 \mathsf{R}_{\vect{z}}\left( P^{\left(Q, \lambda\right)}_{\vect{\Theta}| \vect{Z} = \vect{z}} \right) \leqslant \mathsf{R}_{\vect{z}}\left( Q \right),
\end{IEEEeqnarray}
where, the functional~$\mathsf{R}_{\vect{z}}$ is defined in~\eqref{EqRxy}.
\end{corollary}
 
\section{Variance of the Empirical Risk}\label{SecVariance}

In Lemma~\ref{LemmaContinuousK}, it has been established that if there exists a $\delta > 0$ such that the log-partition function $K_{Q, \vect{z}}$ in~\eqref{EqK} is finite within the open interval $(-\delta, \delta)$,  the log-partition function $K_{Q, \vect{z}}$ is differentiable infinitely many times within the interval $(-\infty, \delta)$. This together with the \emph{mean value theorem} \cite[Theorem $5.10$]{rudin1953bookPrinciples} lead to the following characterization of the differences of the values $K^{(2)}_{Q, \vect{z}}\left( -\frac{1}{t} \right)$ and $K^{(2)}_{Q, \vect{z}}\left( 0 \right)$, with $t > 0$.
\begin{lemma}\label{LemmaDiffVariances}
If the measure  $Q$ in \eqref{EqERMRER} is a probability measure and there exists a $\delta > 0$ such that  the function $K_{Q, \vect{z}}$ in \eqref{EqK} is differentiable within the open interval $(-\delta, \delta)$, then for all $t > 0$,
\begin{IEEEeqnarray}{rcl}
\label{EqDaysBeforeHisBirth}
K^{(2)}_{Q, \vect{z}}\left( -\frac{1}{t} \right) - K^{(2)}_{Q, \vect{z}}\left(  0 \right) & = & -\frac{1}{t}K^{(3)}_{Q, \vect{z}}\left(  -\frac{1}{\beta} \right) < +\infty, \IEEEeqnarraynumspace
\end{IEEEeqnarray}
for some $\beta \in \left( t, + \infty \right)$, where  the functions~$K^{(2)}_{Q, \vect{z}}$, and~$K^{(3)}_{Q, \vect{z}}$ are defined in~\eqref{EqK123}.
\end{lemma}
\begin{IEEEproof}
The proof is an immediate consequence of Lemma~\ref{LemmaContinuousK} and the mean value theorem \cite[Theorem $5.10$]{rudin1953bookPrinciples}. 
\end{IEEEproof}
The relevance of Lemma~\ref{LemmaDiffVariances} lies on the fact that $K^{(2)}_{Q, \vect{z}}\left( -\frac{1}{\lambda} \right)$ and $K^{(2)}_{Q, \vect{z}}\left( 0 \right)$ are the variances of the random variables $W$ in \eqref{EqW} and $V$ in \eqref{EqV}. See Lemma~\ref{CorDerivatives} and Lemma~\ref{CorDerivativesQ}.
Under the assumptions of Lemma~\ref{LemmaDiffVariances}, it follows that the function $K^{(3)}_{Q, \vect{z}}$ is continuous in $(-\infty, \delta)$, where $\delta > 0$. Hence, for all $t>0$, the function $K^{(3)}_{Q, \vect{z}}$ achieves a maximum and a minimum within the interval $\left[ -\frac{1}{t}, 0\right]$.  Such extrema allow providing  lower and upper bounds on the variance $K^{(2)}_{Q, \vect{z}}\left( -\frac{1}{\lambda} \right)$ of the random variable $W$  in terms of the variance $K^{(2)}_{Q, \vect{z}}\left(  0 \right)$ of the random variable $V$, as shown hereunder.
\begin{corollary}\label{CorBoundsOnTheVariance}
If the measure  $Q$ in \eqref{EqERMRER} is a probability measure and there exists a $\delta > 0$ such that  the function $K_{Q, \vect{z}}$ in \eqref{EqK} is differentiable within the open interval $(-\delta, \delta)$, then for all $t > 0$,
\begin{IEEEeqnarray}{rcl}
\label{EqDaysBeforeHisBirthBis}
K^{(2)}_{Q, \vect{z}}\left(  0 \right) - \frac{1}{t}  c_2& \leqslant  K^{(2)}_{Q, \vect{z}}\left( -\frac{1}{t} \right)  \leqslant &  K^{(2)}_{Q, \vect{z}}\left(  0 \right) - \frac{1}{t} c_1 ,
\end{IEEEeqnarray}
where,
\begin{IEEEeqnarray}{rcl}
c_1 & =  & \min_{s \in \left[ -\frac{1}{t}, 0\right]} K^{(3)}_{Q, \vect{z}}\left(  s \right) \mbox{ and }\\
c_2 & =  & \max_{s \in \left[ -\frac{1}{t}, 0\right]} K^{(3)}_{Q, \vect{z}}\left(  s \right),
\end{IEEEeqnarray}
and  the functions~$K^{(2)}_{Q, \vect{z}}$, and~$K^{(3)}_{Q, \vect{z}}$ are defined in~\eqref{EqK123}.
\end{corollary}

The inequality in \eqref{EqDaysBeforeHisBirthBis} reveals that under the assumptions of Corollary~\ref{CorBoundsOnTheVariance}, in the asymptotic regime when $t \to +\infty$, the variances of the random variables $W$ in~\eqref{EqW} and $V$ in~\eqref{EqV} are identical.  
Additionally, unlike the means $K^{(1)}_{Q, \vect{z}}\left( -\frac{1}{\lambda} \right)$ and $K^{(1)}_{Q, \vect{z}}\left( 0 \right)$ of the random variables $W$ and $V$, which satisfy $K^{(1)}_{Q, \vect{z}}\left( -\frac{1}{\lambda} \right) \leqslant K^{(1)}_{Q, \vect{z}}\left( 0 \right)$ (Corollary~\ref{CorSensitivityIneq}), their variances $K^{(2)}_{Q, \vect{z}}\left( -\frac{1}{\lambda} \right)$ and $K^{(2)}_{Q, \vect{z}}\left( 0 \right)$ might satisfy $K^{(2)}_{Q, \vect{z}}\left( -\frac{1}{\lambda} \right) < K^{(2)}_{Q, \vect{z}}\left( 0 \right)$ or $K^{(2)}_{Q, \vect{z}}\left( -\frac{1}{\lambda} \right) \geqslant K^{(2)}_{Q, \vect{z}}\left( 0 \right)$ depending on whether the function $K^{(3)}_{Q, \vect{z}}$ is positive or negative within the interval $[-\frac{1}{\lambda}, 0]$.
 
Lemma~\ref{LemmaDiffVariances} shows that the monotonicity of the expectation of the random variable~$W$ in~\eqref{EqW} with respect to $\lambda$, stated by Theorem~\ref{CorDecreasingAverage}, is not a property exhibited by the  variance nor the third cumulant. The following example highlights this observation.

\begin{example}\label{ExampleDecreasingVariance}
Consider the ERM-RER problem in~\eqref{EqERMRER}, under the assumption that~$Q$ is a probability measure and the empirical risk function~$\mathsf{L}_{\vect{z}}$ in~\eqref{EqLxy} is such that for all $\vect{\theta} \in \set{M}$,
\begin{equation}\label{EqLZeroOne}
\mathsf{L}_{\vect{z}}\left( \vect{\theta} \right) = \left\lbrace 
\begin{array}{lcl}   
0 & \text{ if } &  \vect{\theta} \in \set{A}\\
1 & \text{ if } &  \vect{\theta} \in \set{M}\setminus\set{A},
\end{array}
 \right. 
\end{equation}
where  the sets~$\set{A} \subset \set{M}$ and~$\set{M}\setminus\set{A}$ are nonnegligible with respect to the reference probability measure~$Q$.
In this case, the function~$K_{Q, \vect{z}}$ in~\eqref{EqK} satisfies for all~$\lambda>0$, 
\begin{IEEEeqnarray}{rcl}
\label{EqThisKeko}
K_{Q, \vect{z}}\left( - \frac{1}{\lambda} \right)  
& = & \log\left( Q\left( \set{A} \right)  + \exp\left(- \frac{1}{\lambda}\right) \left( 1-Q\left( \set{A} \right)  \right)  \right). \squeezeequ\IEEEeqnarraynumspace
\end{IEEEeqnarray}
The derivatives~$K^{(1)}_{Q, \vect{z}}$,~$K^{(2)}_{Q, \vect{z}}$,  and~$K^{(3)}_{Q, \vect{z}}$ in~\eqref{EqK123}  of the function~$K_{Q, \vect{z}}$ in~\eqref{EqThisKeko} satisfy for all~$\lambda>0$, 
\begin{IEEEeqnarray}{rcl}
\label{EqDecreasingVariance}
K^{(1)}_{Q, \vect{z}}\left( - \frac{1}{\lambda} \right)  \middlesqueezeequ& = & \frac{ \exp\left(- \frac{1}{\lambda}\right) \left( 1-Q\left( \set{A} \right)  \right) }{ Q\left( \set{A} \right)  + \exp\left(- \frac{1}{\lambda}\right) \left( 1-Q\left( \set{A} \right)  \right)}; \\
K^{(2)}_{Q, \vect{z}}\left( - \frac{1}{\lambda} \right) \middlesqueezeequ & = & \frac{ Q\left( \set{A} \right) \left( 1-Q\left( \set{A} \right)  \right)\exp\left(- \frac{1}{\lambda}\right) }{\left( Q\left( \set{A} \right)  + \exp\left(- \frac{1}{\lambda}\right) \left( 1-Q\left( \set{A} \right)  \right) \right)^2}; \mbox{ and }\IEEEeqnarraynumspace\\
\nonumber
K^{(3)}_{Q, \vect{z}}\left( - \frac{1}{\lambda} \right)  \middlesqueezeequ & = &K^{(2)}_{Q, \vect{z}}\left( - \frac{1}{\lambda} \right)  \\
& & \left(  \frac{Q\left( \set{A} \right)  - \left( 1-Q\left( \set{A} \right)  \right)\exp\left(- \frac{1}{\lambda}\right)  }{ Q\left( \set{A} \right)  + \exp\left(- \frac{1}{\lambda}\right) \left( 1-Q\left( \set{A} \right)  \right) } \right).\middlesqueezeequ\IEEEeqnarraynumspace
\end{IEEEeqnarray}
Note that~$K^{(3)}_{Q, \vect{z}}\left( - \frac{1}{\lambda} \right) > 0$ if and only if  \begin{equation}
\label{EqConditionK3}
Q\left( \set{A} \right)  - \left( 1-Q\left( \set{A} \right)  \right)\exp\left(- \frac{1}{\lambda}\right) >0.
\end{equation}
Assume that~$Q\left( \set{A} \right) \geqslant \frac{1}{2}$. Thus, it holds that for all~$\lambda > 0$, the inequality in~\eqref{EqConditionK3} is always satisfied.  This follows from observing that  for all~$\lambda >0$,
\begin{IEEEeqnarray}{rcl}
\exp\left(- \frac{1}{\lambda}\right) < 1 \leqslant \frac{Q\left( \set{A} \right)}{1-Q\left( \set{A} \right)}.
\end{IEEEeqnarray}
Hence, if~$Q\left( \set{A} \right) \geqslant \frac{1}{2}$, for all decreasing sequences of positive reals~$\lambda_{1} > \lambda_{2} >  \ldots > 0$, it holds that
\begin{IEEEeqnarray}{rcl}
\label{EqDecreasingVariancei7}
\frac{1}{4} \geqslant K^{(2)}_{Q, \vect{z}}\left( - \frac{1}{\lambda_1} \right) > K^{(2)}_{Q, \vect{z}}\left( - \frac{1}{\lambda_2} \right)  > \ldots > 0.
\end{IEEEeqnarray}
Alternatively, assume that~$Q\left( \set{A} \right) < \frac{1}{2}$.  In this case, the inequality in~\eqref{EqConditionK3} is satisfied if and only if 
\begin{equation}
\lambda < \left( \log\left(  \frac{1-Q\left( \set{A} \right)}{Q\left( \set{A} \right)} \right) \right)^{-1}.
\end{equation}
Hence, if~$Q\left( \set{A} \right) < \frac{1}{2}$, then for all decreasing sequences of positive reals~$$\left( \log\left(  \frac{1-Q\left( \set{A} \right)}{Q\left( \set{A} \right)} \right) \right)^{-1} > \lambda_{1} > \lambda_{2} > \ldots > 0,$$ it holds that
\begin{IEEEeqnarray}{rcl}
\label{EqDecreasingVarianceg0}
\frac{1}{4} > K^{(2)}_{Q, \vect{z}}\left( - \frac{1}{\lambda_1} \right) > K^{(2)}_{Q, \vect{z}}\left( - \frac{1}{\lambda_2} \right) >  \ldots > 0.
\end{IEEEeqnarray}
Moreover, for all decreasing sequences of positive reals~$$\lambda_{1} > \lambda_{2}  > \ldots > \left( \log\left(  \frac{1-Q\left( \set{A} \right)}{Q\left( \set{A} \right)} \right) \right)^{-1},$$ it holds that
\begin{IEEEeqnarray}{rcl}
\label{EqDecreasingVarianceg097}
K^{(2)}_{Q, \vect{z}}\left( - \frac{1}{\lambda_1} \right) < K^{(2)}_{Q, \vect{z}}\left( - \frac{1}{\lambda_2} \right)  < \ldots < \frac{1}{4}.
\end{IEEEeqnarray}
The upperbound by~$\frac{1}{4}$ in~\eqref{EqDecreasingVariancei7},~\eqref{EqDecreasingVarianceg0} and~\eqref{EqDecreasingVarianceg097} follows by noticing that  the value $K^{(2)}_{Q, \vect{z}} \left( -\frac{1}{\lambda}\right)$ is maximized when~$\lambda = \left( \log\left(  \frac{1-Q\left( \set{A} \right)}{Q\left( \set{A} \right)} \right) \right)^{-1}$ and~$K^{(2)}_{Q, \vect{z}} \left( -\frac{1}{\lambda}\right) = \frac{1}{4}$.
\end{example}

\begin{figure}[h]
\begin{center}
\includegraphics[width=\linewidth]{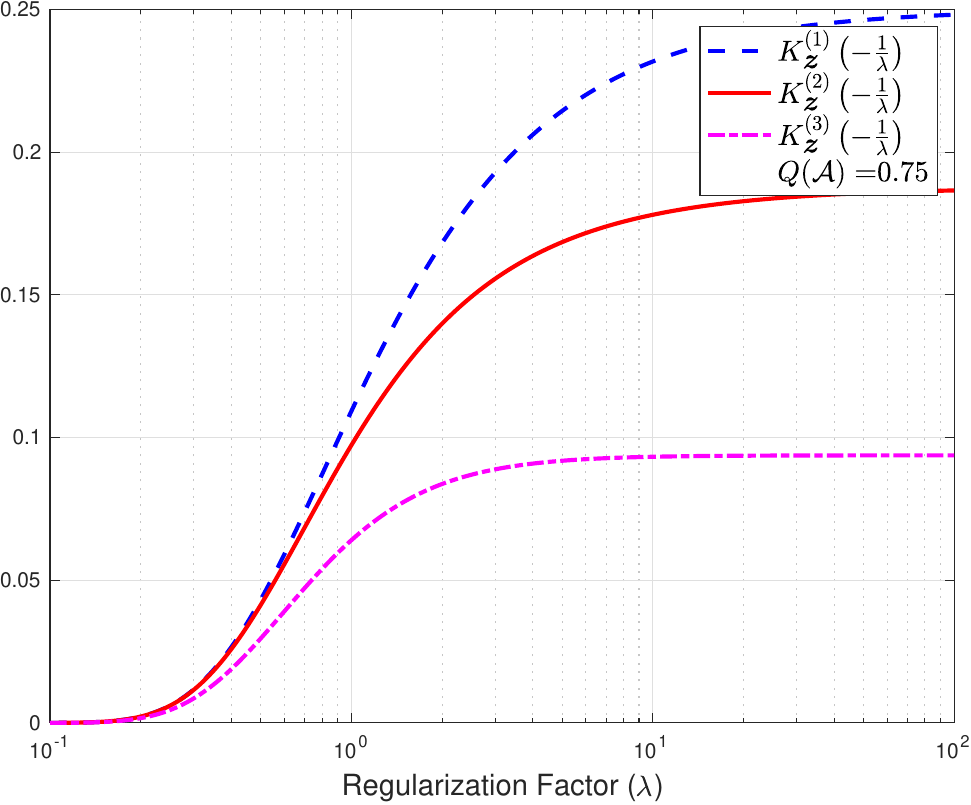}
\end{center}
\caption{Mean $K^{(1)}_{Q, \vect{z}}\left( - \frac{1}{\lambda} \right)$, variance $K^{(2)}_{Q, \vect{z}}\left( - \frac{1}{\lambda} \right)$, and third central moment $K^{(3)}_{Q, \vect{z}}\left( - \frac{1}{\lambda} \right)$ of the empirical risk in Example~\ref{ExampleDecreasingVariance}, with $Q\left( \set{A} \right) = \frac{3}{4}$}
\label{FigQThreequarter}
\end{figure}

\begin{figure}[h!]
\begin{center}
\includegraphics[width=\linewidth]{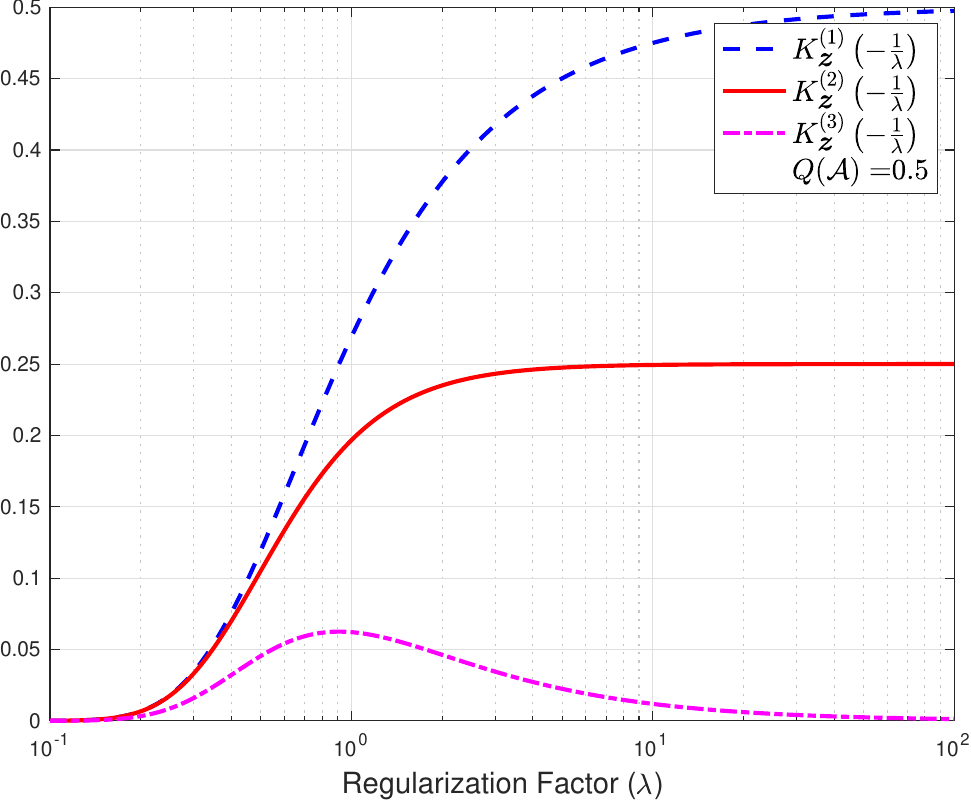}
\end{center}
\caption{Mean $K^{(1)}_{Q, \vect{z}}\left( - \frac{1}{\lambda} \right)$, variance $K^{(2)}_{Q, \vect{z}}\left( - \frac{1}{\lambda} \right)$, and third central moment $K^{(3)}_{Q, \vect{z}}\left( - \frac{1}{\lambda} \right)$ of the empirical risk in Example~\ref{ExampleDecreasingVariance}, with $Q\left( \set{A} \right) = \frac{1}{2}$}
\label{FigQHalf}
\end{figure}

Example~\ref{ExampleDecreasingVariance} provides important insights on the choice of the reference measure~$Q$. Note for instance that when the reference measure assigns a probability to the set of models~$\set{T}\left( \vect{z} \right)$ in~\eqref{EqHatTheta}  that is greater than or equal to the probability of suboptimal models~$\set{M}\setminus\set{T}\left( \vect{z} \right)$, i.e.,~$Q \left( \set{T}\left( \vect{z} \right) \right) \geqslant \frac{1}{2}$, the variance is strictly decreasing to zero when~$
\lambda$ decreases. See for instance, Figure~\ref{FigQThreequarter} and Figure~\ref{FigQHalf}. That is, when the reference measure assigns higher probability to the set of solutions to the ERM problem in~\eqref{EqOriginalOP}, the variance is monotone with respect to the parameter~$\lambda$.

\begin{figure}[h!]
\begin{center}
\includegraphics[width=\linewidth]{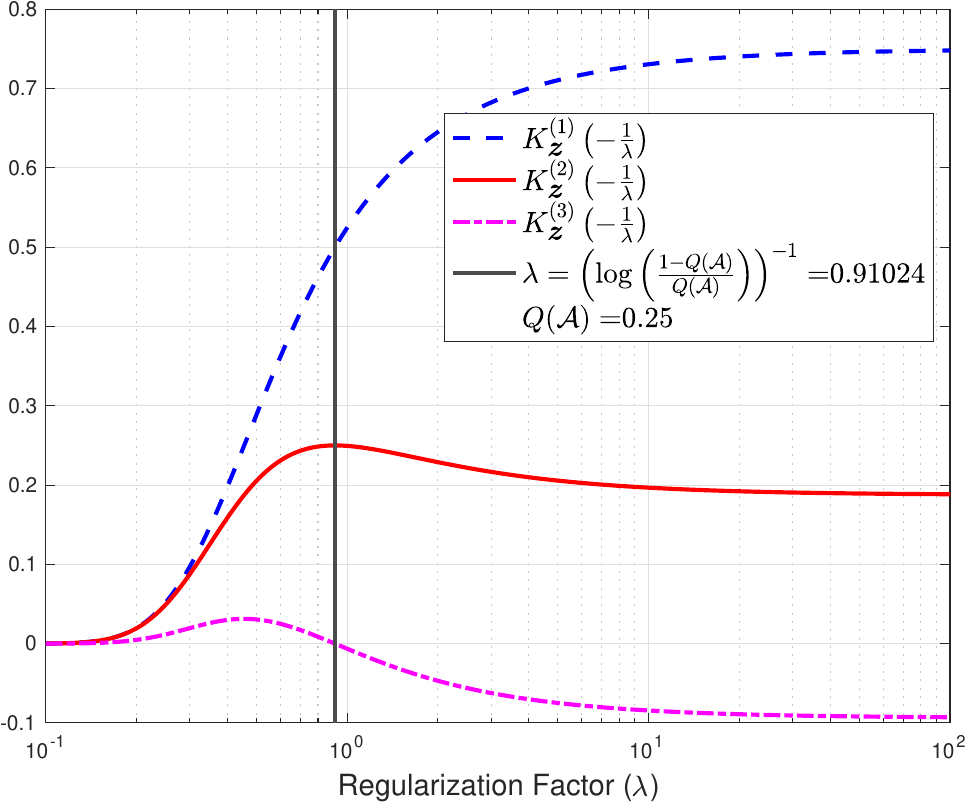}
\end{center}
\caption{Mean $K^{(1)}_{Q, \vect{z}}\left( - \frac{1}{\lambda} \right)$, variance $K^{(2)}_{Q, \vect{z}}\left( - \frac{1}{\lambda} \right)$, and third central moment $K^{(3)}_{Q, \vect{z}}\left( - \frac{1}{\lambda} \right)$ of the empirical risk in Example~\ref{ExampleDecreasingVariance}, with $Q\left( \set{A} \right) = \frac{1}{4}$}
\label{FigQOnequarter}
\end{figure}

Alternatively, when the reference measure assigns a probability to the set~$\set{T}\left( \vect{z} \right)$ that is smaller than the probability of the set~$\set{M}\setminus\set{T}\left( \vect{z} \right)$, i.e.,~$Q \left( \set{T}\left( \vect{z} \right) \right) < \frac{1}{2}$, there exists a critical point for~$\lambda$ at~$\left( \log\left(  \frac{1-Q\left( \set{A} \right)}{Q\left( \set{A} \right)} \right) \right)^{-1}$. See for instance, Figure~\ref{FigQOnequarter}. More importantly, such a critical point can be arbitrarily close to zero depending on  the value~$Q\left( \set{A} \right)$.
The variance strictly decreases when~$\lambda$ decreases beyond the value~$\left( \log\left(  \frac{1-Q\left( \set{A} \right)}{Q\left( \set{A} \right)} \right) \right)^{-1}$. Otherwise, reducing~$\lambda$ above the value~$\left( \log\left(  \frac{1-Q\left( \set{A} \right)}{Q\left( \set{A} \right)} \right) \right)^{-1}$ increases the variance. 

In general, these observations suggest that reference measures~$Q$ that allocate small measures to the sets containing the set~$\set{T}\left( \vect{z} \right)$ might require reducing the value~$\lambda$ beyond a small threshold in order to observe small values of~$K^{(2)}_{Q, \vect{z}}\left( - \frac{1}{\lambda} \right)$, which is the variance of the random variable~$W$, in~\eqref{EqW}.
These observations are central to understanding the concentration of probability that occurs when~$\lambda$ decreases to zero, as discussed in Section~\ref{SecMonotonicConcentration}.


\section{Cumulant Generating Function of the Empirical Risk}\label{SecSubGaussianity}
Consider the transport of the measure~$P^{\left(Q, \lambda \right)}_{\vect{\Theta} | \vect{Z} = \vect{z}}$ in~\eqref{EqGenpdf} from~$\Bormeaspace{\set{M}}$ to~$\Bormeaspace{[0, +\infty]}$ through the function~$\mathsf{L}_{\vect{z}}$ in~\eqref{EqLxy}. 
Denote the resulting probability measure in~$\Bormeaspace{[0, +\infty]}$ by~$P^{(Q, \lambda)}_{W | \vect{Z} = \vect{z}}$. 
That is, for all~$\mathcal{A} \in \mathscr{B}\left([0, +\infty]\right)$,
\begin{equation}
\label{EqDefPZ}
P^{(Q,\lambda)}_{W| \vect{Z} = \vect{z}}\left(\set{A}\right) = P^{\left(Q, \lambda \right)}_{\vect{\Theta} | \vect{Z} = \vect{z}}\left(  \mathsf{L}_{\vect{z}}^{-1}\left( \set{A} \right) \right),
\end{equation}
where the term~$\mathsf{L}_{\vect{z}}^{-1}\left( \set{A} \right)$  represents the set
\begin{IEEEeqnarray}{rcl}
\mathsf{L}_{\vect{z}}^{-1}\left( \set{A} \right) &\triangleq&\left\lbrace \vect{\nu} \in \set{M}:  \mathsf{L}_{\vect{z}}(\vect{\nu}) \in \set{A} \right\rbrace.
\end{IEEEeqnarray}
Note that the random variable~$W$ in~\eqref{EqW} induces the probability measure~$P^{(Q,\lambda)}_{W| \vect{Z} = \vect{z}}$ in~$\Bormeaspace{[0, +\infty]}$. 
The objective of this section is to study the properties of the cumulant generating function of the probability measure~$P^{(Q,\lambda)}_{W| \vect{Z} = \vect{z}}$, denoted by~$J_{\vect{z}, Q, \lambda}: \reals \rightarrow \reals \cup \lbrace + \infty \rbrace$,  which satisfies for all~$t \in \reals$,
\begin{IEEEeqnarray}{rCl}
\label{EqJ}
J_{\vect{z}, Q, \lambda} (t) & = & \log\left( \int \exp\left( t w \right) \mathrm{d}P^{(Q, \lambda)}_{W| \vect{Z} = \vect{z}}(w) \right) 
\\
\label{EqJflow}
&= & \log\left( \int \exp\left( t \, \mathsf{L}_{\vect{z}}\left(\vect{\theta}\right)  \right) \mathrm{d}P^{\left(Q, \lambda\right)}_{\vect{\Theta} | \vect{Z} = \vect{z}}(\vect{\theta}) \right),\IEEEeqnarraynumspace
\end{IEEEeqnarray}
where the equality in~\eqref{EqJflow} follows from \cite[Theorem~$1.6.12$]{ash2000probability}.

The following lemma provides an expression for~$J_{\vect{z}, Q, \lambda}$ in terms of the log-partition function~$K_{Q, \vect{z}}$ in~\eqref{EqK}.
\begin{lemma}\label{LemmaCGFPZ}
If~$\lambda \in \set{K}_{Q, \vect{z}}$, with~$\set{K}_{Q, \vect{z}}$ in~\eqref{EqSetKxy},  then, the  function~$J_{\vect{z}, Q, \lambda}$ in~\eqref{EqJ}, verifies for all~$t \in \reals$,
\begin{IEEEeqnarray}{rCl}
\label{EqHormigasPorTodasPartes}
 J_{\vect{z}, Q, \lambda} (t) \squeezeequ  
& = &  K_{ P^{\left(Q, \lambda\right)}_{\vect{\Theta}| \vect{Z} = \vect{z}},\vect{z}}\left( t \right)\\
\label{EqKP}
 & = &  K_{Q, \vect{z}}\left(t- \frac{1}{\lambda}\right) - K_{Q, \vect{z}}\left( -\frac{1}{\lambda} \right) \\
  \label{EqKPBetaFriday}
  & = & \sum_{m =1}^{+\infty}\frac{t^{m}}{m!}K_{Q,\vect{z}}^{(m)}\left(-\frac{1}{\lambda} \right),
\end{IEEEeqnarray}
with the function~$K_{Q, \vect{z}}$ in~\eqref{EqK} and the function $K_{Q,\vect{z}}^{(m)}$ in~\eqref{EqK123}.
\end{lemma}
\begin{IEEEproof}
The proof of \eqref{EqHormigasPorTodasPartes} follows immediately from \eqref{EqK} and \eqref{EqJflow}. 
The proof of \eqref{EqKP} follows from Lemma~\ref{LemmaSophia101}.
Finally, the proof of \eqref{EqKPBetaFriday} follows by observing that a Taylor expansion of the function $K_{Q, \vect{z}}$ in \eqref{EqK} at the point $- \frac{1}{\lambda}$, yields for all $t \in \left\lbrace \nu \in \reals: K_{Q, \vect{z}} (\nu) < +\infty \right\rbrace$,  
\begin{IEEEeqnarray}{rcl}
\label{EqIzzaIsInParis}
K_{Q, \vect{z}}\left( t \right) & = & K_{Q, \vect{z}}\left(  - \frac{1}{\lambda} \right) 
 + \sum_{s=1}^{+\infty}\frac{K^{(s)}_{Q, \vect{z}} \left( - \frac{1}{\lambda} \right)}{s!} \left( t + \frac{1}{\lambda}\right)^{s}.\IEEEeqnarraynumspace
\end{IEEEeqnarray}
Choosing $\alpha \in \left\lbrace \nu \in \reals: K_{Q, \vect{z}} (\nu - \frac{1}{\lambda}) < +\infty \right\rbrace$ such that $t =  \alpha - \frac{1}{\lambda}$ in~\eqref{EqIzzaIsInParis} yields
\begin{IEEEeqnarray}{rcl}
\label{EqHotStuff}
K_{Q, \vect{z}} \left(\alpha - \frac{1}{\lambda}\right)\squeezeequ & = & K_{Q, \vect{z}} \left( - \frac{1}{\lambda}\right)   + \sum_{s=1}^{+\infty}\frac{\alpha^{s}}{s!}K^{(s)}_{Q, \vect{z}} \left( - \frac{1}{\lambda} \right), \squeezeequ\IEEEeqnarraynumspace
\end{IEEEeqnarray}
which implies that for all $t \in \left\lbrace \nu \in \reals: K_{Q, \vect{z}} (\nu - \frac{1}{\lambda}) < +\infty \right\rbrace$,  
\begin{IEEEeqnarray}{rcl}
\label{EqHotStuffColdish}
K_{Q, \vect{z}} \left(t - \frac{1}{\lambda}\right) -  K_{Q, \vect{z}} \left( - \frac{1}{\lambda}\right)   \squeezeequ & = & \sum_{s=1}^{+\infty}\frac{t^{s}}{s!}K^{(s)}_{Q, \vect{z}} \left( - \frac{1}{\lambda} \right).\squeezeequ\IEEEeqnarraynumspace
\end{IEEEeqnarray}

Let $s^{\star} \in \reals \cup \lbrace +\infty \rbrace$ be defined by
\begin{IEEEeqnarray}{rcl}
s^{\star} & \triangleq & \sup \left\lbrace \nu \in \reals: K_{ Q ,\vect{z}}\left(\nu- \frac{1}{\lambda}\right) < \infty \right\rbrace.
\end{IEEEeqnarray}
If $s^{\star} = + \infty$,  then for all $t \in \reals$, $K_{Q, \vect{z}} \left(t - \frac{1}{\lambda}\right) -  K_{Q, \vect{z}} \left( - \frac{1}{\lambda}\right) < + \infty$, and thus,
\begin{IEEEeqnarray}{rCl}
\label{EqQueAburrimiento1}
+\infty > J_{\vect{z}, Q, \lambda} (t) \squeezeequ  
& = &  K_{ P^{\left(Q, \lambda\right)}_{\vect{\Theta}| \vect{Z} = \vect{z}},\vect{z}}\left( t \right)\\
\label{EqQueAburrimiento2}
 & = &  K_{Q, \vect{z}}\left(t- \frac{1}{\lambda}\right) - K_{Q, \vect{z}}\left( -\frac{1}{\lambda} \right) \\
  \label{EqQueAburrimiento3}
  & = & \sum_{m =1}^{+\infty}\frac{t^{m}}{m!}K_{Q,\vect{z}}^{(m)}\left(-\frac{1}{\lambda} \right).
\end{IEEEeqnarray}
Alternatively, if $s^{\star}  < +\infty$, it follows that for all $t > s^{\star}$, $K_{ Q ,\vect{z}}\left(t - \frac{1}{\lambda} \right) = + \infty$. From the fact that the function $K_{Q,\vect{z}}$ is continuous (Lemma~\ref{LemmaContinuousK}) and $ K_{Q,\vect{z}}\left(- \frac{1}{\lambda} \right) < \infty$ (due to the fact that $\lambda \in \set{K}_{Q, \vect{z}}$ in~\eqref{EqSetKxy}), it follows that
\begin{IEEEeqnarray}{rCl}
\label{EqQueAburrimiento4}
+\infty = J_{\vect{z}, Q, \lambda} (t) \squeezeequ 
 & = &  K_{Q, \vect{z}}\left(t- \frac{1}{\lambda}\right) - K_{Q, \vect{z}}\left( -\frac{1}{\lambda} \right) \\
  \label{EqQueAburrimiento5}
  & = & \sum_{m =1}^{+\infty}\frac{t^{m}}{m!}K_{Q,\vect{z}}^{(m)}\left(-\frac{1}{\lambda} \right),
\end{IEEEeqnarray}
which implies that $\sum_{m =1}^{+\infty}\frac{t^{m}}{m!}K_{Q,\vect{z}}^{(m)}\left(-\frac{1}{\lambda} \right) = +\infty$. 
Hence, in this case, the equality in \eqref{EqKPBetaFriday} is of the form $+\infty = +\infty$.
This completes the proof.
 \end{IEEEproof}

Alternative expressions  for~$J_{\vect{z}, Q, \lambda}$ in~\eqref{EqJ} are provided hereunder.
\begin{lemma}\label{LemmaCGFPZjed}
If~$\lambda \in \set{K}_{Q, \vect{z}}$, with~$\set{K}_{Q, \vect{z}}$ in~\eqref{EqSetKxy},  then, the  function~$J_{\vect{z}, Q, \lambda}$ in~\eqref{EqJ}, verifies for all~$t \in \left( 0 , +\infty\right)$,
\begin{IEEEeqnarray}{rCl}
\nonumber
& &  J_{\vect{z}, Q, \lambda} \left(-\frac{1}{t} \right) \squeezeequ \\
 \label{EqKPBeta}
 & = &-\frac{1}{t} \mathsf{R}_{\vect{z}} \left(P^{\left(Q, \frac{1}{\frac{1}{\lambda} + \frac{1}{t}} \right)}_{\vect{\Theta}| \vect{Z} = \vect{z}}\right) -  \KL{P^{\left(Q, \frac{1}{\frac{1}{\lambda} + \frac{1}{t}} \right)}_{\vect{\Theta}| \vect{Z} = \vect{z}} }{ P^{\left(Q, \lambda \right)}_{\vect{\Theta}| \vect{Z} = \vect{z}}} \squeezeequ \IEEEeqnarraynumspace \\
\label{EqKPBetaJed}
 & = &-\frac{1}{t} \mathsf{R}_{\vect{z}} \left( P^{\left(Q, \lambda \right)}_{\vect{\Theta}| \vect{Z} = \vect{z}}\right) +  \KL{ P^{\left(Q, \lambda \right)}_{\vect{\Theta}| \vect{Z} = \vect{z}}}{P^{\left(Q, \frac{1}{\frac{1}{\lambda} + \frac{1}{t}} \right)}_{\vect{\Theta}| \vect{Z} = \vect{z}} } \IEEEeqnarraynumspace\\
 & \leqslant & 0,
\end{IEEEeqnarray}
where the functional~$\mathsf{R}_{\vect{z}}$ is in~\eqref{EqRxy}; the function~$K_{ P^{\left(Q, \lambda\right)}_{\vect{\Theta}| \vect{Z} = \vect{z}},\vect{z}}$ is in~\eqref{EqGenpdfb}; and the probability measures $P^{\left(Q, \lambda \right)}_{\vect{\Theta}| \vect{Z} = \vect{z}}$ and $P^{\left(Q, \frac{1}{\frac{1}{\lambda} + \frac{1}{t}} \right)}_{\vect{\Theta}| \vect{Z} = \vect{z}}$ are respectively in~\eqref{EqGenpdf} and~\eqref{EqGenpdfc}.
\end{lemma}
\begin{IEEEproof}
The proof of~\eqref{EqKPBeta} follows from~\eqref{EqAgadirEarlyMorningA} in Lemma~\ref{LemmaAgadir} by observing that for all $t \in (0, +\infty)$,
\begin{IEEEeqnarray}{rcl}
\nonumber
& & - t K_{P^{\left(Q, \lambda\right)}_{\vect{\Theta}| \vect{Z} = \vect{z}},\vect{z}}\left(- \frac{1}{t} \right)\\
& = & \mathsf{R}_{\vect{z}}\left( P^{\left(P^{\left(Q, \lambda\right)}_{\vect{\Theta}| \vect{Z} = \vect{z}}, t \right)}_{\vect{\Theta}| \vect{Z} = \vect{z}} \right)  + t  \KL{P^{\left(P^{\left(Q, \lambda\right)}_{\vect{\Theta}| \vect{Z} = \vect{z}}, t \right)}_{\vect{\Theta}| \vect{Z} = \vect{z}}}{P^{\left(Q, \lambda\right)}_{\vect{\Theta}| \vect{Z} = \vect{z}}} \middlesqueezeequ \\
\label{EqCrazyThor}
& = & \mathsf{R}_{\vect{z}}\left( P^{\left(Q, \frac{1}{\frac{1}{\lambda} + \frac{1}{t}} \right)}_{\vect{\Theta}| \vect{Z} = \vect{z}} \right)  + t  \KL{P^{\left(Q, \frac{1}{\frac{1}{\lambda} + \frac{1}{t}} \right)}_{\vect{\Theta}| \vect{Z} = \vect{z}}}{P^{\left(Q, \lambda\right)}_{\vect{\Theta}| \vect{Z} = \vect{z}}} \middlesqueezeequ,
\end{IEEEeqnarray}
where the equality in~\eqref{EqCrazyThor} follows from Lemma~\ref{LemmaSophia101A}.
The proof of~\eqref{EqKPBetaJed} follows from~\eqref{EqAgadirEarlyMorningB} in  Lemma~\ref{LemmaAgadir} by observing that  for all $t \in (0, +\infty)$,
\begin{IEEEeqnarray}{rcl}
\nonumber
& & - t K_{P^{\left(Q, \lambda\right)}_{\vect{\Theta}| \vect{Z} = \vect{z}},\vect{z}}\left(- \frac{1}{t} \right)\\
& = & \mathsf{R}_{\vect{z}}\left( P^{\left(Q, \lambda\right)}_{\vect{\Theta}| \vect{Z} = \vect{z}} \right)  - t  \KL{P^{\left(Q, \lambda\right)}_{\vect{\Theta}| \vect{Z} = \vect{z}}}{P^{\left(P^{\left(Q, \lambda\right)}_{\vect{\Theta}| \vect{Z} = \vect{z}}, t \right)}_{\vect{\Theta}| \vect{Z} = \vect{z}}} \middlesqueezeequ \\
\label{EqCrazyThorcito}
& = & \mathsf{R}_{\vect{z}}\left( P^{\left(Q, \lambda\right)}_{\vect{\Theta}| \vect{Z} = \vect{z}} \right)  - t \KL{P^{\left(Q, \lambda\right)}_{\vect{\Theta}| \vect{Z} = \vect{z}}}{P^{\left(Q, \frac{1}{\frac{1}{\lambda} + \frac{1}{t}} \right)}_{\vect{\Theta}| \vect{Z} = \vect{z}}} \middlesqueezeequ,
\end{IEEEeqnarray}
where the equality in~\eqref{EqCrazyThorcito} follows from Lemma~\ref{LemmaSophia101A}, 
which completes the proof.
 \end{IEEEproof}

From Lemma~\ref{LemmaContinuousK} and Lemma~\ref{LemmaCGFPZ}, it follows that the function $J_{\vect{z}, Q, \lambda}$ in~\eqref{EqJ} is increasing and differentiable infinitely many times in the interior of $\left\lbrace t \in \reals : K_{Q, \vect{z}}\left(t - \frac{1}{\lambda}\right) < +\infty \right\rbrace$. Moreover, note that $\left( -\infty, \frac{1}{\lambda} \right] \subset \left\lbrace t \in \reals : K_{Q, \vect{z}}\left(t - \frac{1}{\lambda}\right) < +\infty \right\rbrace$.
Denote by~$J^{(m)}_{ \vect{z},Q,\lambda}: \reals \rightarrow \reals\cup \lbrace + \infty \rbrace$, with~$m \in \ints$, the~$m$-th derivative of the function~$J_{\vect{z}, Q, \lambda}$ in~\eqref{EqJ}. That is, for all~$s \in \reals$,
\begin{equation}
J^{(m)}_{ \vect{z},Q,\lambda} (s) = \frac{\mathrm{d}^{m}}{\mathrm{d}t^{m}} J_{\vect{z}, Q, \lambda} (t) \Bigr|_{t = s}. 
\end{equation}
From Lemma~\ref{LemmaCGFPZ}, it follows that for all~$m \in \ints$, and for all~$\alpha  \in \reals$, the following holds,
\begin{equation}
\label{EqRelationJK}
J^{(m)}_{ \vect{z},Q,\lambda} (\alpha)  = K^{(m)}_{Q, \vect{z}}\left(\alpha - \frac{1}{\lambda}\right),
\end{equation}
where the function~$K^{(m)}_{Q, \vect{z}}$ denotes the $m$-th derivative of the  function~$K_{Q, \vect{z}}$ in~\eqref{EqK}. See for instance, Lemma~\ref{CorDerivatives}.
The equality in~\eqref{EqRelationJK} establishes a relation between the cumulant generating function~$J_{\vect{z}, Q, \lambda}$ and the function~$K_{Q, \vect{z}}$.  This observation becomes an alternative proof to Lemma~\ref{CorDerivatives}.

The following theorem presents the relation between the cumulant generating function~$J_{\vect{z}, Q, \lambda}$  and the functions~$K^{(1)}_{Q, \vect{z}}$ and ~$K^{(2)}_{Q, \vect{z}}$ in~\eqref{EqNiceK1} and~\eqref{EqNiceK2}.
 \begin{theorem}\label{LemmaGaussGen}
For all~$\alpha  \in \reals$, the function~$J_{\vect{z}, Q, \lambda}$ in~\eqref{EqJ}  verifies the following equality 
\begin{equation}
\label{EqJGoldenGaussian}
J_{\vect{z}, Q, \lambda} (\alpha) = \alpha  K^{(1)}_{Q, \vect{z}}\left( - \frac{1}{\lambda}\right)  + \frac{1}{2}\alpha^2 K^{(2)}_{Q, \vect{z}}\left( \xi_{\alpha} \right) 
\end{equation}
with
\begin{IEEEeqnarray}{rcl}
\label{EqEatingEarlyToday}
 \xi_{\alpha} \in  \left( \min \left\lbrace -\frac{1}{\lambda}, \alpha -\frac{1}{\lambda} \right\rbrace, \max \left\lbrace -\frac{1}{\lambda}, \alpha -\frac{1}{\lambda}  \right\rbrace \right),
\end{IEEEeqnarray}
where the functions~$K^{(1)}_{Q, \vect{z}}$ and ~$K^{(2)}_{Q, \vect{z}}$ are defined in~\eqref{EqNiceK1} and~\eqref{EqNiceK2}, respectively.
\end{theorem}
\begin{IEEEproof}
From Lemma~\ref{LemmaContinuousK}, it follows that the function $K_{Q,\vect{z}}$ is differentiable infinitely many times in the interior of $\left\lbrace t \in \reals : K_{Q, \vect{z}}(t) < +\infty \right\rbrace$.
Then, a Taylor expansion of the function $K_{Q, \vect{z}}$ in~\eqref{EqK} at the point $-\frac{1}{\lambda}$ yields for all $t  \in \left\lbrace \nu \in \reals: K_{Q, \vect{z}} (\nu ) < +\infty \right\rbrace$,
\begin{IEEEeqnarray}{rcl}
\nonumber
K_{Q, \vect{z}}\left( t \right) & = & K_{Q, \vect{z}}\left( - \frac{1}{\lambda} \right) \\
\label{EqFindAWayToMakeItMagic}
& & + \sum_{s=1}^{+\infty}\frac{1}{s!} \left( t + \frac{1}{\lambda}\right)^{s}K^{(s)}_{Q, \vect{z}} \left( - \frac{1}{\lambda} \right).
\end{IEEEeqnarray}
Choosing $t = \alpha - \frac{1}{\lambda}$, with $\alpha  \in \left\lbrace \nu \in \reals: K_{Q, \vect{z}} (\nu - \frac{1}{\lambda}) < +\infty \right\rbrace$ in~\eqref{EqFindAWayToMakeItMagic}, it holds from the Taylor-Lagrange theorem \cite[Theorem~$2.5.4$]{trench2003BookIntroRealAnalysis} that
\begin{IEEEeqnarray}{rcl}
\nonumber
K_{Q, \vect{z}}\left(\alpha -\frac{1}{\lambda} \right) & = &  K_{Q, \vect{z}}\left(- \frac{1}{\lambda} \right) + \alpha K^{(1)}_{Q, \vect{z}} \left( - \frac{1}{\lambda} \right)\\
\label{EqTeEnvioPoemas}
& &+  \frac{1}{2} \alpha^{2} K^{(2)}_{Q, \vect{z}} \left( \xi \right), 
\end{IEEEeqnarray}
where $\xi \in \left( \min\lbrace -\frac{1}{\lambda}, \alpha -\frac{1}{\lambda} \rbrace, \max\lbrace -\frac{1}{\lambda}, \alpha -\frac{1}{\lambda}  \rbrace, \right)$. The proof is completed by noticing that from Lemma~\ref{LemmaCGFPZ}, it holds that $J_{\vect{z}, Q, \lambda} (\alpha) = K_{Q, \vect{z}} \left( \alpha - \frac{1}{\lambda}\right) -  K_{Q, \vect{z}} \left( - \frac{1}{\lambda}\right)$.
\end{IEEEproof}

In \eqref{EqJGoldenGaussian}, the term $\xi_{\alpha}$ depends on $\alpha$ via \eqref{EqEatingEarlyToday}. The focus is now on the term $K^{(2)}_{Q, \vect{z}}\left( \xi_{\alpha} \right)$, when $\alpha \in \left\lbrace t \in \reals : J_{\vect{z}, Q, \lambda}(t) < +\infty \right\rbrace$.
\begin{theorem}\label{LemmaGaussDalilo}
The function~$J_{\vect{z}, Q, \lambda}$ in~\eqref{EqJ} verifies the following inequality for all~$\alpha  \in  \left\lbrace t \in \reals : J_{\vect{z}, Q, \lambda}(t) < +\infty \right\rbrace$, 
\begin{equation}
\label{EqDalianGaussian}
J_{\vect{z}, Q, \lambda} (\alpha) \leqslant \alpha  K^{(1)}_{Q, \vect{z}}\left( - \frac{1}{\lambda}\right)  + \frac{1}{2}\alpha^2 \beta_{Q, \vect{z}}^2
\end{equation}
where $\beta_{Q,\vect{z}}$ satisfies
\begin{IEEEeqnarray}{rcl}
\label{EqDalianLeftUpset}
\beta_{Q,\vect{z}} & = & \sup \left\lbrace \sqrt{K^{(2)}_{Q, \vect{z}} \left( \alpha \right)}: \alpha \in \left(-\infty, b - \frac{1}{\lambda}\right) \right\rbrace,\squeezeequ\IEEEeqnarraynumspace
\end{IEEEeqnarray}%
with 
\begin{IEEEeqnarray}{rCl}
\label{EqDamarOnThePhone}
b & \triangleq &\sup \left\lbrace t \in \reals : J_{\vect{z}, Q, \lambda}(t) < +\infty \right\rbrace,
\end{IEEEeqnarray}
and the functions~$K^{(1)}_{Q, \vect{z}}$ and ~$K^{(2)}_{Q, \vect{z}}$ are defined in~\eqref{EqNiceK1} and~\eqref{EqNiceK2}, respectively.
\end{theorem}
\begin{IEEEproof}
From  Lemma~\ref{LemmaCGFPZ}, it holds that
\begin{IEEEeqnarray}{rcl}
\nonumber
\left\lbrace t \in \reals : J_{\vect{z}, Q, \lambda}(t) < +\infty \right\rbrace & = & \left\lbrace t \in \reals : K_{Q, \vect{z}}\left(t- \frac{1}{\lambda}\right) < +\infty \right\rbrace \squeezeequ,  
\end{IEEEeqnarray}
which 
implies that the set $\left\lbrace t \in \reals : J_{\vect{z}, Q, \lambda}(t) < +\infty \right\rbrace$ contains an interval of the form $\left( -\infty, b \right)$ and is contained within an interval of the form $\left( -\infty, b \right]$, with $b$ in \eqref{EqDamarOnThePhone}. This follows from the fact that the function $K_{Q, \vect{z}}$ is continuous and nondecreasing (Lemma~\ref{LemmaContinuousK}). 
Note also that $K_{Q, \vect{z}}(0) = 0$, and thus, $\frac{1}{\lambda} \in \left\lbrace t \in \reals : J_{\vect{z}, Q, \lambda}(t) < +\infty \right\rbrace$, which implies that $b >\frac{1}{\lambda}$. 

From Theorem~\ref{LemmaGaussGen}, it follows that for all $\alpha \in \left(-\infty, b \right)$, 
\begin{IEEEeqnarray}{rcl}
\label{EqJGoldenGaussianRemade}
J_{\vect{z}, Q, \lambda} (\alpha) & \leqslant & \alpha  K^{(1)}_{Q, \vect{z}}\left( - \frac{1}{\lambda}\right)  + \frac{1}{2}\alpha^2 \sigma^2_{\alpha}, 
\end{IEEEeqnarray}
where $\sigma_{\alpha}$ is 
\begin{IEEEeqnarray}{rcl}
\nonumber 
& & \sigma_{\alpha} \\
\nonumber
& \triangleq & \sup \left\lbrace \sqrt{K^{(2)}_{Q, \vect{z}}\left( \xi \right)}: \xi \in \left( \min \left\lbrace -\, \frac{1}{\lambda}, \alpha -\, \frac{1}{\lambda} \right\rbrace, \max \left\lbrace -\, \frac{1}{\lambda}, \alpha -\, \frac{1}{\lambda}  \right\rbrace \right) \right\rbrace. \Tsupersqueezeequ
\end{IEEEeqnarray}
In the asymptotic regime, when $\alpha \to -\infty$, it holds that
\begin{IEEEeqnarray}{rcl}
\nonumber 
\lim_{\alpha \to - \infty} \sigma_{\alpha} & \triangleq & \sup \left\lbrace \sqrt{K^{(2)}_{Q, \vect{z}}\left( \xi \right)}: \xi \in \left(-\infty,  -\frac{1}{\lambda} \right) \right\rbrace,
\end{IEEEeqnarray}
and when $\alpha \to b$, it holds that
\begin{IEEEeqnarray}{rcl}
\nonumber 
\lim_{\alpha \to b} \sigma_{\alpha} & \triangleq & \sup \left\lbrace \sqrt{K^{(2)}_{Q, \vect{z}}\left( \xi \right)}: \xi \in \left(-\frac{1}{\lambda}, b -\frac{1}{\lambda}\right) \right\rbrace.
\end{IEEEeqnarray}
Thus, for all $\alpha \in \left(-\infty, b \right)$, 
\begin{IEEEeqnarray}{rcl}
\nonumber 
\sigma_{\alpha} & \leqslant & \beta_{Q,\vect{z}}, 
\end{IEEEeqnarray}
with $\beta_{Q,\vect{z}}$ defined in \eqref{EqDalianLeftUpset}, which completes the proof.
\end{IEEEproof}
The main implication of Theorem~\ref{LemmaGaussDalilo} is that if $\beta_{Q,\vect{z}}$  in \eqref{EqDalianLeftUpset} is finite,  the random variable~$W$ in~\eqref{EqW} is a sub-Gaussian random variable with sub-Gaussianity parameter $\beta_{Q,\vect{z}}$ \cite[Section~$2.3$]{boucheron2013book}.
This follows by noticing that the function $J_{\vect{z}, Q, \lambda}$ in \eqref{EqJflow} is the cumulant generating function of the random variable $W$. Hence, whenever $\beta_{Q,\vect{z}}$ is finite, when the models are sampled from the ERM-RER optimal measure~$P^{\left(Q, \lambda \right)}_{\vect{\Theta} | \vect{Z} = \vect{z}}$ in~\eqref{EqGenpdf}, the empirical risk with respect to the dataset $\vect{z}$ is a sub-Gaussian random variable with sub-Gaussianity parameter $\beta_{Q,\vect{z}}$.   
The following corollary of Theorem~\ref{LemmaGaussDalilo} highlights this observation.
\begin{corollary}\label{CorTheLastOne}
If $\beta_{Q,\vect{z}}$  in \eqref{EqDalianLeftUpset} is finite, the random variable~$W$ in~\eqref{EqW} is a sub-Gaussian random variable with sub-Gaussianity parameter $\beta_{Q,\vect{z}}$ in \eqref{EqDalianLeftUpset}.
\end{corollary}

\section{Concentration of Probability}\label{SecMonotonicConcentration}

Consider the following set, 
\begin{IEEEeqnarray}{rcl}
\label{EqSetN}
\set{N}_{Q, \vect{z}}(\lambda) & \triangleq &\left\lbrace \vect{\theta} \in \set{M}: \mathsf{L}_{\vect{z}}\left( \vect{\theta}\right)  \leqslant   \mathsf{R}_{\vect{z}}\left(P^{\left(Q, \lambda\right)}_{\vect{\Theta} | \vect{Z} = \vect{z}} \right) \right \rbrace,
\end{IEEEeqnarray}
where the function~$\mathsf{L}_{\vect{z}}$ is defined by~\eqref{EqLxy};  the functional~$\mathsf{R}_{\vect{z}}$ is defined by~\eqref{EqRxy}; and the probability measure~$P^{\left(Q, \lambda\right)}_{\vect{\Theta} | \vect{Z} = \vect{z}}$ is in~\eqref{EqGenpdf}.  
%
%
This section introduces two results. 
First, in Theorem~\ref{CorDecreasingSet}, it is shown that when~$\lambda$ tends to zero, the set~$\set{N}_{Q, \vect{z}}(\lambda)$ forms an indexed family of sets that is monotonic and decreases to the set 
\begin{equation}\label{EqSetNStar}
\set{N}_{Q, \vect{z}}^{\star} \triangleq \set{L}_{\vect{z}} \left( \delta_{Q,\vect{z}}^{\star} \right),  
\end{equation}
where~$\delta_{Q,\vect{z}}^{\star}$ is defined in~\eqref{EqDeltaStar}; and the set~$\set{L}_{\vect{z}}(\delta_{Q,\vect{z}}^{\star})$ is defined in~\eqref{EqSetL}.
Second, in Theorem~\ref{CorDecreasingProbability}, it is shown that  the probability~$P^{\left(Q, \lambda \right)}_{\vect{\Theta}| \vect{Z} = \vect{z}}(\set{N}_{Q, \vect{z}}(\lambda) )$ strictly increases when~$\lambda$ tends to zero. More importantly, in Theorem~\ref{LemmaLimitProb}, it is shown that the limit of the probability~$P^{\left(Q, \lambda \right)}_{\vect{\Theta}| \vect{Z} = \vect{z}}(\set{N}_{Q, \vect{z}}(\lambda) )$, when $\lambda \to 0$, is  equal to one. These observations justify referring to the set~$\set{N}_{Q, \vect{z}}^{\star}$ as the \emph{limit set}. 
These observations are complementary to those stated in Section~\ref{SecFelicidadeA} and Section~\ref{SecFelicidadeB}.
This section ends by showing that the probability measure~$P^{\left(Q, \lambda\right)}_{\vect{\Theta} | \vect{Z} = \vect{z}}$ concentrates on a specific subset~$\set{L}^{\star}_{Q,\vect{z}}$ in~\eqref{EqSetLStar} of the set~$\set{N}_{Q, \vect{z}}^{\star}$. At the light of this observation, the set~$\set{L}^{\star}_{Q,\vect{z}}$ is referred to as the \emph{nonnegligible limit set}.
Finally, it is shown that when the~$\sigma$-finite measure~$Q$ in~\eqref{EqERMRER} is coherent, the sets~$\set{N}_{Q, \vect{z}}^{\star}$ and 
$\set{L}^{\star}_{Q,\vect{z}}$ are identical. 

\subsection{The Limit Set}

The set~$\set{N}_{Q, \vect{z}}(\lambda)$ in~\eqref{EqSetN}, with~$\lambda \in \set{K}_{Q, \vect{z}}$ and~$\set{K}_{Q, \vect{z}}$ in~\eqref{EqSetKxy}, contains all the models that induce an empirical risk that is smaller than or equal to~$\mathsf{R}_{\vect{z}}\left(P^{\left(Q, \lambda\right)}_{\vect{\Theta} | \vect{Z} = \vect{z}} \right)$, i.e., the ERM-RER-optimal expected empirical risk in~\eqref{EqRK}. This observation unveils the existence of a  relation between the set~$\set{N}_{Q, \vect{z}}^{\star}$ in~\eqref{EqSetNStar} and the set~$\set{T}\left( \vect{z} \right)$ in~\eqref{EqHatTheta}, as shown by the following lemma. 

\begin{lemma}\label{LemNstarT}
The set~$\set{N}_{Q, \vect{z}}^{\star}$ in~\eqref{EqSetNStar} satisfies
\begin{equation}\label{EqCubanInclusion}
\set{T}\left( \vect{z} \right)  \subseteq \set{N}_{Q, \vect{z}}^{\star},
\end{equation}
where the set~$\set{T}\left( \vect{z} \right)$ is in~\eqref{EqHatTheta}. 
Moreover, 
\begin{equation}\label{EqCubanFreedom}
\set{T}\left( \vect{z} \right)  = \set{N}_{Q, \vect{z}}^{\star},
\end{equation}
 if and only if $(a)$ the ERM problem in~\eqref{EqOriginalOP} possesses a solution; and $(b)$ the reference measure~$Q$ in~\eqref{EqERMRER} is coherent.
\end{lemma}
\begin{IEEEproof}
If the set $\set{T}\left( \vect{z} \right)$ in~\eqref{EqHatTheta} is empty, the inclusion in~\eqref{EqCubanInclusion} is trivially true. 
Assume that $\abs{\set{T}\left( \vect{z} \right)} > 0$. Hence, the proof of the inclusion in~\eqref{EqCubanInclusion} follows from observing that for all~$\vect{\theta} \in \set{T}\left( \vect{z} \right)$, it holds that~$\mathsf{L}_{\vect{z}}\left( \vect{\theta}\right) = \rho^{\star}  \leqslant \delta_{Q,\vect{z}}^{\star}$, with~$\delta_{Q,\vect{z}}^{\star}$  in~\eqref{EqDeltaStar} and~$\rho^{\star}$ in~\eqref{EqRhoStar}. Hence,~$\vect{\theta} \in \set{N}_{Q, \vect{z}}^{\star}$. This completes the proof of the inclusion in~\eqref{EqCubanInclusion}.

The proof of the equality in~\eqref{EqCubanFreedom} is presented in two parts. In the first part, it is proved that if~\eqref{EqCubanFreedom} holds, then the ERM problem in~\eqref{EqOriginalOP} possesses a solution and the measure~$Q$ is coherent. 
The second part proves the converse.
The proof of the first part is as follows. Under the assumption that~$\set{T}\left( \vect{z} \right)  = \set{N}_{Q, \vect{z}}^{\star}$ holds, it follows that~$\delta_{Q,\vect{z}}^{\star} = \rho^{\star}$, with $\rho^{\star}$ in~\eqref{EqRhoStar}, which implies that the ERM problem in~\eqref{EqOriginalOP} possesses a solution. Moreover, for all~$\delta \in \left(\rho^{\star}, +\infty \right)$, it holds that~$Q \left( \set{L}_{\vect{z}} \left( \delta \right)  \right) > 0$, which verifies that the measure~$Q$ is coherent and completes the proof of the first part. 
The proof of the second part is as follows. Under the assumption that the ERM problem in~\eqref{EqOriginalOP} possesses a solution and the measure~$Q$ is coherent, it follows that~$\delta_{Q,\vect{z}}^{\star} = \rho^{\star}$. Hence,~$\set{T}\left( \vect{z} \right)  = \set{N}_{Q, \vect{z}}^{\star}$, which completes the proof of the second part.
\end{IEEEproof}

The following theorem highlights that the set~$\set{N}_{Q, \vect{z}}\left( \lambda \right)$ is decreasing with $\lambda$.

\begin{theorem}\label{CorDecreasingSet}
For all~$(\lambda_1,\lambda_2)  \in \set{K}_{Q, \vect{z}}\times\set{K}_{Q, \vect{z}}$, with~$\set{K}_{Q, \vect{z}}$ in~\eqref{EqSetKxy} and~$\lambda_{1} > \lambda_{2}$, the sets~$\set{N}_{Q, \vect{z}}\left( \lambda_1 \right)$ and~$\set{N}_{Q, \vect{z}}\left( \lambda_2 \right)$ in~\eqref{EqSetN} satisfy
\begin{IEEEeqnarray}{rcl}
\label{EqDecreasingsets}
\set{M} \supseteq \set{N}_{Q, \vect{z}}(\lambda_1) \supseteq \set{N}_{Q, \vect{z}}(\lambda_2)  \supseteq \set{N}_{Q, \vect{z}}^{\star} ,
\end{IEEEeqnarray}
with~$\set{N}_{Q, \vect{z}}^{\star}$ being the set defined in~\eqref{EqSetNStar}.
Moreover, if the empirical risk function~$\mathsf{L}_{\vect{z}}$ in~\eqref{EqLxy} is continuous on~$\set{M}$ and separable with respect to the measure~$Q$ in~\eqref{EqERMRER}, then,  
\begin{IEEEeqnarray}{rcl}
\label{EqNavir2315}
\set{M} \supset \set{N}_{Q, \vect{z}}(\lambda_1) \supset \set{N}_{Q, \vect{z}}(\lambda_2) \supset   \set{N}_{Q, \vect{z}}^{\star}.
\end{IEEEeqnarray}
\end{theorem}
\begin{IEEEproof}
The proof is presented in Appendix~\ref{AppProofCorDecreasingSet}.
\end{IEEEproof}
An interesting observation is that for all~$\lambda \in \set{K}_{Q, \vect{z}}$, with~$\set{K}_{Q, \vect{z}}$ in~\eqref{EqSetKxy}, only a subset of~$\set{N}_{Q, \vect{z}}\left( \lambda \right)$ might exhibit nonzero probability with respect to the measure~$P^{\left(Q, \lambda \right)}_{\vect{\Theta}| \vect{Z} = \vect{z}}$ in~\eqref{EqGenpdf}. 
Consider for instance that the measure~$Q$ in~\eqref{EqERMRER} is noncoherent (Definition~\ref{DefCoherentMess}). That is,~$\delta_{Q,\vect{z}}^{\star} > \rho^{\star}$, with~$\delta_{Q,\vect{z}}^{\star}$  in~\eqref{EqDeltaStar} and $\rho^{\star}$ in~\eqref{EqRhoStar}.  Thus, for all~$\gamma \in \left(\rho^{\star},\delta_{Q,\vect{z}}^{\star} \right)$, it holds that~$Q\left( \set{L}_{\vect{z}}\left( \gamma \right) \right)= 0$, with the set~$\set{L}_{\vect{z}}(\cdot)$  in~\eqref{EqSetL}. 
From Lemma~\ref{LemmaMutualAC}, this implies that for all~$\gamma \in \left(\rho^{\star},\delta_{Q,\vect{z}}^{\star} \right)$, the measure $P^{\left(Q, \lambda\right)}_{\vect{\Theta} | \vect{Z} = \vect{z}}$ in~\eqref{EqGenpdf} satisfies $P^{\left(Q, \lambda\right)}_{\vect{\Theta} | \vect{Z} = \vect{z}}\left( \set{L}_{\vect{z}}\left( \gamma \right) \right)= 0$, while verifying that ~$\set{L}_{\vect{z}}\left( \gamma \right) \subseteq  \set{N}_{Q, \vect{z}}\left( \lambda \right)$.
 These observations lead to the analysis of the asymptotic concentration of probability in the following section. 

\subsection{The Nonnegligible Limit Set}
 
The first step in the analysis of the asymptotic concentration of the probability measure~$P^{\left(Q, \lambda\right)}_{\vect{\Theta} | \vect{Z} = \vect{z}}$ in~\eqref{EqGenpdf}  
is to show that the probability~$P^{\left(Q, \lambda\right)}_{\vect{\Theta} | \vect{Z} = \vect{z}} \left(\set{N}_{Q, \vect{z}}(\lambda) \right)$ increases when~$\lambda$ tends to zero, as shown by the following theorem.  

\begin{theorem}\label{CorDecreasingProbability}
For all~$(\lambda_1,\lambda_2)  \in \set{K}_{Q, \vect{z}}\times\set{K}_{Q, \vect{z}}$, with~$\set{K}_{Q, \vect{z}}$ in~\eqref{EqSetKxy} and~$\lambda_{1} > \lambda_{2}$, assume that the measures~$P^{\left(Q, \lambda_1 \right)}_{\vect{\Theta}| \vect{Z} = \vect{z}}$ and~$P^{\left(Q, \lambda_2 \right)}_{\vect{\Theta}| \vect{Z} = \vect{z}}$ satisfy~\eqref{EqGenpdf} with~$\lambda = \lambda_1$ and~$\lambda = \lambda_2$, respectively. Then, the set~$\set{N}_{Q, \vect{z}}\left( \lambda_2 \right)$ in~\eqref{EqSetN} satisfies  
\begin{IEEEeqnarray}{rcl}
\label{EqIncreasingProbabilityA9}
0 < P^{\left(Q, \lambda_1 \right)}_{\vect{\Theta}| \vect{Z} = \vect{z}}(\set{N}_{Q, \vect{z}}(\lambda_{2}) )  \leqslant P^{\left(Q, \lambda_2 \right)}_{\vect{\Theta}| \vect{Z} = \vect{z}}(\set{N}_{Q, \vect{z}}(\lambda_{2}) ),
\end{IEEEeqnarray}
where strict inequality holds if and only if the function~$\mathsf{L}_{\vect{z}}$ is separable with respect to the~$\sigma$-finite measure~$Q$.
\end{theorem}
\begin{IEEEproof}
The proof  is presented in Appendix~\ref{AppProofCorDecreasingProbability}.
\end{IEEEproof}

The following lemma highlights a case in which a stronger  concentration of probability is observed.

\begin{lemma}\label{CorStrongStrinctDecreasingProbability}
Let the function~$\mathsf{L}_{\vect{z}}$ in~\eqref{EqLxy} be separable  with respect to the~$\sigma$-finite measure~$Q$ in~\eqref{EqERMRER}. Let also~$\left( \lambda_{1} , \lambda_{2}\right)  \in \set{K}_{Q, \vect{z}}\times\set{K}_{Q, \vect{z}}$, with~$\set{K}_{Q, \vect{z}}$ in~\eqref{EqSetKxy}, be  two positive reals such that~$\lambda_1 > \lambda_2$ and
\begin{equation}
\label{EqKeyConditionForStrong}
Q \bigg( \set{N}_{Q, \vect{z}}\left( \lambda_1 \right)  \cap \left(\set{N}_{Q, \vect{z}}\left( \lambda_2 \right)\right)^\sfc \bigg) = 0,
\end{equation}
with the complement with respect to the set of models $\set{M}$.
Then, two measures~$P^{\left(Q, \lambda_1 \right)}_{\vect{\Theta}| \vect{Z} = \vect{z}}$ and~$P^{\left(Q, \lambda_2 \right)}_{\vect{\Theta}| \vect{Z} = \vect{z}}$ that respectively satisfy~\eqref{EqGenpdf} with~$\lambda = \lambda_1$ and~$\lambda = \lambda_2$ verify that
\begin{IEEEeqnarray}{rcl}
\label{EqIncreasingProbability245}
P^{\left(Q, \lambda_1 \right)}_{\vect{\Theta}| \vect{Z} = \vect{z}}(\set{N}_{Q, \vect{z}}(\lambda_1) )  < P^{\left(Q, \lambda_2 \right)}_{\vect{\Theta}| \vect{Z} = \vect{z}}(\set{N}_{Q, \vect{z}}(\lambda_{2}) ),
\end{IEEEeqnarray}
where, the set~$\set{N}_{Q, \vect{z}}\left( \cdot \right)$ is defined in~\eqref{EqSetN}.  
\end{lemma}
\begin{IEEEproof}
The proof is presented in Appendix~\ref{AppProofCorStrongStrinctDecreasingProbability}.
\end{IEEEproof}

The following example shows the relevance of Lemma~\ref{CorStrongStrinctDecreasingProbability} in the case in which the empirical risk function~$\mathsf{L}_{\vect{z}}$ in~\eqref{EqLxy} is a simple function and separable  with respect to the~$\sigma$-finite measure~$Q$ in~\eqref{EqERMRER}.

\begin{example}\label{ExSimpleFunction}
Consider Example~\ref{ExampleDecreasingVariance}.  Note that, for all~$\lambda> 0$,
\begin{equation}\label{EqRiceMaker}
0< \mathsf{R}_{\vect{z}}\left(P^{\left(Q, \lambda\right)}_{\vect{\Theta} | \vect{Z} = \vect{z}} \right) < 1,
\end{equation}
where~$\mathsf{R}_{\vect{z}}\left(P^{\left(Q, \lambda\right)}_{\vect{\Theta} | \vect{Z} = \vect{z}} \right)$ is the ERM-RER-optimal expected empirical risk in~\eqref{EqRK}.  
The equality in~\eqref{EqRiceMaker} implies that given two reals~$\lambda_{1}$ and~$\lambda_{2}$ such that~$\lambda_{1} > \lambda_{2} > 0$, it holds that, 
\begin{IEEEeqnarray}{rcl}
\nonumber
& & \set{N}_{Q, \vect{z}}\left( \lambda_1 \right)  \cap \left(\set{N}_{Q, \vect{z}}\left( \lambda_2 \right)\right)^\sfc \\
& = &\left\lbrace \vect{\nu} \in \set{M}: \mathsf{R}_{\vect{z}}\left(P^{\left(Q, \lambda_{2}\right)}_{\vect{\Theta} | \vect{Z} = \vect{z}} \right)
 <  \mathsf{L}_{\vect{z}}\left( \vect{\nu}\right)  \leqslant   \mathsf{R}_{\vect{z}}\left(P^{\left(Q, \lambda_{1}\right)}_{\vect{\Theta} | \vect{Z} = \vect{z}} \right)\right \rbrace \squeezeequ\IEEEeqnarraynumspace\\
 & =&  \emptyset,
\end{IEEEeqnarray}
and moreover,~$\set{N}_{Q, \vect{z}}(\lambda_1) = \set{N}_{Q, \vect{z}}(\lambda_{2})$. Finally,  from Lemma~\ref{CorStrongStrinctDecreasingProbability}, 
\begin{IEEEeqnarray}{rcl}
P^{\left(Q, \lambda_1 \right)}_{\vect{\Theta}| \vect{Z} = \vect{z}}(\set{N}_{Q, \vect{z}}(\lambda_1) )  < P^{\left(Q, \lambda_2 \right)}_{\vect{\Theta}| \vect{Z} = \vect{z}}(\set{N}_{Q, \vect{z}}(\lambda_{2}) ).
\end{IEEEeqnarray}
\end{example}

The main result of this section is presented by the following theorem. 

\begin{theorem}\label{LemmaLimitProb}
The probability measure~$P^{\left(Q, \lambda\right)}_{\vect{\Theta} | \vect{Z} = \vect{z}}$ in~\eqref{EqGenpdf} satisfies
 \begin{equation}
 \lim_{\lambda \rightarrow 0^{+}} P^{\left(Q, \lambda\right)}_{\vect{\Theta} | \vect{Z} = \vect{z}} \left(\set{N}_{Q, \vect{z}}\left( \lambda \right) \right) =  
 1,
 \end{equation}
where, the set~$\set{N}_{Q, \vect{z}}\left( \lambda \right)$ is defined in~\eqref{EqSetN}.
\end{theorem}
\begin{IEEEproof}
The proof follows immediately from Lemma~\ref{TheoSetLStarConcentration} and by noticing that for all $\lambda \in \set{K}_{Q, \vect{z}}$, with~$\set{K}_{Q, \vect{z}}$ in~\eqref{EqSetKxy}, the sets~$\set{L}^{\star}_{Q,\vect{z}}$ in~\eqref{EqSetLStar}  and $\set{N}_{Q, \vect{z}}\left( \lambda \right)$ in~\eqref{EqSetN} satisfies $\set{L}^{\star}_{Q,\vect{z}} \subseteq \set{N}_{Q, \vect{z}}\left( \lambda \right)$.
 \end{IEEEproof}
  
Note that Theorem~\ref{LemmaLimitProb} and Lemma~\ref{TheoSetLStarConcentration} lead to the following conclusion
\begin{equation}
 \lim_{\lambda \rightarrow 0^{+}} P^{\left(Q, \lambda\right)}_{\vect{\Theta} | \vect{Z} = \vect{z}} \left(\set{N}_{Q, \vect{z}}\left( \lambda \right) \setminus  \set{L}^{\star}_{Q,\vect{z}} \right) =  0,
 \end{equation}  
which follows from the fact that~$\set{L}^{\star}_{Q,\vect{z}} \subset \set{N}_{Q, \vect{z}}\left( \lambda \right)$, with~$\set{L}^{\star}_{Q,\vect{z}}$ in~\eqref{EqSetLStar}. This justifies referring to the set~$\set{L}^{\star}_{Q,\vect{z}}$ as the nonnegligible limit set.

\section{$(\delta, \epsilon)$-Optimality}\label{SecDeltaEpsilonOptimality}

This section introduces a PAC guarantee of optimality for the models that are sampled from the probability measure~$P^{\left(Q, \lambda\right)}_{\vect{\Theta} | \vect{Z} = \vect{z}}$ in~\eqref{EqGenpdf} with respect to the ERM problem in~\eqref{EqOriginalOP}. Such guarantee is defined as follows. 
\begin{definition}[$(\delta, \epsilon)$-Optimality] \label{DefDeltaEpsilonOptimal}
Given a pair of positive reals~$(\delta, \epsilon)$, with~$\epsilon < 1$, the probability measure~$P^{\left(Q, \lambda\right)}_{\vect{\Theta} | \vect{Z} = \vect{z}}$ in~\eqref{EqGenpdf} is said to be~$(\delta, \epsilon)$-optimal, if the set~$\set{L}_{\vect{z}} \left( \delta \right)$ in~\eqref{EqSetL} satisfies
\begin{equation}
\label{EqOptimalConcentration}
P^{\left(Q, \lambda\right)}_{\vect{\Theta} | \vect{Z} = \vect{z}} \left(\set{L}_{\vect{z} } \left( \delta \right) \right) > 1 - \epsilon.
\end{equation}
\end{definition}
If the probability measure~$P^{\left(Q, \lambda\right)}_{\vect{\Theta} | \vect{Z} = \vect{z}}$ in~\eqref{EqGenpdf}  is~$(\delta, \epsilon)$-optimal, then it assigns a probability that is always greater than~$1-\epsilon$ to a set that contains models that induce an empirical risk that is smaller than~$\delta$.  From this perspective, particular interest is given to the smallest~$\delta$ and~$\epsilon$ for which~$P^{\left(Q, \lambda\right)}_{\vect{\Theta} | \vect{Z} = \vect{z}}$ is~$\left(\delta , \epsilon \right)$-optimal.

The main result of this section is presented by the following theorem.

\begin{theorem}\label{TheoGibbsDeltaEps} 
For all~$(\delta, \epsilon) \in (\delta_{Q,\vect{z}}^{\star},+\infty) \times (0,1)$, with~$\delta_{Q,\vect{z}}^{\star}$ in~\eqref{EqDeltaStar}, there exists a real~$\lambda \in \set{K}_{Q, \vect{z}}$, with~$\set{K}_{Q, \vect{z}}$ in~\eqref{EqSetKxy}, such that the probability measure~$P^{\left(Q, \lambda\right)}_{\vect{\Theta} | \vect{Z} = \vect{z}}$ is~$\left(\delta , \epsilon \right)$-optimal.
\end{theorem}
\begin{IEEEproof}
Let~$\delta$ be a real in~$\left( \delta_{Q,\vect{z}}^{\star}, +\infty \right)$, with~$\delta_{Q,\vect{z}}^{\star}$ in~\eqref{EqDeltaStar}.  Let also~$\gamma \in \set{K}_{Q, \vect{z}}$ satisfy the following equality:
\begin{equation}
\label{EqProofExistence1}
K^{(1)}_{Q, \vect{z}}\left( - \frac{1}{\gamma} \right)  \leqslant \delta.
\end{equation}
Note that from Lemma~\ref{LemmaContinuousK}, it follows that the function~$K^{(1)}_{Q, \vect{z}}$ is continuous. Moreover, from Theorem~\ref{CorAssymptoticMean}, it follows that such a~$\gamma$ in~\eqref{EqProofExistence1} always exists.  
From~\eqref{EqSetL} and~\eqref{EqSetN}, it holds that
 \begin{equation}
 \label{EqProofExistence03}
\set{N}_{Q, \vect{z}}(\gamma) \subseteq \set{L}_{\vect{z}} \left( \delta \right),
 \end{equation}
 and thus,
 \begin{equation}
  \label{EqProofExistence403}
P^{\left(Q, \gamma\right)}_{\vect{\Theta} | \vect{Z} = \vect{z}}\left(\set{L}_{\vect{z}} \left( \delta \right)\right)  \geqslant  P^{\left(Q, \gamma \right)}_{\vect{\Theta} | \vect{Z} = \vect{z}}\left(\set{N}_{Q, \vect{z}}(\gamma) \right).
\end{equation}
Let~$\lambda$ be a  positive real such that~$\lambda \leqslant \gamma$ and 
\begin{equation}
\label{EqProofExistence13r12}
P^{\left(Q, \lambda \right)}_{\vect{\Theta} | \vect{Z} = \vect{z}}(\set{N}_{Q, \vect{z}}(\lambda) ) > 1 - \epsilon.
\end{equation}
The existence of such a positive real~$\lambda$ follows from Theorem~\ref{LemmaLimitProb}. 
Hence, from~\eqref{EqProofExistence13r12}, it holds that,
\begin{IEEEeqnarray}{rcl}
\label{EqProofExistence634}
1- \epsilon & < &P^{\left(Q, \lambda \right)}_{\vect{\Theta} | \vect{Z} = \vect{z}}(\set{N}_{Q, \vect{z}}(\lambda))  \\
\label{EqProofExistence64}
& \leqslant & P^{\left(Q, \lambda \right)}_{\vect{\Theta} | \vect{Z} = \vect{z}}\left(\set{L}_{\vect{z}}\left( \delta \right) \right), 
\end{IEEEeqnarray}
where 
the inequality in~\eqref{EqProofExistence64} follows from the fact that~$\set{N}_{Q, \vect{z}}(\lambda) \subseteq \set{N}_{Q, \vect{z}}(\gamma) \subseteq \set{L}_{\vect{z}} \left( \delta \right)$.
Finally,  the inequality in~\eqref{EqProofExistence64} implies that the probability measure $P^{\left(Q, \lambda\right)}_{\vect{\Theta} | \vect{Z} = \vect{z}}$ is~$\left(\delta , \epsilon \right)$-optimal (Definition~\ref{DefDeltaEpsilonOptimal}). 
This completes the proof.
\end{IEEEproof}        

A stronger optimality claim can be stated when the reference measure is coherent.
\begin{theorem}\label{TheoGibbsCoherent}
For all~$(\delta, \epsilon) \in (\rho^{\star},+\infty) \times (0,1)$, with $\rho^{\star}$ in~\eqref{EqRhoStar}, there always exists a~$\lambda \in \set{K}_{Q, \vect{z}}$, with~$\set{K}_{Q, \vect{z}}$ in~\eqref{EqSetKxy}, such that the probability measure~$P^{\left(Q, \lambda\right)}_{\vect{\Theta} | \vect{Z} = \vect{z}}$ is~$\left(\delta , \epsilon \right)$-optimal if and only if the reference measure~$Q$ is coherent.
\end{theorem}
\begin{IEEEproof}
The proof is divided into two parts. The first part shows that if for all~$(\delta, \epsilon) \in (\rho^{\star},+\infty) \times (0,1)$, there always exists a~$\lambda \in \set{K}_{Q, \vect{z}}$, with~$\set{K}_{Q, \vect{z}}$ in~\eqref{EqSetKxy}, such that the probability measure~$P^{\left(Q, \lambda\right)}_{\vect{\Theta} | \vect{Z} = \vect{z}}$  in~\eqref{EqGenpdf} is~$(\delta, \epsilon)$-optimal, then, the measure~$Q$ is coherent. 
The second part deals with the converse.

The first part is as follows. Let~$\gamma \in \set{K}_{Q, \vect{z}}$ be such that
\begin{IEEEeqnarray}{rcl}
\label{EqFiestaPromenade634}
P^{\left(Q, \gamma \right)}_{\vect{\Theta} | \vect{Z} = \vect{z}}\left(\set{L}_{\vect{z}}\left( \delta \right) \right)  & > &1- \epsilon,
\end{IEEEeqnarray}
then, for all measurable subsets $\set{A}$ of $\set{L}_{\vect{z}}\left( \delta \right)$, it holds that
\begin{IEEEeqnarray}{rCl}
1- \epsilon \middlesqueezeequ& < & P^{\left(Q, \gamma \right)}_{\vect{\Theta} | \vect{Z} = \vect{z}}\left(\set{L}_{\vect{z}}\left( \delta \right) \right) \\
\nonumber
& = & \int_{\set{A}} \frac{\mathrm{d} P^{\left(Q, \gamma \right)}_{\vect{\Theta} | \vect{Z} = \vect{z}}}{\mathrm{d} Q}\left( \vect{\nu} \right) \mathrm{d} Q \left( \vect{\nu} \right) +  \int_{\set{L}_{\vect{z}}\left( \delta \right)\setminus \set{A}} \frac{\mathrm{d} P^{\left(Q, \gamma \right)}_{\vect{\Theta} | \vect{Z} = \vect{z}}}{\mathrm{d} Q}\left( \vect{\nu} \right) \mathrm{d} Q \left( \vect{\nu} \right), \middlesqueezeequ
\end{IEEEeqnarray}
which, together with Lemma~\ref{CorPositive}, implies that there exists at least one measurable subset $\set{A}$ for which $Q \left( \set{A} \right)  > 0$, and thus, 
\begin{IEEEeqnarray}{rcl}
\label{EqSister2}
Q \left(\set{L}_{\vect{z}}\left( \delta \right) \right) & > & Q \left(\set{A} \right) >  0,
\end{IEEEeqnarray}
which implies that the measure~$Q$ is coherent. This completes the first part of the proof.

The second part of the proof is as follows. Under the assumption that the measure~$Q$ is coherent, it follows that~$\delta_{Q,\vect{z}}^{\star} = \rho^{\star}$. Then, from Theorem~\ref{TheoGibbsDeltaEps}, it follows that for all~$(\delta, \epsilon) \in (\delta_{Q,\vect{z}}^{\star},+\infty) \times (0,1)$,  there always exists a~$\lambda \in \set{K}_{Q, \vect{z}}$, with~$\set{K}_{Q, \vect{z}}$ in~\eqref{EqSetKxy}, such that the probability measure~$P^{\left(Q, \lambda\right)}_{\vect{\Theta} | \vect{Z} = \vect{z}}$  is~$(\delta, \epsilon)$-optimal. This completes the second part of the proof. 
\end{IEEEproof}

\section{Sensitivity  and Generalization}\label{SecSubSensitivity}
This section introduces the notion of sensitivity and establishes its connections with the notion of generalization error of the Gibbs algorithm, cf.~\cite{aminian2021exact}.

\subsection{Sensitivity}
The sensitivity of the expected empirical risk~$\mathsf{R}_{\vect{z}}$ in~\eqref{EqRxy}   to  deviations from the probability measure~$P^{\left(Q, \lambda\right)}_{\vect{\Theta}| \vect{Z} = \vect{z}}$ in~\eqref{EqGenpdf} towards an alternative probability measure~$P\in \triangle\Bormeaspace{\set{M}}$ is introduced as a novel metric to evaluate the generalization capabilities of the ERM-RER-optimal measure $P^{\left(Q, \lambda\right)}_{\vect{\Theta}| \vect{Z} = \vect{z}}$. 
Deviations from the probability measure~$P^{\left(Q, \lambda\right)}_{\vect{\Theta}| \vect{Z} = \vect{z}}$ towards an alternative probability measure~$P$ would allow comparing the ERM-RER-optimal measure with alternative measures (or algorithms). For instance, if new datasets become available, a new ERM-RER problem can be formulated using a larger dataset obtained by aggregating the old and the new datasets, cf. \cite{Perlaza-ISIT2023b} and~\cite{InriaRR9474}. Intuitively, the ERM-RER-optimal measure obtained after the aggregation of datasets might exhibit better generalization capabilities, see for instance~\cite{Perlaza-ISIT2023b}. This analysis is the motivation of the sensitivity, which is defined as follows.
\begin{definition}[Sensitivity]
Given the~$\sigma$-finite measure~$Q$ and the positive real~$\lambda > 0$ in~\eqref{EqERMRER}, let~$\mathsf{S}_{Q, \lambda}: \left( \set{X} \times \set{Y} \right)^n \times \triangle_{Q}\Bormeaspace{\set{M}}\rightarrow \left( - \infty, +\infty \right]$ be a functional such that  
\begin{IEEEeqnarray}{l}
\label{EqDefSensitivity}
\mathsf{S}_{Q, \lambda}\left( \vect{z}, P \right)  =\squeezeequ
\left\lbrace
\begin{array}{cl}
\mathsf{R}_{\vect{z}}\left( P \right)  - \mathsf{R}_{\vect{z}}\left( P^{\left(Q, \lambda\right)}_{\vect{\Theta}| \vect{Z} = \vect{z}} \right) \squeezeequ & \text{ if } \lambda \in \set{K}_{Q,\vect{z}} \squeezeequ \\
+\infty & \text{ otherwise,}
\end{array}
\right. \IEEEeqnarraynumspace
\end{IEEEeqnarray}
where the functional~$\mathsf{R}_{\vect{z}}$ is defined in~\eqref{EqRxy} and the probability measure~$P^{\left(Q, \lambda\right)}_{\vect{\Theta}| \vect{Z} = \vect{z}}$ is in~\eqref{EqGenpdf}.
The sensitivity of the expected empirical risk~$\mathsf{R}_{\vect{z}}$ due to a deviation from~$P^{\left(Q, \lambda\right)}_{\vect{\Theta}| \vect{Z} = \vect{z}}$ to~$P$ is~$\mathsf{S}_{Q, \lambda}\left( \vect{z}, P\right)$.
\end{definition}
Recently, the following exact expression for the sensitivity $\mathsf{S}_{Q, \lambda}\left( \vect{z}, P\right)$ in~\eqref{EqDefSensitivity} was introduced in \cite{Perlaza-ISIT2023b}.
\begin{theorem}[Theorem~1 in \cite{Perlaza-ISIT2023b}]\label{TheoremSensitivityEqual}
The sensitivity $\mathsf{S}_{Q, \lambda}\left( \vect{z}, P \right)$ in~\eqref{EqDefSensitivity} satisfies
\begin{IEEEeqnarray}{rcl}
\nonumber
& & \mathsf{S}_{Q, \lambda}\left( \vect{z}, P \right) \Tsupersqueezeequ \\
\label{EqSensitivityEqual}
& = &  \lambda\Big( \KL{P^{\left(Q, \lambda\right)}_{\vect{\Theta}| \vect{Z} = \vect{z}}}{Q} + \KL{P}{P^{\left(Q, \lambda\right)}_{\vect{\Theta}| \vect{Z} = \vect{z}}} - \KL{P}{Q} \Big),\Tsupersqueezeequ \IEEEeqnarraynumspace
\end{IEEEeqnarray}
where  the probability measure~$P^{\left(Q, \lambda\right)}_{\vect{\Theta}| \vect{Z} = \vect{z}}$ is in~\eqref{EqGenpdf}.
\end{theorem}

 
The following theorem introduces an upper bound on the absolute value of the sensitivity~$\mathsf{S}_{Q, \lambda}\left( \vect{z}, P \right)$ in~\eqref{EqDefSensitivity}, which requires the calculation of only one of the relative entropies in Theorem~\ref{TheoremSensitivityEqual}.

\begin{theorem}\label{TheoremSensitivityA}
For all~$P \in \triangle_{Q}\Bormeaspace{\set{M}}$, the sensitivity $\mathsf{S}_{Q, \lambda}\left( \vect{z}, P \right)$ in~\eqref{EqDefSensitivity} satisfies 
\begin{IEEEeqnarray}{rcl}
\label{EqMinaIronsAt5am}
\abs{\mathsf{S}_{Q, \lambda}\left( \vect{z}, P \right)  }& \leqslant &  \sqrt{2  \beta_{Q,\vect{z}}^2 \KL{P }{P^{\left(Q, \lambda\right)}_{\vect{\Theta}| \vect{Z} = \vect{z}} }},
\end{IEEEeqnarray}
where  the constant~$\beta_{Q,\vect{z}}$ is defined in~\eqref{EqDalianLeftUpset}.
\end{theorem}
\begin{IEEEproof}
The proof is presented in Appendix~\ref{AppProofTheoremSensitivityA}.
\end{IEEEproof}

Note that equality holds in~\eqref{EqMinaIronsAt5am} in the trivial case in which the empirical risk function is not separable with respect to $Q$ (Definition~\ref{DefSeparableLxy}). In such case, for all $P \in \triangle_{Q}\Bormeaspace{\set{M}}$, it holds that $\mathsf{S}_{Q, \lambda}\left( \vect{z}, P \right) = 0$ and $\beta_{Q,\vect{z}} = 0$.

Theorem~\ref{TheoremSensitivityA} establishes an upper and a lower bound on the increase and decrease of the expected empirical risk  that can be obtained by deviating from the optimal solution of the ERM-RER problem in~\eqref{EqERMRER}. More specifically, note that for all  probability measures~$P \in \triangle_{Q}\Bormeaspace{\set{M}}$, it holds that,
\begin{IEEEeqnarray}{rcl}
\mathsf{R}_{\vect{z}}\left( P \right)   \squeezeequ & \geqslant & 
\mathsf{R}_{\vect{z}}\left( P^{\left(Q, \lambda\right)}_{\vect{\Theta} | \vect{Z} = \vect{z}} \right) - \sqrt{2 \beta_{Q,\vect{z}}^2  \KL{P}{P^{\left(Q, \lambda\right)}_{\vect{\Theta} | \vect{Z} = \vect{z}} } } \mbox{ and } \squeezeequ\IEEEeqnarraynumspace\\
\mathsf{R}_{\vect{z}}\left( P \right)   & \leqslant &     \mathsf{R}_{\vect{z}}\left( P^{\left(Q, \lambda\right)}_{\vect{\Theta} | \vect{Z} = \vect{z}} \right) + \sqrt{2\beta_{Q,\vect{z}}^2  \KL{P}{P^{\left(Q, \lambda\right)}_{\vect{\Theta} | \vect{Z} = \vect{z}} } } .
\end{IEEEeqnarray}

\subsection{Generalization Error}%

This section unveils the interesting connection between the notion of sensitivity and the notion of generalization error of the Gibbs algorithm, cf. \cite{aminian2021exact}. The Generalization error is defined under the assumption that datasets are sampled from a probability measure 
\begin{IEEEeqnarray}{rCl}
\label{EqProdPXY}
P_{\vect{Z}} \in \triangle\left( \left( \set{X} \times \set{Y} \right)^{n}, \mathscr{F} \right),
\end{IEEEeqnarray}
where $\mathscr{F}$ denotes a given $\sigma$-field on the set $\left(\mathcal{X}\times\mathcal{Y}\right)^n$.
For such a probability measure $P_{\vect{Z}}$ in~\eqref{EqProdPXY}, let the set~$\set{K}_{Q,P_{\vect{Z}}} \subset \reals$ be 
\begin{equation}
\label{EqSetKQ}
\set{K}_{Q,P_{\vect{Z}}} =  \bigcap_{\vect{z} \in \supp P_{\vect{Z}} } \set{K}_{Q, \vect{z}},
\end{equation}
where the $\sigma$-finite measure~$Q$ is in~\eqref{EqERMRER}. 
The set~$\set{K}_{Q,P_{\vect{Z}}}$ in~\eqref{EqSetKQ} can be empty for some choices of the~$\sigma$-finite measure~$Q$. Nonetheless, from Lemma~\ref{LemmaSubSetK}, it follows that if~$Q$ is a probability measure, then,
\begin{equation}
\label{EqLlamarAPerrineEstaSemana}
\set{K}_{Q,P_{\vect{Z}}}= \left( 0, +\infty \right).
\end{equation}
Under the assumption that datasets are sampled from $P_{\vect{Z}}$ in~\eqref{EqProdPXY}, the generalization error of the Gibbs algorithm with parameters $Q$ and $\lambda$, is defined as the expectation with respect to the product measure $P^{\left(Q, \lambda\right)}_{\vect{\Theta}| \vect{Z}} \cdot P_{\vect{Z}}$, with $P^{\left(Q, \lambda\right)}_{\vect{\Theta}| \vect{Z}}$ in \eqref{EqGenpdf},  of difference between: $(a)$ the \emph{population risk} due to a model $\vect{\theta} \in \set{M}$, 
\begin{IEEEeqnarray}{rcl}
\int \mathsf{L}_{ \vect{z} } \left(\vect{\theta} \right) \mathrm{d}P_{\vect{Z}}\left( \vect{z} \right)
\end{IEEEeqnarray}
with the function $\mathsf{L}_{ \vect{z} }$ defined in \eqref{EqLxy}; and $(b)$ The empirical risk induced by the model $\vect{\theta}$ with respect to a training dataset $
\vect{\vect{z}}$, that is,~$\mathsf{L}_{\vect{z}}\left( \vect{\theta} \right)$.
More specifically, the generalization error of the Gibbs algorithm with parameters $Q$ and $\lambda$ is
\begin{IEEEeqnarray}{rcl}
\nonumber
& & \int \int \left( \int \mathsf{L}_{ \vect{z} } \left(\vect{\theta} \right) \mathrm{d}P_{\vect{Z}}\left( \vect{z} \right) - \mathsf{L}_{\vect{\nu}}\left( \vect{\theta} \right) \right) \mathrm{d}P^{\left(Q, \lambda\right)}_{\vect{\Theta}| \vect{Z} = \vect{\nu}}\left(\vect{\theta} \right) \mathrm{d}P_{\vect{Z}}(\vect{\nu})  \\
\nonumber
& = & \int \left( \int \mathsf{L}_{ \vect{z} } \left(\vect{\theta} \right) \mathrm{d}P_{\vect{Z}}\left( \vect{z} \right) \right) \mathrm{d}P^{\left(Q, \lambda\right)}_{\vect{\Theta}}\left(\vect{\theta} \right)  \\
& &- \int \mathsf{L}_{\vect{\nu}}\left( \vect{\theta} \right) \mathrm{d}P^{\left(Q, \lambda\right)}_{\vect{\Theta}| \vect{Z} = \vect{\nu}}\left(\vect{\theta} \right) \mathrm{d}P_{\vect{Z}}(\vect{\nu})  \\
\nonumber
& = & \int \left( \int \mathsf{L}_{ \vect{z} } \left(\vect{\theta} \right) \mathrm{d}P^{\left(Q, \lambda\right)}_{\vect{\Theta}}\left(\vect{\theta} \right) \right) \mathrm{d}P_{\vect{Z}}\left( \vect{z} \right)  \\
& &- \int \mathsf{L}_{\vect{\nu}}\left( \vect{\theta} \right) \mathrm{d}P^{\left(Q, \lambda\right)}_{\vect{\Theta}| \vect{Z} = \vect{\nu}}\left(\vect{\theta} \right) \mathrm{d}P_{\vect{Z}}(\vect{\nu})  \\
\label{EqPainfullMondayOfTheEndOfMay2023}
& =  & \int \left( \mathsf{R}_{\vect{\nu}}\left( P^{\left(Q, \lambda\right)}_{\vect{\Theta}} \right)  - \mathsf{R}_{\vect{\nu}}\left( P^{\left(Q, \lambda\right)}_{\vect{\Theta}| \vect{Z} = \vect{\nu}} \right) \right) \mathrm{d}P_{\vect{Z}}(\vect{\nu}) ,\supersqueezeequ\IEEEeqnarraynumspace
\end{IEEEeqnarray}
where the probability measure $P^{\left(Q, \lambda\right)}_{\vect{\Theta}}$ satisfies for all  sets~$\set{A} \in \BorSigma{\set{M}}$,
\begin{equation}
\label{EqBarPThetaX}
P^{\left(Q, \lambda\right)}_{\vect{\Theta}}\left( \set{A} \right) = \int P^{\left(Q, \lambda\right)}_{\vect{\Theta} | \vect{Z} = \vect{\nu}  } \left( \set{A} \right)  \mathrm{d}P_{\vect{Z}}\left( \vect{\nu} \right),
\end{equation} 
and the functional $\mathsf{R}_{\vect{\nu}}$ is defined in \eqref{EqRxy}.

The following theorem establishes a connection between sensitivity and generalization error in the particular case in which $Q$ in \eqref{EqERMRER} is a probability measure.
\begin{theorem}\label{Theo30May20238h49}
Under the assumption that datasets are sampled from $P_{\vect{Z}}$ in~\eqref{EqProdPXY}, the generalization error of the Gibbs algorithm with parameters $Q$ (a probability measure) and $\lambda > 0$, is 
\begin{IEEEeqnarray}{rcl}
 \int \hspace{-0.8ex} \mathsf{S}_{Q, \lambda}\left( \vect{\nu},P^{\left(Q, \lambda\right)}_{\vect{\Theta}}\right) \mathrm{d}P_{\vect{Z}}(\vect{\nu}),
\end{IEEEeqnarray}
where the functional~$\mathsf{S}_{Q, \lambda}$ is in~\eqref{EqDefSensitivity}; and the probability measure $P^{\left(Q, \lambda\right)}_{\vect{\Theta}}$ is in \eqref{EqBarPThetaX}.
\end{theorem}
\begin{IEEEproof}
The proof uses the fact that under the assumption that $Q$ is a probability measure,  for all $\vect{\nu} \in \supp P_{\vect{Z}}$, it follows from  Lemma~\ref{LemmaSubSetK} that $\set{K}_{Q,\vect{\nu}}= \left( 0, +\infty \right)$. This implies that for all  $\vect{z} \in \supp P_{\vect{Z}}$ and for all $\lambda > 0$, the ERM-RER problem in \eqref{EqERMRER}, always possesses as solution the measure $P^{\left(Q, \lambda\right)}_{\vect{\Theta}| \vect{Z} = \vect{z}}$ in \eqref{EqGenpdf}. 
Thus, the measure $P^{\left(Q, \lambda\right)}_{\vect{\Theta}}$ in \eqref{EqBarPThetaX} is well defined.
Moreover, $\mathsf{S}_{Q, \lambda}\left( \vect{z}, P^{\left(Q, \lambda\right)}_{\vect{\Theta}} \right) =  \mathsf{R}_{\vect{z}}\left( P^{\left(Q, \lambda\right)}_{\vect{\Theta}} \right)  - \mathsf{R}_{\vect{z}}\left( P^{\left(Q, \lambda\right)}_{\vect{\Theta}| \vect{Z} = \vect{z}} \right)$ and the integral in~\eqref{EqPainfullMondayOfTheEndOfMay2023} is also well defined, which completes the proof.
\end{IEEEproof}

Theorem~\ref{Theo30May20238h49} provides an interesting viewpoint of the generalization error. 
For instance, the probability measure $P^{\left(Q, \lambda\right)}_{\vect{\Theta}}$  in~\eqref{EqProdPXY} can be understood as the barycenter of a subset of $\triangle\Bormeaspace{\set{M}}$ containing the solutions to ERM-RER problems of the form in~\eqref{EqERMRER},  with $\vect{z} \in \supp P_{\vect{Z}}$ in~\eqref{EqProdPXY}. 
Hence, the generalization error of the Gibbs algorithm is the expectation (with respect to $P_{\vect{Z}}$) of the sensitivity of the expected empirical risks $ \mathsf{R}_{\vect{z}}$ in~\eqref{EqRxy} to variations from the ERM-RER-optimal measure $P^{\left(Q, \lambda\right)}_{\vect{\Theta}| \vect{Z} = \vect{z}}$  towards the barycenter, i.e., the measure $P^{\left(Q, \lambda\right)}_{\vect{\Theta}}$.  

The following definition extends the notion of generalization error to Gibbs algorithms obtained by assuming that the reference measure~$Q$ in \eqref{EqERMRER} is a $\sigma$-finite measure. This definition also exploits the relation between  the notions of sensitivity and generalization error introduced by Theorem~\ref{Theo30May20238h49}.

\begin{definition}[Generalization Error of the Gibbs Algorithm]\label{DefGenErr}
Given a $\sigma$-finite measure $Q \in \triangle\Bormeaspace{\set{M}}$ and a real $\lambda>0$, let the functional $\mathsf{G}_{Q, \lambda}:  \triangle\left( \left( \set{X} \times \set{Y} \right)^{n}, \mathscr{F} \right) \to   \left( - \infty, +\infty \right]$ be such that 
\begin{IEEEeqnarray}{rcl}
\label{EqInsop}
\mathsf{G}_{Q, \lambda}\left(P_{\vect{Z}} \right) \Tsupersqueezeequ &=& 
\left\lbrace 
\begin{array}{cl}
\displaystyle \hspace{-1.4ex} \int  \mathsf{S}_{Q, \lambda}\left( \vect{\nu},P^{\left(Q, \lambda\right)}_{\vect{\Theta}}\right) \mathrm{d}P_{\vect{Z}}(\vect{\nu}) \Tsupersqueezeequ& \hspace{-2ex} \text{ if } \lambda \in \set{K}_{Q,P_{\vect{Z}}}\Tsupersqueezeequ\\
+\infty & \hspace{-2ex}  \text{ otherwise,}\Tsupersqueezeequ
\end{array}
\right.\IEEEeqnarraynumspace
\end{IEEEeqnarray} 
where the functional~$\mathsf{S}_{Q, \lambda}$ is in~\eqref{EqDefSensitivity}; the set $\set{K}_{Q,P_{\vect{Z}}}$ is  in~\eqref{EqSetKQ};  and the probability measure $P^{\left(Q, \lambda\right)}_{\vect{\Theta}}$ is in \eqref{EqBarPThetaX}. The generalization error induced by the Gibbs algorithm with parameters $Q$ and $\lambda$ under the assumption that datasets are sampled from the probability measure $P_{\vect{Z}}$, is $\mathsf{G}_{Q, \lambda}\left(P_{\vect{Z}} \right)$.
\end{definition}

The main difficulty for extending the notion of generalization error to Gibbs algorithms obtained under the assumption that the reference measure is not a probability measure, but a $\sigma$-finite measure, is that the integrals in~\eqref{EqPainfullMondayOfTheEndOfMay2023} and~\eqref{EqBarPThetaX} might not be well defined. 
This is essentially due to the fact that, while the ERM-RER problem in \eqref{EqERMRER} always possesses a solution when $Q$ is a probability measure, the existence of a solution when $Q$ is not a probability measure is subject to the condition that for all $\vect{z} \in \supp P_{\vect{Z}}$, $\lambda \in \set{K}_{Q, \vect{z}}$, with $\set{K}_{Q, \vect{z}}$ in~\eqref{EqSetKxy}. This leads to the condition that $\lambda \in \set{K}_{Q,P_{\vect{Z}}}$, with the set $\set{K}_{Q,P_{\vect{Z}}}$ in~\eqref{EqSetKQ}. When such a condition is not met, the definition of sensitivity is void.

The following theorem provides a closed-form expression for the generalization error of the Gibbs algorithm in the general case in which the reference measure $Q$ in \eqref{EqERMRER} is a $\sigma$-finite measure.
\begin{theorem}\label{TheoAminini}
If $\lambda \in \set{K}_{Q,P_{\vect{Z}}}$, with $\set{K}_{Q,P_{\vect{Z}}}$ in~\eqref{EqSetKQ},
the generalization error $\mathsf{G}_{Q, \lambda}\left(P_{\vect{Z}} \right)$ in~\eqref{EqInsop} satisfies 
\begin{IEEEeqnarray}{rCl}
\nonumber
\mathsf{G}_{Q, \lambda}\left(P_{\vect{Z}} \right) & =  & \lambda\Bigg( \int  \KL{P^{\left(Q, \lambda\right)}_{\vect{\Theta}| \vect{Z} = \vect{\nu}} }{P^{\left(Q, \lambda\right)}_{\vect{\Theta}}}  \mathrm{d} P_{\vect{Z}}(\vect{\nu})  \Dsupersqueezeequ\\
\label{EqPulgosos}
&   & + \int \KL{P^{\left(Q, \lambda\right)}_{\vect{\Theta}}}{P^{\left(Q, \lambda\right)}_{\vect{\Theta}| \vect{Z} = \vect{\nu}} }  \mathrm{d} P_{\vect{Z}}(\vect{\nu})  \Bigg)  ,\Dsupersqueezeequ \IEEEeqnarraynumspace
\end{IEEEeqnarray}
where for all $\vect{z} \in \supp P_{\vect{Z}}$, the probability measure~$P^{\left(Q,\lambda\right)}_{\vect{\Theta}| \vect{Z} = \vect{z}}$ is in~\eqref{EqGenpdf}; and  
the probability measure~$P^{\left(Q, \lambda \right)}_{\vect{\Theta}}$ is defined in~\eqref{EqBarPThetaX}.
\end{theorem}
\begin{IEEEproof}
The proof is presented in Appendix~\ref{AppProofOfTheoAminini}.
\end{IEEEproof}

The terms $\int  \KL{P^{\left(Q, \lambda\right)}_{\vect{\Theta}| \vect{Z} = \vect{\nu}}}{P^{\left(Q, \lambda\right)}_{\vect{\Theta}}} \mathrm{d} P_{\vect{Z}}(\vect{\nu})$ and $\int  \KL{P^{\left(Q, \lambda\right)}_{\vect{\Theta}}}{P^{\left(Q, \lambda\right)}_{\vect{\Theta}| \vect{Z} = \vect{\nu}} } \mathrm{d} P_{\vect{Z}}(\vect{\nu})$ in  the right-hand side of~\eqref{EqPulgosos} are respectively the mutual and the lautum information~\cite{palomar2008lautum} induced by a joint probability measure $P_{\vect{\Theta}, \vect{Z}}$ whose marginals are $P_{\vect{Z}}$ in~\eqref{EqProdPXY} and $P^{\left(Q, \lambda\right)}_{\vect{\Theta}}$ in~\eqref{EqBarPThetaX}.
When the reference measure $Q$ in~\eqref{EqERMRER} is a probability measure, Theorem~\ref{TheoAminini} reduces to \cite[Theorem~$1$]{aminian2021exact}. 
Interestingly, independently of whether the reference measure $Q$ in~\eqref{EqERMRER} is a probability measure, or whether the $n$ data points in the datasets are independent and identically distributed, the generalization error $\mathsf{G}_{Q, \lambda}\left(P_{\vect{Z}} \right)$ in~\eqref{EqInsop} is always a factor of the sum of the mutual and lautum information induced by the joint probability measure $P_{\vect{\Theta}, \vect{Z}}$ mentioned above. 

Theorem~\ref{TheoAminini} also provides an alternative interpretation of the generalization error  $\mathsf{G}_{Q, \lambda}\left(P_{\vect{Z}} \right)$ in~\eqref{EqInsop}. Note that by writing one of the factors in the right-hand side of ~\eqref{EqPulgosos} as 
\begin{IEEEeqnarray}{rcl}
\nonumber
\int \left( \KL{P^{\left(Q, \lambda\right)}_{\vect{\Theta}| \vect{Z} = \vect{\nu}} }{P^{\left(Q, \lambda\right)}_{\vect{\Theta}}} + \KL{P^{\left(Q, \lambda\right)}_{\vect{\Theta}}}{P^{\left(Q, \lambda\right)}_{\vect{\Theta}| \vect{Z} = \vect{\nu}} } \right)  \mathrm{d} P_{\vect{Z}}(\vect{\nu}),  \Dsupersqueezeequ 
\end{IEEEeqnarray}
it becomes clear that $\mathsf{G}_{Q, \lambda}\left(P_{\vect{Z}} \right)$ is the expectation with respect to $P_{\vect{Z}}$ of the symmetrized Kullback-Leibler divergence, also known as Jeffrey's divergence~\cite{jeffreys1946invariant}, of the probability measures $P^{\left(Q, \lambda\right)}_{\vect{\Theta}| \vect{Z} = \vect{z}}$ and $P^{\left(Q, \lambda\right)}_{\vect{\Theta}}$. That is, the solution to the ERM-RER problem in~\eqref{EqERMRER} and the barycenter induced by $P_{\vect{Z}}$.

The following theorem provides an upper-bound on the generalization error of the Gibbs algorithm only in terms of the lautum information induced by such a joint probability measure~$P_{\vect{\Theta}, \vect{Z}}$. 
\begin{theorem}\label{LemmaBoomBoomDown}
The generalization error $\mathsf{G}_{Q, \lambda}\left(P_{\vect{Z}} \right)$ in~\eqref{EqInsop} satisfies for all  $\lambda \in \set{K}_{Q,P_{\vect{Z}}}$,
\begin{IEEEeqnarray}{rCl}
0 \leqslant & \mathsf{G}_{Q, \lambda}\left(P_{\vect{Z}} \right)\Dsupersqueezeequ
\label{EqFreeOfLace}
&  \leqslant  \sqrt{2 \sigma^{2}_{Q}  \int  \KL{P^{\left(Q, \lambda\right)}_{\vect{\Theta}}}{P^{\left(Q, \lambda\right)}_{\vect{\Theta}| \vect{Z} = \vect{\nu}} } \mathrm{d} P_{\vect{Z}}(\vect{\nu}) },\Dsupersqueezeequ \IEEEeqnarraynumspace
\end{IEEEeqnarray}
where for all $\vect{z} \in \supp P_{\vect{Z}}$, the probability measure~$P^{\left(Q,\lambda\right)}_{\vect{\Theta}| \vect{Z} = \vect{z}}$ is in~\eqref{EqGenpdf}; 
the probability measure~$P^{\left(Q, \lambda \right)}_{\vect{\Theta}}$ is defined in~\eqref{EqBarPThetaX}; and
\begin{IEEEeqnarray}{rcl}
\label{EqListeningItalians}
\sigma_{Q} & = & \sup \left\lbrace  \beta_{Q,\vect{z}}: \vect{z} \in \left( \set{X} \times \set{Y} \right)^{n}  \right\rbrace,
\end{IEEEeqnarray}
with $\beta_{Q,\vect{z}}$ in~\eqref{EqDalianLeftUpset}.
\end{theorem}
\begin{IEEEproof}
The proof of the inequality $\mathsf{G}_{Q, \lambda}\left(P_{\vect{Z}} \right) \geqslant 0$ follows from observing that for all $\vect{\nu} \in \left( \set{X} \times \set{Y} \right)^n$, the terms  
$\KL{P^{\left(Q, \lambda\right)}_{\vect{\Theta}| \vect{Z} = \vect{\nu}} }{P^{\left(Q, \lambda\right)}_{\vect{\Theta}}} $ and $\KL{P^{\left(Q, \lambda\right)}_{\vect{\Theta}}}{P^{\left(Q, \lambda\right)}_{\vect{\Theta}| \vect{Z} = \vect{\nu}} }$ in \eqref{EqPulgosos} are nonnegative (Theorem~\ref{LemmaREA}).
The proof of the remaining inequality follows from~\eqref{EqInsop} and the following inequalities:
\begin{IEEEeqnarray}{rcl}
\mathsf{G}_{Q, \lambda}\left(P_{\vect{Z}} \right) 
\label{EqDesaparecida1}
 &=& \abs{ \int  \mathsf{S}_{Q, \lambda}\left( \vect{\nu},P^{\left(Q, \lambda\right)}_{\vect{\Theta}}\right) \mathrm{d}P_{\vect{Z}}(\vect{\nu})}\\
 \label{EqDesaparecida2}
&\leqslant &\int \abs{\mathsf{S}_{Q, \lambda}\left( \vect{\nu}, P^{\left(Q, \lambda\right)}_{\vect{\Theta}} \right)  } \mathrm{d} P_{\vect{Z}}(\vect{\nu}) \\
\label{EqEhihAou}
& \leqslant &   \int\sqrt{2 \beta_{Q,\vect{\nu}} \KL{P^{\left(Q, \lambda\right)}_{\vect{\Theta}}}{P^{\left(Q, \lambda\right)}_{\vect{\Theta}| \vect{Z} = \vect{\nu}} } }\mathrm{d} P_{\vect{Z}}(\vect{\nu}) ,\squeezeequ \IEEEeqnarraynumspace\\
\label{EqLaFuerza}
& \leqslant &  \int \hspace{-1ex} \sqrt{2 \sigma_{Q}^2\KL{P^{\left(Q, \lambda\right)}_{\vect{\Theta}}}{P^{\left(Q, \lambda\right)}_{\vect{\Theta}| \vect{Z} = \vect{\nu}} } }\mathrm{d} P_{\vect{Z}}(\vect{\nu})\\
\label{EqLaFuerza351}
& \leqslant &  \sqrt{2 \sigma_{Q}^2 \int  \KL{P^{\left(Q, \lambda\right)}_{\vect{\Theta}}}{P^{\left(Q, \lambda\right)}_{\vect{\Theta}| \vect{Z} = \vect{\nu}} } \mathrm{d} P_{\vect{Z}}(\vect{\nu}) },
\end{IEEEeqnarray}
where
 the equality in~\eqref{EqDesaparecida1} follows from~\eqref{EqInsop};
the inequality in~\eqref{EqDesaparecida2} follows from \cite[Theorem~$1.5.9$(c)]{ash2000probability};
the inequality in~\eqref{EqEhihAou} follows from Theorem~\ref{TheoremSensitivityA};
 the inequality in~\eqref{EqLaFuerza} follows from~\eqref{EqListeningItalians}; and 
the inequality in~\eqref{EqLaFuerza351} follows from Jensen's inequality \cite[Section~$6.3.5$]{ash2000probability}.
This completes the proof.
\end{IEEEproof}
In a nutshell, the generalization error $\mathsf{G}_{Q, \lambda}\left(P_{\vect{Z}} \right)$ in~\eqref{EqInsop} is upper bounded up to a constant factor  by the square root of the lautum information induced by the joint probability measure $P_{\vect{\Theta}, \vect{Z}}$ mentioned above. 
Theorem~\ref{LemmaBoomBoomDown} is reminiscent of \cite[Theorem~$1$]{xu2017information}, which provides a similar upper-bound on $\mathsf{G}_{Q, \lambda}\left(P_{\vect{Z}} \right)$ using the mutual information  instead of the lautum information induced by the joint probability measure $P_{\vect{\Theta}, \vect{Z}}$.
The interest in Theorem~\ref{LemmaBoomBoomDown} for the specific case of the Gibbs algorithm,  lies on the fact that it holds under milder conditions than those in~\cite[Theorem~$1$]{xu2017information}. For instance, no additional conditions on the loss function $\ell$ in~\eqref{EqEll} concerning sub-Gaussianity are assumed. Moreover, the probability measure $P_{\vect{Z}}$ from which datasets are sampled is not necessarily a product measure. 
\section{Conclusions and Final Remarks}\label{SecDiscussion}
 
The classical ERM-RER problem in~\eqref{EqERMRER} has been studied under the assumption that the reference measure~$Q$ is a~$\sigma$-finite measure, instead of a probability measure, which leads to a more general problem that includes the ERM problem with (discrete or differential) entropy regularization and the information-risk minimization problem. While in the case in which the reference measure is a probability measure the solution to the ERM-RER problem always exists, in this general case, the existence of a solution is subject to a condition that depends on the loss function, the reference measure, the regularization factor, and the training dataset. When a solution exists, it has been proved that it is unique. Additionally, if it exists, such a solution and the reference measure are mutually absolutely continuous. Interestingly, the empirical risk observed when models are sampled from the ERM-RER-optimal probability measure is a sub-Gaussian random variable that exhibits a PAC guarantee for the ERM problem. That is, for some positive~$\delta$ and~$\epsilon$, it is shown that there always exist some parameters for the ERM-RER problem such that the set of models that induce an empirical risk smaller than~$\delta$ exhibits a probability that is not smaller that~$1-\epsilon$. Interestingly, none of these results relies on statistical assumptions on the datasets.  

The sensitivity of the expected empirical risk to deviations from the ERM-RER-optimal measure to alternative measures is introduced as a new performance metric to evaluate the generalization capabilities of the Gibbs algorithm. In particular, an upper bound on the absolute value of the sensitivity, which depends on the training dataset, is presented. This bound is formed by a constant factor and the square root of the relative entropy of the alternative measure (the deviation) with respect to the ERM-RER solution. 
Finally, it is shown that the expectation of the sensitivity (with respect to the datasets) to deviations towards a particular measure is equivalent to the generalization error of the Gibbs algorithm. Equipped with this observation, the generalization error is shown to be in the most general case, up to a constant factor, the sum of the mutual and lautum information between the models and the datasets, which was a result known exclusively for the case in which the reference is a probability measure, cf.~\cite{aminian2021exact}.  
From this perspective, it is argued that the study of the generalization capabilities of the Gibbs algorithm based on generalization error is a significantly narrow view. This is essentially because it is looking at an expectation of the sensitivity to deviations to a particular measure, i.e., the barycenter of the set of ERM-RER solutions induced by a prior on the datasets. 
A broader view is offered by the study of the sensitivity to deviations towards other measures, i.e., ERM-RER-optimal measures obtained with different training data sets. This approach has lead already to a few initial results in \cite{Perlaza-ISIT2023b} that highlight the connections to sensitivity, training error, and test error. Nonetheless, the study of the sensitivity in the aim of describing the generalization capabilities of learning algorithms remains by now as an open problem. 

\begin{appendices}

\section{Proof of Theorem~\ref{TheoremCasaBlanca}}\label{AppProofOfTheoremCasaBlanca} 

Consider the function $f:[0, +\infty) \to \reals$ such that 
\begin{IEEEeqnarray}{rcl}
\label{EqVousEtesVikings}
f(x) = \left\lbrace
\begin{array}{rcl}
x \log(x) & \text{ if }& x > 0\\
0 & \text{ if }& x = 0,
\end{array}\right.
\end{IEEEeqnarray}
 and note that it is strictly convex.  
From the assumption that for all $i \in \lbrace 1,2 \rbrace$,  $P_{i}$ and $Q_{i}$ are both measures on the same measurable space $\left( \Omega , \mathscr{F} \right)$, with $P_i$ absolutely continuous with respect to $Q_i$, let $g: \Omega \to [0, \infty)$ be the function  
\begin{IEEEeqnarray}{rcl}
g(x)  & = & \frac{\mathrm{d} \left( \lambda P_1 + (1-\lambda) P_2 \right)}{\mathrm{d}\left( \lambda Q_1 + (1-\lambda) Q_2\right)}\left( x \right),
\end{IEEEeqnarray}
where $\frac{\mathrm{d} \left( \lambda P_1 + (1-\lambda) P_2 \right)}{\mathrm{d}\left( \lambda Q_1 + (1-\lambda) Q_2\right)}$ is the Radon-Nikodym derivative of the measure $ \lambda P_1 + (1-\lambda) P_2$ with respect to $ \lambda Q_1 + (1-\lambda) Q_2$.
Using this notation, for all $\lambda \in (0,1)$, 
\begin{IEEEeqnarray}{rcl}
\nonumber
& & \KL{\lambda P_1 + (1-\lambda) P_2}{\lambda Q_1 + (1-\lambda) Q_2} \\
& & - \lambda \KL{ P_1 }{Q_1} + (1 - \lambda) \KL{ P_2}{Q_2}\\
\nonumber
& = & \int \log\left( g\left( x \right) \right) \mathrm{d} \left( \lambda P_1 + (1-\lambda) P_2 \right) (x) \Dsupersqueezeequ \\
& &\nonumber 
- \lambda \int \log \left(\frac{\mathrm{d} P_1}{\mathrm{d} Q_1}(x)\right) \mathrm{d} P_1(x)
- (1 - \lambda) \int \log \left(\frac{\mathrm{d} P_2}{\mathrm{d} Q_2}(x)\right) \mathrm{d} P_2(x) \Dsupersqueezeequ \\
\nonumber
& = &\lambda \int \log\left( g\left( x \right) \right) \mathrm{d} P_1 (x)  + (1-\lambda)   \int \log\left( g\left( x \right) \right) \mathrm{d} P_2  (x)\\
\nonumber 
& &
- \lambda \int \log \left(\frac{\mathrm{d} P_1}{\mathrm{d} Q_1}(x)\right) \mathrm{d} P_1(x)
- (1 - \lambda) \int \log \left(\frac{\mathrm{d} P_2}{\mathrm{d} Q_2}(x)\right) \mathrm{d} P_2(x)\Dsupersqueezeequ \\
\nonumber
& = &\lambda \int \log\left( \left( \frac{\mathrm{d} P_1}{\mathrm{d} Q_1}(x) \right)^{-1} g\left( x \right) \right) \mathrm{d} P_1 (x)\squeezeequ \\
\nonumber
&& + (1-\lambda)   \int \log\left(  \left(  \frac{\mathrm{d} P_2}{\mathrm{d} Q_2}(x) \right)^{-1} g\left( x \right) \right) \mathrm{d} P_2  (x)\squeezeequ \\
\nonumber
& = &\lambda \int \frac{\mathrm{d} P_1}{\mathrm{d} Q_1}(x) \log\left( \left( \frac{\mathrm{d} P_1}{\mathrm{d} Q_1}(x) \right)^{-1} g\left( x \right) \right) \mathrm{d} Q_1 (x)\squeezeequ \\
\nonumber
&& + (1-\lambda)   \int \frac{\mathrm{d} P_2}{\mathrm{d} Q_2}(x) \log\left(  \left(  \frac{\mathrm{d} P_2}{\mathrm{d} Q_2}(x) \right)^{-1} g\left( x \right) \right) \mathrm{d} Q_2  (x)\squeezeequ \\
\nonumber
& = &\lambda \int \frac{g(x) \frac{\mathrm{d} P_1}{\mathrm{d} Q_1}(x)}{g(x)}\log\left( \left( \frac{\mathrm{d} P_1}{\mathrm{d} Q_1}(x) \right)^{-1} g\left( x \right) \right) \mathrm{d} Q_1 (x)\squeezeequ \\
\nonumber
&& + (1-\lambda)   \int \frac{g(x) \frac{\mathrm{d} P_2}{\mathrm{d} Q_2}(x)}{g(x)}\log\left(  \left(  \frac{\mathrm{d} P_2}{\mathrm{d} Q_2}(x) \right)^{-1} g\left( x \right) \right) \mathrm{d} Q_2  (x)\squeezeequ \\
\nonumber
& = & - \lambda \int g(x) f \left( \frac{\mathrm{d} P_1}{\mathrm{d} Q_1}(x) \left( g(x) \right)^{-1} \right) \mathrm{d} Q_1 (x)\squeezeequ \\
\label{EqValhala101}
&& - (1-\lambda)   \int g(x) f \left( \frac{\mathrm{d} P_2}{\mathrm{d} Q_2}(x) \left( g(x) \right)^{-1} \right)  \mathrm{d} Q_2  (x),\squeezeequ 
\end{IEEEeqnarray}
where the function $f$ is defined in~\eqref{EqVousEtesVikings}.
Let $\beta_1$ and $\beta_2$ be the following constants:
\begin{IEEEeqnarray}{rcl}
\label{EqWishpering}
\beta_1  \triangleq  \int g(\nu) \mathrm{d}Q_{1}(\nu) \mbox{ and } \beta_2  \triangleq  \int g(\nu) \mathrm{d}Q_{2}(\nu).
\end{IEEEeqnarray}
From~\eqref{EqValhala101} and~\eqref{EqWishpering}, it follows that for all $\lambda \in (0,1)$, 
\begin{IEEEeqnarray}{rcl}
\nonumber
& & \KL{\lambda P_1 + (1-\lambda) P_2}{\lambda Q_1 + (1-\lambda) Q_2} \\
\nonumber
& & - \lambda \KL{ P_1 }{Q_1} + (1 - \lambda) \KL{ P_2}{Q_2}\\
\nonumber
& = & - \lambda \beta_1 \int \frac{g(x)}{\beta_1} f \left( \frac{\mathrm{d} P_1}{\mathrm{d} Q_1}(x) \left( g(x) \right)^{-1} \right) \mathrm{d} Q_1 (x)\squeezeequ \\
\nonumber
&& - (1-\lambda)\beta_2   \int \frac{g(x)}{\beta_2} f \left( \frac{\mathrm{d} P_2}{\mathrm{d} Q_2}(x) \left( g(x) \right)^{-1} \right)  \mathrm{d} Q_2  (x)\squeezeequ \\
\label{EqAfternoonInAgadir1}
& \leqslant & - \lambda \beta_1 f \left(\int \frac{g(x)}{\beta_1}   \frac{\mathrm{d} P_1}{\mathrm{d} Q_1}(x) \left( g(x) \right)^{-1} \mathrm{d} Q_1 (x) \right) \squeezeequ \\
\nonumber
&& - (1-\lambda)\beta_2    f \left( \int \frac{g(x)}{\beta_2} \frac{\mathrm{d} P_2}{\mathrm{d} Q_2}(x) \left( g(x) \right)^{-1} \mathrm{d} Q_2  (x)\right)  \squeezeequ \\
& = & - \lambda \beta_1 f \left( \frac{1}{\beta_1} \int\mathrm{d} P_1 (x) \right)  - (1-\lambda)\beta_2    f \left(\frac{1}{\beta_2} \int\mathrm{d} P_2  (x)\right)  \squeezeequ \IEEEeqnarraynumspace
\end{IEEEeqnarray}
\begin{IEEEeqnarray}{rcl}
& = & - \lambda \beta_1 f \left( \frac{1}{\beta_1} \right)  - (1-\lambda)\beta_2    f \left(\frac{1}{\beta_2} \right)  \squeezeequ \\
\label{EqAfternoonInAgadir2}
& \leqslant & -  f \left( \lambda \beta_1 \frac{1}{\beta_1} +  (1-\lambda)\beta_2   \frac{1}{\beta_2} \right) \\
& = & -  f \left( 1 \right) \\
&=&  0,
\end{IEEEeqnarray}
where the inequalities in~\eqref{EqAfternoonInAgadir1} and ~\eqref{EqAfternoonInAgadir2}  follow from Jensen's inequality \cite[Section~$6.3.5$]{ash2000probability} and  the fact that the function $f$ in~\eqref{EqValhala101} is strictly concave. Note that from~\eqref{EqWishpering}, in~\eqref{EqAfternoonInAgadir1}, for all $i \inCountTwo$, $ \int \frac{g(x)}{\beta_i} \mathrm{d} Q_i (x) = 1$; while in \eqref{EqAfternoonInAgadir2},
\begin{IEEEeqnarray}{rcl}
\lambda \beta_1 + (1 - \lambda) \beta_2 & = & \int g(\nu) \mathrm{d}\left(\lambda Q_1 + (1 - \lambda) Q_{2} \right)(\nu)\squeezeequ\\
& = & \int \mathrm{d}\left(\lambda P_1 + (1 - \lambda) P_{2} \right)(\nu)\\
& = & \lambda\int  \mathrm{d} P_1(\nu)  + (1 - \lambda) \int \mathrm{d}P_{2} (\nu) \squeezeequ\\
& = & 1.
\end{IEEEeqnarray}
Given the strict convexity of the function $f$ in \eqref{EqVousEtesVikings},  equality in \eqref{EqAfternoonInAgadir1} and \eqref{EqAfternoonInAgadir2} hold if and only if $P_{1} = P_{2}$ and $Q_{1} = Q_{2}$.
This completes the proof.

\section{Proof of Lemma~\ref{LemmaSubSetK}}\label{AppProofLemmaSubSetK} 
The proof is divided into two parts.  
The first part is as follows. Under the assumption that the set~$\set{K}_{Q, \vect{z}}$ in~\eqref{EqSetKxy} is empty, there is nothing to prove.  Alternatively, under the assumption that the set~$\set{K}_{Q, \vect{z}}$ is not empty, there always exists a real~$b \in \set{K}_{Q, \vect{z}}$, such that~$K_{Q, \vect{z}}\left(-\frac{1}{b} \right)  < +\infty$.
Note that for all~$\vect{\theta} \in \set{M}$,
\begin{IEEEeqnarray}{rCl}
\frac{\mathrm{d}}{\mathrm{d} t} \exp\left( -\frac{1}{t} \; \mathsf{L}_{\vect{z}} \left( \vect{\theta} \right) \right)  & = & \frac{1}{t^2} \mathsf{L}_{\vect{z}}\left( \vect{\theta} \right)\exp\left( -\frac{1}{t} \; \mathsf{L}_{\vect{z}}\left( \vect{\theta} \right) \right) \geqslant 0, \squeezeequ \IEEEeqnarraynumspace
\end{IEEEeqnarray}
with~$\mathsf{L}_{\vect{z}}$ in~\eqref{EqLxy}. 
Thus, from~\eqref{EqK}, it follows that~$K_{Q, \vect{z}}\left(-\frac{1}{b} \right)$ is nondecreasing with~$b$.
This implies that~$(0, b] \subseteq \set{K}_{Q, \vect{z}}$. 

Let~$b^{\star} \in (0,+\infty]$ be
\begin{equation}
b^{\star} = \sup \set{K}_{Q, \vect{z}}.
\end{equation}
Hence, if~$b^{\star} = +\infty$, it follows from~\eqref{EqSetKxy} that 
\begin{equation}\label{EqLemmaSubSetK1}
\set{K}_{Q, \vect{z}} = (0, +\infty ).
\end{equation}
Alternatively, if~$b^{\star} < +\infty$, it holds that
\begin{equation}\label{EqLemmaSubSetK01}
(0, b^{\star} ) \subseteq \set{K}_{Q, \vect{z}} \subseteq (0, b^{\star}].
\end{equation}
In either case,  it follows that $\set{K}_{Q, \vect{z}}$ is a convex set. This completes the first part of the proof. 

The second part of the proof is under the assumption that $Q$ is a probability measure. Under this assumption, for all~$\vect{\theta} \in \set{M}$ and for all for all~$t > 0$, it follows that
\begin{equation}
\exp\left( -\frac{1}{t} \; \mathsf{L}_{\vect{z}} \left( \vect{\theta} \right) \right) \leqslant1,
\end{equation}
with~$\mathsf{L}_{\vect{z}}$ in~\eqref{EqLxy}. 
Thus, 
\begin{IEEEeqnarray}{rcl}
K_{Q, \vect{z}}\left(-\frac{1}{t} \right) & = &  \log\left( \int \exp\left( -\frac{1}{t} \; \mathsf{L}_{\vect{z}}\left(\vect{\theta}\right)  \right) \mathrm{d}Q(\vect{\theta}) \right) \\
& \leqslant & \log\left( \int  \mathrm{d}Q(\vect{\theta}) \right) \\
& = & 0,
\end{IEEEeqnarray} 
which implies that~$(0, +\infty) \subseteq \set{K}_{Q, \vect{z}}$. Thus, if $Q$ is a probability measure, from~\eqref{EqSetKxy}, it holds that~$\set{K}_{Q, \vect{z}} = (0, +\infty)$, which completes the proof.

\section{Proof of Theorem~\ref{TheoremOptimalModel}}\label{AppProofTheoremOptimalModel}

The optimization problem in~\eqref{EqERMRER} can be re-written in terms of the Radon-Nikodym derivative of the optimization measure $P$ with respect to the measure $Q$, denoted by $\frac{\mathrm{d} P}{\mathrm{d} Q}: \set{M} \to [0, \infty)$, which yields: 
\begin{subequations}\label{EqTheKeyOP}
\begin{IEEEeqnarray}{ccl}
\nonumber
 \min_{P \in \triangle_{Q}\Bormeaspace{\set{M}}} & & \int \mathsf{L}_{\vect{z}} \left(\vect{\theta} \right)\frac{\mathrm{d}P}{\mathrm{d}Q}(\vect{\theta}) \mathrm{d}Q(\vect{\theta}) \\
 & &  + \lambda  \int  \frac{\mathrm{d}P}{\mathrm{d}Q}(\vect{\theta}) \log\left( \frac{\mathrm{d}P}{\mathrm{d}Q}(\vect{\theta})\right)  \mathrm{d}Q\left( \vect{\theta} \right) \Dsupersqueezeequ
\IEEEeqnarraynumspace\\
   \label{EqTheKeyOPConstraint}
    \mathrm{s.~t.} & & \int \frac{\mathrm{d}P}{\mathrm{d}Q}(\vect{\theta})  \d Q \left( \vect{ \theta} \right) = 1.
\end{IEEEeqnarray}
\end{subequations}
The remainder of the proof focuses on the problem in which the optimization is over  the function  $\frac{\mathrm{d} P}{\mathrm{d} Q}$ instead of the measure $P$. This is due to the fact that for all $P \in \bigtriangleup_{Q}\left( \set{M} \right)$, the Radon-Nikodym derivate $\frac{\mathrm{d} P}{\mathrm{d} Q}$ is unique up to sets of zero measure with respect to the measure $Q$.
Let $\mathscr{M}$ be the set of measurable functions $\set{M} \to \reals$ with respect to the measurable spaces $\left( \set{M}, \BorSigma{\set{M}} \right)$ and $\left( \reals, \BorSigma{\reals} \right)$ that are absolutely integrable with respect to $Q$.
That is, for all $\hat{g} \in \mathscr{M}$, it holds that
\begin{IEEEeqnarray}{rcl}
\label{EqValioLaPenaC}
\int \abs{ \hat{g} (\vect{\theta})} \mathrm{d} Q(\vect{\theta}) & < & \infty.
\end{IEEEeqnarray}

Hence, the optimization problem of interest is:
 \begin{subequations}\label{EqTheKeyOPInTheMaternity}
\begin{IEEEeqnarray}{ccl}
 \min_{g \in \mathscr{M}} & & \int \mathsf{L}_{\vect{z}} \left(\vect{\theta} \right)g(\vect{\theta}) \mathrm{d}Q(\vect{\theta}) + \lambda  \int  g(\vect{\theta}) \log\left( g(\vect{\theta})\right)  \mathrm{d}Q\left( \vect{\theta} \right) \Dsupersqueezeequ
\IEEEeqnarraynumspace\\
   \label{EqTheKeyOPConstraintInTheMaternity}
    \mathrm{s.~t.} & & \int g(\vect{\theta}) \d Q \left( \vect{ \theta} \right) = 1.
\end{IEEEeqnarray}
\end{subequations}

Let the Lagrangian of the optimization problem in~\eqref{EqTheKeyOPInTheMaternity} be the functional $L: \mathscr{M}\times \reals \rightarrow \reals$ such that
\begin{IEEEeqnarray}{rcl}
\nonumber
 L\left(g, \beta \right)  
&=& \displaystyle\int \mathsf{L}_{\vect{z}}\left( \vect{\nu}\right)g \left( \vect{\nu} \right)  \mathrm{d}Q\left( \vect{\nu} \right) 
   + \lambda \displaystyle\int g\left( \vect{\nu} \right)  \log \left(g\left( \vect{\nu} \right) \right)  \mathrm{d}Q\left( \vect{\nu} \right)\squeezeequ \\
\label{EqFunctionalL}
   & &    + \beta \left(\displaystyle\int g\left( \vect{\nu} \right) \mathrm{d}Q\left( \vect{\nu} \right)   -1 \right),
\end{IEEEeqnarray}
where~$\beta$ is a real that acts as a Lagrangian multiplier due to the constraint~\eqref{EqTheKeyOPConstraintInTheMaternity}.  
Let $\hat{g}: \set{M} \rightarrow \reals$ be a function in $\mathscr{M}$. 
The Gateaux differential of the functional $L$ in~\eqref{EqFunctionalL} at $\left(g, \beta\right) \in \mathscr{M}\times \reals$ in the direction of $\hat{g}$, if it exists, is
\begin{IEEEeqnarray}{rcl}
\label{EqNecessaryCondtion5321}
    \partial L(g, \beta; \hat{g} ) & \triangleq & \left.\frac{\mathrm{d}}{\mathrm{d} \gamma}  L(g + \gamma \hat{g}, \beta) \right|_{\gamma = 0}.
\end{IEEEeqnarray}
The proof continues under the assumption that the functions $g$ and $\hat{g}$ are such that the Gateaux differential in \eqref{EqNecessaryCondtion5321} exists. 
Under such an assumption, let the function $r: \reals \rightarrow \reals$ satisfy for all $\alpha \in (-\epsilon, \epsilon)$, with $\epsilon$ arbitrarily small, that
\begin{IEEEeqnarray}{lcl}
\nonumber
r(\alpha) 
& = & \displaystyle\int \mathsf{L}_{\vect{z}}\left( \vect{\nu}\right) \left(g \left( \vect{\nu} \right)  + \alpha \hat{g}\left( \vect{\nu} \right) \right) \mathrm{d}Q\left( \vect{\nu} \right) \\
\nonumber
& & +   \beta \left(\displaystyle\int\left(g\left( \vect{\nu} \right)  + \alpha \hat{g}\left( \vect{\nu} \right) \right)\mathrm{d}Q\left( \vect{\nu} \right)   -1 \right)\\
  &  &+  \lambda  \int  \left(\hat{g}\left( \vect{\nu} \right)  + \alpha \hat{g}\left( \vect{\nu} \right) \right) \log \left(g\left( \vect{\nu} \right)  + \alpha \hat{g}\left( \vect{\nu} \right) \right) \mathrm{d}Q\left( \vect{\nu} \right) \squeezeequ
\IEEEeqnarraynumspace\\
\nonumber
& = & \displaystyle\int g \left( \vect{\nu} \right) \left( \mathsf{L}_{\vect{z}}\left( \vect{\nu}\right)    + \beta \right) \mathrm{d}Q\left( \vect{\nu} \right) - \beta \\
\nonumber
& & +   \alpha \left(\displaystyle\int \hat{g}\left( \vect{\nu} \right) \left( \mathsf{L}_{\vect{z}}\left( \vect{\nu}\right) + \beta \right)\mathrm{d}Q\left( \vect{\nu} \right)  \right)\\
\label{Eqr}    
  &  &+  \lambda  \int  \left(g\left( \vect{\nu} \right)  + \alpha \hat{g}\left( \vect{\nu} \right) \right) \log \left(g\left( \vect{\nu} \right)  + \alpha \hat{g}\left( \vect{\nu} \right) \right) \mathrm{d}Q\left( \vect{\nu} \right), \squeezeequ
\IEEEeqnarraynumspace
\end{IEEEeqnarray}
where the last equality is simply an algebraic re-arrangement of terms. 
From the assumption that the functions $g$ and $\hat{g}$ are such that the Gateaux differential in \eqref{EqNecessaryCondtion5321} exists, it follows that  the function $r$ in \eqref{Eqr} is differentiable at zero. 
Note that the first two terms in \eqref{Eqr} are independent of $\alpha$; the third term is linear with $\alpha$; and the fourth term can be written using the function $\hat{r}: \reals \to \reals$ such that for all $\alpha \in (-\epsilon, \epsilon)$, with $\epsilon$ arbitrarily small, satisfies
\begin{IEEEeqnarray}{rcl}
\label{EqrHatForDaunas}
\hat{r}(\alpha) & = &   \lambda  \int  \left(g\left( \vect{\nu} \right)  + \alpha \hat{g}\left( \vect{\nu} \right) \right) \log \left(g\left( \vect{\nu} \right)  + \alpha \hat{g}\left( \vect{\nu} \right) \right) \mathrm{d}Q\left( \vect{\nu} \right) \squeezeequ\IEEEeqnarraynumspace\\
\label{EqNextToPequenoDali}
& = &  \lambda  \int  f\left(g\left( \vect{\nu} \right)  + \alpha \hat{g}\left( \vect{\nu} \right) \right)\mathrm{d}Q\left( \vect{\nu} \right),
\end{IEEEeqnarray}
where $f: (0, +\infty) \to \reals$ is such that $f(t) = t \log(t)$.
Under the same assumption,  it follows that  the function $\hat{r}$ in \eqref{EqrHatForDaunas} is differentiable at zero. That is, the limit 
\begin{IEEEeqnarray}{rcl}
\label{EqNeverDoingThisAgain}
\lim_{\delta \to 0}\frac{1}{\delta} ( \hat{r}(\gamma + \delta) - \hat{r}(\gamma) ) 
\end{IEEEeqnarray}
exists for all $\gamma \in (-\epsilon, \epsilon)$, with $\epsilon$ arbitrarily small.
Note that the function $f$ in \eqref{EqNextToPequenoDali} is continuous and differentiable (with finite derivate) in $(0, +\infty)$. 
Thus, the function $f$ is also Lipschitz continuous. 
This implies that for all $\vect{\theta} \in\supp Q$, and for all $\gamma \in (-\epsilon, \epsilon)$, with $\epsilon > 0$ arbitrarily small, it holds that
\begin{IEEEeqnarray}{rcl}
\abs{f(g(\vect{\theta}) + (\gamma + \delta) \hat{g}(\vect{\theta}))   - f(g(\vect{\theta}) + \gamma \hat{g}(\vect{\theta}))} & \leqslant & c \abs{\hat{g}(\vect{\theta})} \abs{\delta},\squeezeequ\IEEEeqnarraynumspace
\end{IEEEeqnarray}
with $\delta >0$, for some constant $c$ positive and finite. This implies that
\begin{IEEEeqnarray}{rcl}
\label{EqCoffeeTimeIsPast}
\abs{\frac{f(g(\vect{\theta}) + (\gamma + \delta) \hat{g}(\vect{\theta}))   - f(g(\vect{\theta}) + \gamma \hat{g}(\vect{\theta}))}{\delta}} & \leqslant &  c \abs{\hat{g}(\vect{\theta})}.\squeezeequ\IEEEeqnarraynumspace
\end{IEEEeqnarray}
Using these arguments, the limit in \eqref{EqNeverDoingThisAgain} satisfies for all $\gamma  \in (-\epsilon, \epsilon)$, with $\epsilon > 0$ arbitrarily small, that
\begin{IEEEeqnarray}{rcl}
\nonumber
& & \lim_{\delta \to 0}\frac{1}{\delta} ( \hat{r}(\gamma + \delta) - \hat{r}(\gamma) ) \\
\nonumber
& = &\lambda \lim_{\delta \to 0}  \int \frac{f(g(\vect{\theta}) + (\gamma + \delta) \hat{g}(\vect{\theta}))   - f(g(\vect{\theta}) + \gamma \hat{g}(\vect{\theta})) }{\delta} \,\mathrm{d} Q(\vect{\theta}) \\
\label{EqTheKeyEqualityOfThisS}
& = & \lambda \int  \dot{f}(g(\vect{\theta}) + \gamma \hat{g}(\vect{\theta}))  \hat{g}(\vect{\theta})\,\mathrm{d} Q(\vect{\theta})\\
\label{EqTheKeyEqualityOfThisSh}
& < & \infty,
\end{IEEEeqnarray}
where the function $\dot{f}: \left( 0, +\infty \right) \to \reals$ is the derivative of $f$. That is, $\dot{f}(t) = 1 +\log(t)$.
The equality in \eqref{EqTheKeyEqualityOfThisS} and the inequality in \eqref{EqTheKeyEqualityOfThisSh} follow from noticing that the conditions for the dominated convergence theorem hold \cite[Theorem~$1.6.9$]{ash2000probability}, namely:
\begin{itemize}
\item For all $\gamma  \in (-\epsilon, \epsilon)$, with $\epsilon > 0$, the inequality in \eqref{EqCoffeeTimeIsPast} holds;
\item The function $\hat{g}$   in \eqref{EqCoffeeTimeIsPast} satisfies the inequality in \eqref{EqValioLaPenaC}; and 
\item For all $\vect{\theta} \in\supp Q$ and for all $\gamma \in (-\epsilon, \epsilon)$, with $\epsilon > 0$ arbitrarily small, it holds  that
\begin{IEEEeqnarray}{rcl}
\nonumber
& &\lim_{\delta \to 0} \frac{f(g(\vect{\theta}) + (\gamma + \delta) \hat{g}(\vect{\theta}))   - f(g(\vect{\theta}) + \gamma \hat{g}(\vect{\theta})) }{\delta} \\
& =&  \frac{\mathrm{d}}{\mathrm{d} \gamma} f(g(\vect{\theta}) + \gamma \hat{g}(\vect{\theta})) 
\IEEEeqnarraynumspace \\
\label{EqCarryingDali}
& = &  \dot{f}(g(\vect{\theta}) + \gamma \hat{g}(\vect{\theta})) \hat{g}(\vect{\theta}).
\end{IEEEeqnarray}
\end{itemize}

Hence, the derivative of the real function~$r$ in~\eqref{Eqr} is
\begin{IEEEeqnarray}{rcl}
\nonumber
\frac{\mathrm{d}}{\mathrm{d}\alpha} r(\alpha) \squeezeequ
& = & \displaystyle\int  \mathsf{L}_{\vect{z}}\left( \vect{\nu}\right) \hat{g}\left( \vect{\nu} \right) \mathrm{d}Q\left( \vect{\nu} \right) + \beta\displaystyle\int  \hat{g}\left( \vect{\nu} \right) \mathrm{d}Q\left( \vect{\nu} \right)\Dsupersqueezeequ \\
\label{EqDerr} 
&+ & \lambda \displaystyle\int   \hat{g}\left( \vect{\nu} \right) \left( 1 +  \log \left(g \left( \vect{\nu} \right)  + \alpha \hat{g}\left( \vect{\nu} \right) \right) \right) \mathrm{d}Q\left( \vect{\nu} \right). \Dsupersqueezeequ
\IEEEeqnarraynumspace
\end{IEEEeqnarray}
From~\eqref{EqNecessaryCondtion5321} and~\eqref{EqDerr}, it follows that
\begin{IEEEeqnarray}{rcl}
\nonumber
& & \partial L\left(g, \beta; \hat{g} \right) \\
\label{EqGateauxDiff}
& =&  \displaystyle\int   \hat{g}\left( \vect{\nu} \right) \left(\mathsf{L}_{\vect{z}}\left( \vect{\nu}\right) 
    + \lambda \left( 1 +  \log\left(g \left( \vect{\nu} \right)  \right)  \right)  
    + \beta  \right) \mathrm{d}Q\left( \vect{\nu} \right). \Dsupersqueezeequ
\IEEEeqnarraynumspace
\end{IEEEeqnarray}

The relevance of the Gateaux differential in~\eqref{EqGateauxDiff} stems from \cite[Theorem~$1$, page~$178$]{luenberger1997bookOptimization}, which unveils the fact that a necessary condition for the functional~$L$ in~\eqref{EqFunctionalL} to have a minimum at~$\left( \frac{\mathrm{d}P^{\left(Q, \lambda\right)}_{\vect{\Theta} | \vect{Z} = \vect{z}}}{\mathrm{d}Q}, \beta \right) \in \mathscr{M} \times \reals$ is that for all functions~$\hat{g} \in \mathscr{M}$,
\begin{equation}\label{EqConditionhi}
\partial L\left(\frac{\mathrm{d}P^{\left(Q, \lambda\right)}_{\vect{\Theta} | \vect{Z} = \vect{z}}}{\mathrm{d}Q},\beta; \hat{g} \right)  = 0.    
\end{equation}

From~\eqref{EqConditionhi}, it follows that~$\frac{\mathrm{d}P^{\left(Q, \lambda\right)}_{\vect{\Theta} | \vect{Z} = \vect{z}}}{\mathrm{d}Q}~$ must satisfy for all functions~$\hat{g}$ in~$\mathscr{M}$ that
\begin{IEEEeqnarray}{l}
\nonumber
0 =\\ 
\nonumber
\displaystyle\int   \hat{g}\left( \vect{\nu} \right) \left(\mathsf{L}_{\vect{z}}\left( \vect{\nu}\right) 
    + \lambda \left( 1 +  \log\left( \frac{\mathrm{d}P^{\left(Q, \lambda\right)}_{\vect{\Theta} | \vect{Z} = \vect{z}}}{\mathrm{d}Q}\left( \vect{\nu} \right) \right) \right)  + \beta  \right) \mathrm{d}Q\left( \vect{\nu} \right),
    \Dsupersqueezeequ
\IEEEeqnarraynumspace
\end{IEEEeqnarray}
which implies that for all~$\vect{\nu} \in \supp Q$,
\begin{equation}
\mathsf{L}_{\vect{z}}\left( \vect{\nu}\right) 
    + \lambda \left( 1 +  \log\left( \frac{\mathrm{d}P^{\left(Q, \lambda\right)}_{\vect{\Theta} | \vect{Z} = \vect{z}}}{\mathrm{d}Q}\left( \vect{\nu} \right) \right) \right)  + \beta   = 0,
\end{equation}
and thus,
\begin{equation}
\label{EqKeyNorthStarx}
  \frac{\mathrm{d}P^{\left(Q, \lambda\right)}_{\vect{\Theta} | \vect{Z} = \vect{z}}}{\mathrm{d}Q}\left( \vect{\nu} \right) = \exp\left( - \frac{\beta + \lambda}{\lambda}\right) \exp\left( -\frac{\mathsf{L}_{\vect{z}}\left( \vect{\nu}\right) }{\lambda}\right),
\end{equation}
with~$\beta$ chosen to satisfy~\eqref{EqTheKeyOPConstraint}.
That is,
\begin{IEEEeqnarray}{rcl}
\label{EqKeyNorthStarHI}
 \frac{\mathrm{d}P^{\left(Q, \lambda\right)}_{\vect{\Theta} | \vect{Z} = \vect{z}}}{\mathrm{d}Q}\left( \vect{\nu} \right)  & = & \frac{\exp\left( -\frac{\mathsf{L}_{\vect{z}}\left( \vect{\nu}\right) }{\lambda}\right)}{\displaystyle\int \exp\left( -\frac{\mathsf{L}_{\vect{z}}\left( \vect{\theta}\right) }{\lambda}\right) \mathrm{d} Q\left( \vect{\theta}\right)}\\
 \label{EqKeyNorthStarHI65}
 & = &  \exp\left( - K_{Q, \vect{z}}\left(- \frac{1}{\lambda} \right) - \frac{1}{\lambda} \mathsf{L}_{\vect{z}}\left( \vect{\nu}\right)\right). 
\end{IEEEeqnarray}
The proof continues by verifying that the measure $P^{\left(Q, \lambda\right)}_{\vect{\Theta} | \vect{Z} = \vect{z}}$ that satisfies~\eqref{EqKeyNorthStarHI} is the unique solution to the ERM-RER problem in \eqref{EqERMRER}. Such verification is done by showing that the objective function in \eqref{EqERMRER} is strictly convex with the optimization variable. Let $P_1$ and $P_2$ be two different probability measures in $\Bormeaspace{\set{M}}$ and let $\alpha$ be in $(0,1)$. Hence,
\begin{IEEEeqnarray}{rcl}
\nonumber
& & \mathsf{R}_{\vect{z}}\left( \alpha P_1 + (1-\alpha) P_2 \right) + \lambda \KL{\alpha P_1 + (1-\alpha) P_2}{Q} \\
\nonumber
& = & \alpha  \mathsf{R}_{\vect{z}}\left(P_1\right) + (1-\alpha) \mathsf{R}_{\vect{z}}\left(P_2\right) + \lambda \KL{\alpha P_1 + (1-\alpha) P_2}{Q}  \squeezeequ\IEEEeqnarraynumspace\\
\nonumber
& > & \alpha  \left( \mathsf{R}_{\vect{z}}\left(P_1\right)  + \lambda \KL{ P_1}{Q} \right)
+ (1-\alpha) \left( \mathsf{R}_{\vect{z}}\left(P_2\right) + \lambda   \KL{ P_2}{Q} \right) \squeezeequ
\end{IEEEeqnarray}
where the functional~$\mathsf{R}_{\vect{z}}$ is defined in~\eqref{EqRxy}. The equality above follows from the properties of the Lebesgue integral, while  the inequality follows from Theorem~\ref{TheoremCasaBlanca}. This proves that the solution is unique due to the strict concavity of the objective function, which completes the proof.

\section{Proof of Lemma~\ref{CorPositive}}\label{AppProofCorPositive}

From Theorem~\ref{TheoremOptimalModel}, it follows that for all~$\vect{\theta} \in  \supp Q$,  %
\begin{IEEEeqnarray}{rcl} 
 \frac{\mathrm{d}P^{\left(Q, \lambda\right)}_{\vect{\Theta} | \vect{Z} = \vect{z}}}{\mathrm{d}Q}\left( \vect{\theta}\right)  
\label{EqBoolingLu}
& = & \exp\left( - K_{Q, \vect{z}}\left(- \frac{1}{\lambda} \right) - \frac{1}{\lambda} \mathsf{L}_{\vect{z}}\left( \vect{\theta}\right)\right)\\
\label{EqBoolingLun}
& \leqslant & \exp\left( - K_{Q, \vect{z}}\left(- \frac{1}{\lambda} \right) \right)\\
\label{EqBoolingLung}
& < & +\infty,
\end{IEEEeqnarray}
where the inequality in~\eqref{EqBoolingLun} follows from the fact that the function~$\mathsf{L}_{\vect{z}}$ is nonnegative; and the equality in~\eqref{EqBoolingLung} follows from the fact that~$\lambda \in \set{K}_{Q, \vect{z}}$.
This completes the proof of finiteness. 

The proof of positivity follows from observing that~$\lambda \in \set{K}_{Q, \vect{z}}$ and thus,~$K_{Q, \vect{z}}\left(- \frac{1}{\lambda} \right)  <  +\infty$.
Moreover,  for all~$\vect{\theta} \in  \supp Q$, it holds that~$\mathsf{L}_{\vect{z}}\left( \vect{\theta}\right) \leqslant +\infty$, with equality if and only if~$\mathsf{L}_{\vect{z}}\left( \vect{\theta}\right) = +\infty$.
These two observations put together yield
\begin{IEEEeqnarray}{rcl}
\frac{\mathrm{d}P^{\left(Q, \lambda\right)}_{\vect{\Theta} | \vect{Z} = \vect{z}}}{\mathrm{d}Q}\left( \vect{\theta}\right)  
& = & \exp\left( - K_{Q, \vect{z}}\left(- \frac{1}{\lambda} \right) - \frac{1}{\lambda} \mathsf{L}_{\vect{z}}\left( \vect{\theta}\right)\right)\\
& \geqslant & 0,
\end{IEEEeqnarray}
with equality if and only if~$\mathsf{L}_{\vect{z}}\left( \vect{\theta}\right) = +\infty$.
The proof continues  by showing that $Q\left( \left\lbrace \vect{\theta} \in \set{M}:  \mathsf{L}_{\vect{z}}\left( \vect{\theta}\right) = +\infty \right\rbrace  \right) =0$. From the assumption that~$\lambda \in \set{K}_{Q, \vect{z}}$, it follows that
\begin{IEEEeqnarray}{rCl}
+ \infty & > & \log\left( \int \exp\left( -\frac{1}{\lambda} \mathsf{L}_{\vect{z}} \left( \vect{\theta} \right)\right) \mathrm{d}Q\left( \vect{\theta} \right) \right)\\
\label{EqInsufficientPerformance}
& \geqslant &  \log\left( \int \left(1 -\frac{1}{\lambda} \mathsf{L}_{\vect{z}}\left( \vect{\theta} \right) \right) \mathrm{d}Q\left( \vect{\theta} \right) \right)\\
\label{EqInsufficientPerformanceA}
& = & \log\left(Q\left( \set{M} \right) - \frac{1}{\lambda} \int \mathsf{L}_{\vect{z}}\left( \vect{\theta} \right)  \mathrm{d}Q\left( \vect{\theta} \right) \right), 
\end{IEEEeqnarray}
where the inequality in \eqref{EqInsufficientPerformance} holds from the fact that for all $x \geqslant 0$, $\exp(x) > 1 + x$. 
Hence, the inequality in \eqref{EqInsufficientPerformanceA}  implies that
\begin{IEEEeqnarray}{rCl}
Q\left( \set{M} \right) & < & +\infty \mbox{ and }\\
\label{EqISITinPeril}
 \int \mathsf{L}_{\vect{z}}\left( \vect{\theta} \right)  \mathrm{d}Q\left( \vect{\theta} \right) & < & +\infty.
\end{IEEEeqnarray}
Finally, from \cite[Lemma~$1.6.6$]{ash2000probability} and the inequality in \eqref{EqISITinPeril} it follows that $Q\left( \left\lbrace \vect{\theta} \in \set{M}:  \mathsf{L}_{\vect{z}}\left( \vect{\theta}\right) = +\infty \right\rbrace  \right) =0$. Hence,
\begin{IEEEeqnarray}{rcl}
\frac{\mathrm{d}P^{\left(Q, \lambda\right)}_{\vect{\Theta} | \vect{Z} = \vect{z}}}{\mathrm{d}Q}\left( \vect{\theta}\right)  
& > &  0,
\end{IEEEeqnarray}
almost surely with respect to the measure $Q$. And this completes the proof.

\section{Proof of Lemma~\ref{LemmaMutualAC}} \label{AppProofLemmaMutualAC} 
The probability measure~$P^{\left(Q, \lambda\right)}_{\vect{\Theta} | \vect{Z} = \vect{z}}$ in~\eqref{EqGenpdf} satisfies for all~$\set{C} \in \BorSigma{\set{M}}$,
\begin{IEEEeqnarray}{rcl}
P^{\left(Q, \lambda\right)}_{\vect{\Theta} | \vect{Z} = \vect{z}}\left( \set{C} \right) & = & \int_{\set{C}} \frac{\mathrm{d}P^{\left(Q, \lambda\right)}_{\vect{\Theta} | \vect{Z} = \vect{z}}}{\mathrm{d}Q} \left( \vect{\theta} \right) \mathrm{d}Q \left( \vect{\theta} \right), 
\end{IEEEeqnarray}
and thus, if~$Q \left( \set{C} \right) = 0$, then
\begin{IEEEeqnarray}{rcl}
P^{\left(Q, \lambda\right)}_{\vect{\Theta} | \vect{Z} = \vect{z}}\left( \set{C} \right) & = & 0,
\end{IEEEeqnarray}
which implies the absolute continuity of~$P^{\left(Q, \lambda\right)}_{\vect{\Theta} | \vect{Z} = \vect{z}}$ with respect to~$Q$. 

Alternatively,  given a set~$\set{C}\in \BorSigma{\set{M}}$, assume now that~$P^{\left(Q, \lambda\right)}_{\vect{\Theta} | \vect{Z} = \vect{z}}\left( \set{C} \right) = 0$. Hence,  it follows that 
\begin{IEEEeqnarray}{rcl}
\label{EqRikitea756}
0 & = & P^{\left(Q, \lambda\right)}_{\vect{\Theta} | \vect{Z} = \vect{z}}\left( \set{C} \right)\\
\label{EqRikitea757}
& = &  \int_{\set{C}} \frac{\mathrm{d}P^{\left(Q, \lambda\right)}_{\vect{\Theta} | \vect{Z} = \vect{z}}}{\mathrm{d}Q} \left( \vect{\theta} \right) \mathrm{d}Q \left( \vect{\theta} \right).
\end{IEEEeqnarray}
From Lemma~\ref{CorPositive}, it follows that~$Q\left( \left\lbrace \vect{\theta} \in \set{M}:  \frac{\mathrm{d}P^{\left(Q, \lambda\right)}_{\vect{\Theta} | \vect{Z} = \vect{z}}}{\mathrm{d}Q} \left( \vect{\theta} \right) = 0 \right\rbrace  \right) =0$. Hence,  
\begin{IEEEeqnarray}{rcl}
\int_{\set{C}} \frac{\mathrm{d}P^{\left(Q, \lambda\right)}_{\vect{\Theta} | \vect{Z} = \vect{z}}}{\mathrm{d}Q} \left( \vect{\theta} \right) \mathrm{d}Q \left( \vect{\theta} \right) & = & 0,
\end{IEEEeqnarray}
if and only if~$Q\left( \set{C} \right) = 0$. This verifies the absolute continuity of~$Q$ with respect to~$P^{\left(Q, \lambda\right)}_{\vect{\Theta} | \vect{Z} = \vect{z}}$, and completes the proof. 
  
\section{Proof of Lemma~\ref{LemmaMutualAlphaBeta}}\label{AppProofLemmaMutualAlphaBeta}

Consider the function~$g: \set{M} \rightarrow [0, +\infty)$, 
\begin{equation}
g(\vect{\theta}) =  
\frac{\mathrm{d}P^{(Q,\alpha)}_{\vect{\Theta}| \vect{Z} = \vect{z}}}{\mathrm{d}Q}\left( \vect{\theta} \right) \left( \frac{\mathrm{d}P^{(Q,\beta)}_{\vect{\Theta}| \vect{Z} = \vect{z}}}{\mathrm{d}Q}\left( \vect{\theta} \right)  \right)^{-1},
\end{equation}
and note that for all~$\vect{\theta} \in \supp Q$ $\setminus$ $\big\lbrace \vect{\nu} \in \set{M}$:~$\mathsf{L}_{\vect{z}}\left( \vect{\nu}\right) = +\infty \big\rbrace$,~$g\left( \vect{\theta} \right) > 0$.
Alternatively, for all~$\vect{\theta} \in  \left\lbrace \vect{\nu} \in \set{M}:  \mathsf{L}_{\vect{z}}\left( \vect{\nu}\right) = +\infty \right\rbrace~$,~$g\left( \vect{\theta} \right) = 0$, which follows from the assumption~$0 \cdot \frac{1}{0} = 0$. 

Consider a measure~$P$ on~$\Bormeaspace{\set{M}}$, such that for all  sets~$\set{A} \in \BorSigma{\set{M}}$,
\begin{equation}
\label{EqTelephonePerdu}
P\left( \set{A}\right) = \int_{\set{A}} g(\vect{\theta}) \d P^{(Q,\beta)}_{\vect{\Theta}| \vect{Z} = \vect{z}} (\vect{\theta}),
\end{equation}
and note that if~$P^{(Q,\beta)}_{\vect{\Theta}| \vect{Z} = \vect{z}} (\set{A}) = 0$, then~$P\left( \set{A}\right) = 0$. This implies that~$P$ is absolutely continuous with respect to~$P^{(Q,\beta)}_{\vect{\Theta}| \vect{Z} = \vect{z}} (\set{A})$.
Moreover, from~\eqref{EqTelephonePerdu}, it follows that
\begin{IEEEeqnarray}{rCl}
\nonumber
& &P\left( \set{A}\right) \\
& = & \int_{\set{A}} \frac{\mathrm{d}P^{(Q,\alpha)}_{\vect{\Theta}| \vect{Z} = \vect{z}}}{\mathrm{d}Q}\left( \vect{\theta} \right) \left( \frac{\mathrm{d}P^{(Q,\beta)}_{\vect{\Theta}| \vect{Z} = \vect{z}}}{\mathrm{d}Q}\left( \vect{\theta} \right)  \right)^{-1}\d P^{(Q,\beta)}_{\vect{\Theta}| \vect{Z} = \vect{z}} (\vect{\theta})\supersqueezeequ
\IEEEeqnarraynumspace\\
& = &\hspace{-1ex} \int_{\set{A}} \frac{\mathrm{d}P^{(Q,\alpha)}_{\vect{\Theta}| \vect{Z} = \vect{z}}}{\mathrm{d}Q}\left( \vect{\theta} \right) \left( \frac{\mathrm{d}P^{(Q,\beta)}_{\vect{\Theta}| \vect{Z} = \vect{z}}}{\mathrm{d}Q}\left( \vect{\theta} \right)  \right)^{-1}\frac{\mathrm{d}P^{(Q,\beta)}_{\vect{\Theta}| \vect{Z} = \vect{z}}}{\mathrm{d}Q}\left( \vect{\theta} \right)  \d Q (\vect{\theta}) \supersqueezeequ
\IEEEeqnarraynumspace\\
& = & \int_{\set{A}} \frac{\mathrm{d}P^{(Q,\alpha)}_{\vect{\Theta}| \vect{Z} = \vect{z}}}{\mathrm{d}Q}\left( \vect{\theta} \right)  \d Q (\vect{\theta})\\
& = & \int_{\set{A}}   \d P^{(Q,\alpha)}_{\vect{\Theta}| \vect{Z} = \vect{z}} (\vect{\theta})\\
& = & P^{(Q,\alpha)}_{\vect{\Theta}| \vect{Z} = \vect{z}} ( \set{A} ),
\end{IEEEeqnarray}
which implies that the probability measures $P$ in \eqref{EqTelephonePerdu} and $P^{(Q,\alpha)}_{\vect{\Theta}| \vect{Z} = \vect{z}}$ are identical. 
Thus,~$P^{(Q,\alpha)}_{\vect{\Theta}| \vect{Z} = \vect{z}}$ is absolutely continuous with respect to~$P^{(Q,\beta)}_{\vect{\Theta}| \vect{Z} = \vect{z}}$.
The proof that~$P^{(Q,\beta)}_{\vect{\Theta}| \vect{Z} = \vect{z}}$ is absolutely continuous with respect to~$P^{(Q,\alpha)}_{\vect{\Theta}| \vect{Z} = \vect{z}}$ follows the same argument. This completes the proof.

 \section{Proof of Lemma~\ref{CorAssymptoticZero}} \label{AppProofCorAssymptoticZero}
 
From Theorem~\ref{TheoremOptimalModel}, the probability measure~$P^{\left(Q, \lambda\right)}_{\vect{\Theta} | \vect{Z} = \vect{z}}$ in~\eqref{EqGenpdf} satisfies for all~$\vect{\theta} \in \supp Q$,
\begin{IEEEeqnarray}{rcl} 
\frac{\mathrm{d}P^{\left(Q, \lambda\right)}_{\vect{\Theta} | \vect{Z} = \vect{z}}}{\mathrm{d}Q} \left( \vect{\theta} \right) \supersqueezeequ
& = & \frac{\exp\left( -\frac{\mathsf{L}_{\vect{z}}\left( \vect{\theta}\right)}{\lambda}\right)}{\displaystyle\int\exp\left( -\frac{\mathsf{L}_{\vect{z}}\left( \vect{\nu} \right)}{\lambda} \right)  \mathrm{d} Q\left(\vect{\nu} \right)}\\
  & =& \left(\exp\left( \frac{\mathsf{L}_{\vect{z}}\left( \vect{\theta}\right)}{\lambda}\right) \displaystyle\int\exp\left( -\frac{\mathsf{L}_{\vect{z}}\left( \vect{\nu} \right)}{\lambda} \right)  \mathrm{d} Q\left(\vect{\nu} \right)  \right)^{-1}\supersqueezeequ
\IEEEeqnarraynumspace\\
  \label{EqSophia}
    & =& \left( \displaystyle\int\exp\left( \frac{1}{\lambda}\left(\mathsf{L}_{\vect{z}}\left( \vect{\theta}\right)-  \mathsf{L}_{\vect{z}}\left( \vect{\nu} \right)  \right)   \right)   \mathrm{d} Q\left(\vect{\nu} \right) \right)^{-1}.\supersqueezeequ
\IEEEeqnarraynumspace
\end{IEEEeqnarray}
Given~$\vect{\theta} \in \supp Q$, consider the partition of~$\supp Q$ formed by the sets~$\set{A}_{0}\left( \vect{\theta} \right)$, $\set{A}_{1}\left( \vect{\theta} \right)$, and~$\set{A}_{2}\left( \vect{\theta} \right)$, which satisfy the following:
\begin{subequations}\label{EqKeyEquationsTaravao}
\begin{IEEEeqnarray}{rcl}
\set{A}_{0}\left( \vect{\theta} \right) & \triangleq & \left\lbrace \vect{\nu} \in \supp Q:  \mathsf{L}_{\vect{z}}\left( \vect{\theta}\right) -  \mathsf{L}_{\vect{z}}\left( \vect{\nu}\right) = 0 \right\rbrace, \\
\set{A}_{1}\left( \vect{\theta} \right) & \triangleq & \left\lbrace \vect{\nu} \in \supp Q:  \mathsf{L}_{\vect{z}}\left( \vect{\theta}\right) -  \mathsf{L}_{\vect{z}}\left( \vect{\nu}\right) < 0 \right\rbrace,  \mbox{ and }
\IEEEeqnarraynumspace\\
\set{A}_{2}\left( \vect{\theta} \right) & \triangleq & \left\lbrace \vect{\nu} \in \supp Q:  \mathsf{L}_{\vect{z}}\left( \vect{\theta}\right) -  \mathsf{L}_{\vect{z}}\left( \vect{\nu}\right) > 0 \right\rbrace.
\end{IEEEeqnarray}
\end{subequations}
Using the sets~$\set{A}_{0}\left( \vect{\theta} \right)$,~$\set{A}_{1}\left( \vect{\theta} \right)$, and~$\set{A}_{2}\left( \vect{\theta} \right)$ in~\eqref{EqSophia}, the following holds for all~$\vect{\theta} \in \supp Q$,
\begin{IEEEeqnarray}{rcl} 
\nonumber
& & \frac{\mathrm{d}P^{\left(Q, \lambda\right)}_{\vect{\Theta} | \vect{Z} = \vect{z}}}{\mathrm{d}Q} \left( \vect{\theta} \right)      \supersqueezeequ\\ 
\nonumber
& =& \Bigg( \displaystyle\int_{\set{A}_{0}\left( \vect{\theta} \right)}\exp\left( \frac{1}{\lambda}\left(\mathsf{L}_{\vect{z}}\left( \vect{\theta}\right)-  \mathsf{L}_{\vect{z}}\left( \vect{\nu} \right)  \right)   \right)\mathrm{d} Q\left(\vect{\nu} \right)\squeezeequ \\
\nonumber
    & +&  \displaystyle\int_{\set{A}_{1}\left( \vect{\theta} \right)}\exp\left( \frac{1}{\lambda}\left(\mathsf{L}_{\vect{z}}\left( \vect{\theta}\right)-  \mathsf{L}_{\vect{z}}\left( \vect{\nu} \right)  \right)   \right) \mathrm{d} Q\left(\vect{\nu} \right)  \squeezeequ\\
    & + &  \displaystyle\int_{\set{A}_{2}\left( \vect{\theta} \right)}\exp\left( \frac{1}{\lambda}\left(\mathsf{L}_{\vect{z}}\left( \vect{\theta}\right)-  \mathsf{L}_{\vect{z}}\left( \vect{\nu} \right)  \right)   \right) \mathrm{d} Q\left(\vect{\nu} \right)  \Bigg)^{-1}\supersqueezeequ\\
    \nonumber
    & =& \Bigg( Q\left( \set{A}_{0}\left( \vect{\theta} \right) \right)   +  \displaystyle\int_{\set{A}_{1}\left( \vect{\theta} \right)}\exp\left( \frac{1}{\lambda}\left(\mathsf{L}_{\vect{z}}\left( \vect{\theta}\right)-  \mathsf{L}_{\vect{z}}\left( \vect{\nu} \right)  \right)   \right) \mathrm{d} Q\left(\vect{\nu} \right) \supersqueezeequ
\IEEEeqnarraynumspace \\
    \label{EqTitereFue}
    & +&   \displaystyle\int_{\set{A}_{2}\left( \vect{\theta} \right)}\exp\left( \frac{1}{\lambda}\left(\mathsf{L}_{\vect{z}}\left( \vect{\theta}\right)-  \mathsf{L}_{\vect{z}}\left( \vect{\nu} \right)  \right)   \right) \mathrm{d} Q \left(\vect{\nu} \right)  \Bigg)^{-1}.\supersqueezeequ
\end{IEEEeqnarray}

Note that the sets 
\begin{IEEEeqnarray}{l}
\left\lbrace \vect{\nu} \in \supp Q:  \mathsf{L}_{\vect{z}}\left( \vect{\nu}\right) = \delta_{Q,\vect{z}}^{\star} \right\rbrace, \\
\left\lbrace \vect{\nu} \in \supp Q:  \mathsf{L}_{\vect{z}}\left( \vect{\nu}\right) > \delta_{Q,\vect{z}}^{\star} \right\rbrace\mbox{,  and}\\
\left\lbrace \vect{\nu} \in \supp Q:  \mathsf{L}_{\vect{z}}\left( \vect{\nu}\right) < \delta_{Q,\vect{z}}^{\star} \right\rbrace,
\end{IEEEeqnarray}
with~$\delta_{Q,\vect{z}}^{\star}$ in~\eqref{EqDeltaStar},  form a partition of the set~$\supp Q$. Following this observation, the rest of the proof is divided into three parts.  
The first part evaluates~$\lim_{\lambda \rightarrow 0^{+}} \frac{\mathrm{d}P^{\left(Q, \lambda\right)}_{\vect{\Theta} | \vect{Z} = \vect{z}}}{\mathrm{d}Q} \left( \vect{\theta} \right)$, with~$\vect{\theta} \in\left\lbrace \vect{\nu} \in \set{M}:  \mathsf{L}_{\vect{z}}\left( \vect{\nu}\right) = \delta_{Q,\vect{z}}^{\star} \right\rbrace$. 
The second part considers the case in which~$\vect{\theta} \in \left\lbrace \vect{\nu} \in \set{M}:  \mathsf{L}_{\vect{z}}\left( \vect{\nu}\right) > \delta_{Q,\vect{z}}^{\star} \right\rbrace$. 
The third part considers the remaining case.

The first part is as follows. 
Consider that~$\vect{\theta}  \in  \left\lbrace \vect{\nu} \in \set{M}:  \mathsf{L}_{\vect{z}}\left( \vect{\nu}\right) = \delta_{Q,\vect{z}}^{\star} \right\rbrace$ and note that $ \left\lbrace \vect{\nu} \in \set{M}:  \mathsf{L}_{\vect{z}}\left( \vect{\nu}\right) = \delta_{Q,\vect{z}}^{\star} \right\rbrace = \set{L}^{\star}_{Q,\vect{z}}$. Hence, the sets~$\set{A}_{0}\left( \vect{\theta} \right)$,~$\set{A}_{1}\left( \vect{\theta} \right)$, and~$\set{A}_{2}\left( \vect{\theta} \right)$  in~\eqref{EqKeyEquationsTaravao} satisfy the following:
\begin{subequations}\label{EqBeaudeEDF}
\begin{IEEEeqnarray}{rCl} 
\set{A}_{0}\left( \vect{\theta} \right) & = &    \set{L}^{\star}_{Q,\vect{z}},\\
\set{A}_{1}\left( \vect{\theta} \right) & = &   \left\lbrace \vect{\mu} \in \supp Q: \mathsf{L}_{\vect{z}}\left( \vect{\mu}\right) > \delta_{Q,\vect{z}}^{\star}\right\rbrace, \mbox{ and }\\
\set{A}_{2}\left( \vect{\theta} \right) & = &   \left\lbrace \vect{\mu} \in \supp Q: \mathsf{L}_{\vect{z}}\left( \vect{\mu}\right) < \delta_{Q,\vect{z}}^{\star}\right\rbrace.
\end{IEEEeqnarray}
\end{subequations}
From the definition of~$\delta_{Q,\vect{z}}^{\star}$ in~\eqref{EqDeltaStar}, it follows that~$Q\left(  \set{A}_{2}\left( \vect{\theta} \right) \right) = 0$. Plugging the equalities in~\eqref{EqBeaudeEDF} in~\eqref{EqTitereFue} yields for all~$\vect{\theta}  \in  \left\lbrace \vect{\nu} \in \set{M}:  \mathsf{L}_{\vect{z}}\left( \vect{\nu}\right) = \delta_{Q,\vect{z}}^{\star} \right\rbrace$,
\begin{IEEEeqnarray}{rcl} 
\nonumber
& &\frac{\mathrm{d}P^{\left(Q, \lambda\right)}_{\vect{\Theta} | \vect{Z} = \vect{z}}}{\mathrm{d}Q} \left( \vect{\theta} \right) \supersqueezeequ \\
    \label{EqBlanc0153}
& =& \bigg( Q\left(\set{L}^{\star}_{Q,\vect{z}} \right)  +  \displaystyle \int_{\set{A}_{1}\left( \vect{\theta}\right)}\hspace{-2ex}\exp\left( \frac{1}{\lambda}\left(\mathsf{L}_{\vect{z}}\left( \vect{\theta}\right)-  \mathsf{L}_{\vect{z}}\left( \vect{\nu} \right)  \right)   \right) \mathrm{d} Q\left(\vect{\nu} \right)  \bigg)^{-1}.\supersqueezeequ
\IEEEeqnarraynumspace
\end{IEEEeqnarray}
The equality in~\eqref{EqBlanc0153} implies that  for all~$\vect{\theta}  \in  \left\lbrace \vect{\nu} \in \set{M}:  \mathsf{L}_{\vect{z}}\left( \vect{\nu}\right) = \delta_{Q,\vect{z}}^{\star} \right\rbrace~$,
\begin{IEEEeqnarray}{rcl}
\nonumber
& & \lim_{\lambda \rightarrow 0^{+}} \frac{\mathrm{d}P^{\left(Q, \lambda\right)}_{\vect{\Theta} | \vect{Z} = \vect{z}}}{\mathrm{d}Q} \left( \vect{\theta} \right) \\
\nonumber
  & = & \bigg( \lim_{\lambda \rightarrow 0^{+}}  \displaystyle\int_{\set{A}_{1}\left( \vect{\theta}\right) }\exp\left( \frac{1}{\lambda}\left(\mathsf{L}_{\vect{z}}\left( \vect{\theta}\right)-  \mathsf{L}_{\vect{z}}\left( \vect{\nu} \right)  \right)   \right) \mathrm{d} Q\left(\vect{\nu} \right) + \squeezeequ
\IEEEeqnarraynumspace \\
&+ &  Q\left( \set{L}^{\star}_{Q,\vect{z}}  \right)  \bigg)^{-1} \\
\label{EqBusinessClassB}
& = &\left\lbrace
\begin{array}{cl}
+ \infty &\text{ if } Q\left( \set{L}^{\star}_{Q,\vect{z}}\right) = 0\\
  \frac{1}{Q\left( \set{L}^{\star}_{Q,\vect{z}}\right) }  & \text{ otherwise}. 
\end{array}
\right.
\end{IEEEeqnarray}
where 
the equality in~\eqref{EqBusinessClassB} follows from verifying that   the dominated convergence theorem~\cite[Theorem~$2.6.9$]{ash2000probability} holds. That is, \newline
$(a)$ For all~$\vect{\nu} \in \set{A}_{1}\left( \vect{\theta}\right)$, it holds that~$\exp\left( \frac{1}{\lambda}\left(\mathsf{L}_{\vect{z}}\left( \vect{\theta}\right)-  \mathsf{L}_{\vect{z}}\left( \vect{\nu} \right)  \right)   \right)$ $< 1$; and \newline
$(b)$ For all~$\vect{\nu} \in \set{A}_{1}\left( \vect{\theta}\right)$, it holds that 
\begin{equation}
\lim_{\lambda \rightarrow 0^{+}} \exp\left( \frac{1}{\lambda}\left(\mathsf{L}_{\vect{z}}\left( \vect{\theta}\right)-  \mathsf{L}_{\vect{z}}\left( \vect{\nu} \right)  \right)   \right)=0.
\end{equation}
This completes the first part of the proof.

The second part is as follows. For all~$\delta > \delta_{Q,\vect{z}}^{\star}$ and for all~$\vect{\theta}$ $\in$ $\big\lbrace \vect{\nu} \in \supp Q:$ $\mathsf{L}_{\vect{z}}\left( \vect{\nu}\right) = \delta \big\rbrace$, the sets~$\set{A}_{0}\left( \vect{\theta} \right)$,~$\set{A}_{1}\left( \vect{\theta} \right)$, and~$\set{A}_{2}\left( \vect{\theta} \right)$ in~\eqref{EqKeyEquationsTaravao} satisfy the following:
\begin{subequations}\label{EqBeaudeEDFa}
\begin{IEEEeqnarray}{rcl}
\set{A}_{0}\left( \vect{\theta} \right) & = &   \left\lbrace \vect{\mu} \in \supp Q: \mathsf{L}_{\vect{z}}\left( \vect{\mu}\right) = \delta \right\rbrace,\\
\set{A}_{1}\left( \vect{\theta} \right) & = &   \left\lbrace \vect{\mu} \in \supp Q: \mathsf{L}_{\vect{z}}\left( \vect{\mu}\right) > \delta \right\rbrace, \mbox{ and }\\
\set{A}_{2}\left( \vect{\theta} \right) & = &   \left\lbrace \vect{\mu} \in \supp Q: \mathsf{L}_{\vect{z}}\left( \vect{\mu}\right) < \delta \right\rbrace.
\end{IEEEeqnarray}
\end{subequations}

Consider the sets 
\begin{IEEEeqnarray}{rcl}
\set{A}_{2,1}\left( \vect{\theta} \right) & \triangleq & \left\lbrace \vect{\mu} \in \set{A}_{2}\left( \vect{\theta} \right): \mathsf{L}_{\vect{z}}\left( \vect{\mu}\right) < \delta_{Q,\vect{z}}^{\star} \right\rbrace,  \mbox{ and }\\
\set{A}_{2,2}\left( \vect{\theta} \right) & \triangleq & \left\lbrace \vect{\mu} \in \set{A}_{2}\left( \vect{\theta} \right): \delta_{Q,\vect{z}}^{\star} \leqslant \mathsf{L}_{\vect{z}}\left( \vect{\mu}\right) < \delta \right\rbrace, 
\end{IEEEeqnarray}
and note that~$\set{A}_{2,1}\left( \vect{\theta} \right)$ and~$\set{A}_{2,2}\left( \vect{\theta} \right)$ form a partition of~$\set{A}_{2}\left( \vect{\theta} \right)$. Moreover,  from the definition of~$\delta_{Q,\vect{z}}^{\star}$ in~\eqref{EqDeltaStar}, it holds that
\begin{IEEEeqnarray}{rcl}
\label{EqBackHome}
Q\left( \set{A}_{2,1}\left( \vect{\theta} \right) \right) = 0.
\end{IEEEeqnarray}
Hence, plugging the equalities in~\eqref{EqBeaudeEDFa} and~\eqref{EqBackHome} in~\eqref{EqTitereFue} yields,  for all~$\delta > \delta_{Q,\vect{z}}^{\star}$ and for all~$\vect{\theta} \in \left\lbrace \vect{\nu} \in \set{M}:  \mathsf{L}_{\vect{z}}\left( \vect{\nu}\right) = \delta \right\rbrace$, 
\begin{IEEEeqnarray}{rcl} 
\nonumber
& & \frac{\mathrm{d}P^{\left(Q, \lambda\right)}_{\vect{\Theta} | \vect{Z} = \vect{z}}}{\mathrm{d}Q} \left( \vect{\theta} \right) \\
\nonumber
   & =& \bigg( Q\left( \set{A}_{0}\left( \vect{\theta}\right) \right)  +  \displaystyle\int_{\set{A}_{1}\left( \vect{\theta}\right)}\exp\left( \frac{1}{\lambda}\left(\mathsf{L}_{\vect{z}}\left( \vect{\theta}\right) -  \mathsf{L}_{\vect{z}}\left( \vect{\nu} \right)  \right)   \right) \mathrm{d} Q\left(\vect{\nu} \right) \supersqueezeequ
\IEEEeqnarraynumspace\\
   \label{EqBlanc2087}
   & +&   \displaystyle\int_{\set{A}_{2,2}\left( \vect{\theta}\right)}\exp\left( \frac{1}{\lambda}\left(\mathsf{L}_{\vect{z}}\left( \vect{\theta}\right)-  \mathsf{L}_{\vect{z}}\left( \vect{\nu} \right)  \right)   \right) \mathrm{d} Q\left(\vect{\nu} \right) \bigg)^{-1}.
\end{IEEEeqnarray}
The equality in~\eqref{EqBlanc2087} implies that  for all~$\delta > \delta_{Q,\vect{z}}^{\star}$ and for all~$\vect{\theta}$ $\in$ $\big\lbrace \vect{\nu} \in \set{M}$:~$\mathsf{L}_{\vect{z}}\left( \vect{\nu}\right) = \delta \big\rbrace$,
\begin{IEEEeqnarray}{rcl}
\nonumber
 & & \lim_{\lambda \rightarrow 0^{+}} \frac{\mathrm{d}P^{\left(Q, \lambda\right)}_{\vect{\Theta} | \vect{Z} = \vect{z}}}{\mathrm{d}Q} \left( \vect{\theta} \right) \\
\nonumber
& = & \Bigg(\lim_{\lambda \rightarrow 0^{+}}\displaystyle\int_{\set{A}_{1}\left( \vect{\theta} \right)}\exp\left( \frac{1}{\lambda}\left(\mathsf{L}_{\vect{z}}\left( \vect{\theta}\right)-  \mathsf{L}_{\vect{z}}\left( \vect{\nu} \right)  \right)   \right) \mathrm{d} Q\left(\vect{\nu} \right)  \\
\nonumber
    & + &   \lim_{\lambda \rightarrow 0^{+}}\displaystyle\int_{\set{A}_{2,2}\left( \vect{\theta} \right)}\exp\left( \frac{1}{\lambda}\left(\mathsf{L}_{\vect{z}}\left( \vect{\theta}\right)-  \mathsf{L}_{\vect{z}}\left( \vect{\nu} \right)  \right)   \right) \mathrm{d} Q\left(\vect{\nu} \right) \\
    \label{EqNowIntheAir1}    
    &+ & Q\left( \set{A}_{0}\left( \vect{\theta} \right) \right) \Bigg)^{-1}\\
\nonumber    
    & = & \Bigg(  \lim_{\lambda \rightarrow 0^{+}}\displaystyle\int_{\set{A}_{2,2}\left( \vect{\theta} \right)}\exp\left( \frac{1}{\lambda}\left(\mathsf{L}_{\vect{z}}\left( \vect{\theta}\right)-  \mathsf{L}_{\vect{z}}\left( \vect{\nu} \right)  \right)   \right) \mathrm{d} Q\left(\vect{\nu} \right)  \\
\label{EqNowIntheAir2}    
    & + &    Q\left( \set{A}_{0}\left( \vect{\theta} \right) \right) \Bigg)^{-1}\\
    \label{EqNowIntheAir3}    
    & = & \big( Q\left( \set{A}_{0}\left( \vect{\theta} \right) \right)   + \infty \big)^{-1}\\
    \label{EqNowIntheAir4}
    & = & 0,
\end{IEEEeqnarray}
where 
the equality in~\eqref{EqNowIntheAir2} follows by verifying that   the dominated convergence theorem~\cite[Theorem~$2.6.9$]{ash2000probability} holds. That is, \newline
$(a)$ For all~$\vect{\nu} \in \set{A}_{1}\left( \vect{\theta}\right)$, it holds that~$\exp\left( \frac{1}{\lambda}\left(\mathsf{L}_{\vect{z}}\left( \vect{\theta}\right)-  \mathsf{L}_{\vect{z}}\left( \vect{\nu} \right)  \right)   \right)$ $< 1$; and \newline
$(b)$ For all~$\vect{\nu} \in \set{A}_{1}\left( \vect{\theta}\right)$, it holds that 
\begin{equation}
\lim_{\lambda \rightarrow 0^{+}} \exp\left( \frac{1}{\lambda}\left(\mathsf{L}_{\vect{z}}\left( \vect{\theta}\right)-  \mathsf{L}_{\vect{z}}\left( \vect{\nu} \right)  \right)   \right)=0.
\end{equation}
This completes the second part.

The third part of the proof follows by noticing that the set~$\big\lbrace \vect{\nu} \in \supp Q$:~$\mathsf{L}_{\vect{z}}\left( \vect{\nu}\right) < \delta_{Q,\vect{z}}^{\star} \big\rbrace$ is a negligible set with respect to~$Q$ and thus, for all~$\vect{\theta} \in \big\lbrace \vect{\nu} \in \supp Q:$  $\mathsf{L}_{\vect{z}}\left( \vect{\nu}\right) < \delta_{Q,\vect{z}}^{\star} \big\rbrace$, the value~$\frac{\mathrm{d}P^{\left(Q, \lambda\right)}_{\vect{\Theta} | \vect{Z} = \vect{z}}}{\mathrm{d}Q} \left( \vect{\theta} \right)$ is immaterial.  Hence, it is arbitrarily  assumed that  for all~$\vect{\theta} \in \left\lbrace \vect{\nu} \in \supp Q:  \mathsf{L}_{\vect{z}}\left( \vect{\nu}\right) < \delta_{Q,\vect{z}}^{\star} \right\rbrace$, it holds that 
\begin{equation}
\label{EqWaitingToTakeOff}
\frac{\mathrm{d}P^{\left(Q, \lambda\right)}_{\vect{\Theta} | \vect{Z} = \vect{z}}}{\mathrm{d}Q} \left( \vect{\theta} \right) = 0.
\end{equation}

This completes the third part and completes the proof.

 \section{Proof of Lemma~\ref{TheoSetLStarConcentration}}\label{AppProofTheoSetLStarConcentration}


Consider the following partition of the set~$\set{M}$ formed by the sets
\begin{subequations}\label{EqSetsA}
\begin{IEEEeqnarray}{rcl}
\label{EqSetA00987}
\set{A}_0  & \triangleq &\left\lbrace \vect{\theta} \in \set{M}: \mathsf{L}_{\vect{z}}\left( \vect{\theta}\right)  =   \delta_{Q,\vect{z}}^{\star} \right \rbrace,\\
\label{EqSetA10987}
\set{A}_1& \triangleq &\left\lbrace \vect{\theta} \in \set{M}: \mathsf{L}_{\vect{z}}\left( \vect{\theta}\right)  <    \delta_{Q,\vect{z}}^{\star}  \right \rbrace, \mbox{ and }\\
\label{EqSetA20987}
\set{A}_2 & \triangleq & \left\lbrace \vect{\theta} \in \set{M}: \mathsf{L}_{\vect{z}}\left( \vect{\theta}\right)  >  \delta_{Q,\vect{z}}^{\star} \right \rbrace,
\end{IEEEeqnarray}
\end{subequations}
with~$\delta_{Q,\vect{z}}^{\star}$ in~\eqref{EqDeltaStar} and the function~$\mathsf{L}_{\vect{z}}$ in~\eqref{EqLxy}.
Note that~$\set{A}_0 =\set{L}^{\star}_{Q,\vect{z}}$, with~$\set{L}^{\star}_{Q,\vect{z}}$ in~\eqref{EqSetLStar} and 
 \begin{IEEEeqnarray}{rcl}
 \label{EqLeila12543}
1 & = & P^{\left(Q, \lambda\right)}_{\vect{\Theta} | \vect{Z} = \vect{z}}\left( \set{A}_{0} \right)  + P^{\left(Q, \lambda\right)}_{\vect{\Theta} | \vect{Z} = \vect{z}}\left( \set{A}_{1} \right) + P^{\left(Q, \lambda\right)}_{\vect{\Theta} | \vect{Z} = \vect{z}}\left( \set{A}_{2} \right) \\
 \label{EqLeila1253}
 & = & P^{\left(Q, \lambda\right)}_{\vect{\Theta} | \vect{Z} = \vect{z}}\left( \set{A}_{0} \right)  
 + P^{\left(Q, \lambda\right)}_{\vect{\Theta} | \vect{Z} = \vect{z}}\left( \set{A}_{2} \right) \\
 \label{EqLeila143}
 & = & 
P^{\left(Q, \lambda\right)}_{\vect{\Theta} | \vect{Z} = \vect{z}}\left( \set{A}_{0} \right) 
+ \int_{\set{A}_2}  \mathrm{d}P^{\left(Q, \lambda\right)}_{\vect{\Theta} | \vect{Z} = \vect{z}}  (\vect{\theta}),
 \end{IEEEeqnarray}
 where, 
 the equality in~\eqref{EqLeila1253} follows from noticing that 
~$P^{\left(Q, \lambda\right)}_{\vect{\Theta} | \vect{Z} = \vect{z}}  (\set{A}_{1}) = 0$, which follows from the definition of~$\delta_{Q,\vect{z}}^{\star}$ in~\eqref{EqDeltaStar} and the fact that the probability measure~$P^{\left(Q, \lambda\right)}_{\vect{\Theta} | \vect{Z} = \vect{z}}$ is absolutely continuous with respect to the measure~$Q$.

The above implies that 
 \begin{IEEEeqnarray}{rcl}
\label{EqNetflix239875}
1& = & 
 \lim_{\lambda \rightarrow 0^{+}} P^{\left(Q, \lambda\right)}_{\vect{\Theta} | \vect{Z} = \vect{z}}\left( \set{A}_{0} \right)   
+ \lim_{\lambda \rightarrow 0^{+}} \int_{\set{A}_2}   \frac{\mathrm{d}P^{\left(Q, \lambda\right)}_{\vect{\Theta} | \vect{Z} = \vect{z}}}{\mathrm{d}Q} \left( \vect{\theta} \right)   \mathrm{d}Q (\vect{\theta})\supersqueezeequ
\IEEEeqnarraynumspace\\
\label{EqNetflix235}
& = & \lim_{\lambda \rightarrow 0^{+}} P^{\left(Q, \lambda\right)}_{\vect{\Theta} | \vect{Z} = \vect{z}}\left( \set{A}_{0} \right)   
+ \int_{\set{A}_2} \lim_{\lambda \rightarrow 0^{+}}   \frac{\mathrm{d}P^{\left(Q, \lambda\right)}_{\vect{\Theta} | \vect{Z} = \vect{z}}}{\mathrm{d}Q} \left( \vect{\theta} \right)   \mathrm{d}Q (\vect{\theta})\supersqueezeequ
\IEEEeqnarraynumspace\\
\label{EqNetflix2375}
& = & \lim_{\lambda \rightarrow 0^{+}} P^{\left(Q, \lambda\right)}_{\vect{\Theta} | \vect{Z} = \vect{z}}\left( \set{A}_{0} \right),
 \end{IEEEeqnarray}
 where,  
 the equality in~\eqref{EqNetflix235} follows from the dominated convergence theorem~\cite[Theorem~$1.6.9$]{ash2000probability}, given that the Randon-Nikodym derivative~$\frac{\mathrm{d}P^{\left(Q, \lambda\right)}_{\vect{\Theta} | \vect{Z} = \vect{z}}}{\mathrm{d}Q}$ is positive and finite (Lemma~\ref{CorPositive}); and 
the  inequality in~\eqref{EqNetflix2375} holds from the fact that for all~$\vect{\theta} \in \set{A}_2$, it holds that~$ \lim_{\lambda \rightarrow 0^{+}}  \frac{\mathrm{d}P^{\left(Q, \lambda\right)}_{\vect{\Theta} | \vect{Z} = \vect{z}}}{\mathrm{d}Q}\left( \vect{\theta} \right) = 0$ (Lemma~\ref{CorAssymptoticZero}).
Hence, it finally holds that
\begin{equation}
 \lim_{\lambda \rightarrow 0^{+}} P^{\left(Q, \lambda\right)}_{\vect{\Theta} | \vect{Z} = \vect{z}}\left( \set{L}^{\star}_{Q,\vect{z}}  \right) = 1,
\end{equation}
which completes the proof.

\section{Proof of Lemma~\ref{LemmaHotelAmourA}}\label{AppProofLemmaHotelAmourA}

The proof is presented in two parts. The first part shows that if  for all~$\delta \in \left(  \rho^{\star}, +\infty \right)$, the inequality in~\eqref{EqCrazyEight} holds, then,~$Q$ is coherent. The second part shows that if $Q$ is not coherent, then there exists a~$\delta \in \left(  \rho^{\star}, +\infty \right)$ such that 
\begin{IEEEeqnarray}{ccl}
\label{EqCrazyEightZero}
P^{\left(Q, \lambda\right)}_{\vect{\Theta} | \vect{Z} = \vect{z}}\left( \set{L}_{\vect{z}} \left( \delta \right)  \right) & =  & 0. 
\end{IEEEeqnarray}

The first part is as follows. Note that for all~$\delta \in \left(  \rho^{\star}, +\infty \right)$ and for all~$\vect{\theta} \in \set{L}_{\vect{z}} \left( \delta \right) \cap \supp Q$, it holds from Lemma~\ref{CorPositive} that
 \begin{equation}
 \label{EqElCamino}
\frac{\mathrm{d}P^{\left(Q, \lambda\right)}_{\vect{\Theta} | \vect{Z} = \vect{z}}}{\mathrm{d}Q}\left( \vect{\theta}\right) > 0.
    \end{equation}
Hence, if for all~$\delta \in \left(  \rho^{\star}, +\infty \right)$, the inequality in~\eqref{EqCrazyEight} holds, then
\begin{IEEEeqnarray}{CCL}
\label{EqCrazyEightZeroZero}
0 &< &P^{\left(Q, \lambda\right)}_{\vect{\Theta} | \vect{Z} = \vect{z}}\left( \set{L}_{\vect{z}} \left( \delta \right)  \right)  \\
& = & \int_{\set{L}_{\vect{z}} \left( \delta \right) } \d P^{\left(Q, \lambda\right)}_{\vect{\Theta} | \vect{Z} = \vect{z}}\left( \vect{\theta} \right)\\
& = & \int_{\set{L}_{\vect{z}} \left( \delta \right) } \frac{\mathrm{d}P^{\left(Q, \lambda\right)}_{\vect{\Theta} | \vect{Z} = \vect{z}}}{\mathrm{d}Q}\left( \vect{\theta}\right)  \d Q\left( \vect{\theta} \right),
\end{IEEEeqnarray}
which, together with~\eqref{EqElCamino}, implies that for all~$\delta \in \left(  \rho^{\star}, +\infty \right)$,~$Q\left( \set{L}_{\vect{z}} \left( \delta \right)  \right) > 0$. Hence,~$Q$ is coherent.

The second part is as follows. Assume that $Q$ is not coherent. Then, there exists a $\delta \in \left(  \rho^{\star}, +\infty \right)$ such that~$Q\left( \set{L}_{\vect{z}} \left( \delta \right)  \right) = 0$. Hence, from the fact that~$P^{\left(Q, \lambda\right)}_{\vect{\Theta} | \vect{Z} = \vect{z}}$ is absolutely continuous with respect to $Q$, it follows that~$P^{\left(Q, \lambda\right)}_{\vect{\Theta} | \vect{Z} = \vect{z}}\left( \set{L}_{\vect{z}} \left( \delta \right)  \right)$ $=$ $0$.
 This completes the proof.

\section{Proof of Theorem~\ref{TheoremSensitivityB}}\label{AppProofTheoremSensitivityB}

The optimization problem in~\eqref{EqImprovementp235a8} can be re-written in terms of the Radon-Nikodym derivative of the optimization measure $P$ with respect to the measure $P^{\left(Q, \lambda\right)}_{\vect{\Theta} | \vect{Z} = \vect{z}}$, denoted by $\frac{\mathrm{d} P}{\mathrm{d} P^{\left(Q, \lambda\right)}_{\vect{\Theta} | \vect{Z} = \vect{z}}}: \set{M} \to [0, +\infty)$, which yields: 
\begin{subequations}\label{EqOptimRouge}
\begin{IEEEeqnarray}{lcl}
\label{EqOptimRougeA}
& & \min_{P \in \triangle_{Q}\Bormeaspace{\set{M}}}  \displaystyle\int \mathsf{L}_{\vect{z}}\left( \vect{\nu}\right)\frac{\mathrm{d}P}{\mathrm{d}P^{\left(Q, \lambda\right)}_{\vect{\Theta} | \vect{Z} = \vect{z}} } \left( \vect{\nu} \right)  \mathrm{d}P^{\left(Q, \lambda\right)}_{\vect{\Theta} | \vect{Z} = \vect{z}} \left( \vect{\nu} \right), \squeezeequ\\
\nonumber
&  & \text{subject to: }\\
\label{EqOptimRougeB}
& &   \int \frac{\mathrm{d}P}{\mathrm{d}P^{\left(Q, \lambda\right)}_{\vect{\Theta} | \vect{Z} = \vect{z}} } \left( \vect{\nu} \right)  \log \left(\hspace{-1ex}\frac{\mathrm{d}P}{\mathrm{d}P^{\left(Q, \lambda\right)}_{\vect{\Theta} | \vect{Z} = \vect{z}} } \left( \vect{\nu} \right) \hspace{-1ex}\right)  \mathrm{d}P^{\left(Q, \lambda\right)}_{\vect{\Theta} | \vect{Z} = \vect{z}} \left( \vect{\nu} \right)  \leqslant c,  \mbox{ and }\supersqueezeequ\IEEEeqnarraynumspace \\
\label{EqConstrainth}
 &  & \int\frac{\mathrm{d}P}{\mathrm{d}P^{\left(Q, \lambda\right)}_{\vect{\Theta} | \vect{Z} = \vect{z}} } \left( \vect{\theta} \right)  \mathrm{d} P^{\left(Q, \lambda\right)}_{\vect{\Theta} | \vect{Z} = \vect{z}} (\vect{\theta}) = 1.
\end{IEEEeqnarray}
\end{subequations}
The remainder of the proof focuses on the problem in which the optimization is over  the function   $\frac{\mathrm{d} P}{\mathrm{d} P^{\left(Q, \lambda\right)}_{\vect{\Theta} | \vect{Z} = \vect{z}}}$ instead of the measure $P$. This is due to the fact that for all $P \in \bigtriangleup_{Q}\left( \set{M} \right)$, the Radon-Nikodym derivate  $\frac{\mathrm{d} P}{\mathrm{d} P^{\left(Q, \lambda\right)}_{\vect{\Theta} | \vect{Z} = \vect{z}}}$ is unique up to sets of zero measure with respect to the measure $P^{\left(Q, \lambda\right)}_{\vect{\Theta} | \vect{Z} = \vect{z}}$.
Let $\mathscr{M}$ be the set of measurable functions $\set{M} \to \reals$ with respect to the measurable spaces $\left( \set{M}, \BorSigma{\set{M}} \right)$ and $\left( \reals, \BorSigma{\reals} \right)$ that are absolutely integrable with respect to $P^{\left(Q, \lambda\right)}_{\vect{\Theta} | \vect{Z} = \vect{z}}$.
That is, for all $\hat{g} \in \mathscr{M}$, it holds that
\begin{IEEEeqnarray}{rcl}
\label{EqValioLaPenaJ}
\int \abs{ \hat{g} (\vect{\theta})} \mathrm{d} P^{\left(Q, \lambda\right)}_{\vect{\Theta} | \vect{Z} = \vect{z}}(\vect{\theta}) & < & \infty.
\end{IEEEeqnarray}

Hence, the optimization problem of interest is:
\begin{subequations}\label{EqDaliJustArrivedHome}
\begin{IEEEeqnarray}{lCl}
\label{EqDaliJustArrivedHomeA}
& \min_{g \in \mathscr{M}}  & \displaystyle\int \mathsf{L}_{\vect{z}}\left( \vect{\nu}\right)g\left( \vect{\nu} \right)  \mathrm{d}P^{\left(Q, \lambda\right)}_{\vect{\Theta} | \vect{Z} = \vect{z}} \left( \vect{\nu} \right) \\
\label{EqDaliJustArrivedHomeB}
& \text{s.t: } &    \int g \left( \vect{\nu} \right)  \log \left(g \left( \vect{\nu} \right)  \right)  \mathrm{d}P^{\left(Q, \lambda\right)}_{\vect{\Theta} | \vect{Z} = \vect{z}} \left( \vect{\nu} \right)  \leqslant c,  \mbox{ and }\IEEEeqnarraynumspace \\
\label{EqDaliJustArrivedHomeC}
 &  & \int g \left( \vect{\theta} \right)  \mathrm{d} P^{\left(Q, \lambda\right)}_{\vect{\Theta} | \vect{Z} = \vect{z}} (\vect{\theta}) = 1.
\end{IEEEeqnarray}
\end{subequations}
The Lagrangian of the optimization problem in~\eqref{EqDaliJustArrivedHome} is a functional~$L: \mathscr{M} \times [0, +\infty)^2 \rightarrow \reals$ of the form  
\begin{IEEEeqnarray}{rcl}
\nonumber
L\left(g, 
\alpha, \beta \right) &=& \displaystyle\int \mathsf{L}_{\vect{z}}\left( \vect{\nu}\right)g \left( \vect{\nu} \right)  \mathrm{d}P^{\left(Q, \lambda\right)}_{\vect{\Theta} | \vect{Z} = \vect{z}} \left( \vect{\nu} \right) \\
\nonumber
& &+ \alpha \left( \displaystyle\int g \left( \vect{\nu} \right)  \log \left( g \left( \vect{\nu} \right) \right)  \mathrm{d}P^{\left(Q, \lambda\right)}_{\vect{\Theta} | \vect{Z} = \vect{z}} \left( \vect{\nu} \right)  - c \right) \\
\label{EqFunctionalL865d1}
   & &    + \beta \left(\displaystyle\int g \left( \vect{\nu} \right) \mathrm{d}P^{\left(Q, \lambda\right)}_{\vect{\Theta} | \vect{Z} = \vect{z}} \left( \vect{\nu} \right)   -1 \right),
\end{IEEEeqnarray}
where~the reals~$\alpha$ and~$\beta$ are both nonnegative and act as Lagrangian multipliers due to the constraints~\eqref{EqDaliJustArrivedHomeB} and~\eqref{EqDaliJustArrivedHomeC}, respectively.  

Let~$h:\set{M} \rightarrow \reals$ be a function in~$\mathscr{M}$. The Gateaux differential of the functional~$L$ in~\eqref{EqFunctionalL865d1} at~$\left( g, \alpha, \beta \right) \in \mathscr{M}\times [0, +\infty)^2$ in the direction of~$h$, if it exists, is
\begin{equation}
\label{EqNecessaryCondtion976d}
    \partial L\left(g, \alpha, \beta; h \right) \triangleq \left. \frac{\mathrm{d}}{\mathrm{d}\gamma}  r(\gamma) \right|_{\gamma = 0} ,
\end{equation}
where the real function~$r: \reals \rightarrow \reals$ is such that for all~$\gamma \in \left( -\epsilon, \epsilon \right)$, with some $\epsilon > 0$, satisfies
\begin{IEEEeqnarray}{rcl}
\nonumber
& &    r(\gamma)  \\
\nonumber
& = & \displaystyle\int \mathsf{L}_{\vect{z}}\left( \vect{\nu}\right) \left( g \left( \vect{\nu} \right) + \gamma h\left( \vect{\nu} \right) \right) \mathrm{d}P^{\left(Q, \lambda\right)}_{\vect{\Theta} | \vect{Z} = \vect{z}} \left( \vect{\nu} \right)  \\
    \nonumber
    & + & \alpha \left( \displaystyle\int \left( g \left( \vect{\nu} \right) + \gamma h\left( \vect{\nu} \right) \right)  \log \left( g \left( \vect{\nu} \right) + \gamma h\left( \vect{\nu} \right) \right)  \mathrm{d}P^{\left(Q, \lambda\right)}_{\vect{\Theta} | \vect{Z} = \vect{z}} \left( \vect{\nu} \right)  - c \right) \supersqueezeequ\\
\label{EqFunctionalL07w087f}
   & + &   \beta \left(\displaystyle\int \left( g \left( \vect{\nu} \right) + \gamma h\left( \vect{\nu} \right) \right)  \mathrm{d}P^{\left(Q, \lambda\right)}_{\vect{\Theta} | \vect{Z} = \vect{z}} \left( \vect{\nu} \right)   -1 \right).
\end{IEEEeqnarray}
The proof continues under the assumption that the functions $g$ and $h$ are such that the Gateaux differential in \eqref{EqNecessaryCondtion976d} exists. 
That is, the function $r$ in \eqref{EqFunctionalL07w087f} is differentiable in $\left( -\epsilon, \epsilon \right)$, with some $\epsilon > 0$. 
Using the same arguments as in the proof of Theorem~\ref{TheoremOptimalModel}, it follows that the derivative of the real function~$r$ in~\eqref{EqFunctionalL07w087f} is
\begin{IEEEeqnarray}{rcl}
\nonumber
& & \frac{\mathrm{d}}{\mathrm{d}\gamma} r(\gamma) \\
\nonumber
& = &  \displaystyle\int \mathsf{L}_{\vect{z}} h\left( \vect{\nu}\right) \mathrm{d}P^{\left(Q, \lambda\right)}_{\vect{\Theta} | \vect{Z} = \vect{z}} \left( \vect{\nu} \right) + \alpha \displaystyle\int h\left( \vect{\nu} \right) \mathrm{d}P^{\left(Q, \lambda\right)}_{\vect{\Theta} | \vect{Z} = \vect{z}} \left( \vect{\nu} \right)  \\
    \nonumber
    & +&    \alpha \displaystyle\int   h\left( \vect{\nu} \right)  \log \left( g \left( \vect{\nu} \right) + \gamma h\left( \vect{\nu} \right) \right)  \mathrm{d}P^{\left(Q, \lambda\right)}_{\vect{\Theta} | \vect{Z} = \vect{z}} \left( \vect{\nu} \right)   \\
\label{EqFunctionalL07w087fa}
   & + &  \beta \displaystyle\int h\left( \vect{\nu} \right) \mathrm{d}P^{\left(Q, \lambda\right)}_{\vect{\Theta} | \vect{Z} = \vect{z}} \left( \vect{\nu} \right) .
\end{IEEEeqnarray}
From~\eqref{EqNecessaryCondtion976d} and~\eqref{EqFunctionalL07w087fa}, it follows that
\begin{IEEEeqnarray}{rcl}
 \nonumber
 &&    \partial L\left(g, \alpha, \beta; h \right) \\
 \label{EqGateauxDiff52043058}
 & = &  \displaystyle\int h\left( \vect{\nu}\right) \left( \mathsf{L}_{\vect{z}}\left( \vect{\nu}\right) + \alpha \left(1 + \log g \left( \vect{\nu} \right) \right) + \beta \right)     \mathrm{d}P^{\left(Q, \lambda\right)}_{\vect{\Theta} | \vect{Z} = \vect{z}} \left( \vect{\nu} \right). \IEEEeqnarraynumspace
\end{IEEEeqnarray}

From \cite[Theorem~$1$, page~$217$]{luenberger1997bookOptimization}, it holds that a necessary condition for the functional~$L$ in~\eqref{EqFunctionalL865d1} to have a minimum at~$\left( g, \alpha, \beta \right)  \in \mathscr{M} \times [0, +\infty)^2$ is that for all functions~$h~\in~\mathscr{M}$,
\begin{equation}\label{EqConditionh}
\partial L\left(g, \alpha, \beta; h \right)  = 0, 
\end{equation}
which implies that for all~$\vect{\nu} \in \set{M}$,
\begin{equation}
 \mathsf{L}_{\vect{z}}\left( \vect{\nu}\right) + \alpha \left(1 + \log g \left( \vect{\nu} \right) \right) + \beta  = 0.
\end{equation}
Thus,
\begin{equation}
\label{EqKeyNorthStar}
  g\left( \vect{\nu} \right) =  \exp\left( -\frac{\mathsf{L}_{\vect{z}}\left( \vect{\nu}\right) }{\alpha}\right) \exp\left( -\frac{\beta + \alpha}{\alpha}\right),
\end{equation}
where~$\alpha$ and~$\beta$ are chosen to satisfy their corresponding constraints with equality.
Denote by~$P^{\star}$ the solution of the optimization problem in~\eqref{EqImprovementp235a8}. Hence, from~\eqref{EqKeyNorthStar}, it follows that
\begin{equation}
 \frac{\mathrm{d}P^{\star}}{\mathrm{d}P^{\left(Q, \lambda\right)}_{\vect{\Theta} | \vect{Z} = \vect{z}}} \left( \vect{\nu} \right)  = \frac{ \exp\left( -\frac{\mathsf{L}_{\vect{z}}\left( \vect{\nu}\right) }{\alpha}\right) }{\displaystyle\int \exp\left( -\frac{\mathsf{L}_{\vect{z}}\left( \vect{\theta}\right) }{\alpha}\right) \mathrm{d} P^{\left(Q, \lambda\right)}_{\vect{\Theta} | \vect{Z} = \vect{z}} \left( \vect{\theta}\right)},
\end{equation}
where~$\alpha$ is chosen to satisfy 
\begin{IEEEeqnarray}{rcl}
\label{EqRecycleThis}
  D\left( P^{\star} \|P^{\left(Q, \lambda\right)}_{\vect{\Theta} | \vect{Z} = \vect{z}} \right) = c.
\end{IEEEeqnarray}
From Lemma~\ref{LemmaMutualAC}, it follows that the probability measure~$P^{\star}$ and the~$\sigma$-finite measure~$Q$ satisfy, 
 \begin{IEEEeqnarray}{rcl}
 \frac{\mathrm{d}P^{\star}}{\mathrm{d}Q} \left( \vect{\nu} \right)  & = & \frac{\mathrm{d}P^{\star}}{\mathrm{d}P^{\left(Q, \lambda\right)}_{\vect{\Theta} | \vect{Z} = \vect{z}}} \left( \vect{\nu} \right)  \frac{\mathrm{d}P^{\left(Q, \lambda\right)}_{\vect{\Theta} | \vect{Z} = \vect{z}}}{\mathrm{d}Q} \left( \vect{\nu} \right) \\
\nonumber
 & = & \left(\frac{ \exp\left( -\frac{\mathsf{L}_{\vect{z}}\left( \vect{\nu}\right) }{\alpha}\right) }{\displaystyle\int \exp\left( -\frac{\mathsf{L}_{\vect{z}}\left( \vect{\theta}\right) }{\alpha}\right) \mathrm{d} P^{\left(Q, \lambda\right)}_{\vect{\Theta} | \vect{Z} = \vect{z}} \left( \vect{\theta}\right)} \right) \\
 & & \left(\frac{ \exp\left( -\frac{\mathsf{L}_{\vect{z}}\left( \vect{\nu}\right) }{\lambda}\right) }{\displaystyle\int \exp\left( -\frac{\mathsf{L}_{\vect{z}}\left( \vect{\theta}\right) }{\lambda}\right) \mathrm{d} Q \left( \vect{\theta}\right)} \right) \IEEEeqnarraynumspace \\
\nonumber
 & = & \left(\frac{ \exp\left( -\frac{\mathsf{L}_{\vect{z}}\left( \vect{\nu}\right) }{\alpha}\right) }{\displaystyle\int \frac{ \exp\left( -\frac{\mathsf{L}_{\vect{z}}\left( \vect{\theta}\right) }{\alpha}\right) \exp\left( -\frac{\mathsf{L}_{\vect{z}}\left( \vect{\theta}\right) }{\lambda}\right) }{\displaystyle\int \exp\left( -\frac{\mathsf{L}_{\vect{z}}\left( \vect{\alpha}\right) }{\lambda}\right) \mathrm{d} Q \left( \vect{\alpha}\right)}    \mathrm{d} Q \left( \vect{\theta}\right)} \right) \\
 & &  \left(\frac{ \exp\left( -\frac{\mathsf{L}_{\vect{z}}\left( \vect{\nu}\right) }{\lambda}\right) }{\displaystyle\int \exp\left( -\frac{\mathsf{L}_{\vect{z}}\left( \vect{\theta}\right) }{\lambda}\right) \mathrm{d} Q \left( \vect{\theta}\right)} \right) \IEEEeqnarraynumspace\\
  & = & \frac{ \exp\left( -\left( \frac{1}{\alpha} + \frac{1}{\lambda}\right) \mathsf{L}_{\vect{z}}\left( \vect{\nu}\right)  \right) }{\displaystyle\int \exp\left(-\left( \frac{1}{\alpha} + \frac{1}{\lambda}\right) \mathsf{L}_{\vect{z}}\left( \vect{\nu}\right)  \right) \mathrm{d} Q\left( \vect{\theta}\right)},
\end{IEEEeqnarray}
which implies that~$P^{\star}$ is a Gibbs probability measure on~$\Bormeaspace{\set{M}}$, with energy function~$\mathsf{L}_{\vect{z}}$, reference measure~$Q$,  and regularization parameter~$\frac{1}{\frac{1}{\alpha} + \frac{1}{\lambda}}$, where~$\alpha$ is chosen to satisfy~\eqref{EqRecycleThis}.  
Let the positive real~$\omega$ be~$\omega \triangleq \frac{\alpha \lambda}{\alpha + \lambda}$ and note that~$\omega \in \left(0, \lambda\right]$ and satisfies~$D\left( P^{\left(Q,\omega\right)}_{\vect{\Theta}| \vect{Z} = \vect{z}} \left( \vect{\nu} \right) \| P^{\left(Q, \lambda\right)}_{\vect{\Theta} | \vect{Z} = \vect{z}} \right) = c$.
The proof ends by verifying that the objective function in~\eqref{EqFunctionalL865d1} is strictly convex, and thus, the measure~$P^{\left(Q,\omega \right)}_{\vect{\Theta}| \vect{Z} = \vect{z}}$ is the unique minimizer.  
This completes the proof.

\section{Proof of Lemma~\ref{LemmaContinuousK}} \label{AppProofLemmaContinuousK}

Note that for all $(\lambda_1, \lambda_2) \in \left\lbrace x \in \reals : K_{Q, \vect{z}}(x) < +\infty \right\rbrace^2$, such that $\lambda_1 > 
\lambda_2$, it follows that for all $\vect{\theta} \in \supp Q$, the inequality $\exp\left( \lambda_2 \; \mathsf{L}_{\vect{z}}\left(\vect{\theta}\right) \right) \leqslant \exp\left( \lambda_1 \; \mathsf{L}_{\vect{z}}\left(\vect{\theta}\right) \right)$ holds. This implies that $K_{Q, \vect{z}}\left(\lambda_2 \right) \leqslant K_{Q, \vect{z}}\left(\lambda_1\right)  < +\infty $, which proves that the function is nondecreasing. 

The proof of continuity of the function $K_{Q, \vect{z}}$ follows from observing that for all $\alpha \in \left\lbrace x \in \reals : K_{Q, \vect{z}}(x) < +\infty \right\rbrace$, it holds that
\begin{IEEEeqnarray}{rcl}
\label{EqElCabelloDeCelesteA}
\lim_{t \to \alpha} K_{Q, \vect{z}}(t) & = & \lim_{t \to \alpha} \log\left( \int \exp\left( t \; \mathsf{L}_{\vect{z}}\left(\vect{\theta}\right)  \right) \mathrm{d}Q(\vect{\theta}) \right)\IEEEeqnarraynumspace\\
\label{EqElCabelloDeCelesteB}
& = &  \log\left( \lim_{t \to \alpha} \int \exp\left( t \; \mathsf{L}_{\vect{z}}\left(\vect{\theta}\right)  \right) \mathrm{d}Q(\vect{\theta}) \right)\IEEEeqnarraynumspace\\
\label{EqElCabelloDeCelesteC}
& = & \log\left(\int \lim_{t \to \alpha}  \exp\left( t \; \mathsf{L}_{\vect{z}}\left(\vect{\theta}\right)  \right) \mathrm{d}Q(\vect{\theta}) \right)\IEEEeqnarraynumspace\\
\label{EqElCabelloDeCelesteD}
& = & \log\left(\int \exp\left( \alpha \; \mathsf{L}_{\vect{z}}\left(\vect{\theta}\right)  \right) \mathrm{d}Q(\vect{\theta}) \right)\IEEEeqnarraynumspace\\
\label{EqElCabelloDeCelesteE}
 & = &  K_{Q, \vect{z}}(\alpha),
\end{IEEEeqnarray}
where \eqref{EqElCabelloDeCelesteB} and \eqref{EqElCabelloDeCelesteD} follow from the fact that both the logarithmic and exponential functions are continuous; and the equality in \eqref{EqElCabelloDeCelesteC} follows from the monotone convergence theorem~\cite[Theorem~$1.6.2$]{ash2000probability}. This shows that the function $K_{Q, \vect{z}}$ is continuous in $\left\lbrace x \in \reals : K_{Q, \vect{z}}(x) < +\infty \right\rbrace$.

The proof of differentiability follows by considering the transport of the~$\sigma$-finite measure~$Q$ in~\eqref{EqK} from the measure space~$\Bormeaspace{\set{M}}$ to the measure space~$\Bormeaspace{\left[ 0, +\infty \right)}$ through the function~$\mathsf{L}_{\vect{z}}$ in~\eqref{EqLxy}.  Denote the resulting measure in~$\Bormeaspace{\left[ 0, +\infty \right)}$ by~$P$. 
More specifically, for all~$\mathcal{A} \in \mathscr{B}\left( \left[0, +\infty \right) \right)$, it holds that~$P\left( \set{A}\right) = Q\left(\left\lbrace \vect{\theta} \in\set{M}: \mathsf{L}_{\vect{z}}\left(\vect{\theta} \right) \in \set{A}  \right\rbrace \right)$.
Hence, the function~$K_{Q, \vect{z}}$ satisfies for all~$t \in \left\lbrace \nu \in \reals : K_{Q, \vect{z}}(\nu) < +\infty \right\rbrace$,
\begin{IEEEeqnarray}{rcl}
\label{EqKV}
K_{Q, \vect{z}}\left(t \right) & = &  \log\left( \int \exp\left( t \; \mathsf{L}_{\vect{z}}\left(\vect{\theta}\right)  \right) \mathrm{d}Q(\vect{\theta}) \right) \\
\label{EqAnimalMalHerido}
& = & \log\left( \int \exp\left( t \, w \right) \mathrm{d}P(w) \right),
\end{IEEEeqnarray} 
where the equality \eqref{EqAnimalMalHerido} follows from \cite[Theorem~$1.6.12$]{ash2000probability}.
Denote by~$\phi$ the Laplace transform of the measure~$P$. That is, for all~$t \in \left\lbrace x \in \reals : K_{Q, \vect{z}}(x) < +\infty \right\rbrace$,
\begin{equation}
\phi(t) =  \int  \exp\left( t \, v \right) \mathrm{d}P(v).
\end{equation}
Hence, $\phi(t) = \exp\left( K_{Q, \vect{z}}\left( t \right) \right)$. 
From  \cite[Theorem~$1a$ (page 439)]{FellerBookII}, it follows that the function~$\phi$  has derivatives of all orders in~$\left\lbrace x \in \reals : K_{Q, \vect{z}}(x) < +\infty \right\rbrace$, and thus, so does the function~$K_{Q, \vect{z}}$ in the interior of~$\left\lbrace x \in \reals : K_{Q, \vect{z}}(x) < +\infty \right\rbrace$.  This completes the proof.

\section{Proof of Lemma~\ref{LemmaConvexK}}\label{AppProofLemmaConvexK}
Let~$(\gamma_1, \gamma_2) \in \reals^2$, with~$\gamma_1 \neq \gamma_2$ and~$\alpha \in [0,1]$ be fixed. Assume that~$K_{Q, \vect{z}}\left( \gamma_1 \right) < +\infty$ and~$K_{Q, \vect{z}}\left( \gamma_2 \right) < +\infty$.
Then, for all~$\alpha \in (0,1)$, the following holds
\begin{IEEEeqnarray}{rcl}
\nonumber
& & \alpha K_{Q, \vect{z}}\left( \gamma_1 \right) + (1 - \alpha ) K_{Q, \vect{z}}\left( \gamma_2 \right)  \\
\nonumber
& = &  \alpha   \log\left( \int \exp\left( \gamma_1 \, \mathsf{L}_{\vect{z}}\left(\vect{\theta}\right)  \right) \d Q(\vect{\theta}) \right) \\
& & +  (1- \alpha)  \log\left( \int \exp\left( \gamma_2 \, \mathsf{L}_{\vect{z}}\left(\vect{\theta}\right)  \right) \d Q(\vect{\theta}) \right)\\
\nonumber
& = &   \log\left( \left( \int \exp\left( \gamma_1 \, \mathsf{L}_{\vect{z}}\left(\vect{\theta}\right)  \right) \d Q(\vect{\theta}) \right)^{\alpha}\right)  \\
& & + \log\left( \left( \int \exp\left( \gamma_2 \, \mathsf{L}_{\vect{z}}\left(\vect{\theta}\right)  \right) \d Q(\vect{\theta}) \right)^{(1- \alpha) }\right)
\end{IEEEeqnarray}
\begin{IEEEeqnarray}{rcl}
\nonumber
& = &   \log\bigg( \left( \int \exp\left( \gamma_1 \, \mathsf{L}_{\vect{z}}\left(\vect{\theta}\right)  \right) \d Q(\vect{\theta}) \right)^{\alpha}\\
& &   \left( \int \exp\left( \gamma_2 \, \mathsf{L}_{\vect{z}}\left(\vect{\theta}\right)  \right) \d Q(\vect{\theta}) \right)^{(1- \alpha) }\bigg)\supersqueezeequ
\IEEEeqnarraynumspace\\
   \nonumber
& = &    \log\Bigg( \left( \int \exp\left( \gamma_1 \alpha \mathsf{L}_{\vect{z}}\left(\vect{\theta}\right)  \right)^{p} \d Q(\vect{\theta}) \right)^{\frac{1}{p}}\\
  \label{EqProofLemmaCGFPZ1afa}
  & &  \left( \int \exp\left( \gamma_2 (1-\alpha) \mathsf{L}_{\vect{z}}\left(\vect{\theta}\right)  \right)^{q} \d Q(\vect{\theta}) \right)^{\frac{1}{q}}\Bigg)\\
   \label{EqProofLemmaCGFPZ1}
&\geqslant&      \log\left( \int \exp\left(  \gamma_1 \alpha \mathsf{L}_{\vect{z}}\left(\vect{\theta}\right)  \right) \exp\left( \gamma_2 (1-\alpha) \mathsf{L}_{\vect{z}}\left(\vect{\theta}\right)  \right) \d Q(\vect{\theta}) \right) \supersqueezeequ
\IEEEeqnarraynumspace \\
   & = &     \log\left( \int \exp\bigg( \left( \gamma_1\alpha + \gamma_2(1-\alpha) \right) \mathsf{L}_{\vect{z}}\left(\vect{\theta}\right) \bigg) \d Q(\vect{\theta}) \right)\\
&  = &    K_{Q, \vect{z}}\left( \gamma_1\alpha + \gamma_2 (1 - \alpha )  \right),
\end{IEEEeqnarray}
where 
the inequality in~\eqref{EqProofLemmaCGFPZ1afa} follows with~$\alpha \triangleq \frac{1}{p}$ and~$1 - \alpha \triangleq \frac{1}{q}$; 
the inequality in~\eqref{EqProofLemmaCGFPZ1} follows from H\"{o}lder's inequality. Hence,  equality in~\eqref{EqProofLemmaCGFPZ1} holds if and only if there exist two constants~$\beta_1$ and~$\beta_2$, not simultaneously equal to zero, such that the set
\begin{IEEEeqnarray}{rCl}
\nonumber
\set{A} 
& \triangleq & 
\left\lbrace \vect{\theta} \in \set{M}:  \beta_1 \exp\left(  \gamma_1 \mathsf{L}_{\vect{z}}\left(\vect{\theta}\right)  \right) = \beta_2 \exp\left(  \gamma_2 \mathsf{L}_{\vect{z}}\left(\vect{\theta}\right)  \right) \right\rbrace
\IEEEeqnarraynumspace\\
& = & \left\lbrace \vect{\theta} \in \set{M}:  \exp\left(  \left( \gamma_1 - \gamma_2\right) \mathsf{L}_{\vect{z}}\left(\vect{\theta}\right)  \right) = \frac{\beta_2}{\beta_1} \right\rbrace\\
& = & \left\lbrace \vect{\theta} \in \set{M}:   \mathsf{L}_{\vect{z}}\left(\vect{\theta}\right)  = \frac{\log \frac{\beta_2}{\beta_1}}{\left( \gamma_1 - \gamma_2\right)} \right\rbrace,
\IEEEeqnarraynumspace
\end{IEEEeqnarray}
satisfies $Q\left( \set{A} \right) = 1$.
That is, strict inequality in~\eqref{EqProofLemmaCGFPZ1} holds if and only if the function~$\mathsf{L}_{\vect{z}}$ is separable with respect to the~$\sigma$-finite measure~$Q$.
When~$\alpha = 0$ or~$\alpha = 1$, the proof is trivial.  
This completes the proof.  
  
  \section{Proof of Lemma~\ref{CorDerivatives}}\label{ProofCorDerivatives}

For all~$s \in \set{K}_{Q, \vect{z}}$, with~$\set{K}_{Q, \vect{z}}$ in~\eqref{EqSetKxy}, the equality in~\eqref{EqK1} implies the following,
\begin{IEEEeqnarray}{rcl}
\nonumber
& & K^{(1)}_{Q, \vect{z}} \left(-\frac{1}{s} \right)\\
& = & \frac{\mathrm{d}}{\mathrm{d}t} \log\left( \int  \exp\left( t \, \mathsf{L}_{\vect{z}}\left(\vect{\theta}\right)  \right) \d Q(\vect{\theta}) \right) \Biggr|_{t = -\frac{1}{s}} \\
\label{EqProofCorDerivativesA95}
& = &  \int  \frac{\mathsf{L}_{\vect{z}}\left(\vect{\theta}\right)  \exp\left( t \, \mathsf{L}_{\vect{z}}\left(\vect{\theta}\right)  \right) }{\int  \exp\left( t \, \mathsf{L}_{\vect{z}}\left(\vect{v}\right)  \right) \d Q(\vect{v})}  \d Q(\vect{\theta}) \Biggr|_{t = -\frac{1}{s}} \\
\label{EqProofCorDerivativesA96}
& = & \int  \frac{\mathsf{L}_{\vect{z}}\left(\vect{\theta}\right)  \exp\left( -\frac{1}{s} \, \mathsf{L}_{\vect{z}}\left(\vect{\theta}\right)  \right)}{  \int  \exp\left(- \frac{1}{s} \, \mathsf{L}_{\vect{z}}\left(\vect{v}\right)  \right) \d Q(\vect{v}) }    \d Q(\vect{\theta})   \\
\label{EqProofCorDerivativesA97}
& = & \exp \left( -K_{Q, \vect{z}}\left(-\frac{1}{s} \right) \right) \int  \mathsf{L}_{\vect{z}}\left(\vect{\theta}\right)  \exp\left(- \,\frac{1}{s} \mathsf{L}_{\vect{z}}\left(\vect{\theta}\right)  \right) \d Q(\vect{\theta})
\middlesqueezeequ
\IEEEeqnarraynumspace \\
\label{EqProofCorDerivativesA98}
& = &  \int  \mathsf{L}_{\vect{z}}\left(\vect{\theta}\right)  \exp\left(-K_{Q, \vect{z}}\left(-\frac{1}{s} \right) -\frac{1}{s}\, \mathsf{L}_{\vect{z}}\left(\vect{\theta}\right)  \right) \d Q(\vect{\theta})\\
\label{EqProofCorDerivativesA}
& = &  \int \mathsf{L}_{\vect{z}}\left( \vect{\theta} \right)  \d P^{(Q,s)}_{\vect{\Theta}| \vect{Z} = \vect{z}} (\vect{\theta}),
\end{IEEEeqnarray}
where 
the equality in~\eqref{EqProofCorDerivativesA95} holds  from the dominated convergence theorem \cite[Theorem~$1.6.9$]{ash2000probability};
the equality in~\eqref{EqProofCorDerivativesA97} follows from~\eqref{EqK}; and 
the equality in~\eqref{EqProofCorDerivativesA} follows from~\eqref{EqGenpdf}.

For all~$s \in \set{K}_{Q, \vect{z}}$, with~$\set{K}_{Q, \vect{z}}$ in~\eqref{EqSetKxy}, the equalities in~\eqref{EqK123} and~\eqref{EqProofCorDerivativesA98} imply that
\begin{IEEEeqnarray}{rcl}
\nonumber
& & K^{(2)}_{Q, \vect{z}}\left(- \frac{1}{s} \right) \\
& = &  \frac{\mathrm{d}}{\mathrm{d}t} \int  \mathsf{L}_{\vect{z}}\left(\vect{\theta}\right)  \exp\left(-K_{Q, \vect{z}}\left(t \right) + t\, \mathsf{L}_{\vect{z}}\left(\vect{\theta}\right)  \right) \d Q(\vect{\theta}) \Biggr|_{t = -\frac{1}{s}} \\ 
\nonumber
& = & \int  \mathsf{L}_{\vect{z}}\left(\vect{\theta}\right)\left(-K^{(1)}_{Q, \vect{z}}\left(t \right) + \mathsf{L}_{\vect{z}}\left(\vect{\theta}\right) \right) \\ 
\label{EqProofCorDerivativesB82}
& & \exp\left(-K_{Q, \vect{z}}\left(t \right) + t \mathsf{L}_{\vect{z}}\left(\vect{\theta}\right)  \right) \d Q(\vect{\theta}) \Biggr|_{t = -\frac{1}{s}} \\ 
\nonumber
& = & \int  \mathsf{L}_{\vect{z}}\left(\vect{\theta}\right)\left(-K^{(1)}_{Q, \vect{z}}\left(-\frac{1}{s}\right) + \mathsf{L}_{\vect{z}}\left(\vect{\theta}\right) \right) \\
& &  \exp\left(-K_{Q, \vect{z}}\left(-\frac{1}{s}  \right) -\frac{1}{s}  \mathsf{L}_{\vect{z}}\left(\vect{\theta}\right)  \right) \d Q(\vect{\theta})  \\ 
\label{EqProofCorDerivativesB83Texas}
& = & \int  \mathsf{L}_{\vect{z}}\left(\vect{\theta}\right)\left(-K^{(1)}_{Q, \vect{z}}\left(-\frac{1}{s} \right) + \mathsf{L}_{\vect{z}}\left(\vect{\theta}\right) \right)   \mathrm{d}P^{(Q,s)}_{\vect{\Theta}| \vect{Z} = \vect{z}}\left(\vect{\theta}\right)  \\ 
\nonumber
& = & -K^{(1)}_{Q, \vect{z}}\left(-\frac{1}{s} \right)  \int  \mathsf{L}_{\vect{z}}\left(\vect{\theta}\right)     \mathrm{d}P^{(Q,s)}_{\vect{\Theta}| \vect{Z} = \vect{z}}\left(\vect{\theta}\right)  \\
\label{EqProofCorDerivativesB84}
& &  +  \int  \left( \mathsf{L}_{\vect{z}}\left(\vect{\theta}\right) \right)^2 \mathrm{d}P^{(Q,s)}_{\vect{\Theta}| \vect{Z} = \vect{z}}\left(\vect{\theta}\right) \\  
\label{EqProofCorDerivativesB85}
& = & - \left( K^{(1)}_{Q, \vect{z}}\left(-\frac{1}{s} \right)\right)^2  +  \int  \left( \mathsf{L}_{\vect{z}}\left(\vect{\theta}\right) \right)^2   \mathrm{d} P^{(Q,s)}_{\vect{\Theta}| \vect{Z} = \vect{z}}\left(\vect{\theta}\right)    \\  
\label{EqProofCorDerivativesB86}
& = & \int  \left( \mathsf{L}_{\vect{z}}\left(\vect{\theta}\right) - K^{(1)}_{Q, \vect{z}}\left(-\frac{1}{s} \right) \right)^2   \mathrm{d} P^{(Q,s)}_{\vect{\Theta}| \vect{Z} = \vect{z}}\left(\vect{\theta}\right)   , 
\end{IEEEeqnarray}
where 
the equality in~\eqref{EqProofCorDerivativesB82} follows from the dominated convergence theorem \cite[Theorem~$1.6.9$]{ash2000probability}; 
the equality in~\eqref{EqProofCorDerivativesB83Texas} is due to a change of measure through the Radon-Nikodym derivative in~\eqref{EqGenpdf}; and 
the equality in~\eqref{EqProofCorDerivativesB85} follows from~\eqref{EqProofCorDerivativesA}.

For all~$s \in \set{K}_{Q, \vect{z}}$, with~$\set{K}_{Q, \vect{z}}$ in~\eqref{EqSetKxy}, the equalities in~\eqref{EqK123} and~\eqref{EqProofCorDerivativesB85} imply that

\begin{IEEEeqnarray}{lcl}
\nonumber
& & K^{(3)}_{Q, \vect{z}}\left(- \frac{1}{s} \right) \\
& = &  \frac{\mathrm{d}}{\mathrm{d} t} \left( \displaystyle\int \left( \mathsf{L}_{\vect{z}}\left(\vect{\theta}\right) \right)^2 \mathrm{d}P^{\left(Q,-\frac{1}{t}\right)}_{\vect{\Theta}| \vect{Z} = \vect{z}}\left(\vect{\theta}\right) -\left( K^{(1)}_{Q, \vect{z}}\left( t \right) \right)^2 \right)\Biggr|_{t = -\frac{1}{s}} \\ 
\nonumber
& = &  \frac{\mathrm{d}}{\mathrm{d} t} \Bigg(\displaystyle\int \bigg(
\left( \mathsf{L}_{\vect{z}}\left(\vect{\theta}\right) \right)^2 \exp\left( - K_{Q, \vect{z}}\left(t \right) +t \mathsf{L}_{\vect{z}}\left( \vect{\theta}\right) \right) \bigg)\d Q \left( \vect{\theta} \right) \\
\label{EqRurutu}
& & -\left( K^{(1)}_{Q, \vect{z}}\left( t \right) \right)^2 \Bigg) \Biggr|_{t = -\frac{1}{s}}\\
\nonumber
& = &  \displaystyle\int 
\left( \mathsf{L}_{\vect{z}}\left(\vect{\theta}\right) \right)^2 \left( \frac{\mathrm{d}}{\mathrm{d} t} \exp\left( - K_{Q, \vect{z}}\left( t \right) + t \mathsf{L}_{\vect{z}}\left( \vect{\theta}\right)\right)\Biggr|_{t = -\frac{1}{s}} \right) \d Q \left( \vect{\theta} \right)\squeezeequ
 \\
\label{EqTikehau}
& &- 2 K^{(1)}_{Q, \vect{z}}\left(t \right) K^{(2)}_{Q, \vect{z}}\left(t \right) \Biggr|_{t = -\frac{1}{s}}  
\end{IEEEeqnarray}
\begin{IEEEeqnarray}{rcl}
\nonumber
 & = &  \displaystyle\int 
\left( \mathsf{L}_{\vect{z}}\left(\vect{\theta}\right) \right)^2 \\
\nonumber
& & \left(\mathsf{L}_{\vect{z}}\left( \vect{\theta}\right) - K^{(1)}_{Q, \vect{z}}\left( t\right) \right)\exp\left( - K_{Q, \vect{z}}\left( t \right) + t \mathsf{L}_{\vect{z}}\left( \vect{\theta}\right)\right)\Biggr|_{t = -\frac{1}{s}} \d Q \left( \vect{\theta} \right)\supersqueezeequ \\
& & - 2 K^{(1)}_{Q, \vect{z}}\left(t \right) K^{(2)}_{Q, \vect{z}}\left(t \right) \Biggr|_{t = -\frac{1}{s}}  \\
\nonumber
 & = &  \displaystyle\int 
\left( \mathsf{L}_{\vect{z}}\left(\vect{\theta}\right) \right)^2 \left(\mathsf{L}_{\vect{z}}\left( \vect{\theta}\right) - K^{(1)}_{Q, \vect{z}}\left( -\frac{1}{s} \right) \right) \\
\nonumber
& & \exp\left( - K_{Q, \vect{z}}\left( -\frac{1}{s} \right) -\frac{1}{s} \mathsf{L}_{\vect{z}}\left( \vect{\theta}\right)\right)  \d Q \left( \vect{\theta} \right) \\
\label{EqLastCall112}
& & - 2 K^{(1)}_{Q, \vect{z}}\left(-\frac{1}{s} \right) K^{(2)}_{Q, \vect{z}}\left(-\frac{1}{s}\right)\\
 & = &  \displaystyle\int 
\nonumber
\left( \mathsf{L}_{\vect{z}}\left(\vect{\theta}\right) \right)^2 \left(\mathsf{L}_{\vect{z}}\left( \vect{\theta}\right) - K^{(1)}_{Q, \vect{z}}\left( -\frac{1}{s} \right) \right) \d P^{\left(Q,s\right)}_{\vect{\Theta}| \vect{Z} = \vect{z}}\left(\vect{\theta}\right)\\
 \label{EqMarquezas}
 && - 2 K^{(1)}_{Q, \vect{z}}\left(-\frac{1}{s} \right) K^{(2)}_{Q, \vect{z}}\left(-\frac{1}{s}\right) \\
\nonumber 
& = &\int  \left( \mathsf{L}_{\vect{z}}\left(\vect{\theta}\right) \right)^3  \d P^{\left(Q,s\right)}_{\vect{\Theta}| \vect{Z} = \vect{z}}\left(\vect{\theta}\right)  \\
\nonumber
& & - K^{(1)}_{Q, \vect{z}}\left( -\frac{1}{s} \right) \int  \left( \mathsf{L}_{\vect{z}}\left(\vect{\theta}\right) \right)^2 \d P^{\left(Q,s\right)}_{\vect{\Theta}| \vect{Z} = \vect{z}}\left(\vect{\theta}\right)  \\
& & -  2 K^{(1)}_{Q, \vect{z}}\left(- \frac{1}{s} \right) K^{(2)}_{Q, \vect{z}}\left(- \frac{1}{s} \right)\\
\nonumber
& = &  \int  \left( \mathsf{L}_{\vect{z}}\left(\vect{\theta}\right) \right)^3  \d P^{\left(Q,s\right)}_{\vect{\Theta}| \vect{Z} = \vect{z}}\left(\vect{\theta}\right)  \\
\nonumber
 &  & - K^{(1)}_{Q, \vect{z}}\left( -\frac{1}{s} \right)   \left( K^{(2)}_{Q, \vect{z}}\left(- \frac{1}{s} \right) + \left( K^{(1)}_{Q, \vect{z}}\left(- \frac{1}{s} \right) \right)^2 \right)\\
 \label{EqUrkinaona}
& &   -  2 K^{(1)}_{Q, \vect{z}}\left(- \frac{1}{s} \right) K^{(2)}_{Q, \vect{z}}\left(- \frac{1}{s} \right)\\
& = & \int  \left( \mathsf{L}_{\vect{z}}\left(\vect{\theta}\right) \right)^3  \d P^{\left(Q,s\right)}_{\vect{\Theta}| \vect{Z} = \vect{z}}\left(\vect{\theta}\right)  - K^{(1)}_{Q, \vect{z}}\left( -\frac{1}{s} \right)^3 \\
& &-  3 K^{(1)}_{Q, \vect{z}}\left(- \frac{1}{s} \right) K^{(2)}_{Q, \vect{z}}\left(- \frac{1}{s} \right)\\
& = &  \int \left( \mathsf{L}_{\vect{z}}\left(\vect{\theta}\right) - K^{(1)}_{Q, \vect{z}}\left( -\frac{1}{s} \right) \right)^3 \d P^{\left(Q, s\right)}_{\vect{\Theta}| \vect{Z} = \vect{z}}\left(\vect{\theta}\right),
\end{IEEEeqnarray}
where 
the equality in~\eqref{EqRurutu} follows from~\eqref{EqGenpdf}; 
and the equality in~\eqref{EqTikehau} follows from the dominated convergence theorem \cite[Theorem~$1.6.9$]{ash2000probability};
the equality in~\eqref{EqMarquezas} follows from~\eqref{EqGenpdf}; and
the equality in~\eqref{EqUrkinaona} follows from~\eqref{EqProofCorDerivativesB85}.

This completes the proof.

\section{Proof of Theorem~\ref{CorDecreasingAverage}}\label{AppProofCorDecreasingAverage}

The proof is based on the analysis of the derivative of~$K^{(1)}_{Q, \vect{z}} \left(-\frac{1}{\lambda} \right)$ with respect to~$\lambda$ in~$\mathrm{int}\set{K}_{Q, \vect{z}}$. This is due to Corollary~\ref{CorRK}. 
For instance, note that
\begin{IEEEeqnarray}{lcl}
\frac{\mathrm{d}}{\mathrm{d} \lambda} \mathsf{R}_{\vect{z}}\left(P^{\left(Q, \lambda\right)}_{\vect{\Theta} | \vect{Z} = \vect{z}} \right) & = & \frac{\mathrm{d}}{\mathrm{d} \lambda} K^{(1)}_{Q, \vect{z}} \left(-\frac{1}{\lambda} \right) \\
\label{EqRangiroaB8765BackToNice}
& = & \frac{1}{\lambda^2} K^{(2)}_{Q, \vect{z}} \left(-\frac{1}{\lambda} \right)\\
\label{EqRangiroaB8765}
 & \geqslant &  0,
\end{IEEEeqnarray}
where the equality in~\eqref{EqRangiroaB8765BackToNice} follows from Lemma~\ref{CorDerivatives}.
The inequality in~\eqref{EqRangiroaB8765} implies that the expected empirical risk~$\mathsf{R}_{\vect{z}}\left(P^{\left(Q, \lambda\right)}_{\vect{\Theta} | \vect{Z} = \vect{z}} \right)  = K^{(1)}_{Q, \vect{z}}\left(- \frac{1}{\lambda} \right)$ in~\eqref{EqRK} is nondecreasing with respect to~$\lambda$.
The rest of the proof consists in showing that for all~$\alpha \in \set{K}_{Q, \vect{z}}$,  the function~$K^{(2)}_{Q, \vect{z}}$  in~\eqref{EqK123} satisfies~$K^{(2)}_{Q, \vect{z}}\left(-\frac{1}{\alpha} \right)> 0$ if and only if the function~$\mathsf{L}_{\vect{z}}$ in~\eqref{EqLxy} is separable. 
For doing so,  a handful of preliminary results are described in the following subsection. The proof of Theorem~\ref{CorDecreasingAverage} resumes in Subsection~\ref{SecAppProofCorDecreasingAverageResumes}

\subsection{Preliminaries}

Given a positive real~$\lambda \in \set{K}_{Q, \vect{z}}$, with~$\set{K}_{Q, \vect{z}}$ in~\eqref{EqSetKxy},  consider a partition of~$\set{M}$ formed by the sets~$\set{R}_0(\lambda)$,~$\set{R}_1(\lambda)$ and~$\set{R}_2(\lambda)$, such that
\begin{subequations}\label{EqSetR}
\begin{IEEEeqnarray}{rcl}
\label{EqSetR0}
\set{R}_0(\lambda) & \triangleq &\left\lbrace \vect{\nu} \in \set{M}: \mathsf{L}_{\vect{z}}\left( \vect{\nu}\right)  =    \mathsf{R}_{\vect{z}}\left(P^{\left(Q, \lambda\right)}_{\vect{\Theta} | \vect{Z} = \vect{z}} \right) \right \rbrace,\\
\label{EqSetR1}
\set{R}_1(\lambda) & \triangleq &\left\lbrace \vect{\nu} \in \set{M}: \mathsf{L}_{\vect{z}}\left( \vect{\nu}\right)  <   \mathsf{R}_{\vect{z}}\left(P^{\left(Q, \lambda\right)}_{\vect{\Theta} | \vect{Z} = \vect{z}} \right) \right \rbrace, \mbox{ and } \IEEEeqnarraynumspace\\
\label{EqSetR2}
\set{R}_2 (\lambda)& \triangleq & \left\lbrace \vect{\nu} \in \set{M}: \mathsf{L}_{\vect{z}}\left( \vect{\nu}\right)  >  \mathsf{R}_{\vect{z}}\left(P^{\left(Q, \lambda\right)}_{\vect{\Theta} | \vect{Z} = \vect{z}} \right) \right \rbrace,
\end{IEEEeqnarray}
\end{subequations}
where the functional~$\mathsf{R}_{\vect{z}}$ is in~\eqref{EqRxy} and the probability measure~$P^{\left(Q, \lambda\right)}_{\vect{\Theta} | \vect{Z} = \vect{z}}$ is in~\eqref{EqGenpdf}.
The sets in~\eqref{EqSetR} exhibit several properties that are central for proving the main results of this section. 
\begin{lemma}\label{LemmaR1R2}
The probability measure~$P^{\left(Q, \lambda \right)}_{\vect{\Theta}| \vect{Z} = \vect{z}}$ in~\eqref{EqGenpdf}, satisfies  
\begin{equation}
P^{\left(Q, \lambda \right)}_{\vect{\Theta}| \vect{Z} = \vect{z}} \left( \set{R}_1(\lambda)\right) > 0,
\end{equation}
if and only if  
\begin{equation}
P^{\left(Q, \lambda \right)}_{\vect{\Theta}| \vect{Z} = \vect{z}}\left( \set{R}_2(\lambda)\right)> 0, 
\end{equation}
where the sets~$\set{R}_1(\cdot)$ and~$\set{R}_2(\cdot)$ are in~\eqref{EqSetR1} and~\eqref{EqSetR2}, respectively.
\end{lemma}
\begin{IEEEproof}
The proof is divided into two parts. In the first part, given a  real~$\alpha \in \set{K}_{Q, \vect{z}}$, it is proven that  if the set~$\set{R}_1\left( \alpha \right)$ is nonnegligible with respect to~$P^{(Q,\alpha)}_{\vect{\Theta}| \vect{Z} = \vect{z}}~$, then the set~$\set{R}_2\left( \alpha \right)$ is nonnegligible with respect to~$P^{(Q,\alpha)}_{\vect{\Theta}| \vect{Z} = \vect{z}}~$.
The second part proves the converse.

The first part  is proved by contradiction. Assume that set~$\set{R}_2\left( \alpha \right)$ is negligible with respect to~$P^{(Q,\alpha)}_{\vect{\Theta}| \vect{Z} = \vect{z}}~$.  Hence, from Lemma~\ref{CorDerivatives}, it holds that
\begin{IEEEeqnarray}{rcl}
\nonumber
& & K^{(1)}_{Q, \vect{z}}\left(-\frac{1}{\alpha} \right) \\
\nonumber
&=&  \displaystyle\int_{\set{R}_0(\alpha)}  \mathsf{L}_{\vect{z}}\left( \vect{\nu}\right)   \mathrm{d}P^{(Q,\alpha)}_{\vect{\Theta}| \vect{Z} = \vect{z}} \left( \vect{\nu} \right) +   \displaystyle\int_{\set{R}_1(\alpha)}  \mathsf{L}_{\vect{z}}\left( \vect{\nu}\right)   \mathrm{d}P^{(Q,\alpha)}_{\vect{\Theta}| \vect{Z} = \vect{z}} \left( \vect{\nu} \right) \\
\label{EqTahitiNui876235}
& & + \displaystyle\int_{\set{R}_2(\alpha)}  \mathsf{L}_{\vect{z}}\left( \vect{\nu}\right)   \mathrm{d}P^{(Q,\alpha)}_{\vect{\Theta}| \vect{Z} = \vect{z}} \left( \vect{\nu} \right)\\
\label{EqTahitiNui8772}
&=&  \displaystyle\int_{\set{R}_0(\alpha)} \hspace{-2.5ex} \mathsf{L}_{\vect{z}}\left( \vect{\nu}\right)   \mathrm{d}P^{(Q,\alpha)}_{\vect{\Theta}| \vect{Z} = \vect{z}} \left( \vect{\nu} \right) +   \displaystyle\int_{\set{R}_1(\alpha)} \hspace{-2.5ex}  \mathsf{L}_{\vect{z}}\left( \vect{\nu}\right)   \mathrm{d}P^{(Q,\alpha)}_{\vect{\Theta}| \vect{Z} = \vect{z}} \left( \vect{\nu} \right) \\
\nonumber
& = &  K^{(1)}_{Q, \vect{z}}\left(-\frac{1}{\alpha} \right)  P^{(Q,\alpha)}_{\vect{\Theta}| \vect{Z} = \vect{z}} \left( \set{R}_0(\alpha) \right) \\
\label{EqTahitiNui878myw}
& &  +   \displaystyle\int_{\set{R}_1(\alpha)}  \mathsf{L}_{\vect{z}}\left( \vect{\nu}\right)   \mathrm{d}P^{(Q,\alpha)}_{\vect{\Theta}| \vect{Z} = \vect{z}} \left( \vect{\nu} \right) \\ 
\nonumber
& < &  K^{(1)}_{Q, \vect{z}}\left(-\frac{1}{\alpha} \right)  P^{(Q,\alpha)}_{\vect{\Theta}| \vect{Z} = \vect{z}} \left( \set{R}_0(\alpha) \right) \\
\label{EqTahitiNui8791r}
& &+  K^{(1)}_{Q, \vect{z}}\left(-\frac{1}{\alpha} \right)  P^{(Q,\alpha)}_{\vect{\Theta}| \vect{Z} = \vect{z}} \left( \set{R}_1(\alpha) \right) \\ 
\label{EqTahitiNui880}
& = &  K^{(1)}_{Q, \vect{z}}\left(-\frac{1}{\alpha} \right)  \big( P^{(Q,\alpha)}_{\vect{\Theta}| \vect{Z} = \vect{z}} \left( \set{R}_0(\alpha) \right)  +  P^{(Q,\alpha)}_{\vect{\Theta}| \vect{Z} = \vect{z}} \left( \set{R}_1(\alpha) \right) \big) \IEEEeqnarraynumspace\\ 
\label{EqTahitiNui8802r31}
& = &  K^{(1)}_{Q, \vect{z}}\left(-\frac{1}{\alpha} \right),  
\end{IEEEeqnarray}
which is a contradiction.

The second part of the proof follows the same arguments as in the first part. Assume that the set~$\set{R}_1\left( \alpha \right)$ is negligible with respect to~$P^{(Q,\alpha)}_{\vect{\Theta}| \vect{Z} = \vect{z}}~$. Hence, from Lemma~\ref{CorDerivatives}, it holds that
\begin{IEEEeqnarray}{rcl}
\nonumber
& & K^{(1)}_{Q, \vect{z}}\left(-\frac{1}{\alpha} \right)\\
\nonumber &=&  \displaystyle\int_{\set{R}_0(\alpha)}  \mathsf{L}_{\vect{z}}\left( \vect{\nu}\right)   \mathrm{d}P^{(Q,\alpha)}_{\vect{\Theta}| \vect{Z} = \vect{z}} \left( \vect{\nu} \right) +   \displaystyle\int_{\set{R}_1(\alpha)}  \mathsf{L}_{\vect{z}}\left( \vect{\nu}\right)   \mathrm{d}P^{(Q,\alpha)}_{\vect{\Theta}| \vect{Z} = \vect{z}} \left( \vect{\nu} \right) \\
\label{EqTahitiNui876}
& & + \displaystyle\int_{\set{R}_2(\alpha)}  \mathsf{L}_{\vect{z}}\left( \vect{\nu}\right)   \mathrm{d}P^{(Q,\alpha)}_{\vect{\Theta}| \vect{Z} = \vect{z}} \left( \vect{\nu} \right)\\
\label{EqTahitiNui877}
&=&  \displaystyle\int_{\set{R}_0(\alpha)} \hspace{-3ex} \mathsf{L}_{\vect{z}}\left( \vect{\nu}\right)   \mathrm{d}P^{(Q,\alpha)}_{\vect{\Theta}| \vect{Z} = \vect{z}} \left( \vect{\nu} \right) +   \displaystyle\int_{\set{R}_2(\alpha)} \hspace{-3ex} \mathsf{L}_{\vect{z}}\left( \vect{\nu}\right)   \mathrm{d}P^{(Q,\alpha)}_{\vect{\Theta}| \vect{Z} = \vect{z}} \left( \vect{\nu} \right) \\
\nonumber
& = &  K^{(1)}_{Q, \vect{z}}\left(-\frac{1}{\alpha} \right)  P^{(Q,\alpha)}_{\vect{\Theta}| \vect{Z} = \vect{z}} \left( \set{R}_0(\alpha) \right)   \\
\label{EqTahitiNui878}
& &+   \displaystyle\int_{\set{R}_2(\alpha)} \mathsf{L}_{\vect{z}}\left( \vect{\nu}\right)   \mathrm{d}P^{(Q,\alpha)}_{\vect{\Theta}| \vect{Z} = \vect{z}} \left( \vect{\nu} \right)  \IEEEeqnarraynumspace\\ 
\nonumber
& > &  K^{(1)}_{Q, \vect{z}}\left(-\frac{1}{\alpha} \right)  P^{(Q,\alpha)}_{\vect{\Theta}| \vect{Z} = \vect{z}} \left( \set{R}_0(\alpha) \right) \\
\label{EqTahitiNui879}
& &+  K^{(1)}_{Q, \vect{z}}\left(-\frac{1}{\alpha} \right)  P^{(Q,\alpha)}_{\vect{\Theta}| \vect{Z} = \vect{z}} \left( \set{R}_2(\alpha) \right) \\ 
\label{EqTahitiNui880r3}
& = &  K^{(1)}_{Q, \vect{z}}\left(-\frac{1}{\alpha} \right)  \left( P^{(Q,\alpha)}_{\vect{\Theta}| \vect{Z} = \vect{z}} \left( \set{R}_0(\alpha) \right)  +  P^{(Q,\alpha)}_{\vect{\Theta}| \vect{Z} = \vect{z}} \left( \set{R}_2(\alpha) \right) \right) \\ 
\label{EqTahitiNui880234}
& = &  K^{(1)}_{Q, \vect{z}}\left(-\frac{1}{\alpha} \right),  
\end{IEEEeqnarray}
which is also a contradiction. This completes the proof.
\end{IEEEproof}

A more general result can be immediately obtained by combining Lemma~\ref{LemmaMutualAlphaBeta} and  Lemma~\ref{LemmaR1R2}.

\begin{lemma}\label{LemmaGenR1R2}
For all~$\alpha \in \set{K}_{Q, \vect{z}}$, with~$\set{K}_{Q, \vect{z}}$ in~\eqref{EqSetKxy}, the measure~$P^{\left(Q, \lambda \right)}_{\vect{\Theta}| \vect{Z} = \vect{z}}$ in~\eqref{EqGenpdf}, satisfies  
\begin{equation}
P^{\left(Q, \lambda \right)}_{\vect{\Theta}| \vect{Z} = \vect{z}} \left( \set{R}_1(\alpha)\right) > 0,
\end{equation}
if and only if  
\begin{equation}
P^{\left(Q, \lambda \right)}_{\vect{\Theta}| \vect{Z} = \vect{z}}\left( \set{R}_2(\alpha)\right)> 0, 
\end{equation}
where the sets~$\set{R}_1(\alpha)$ and~$\set{R}_2(\alpha)$ are in~\eqref{EqSetR1} and~\eqref{EqSetR2}, respectively.
\end{lemma}

\subsection{The proof}\label{SecAppProofCorDecreasingAverageResumes}

The rest of the proof of Theorem~\ref{CorDecreasingAverage} is divided into two parts. 
In the first part, it is shown that if for all~$\alpha \in \set{K}_{Q, \vect{z}}$,~$K^{(2)}_{Q, \vect{z}}\left(-\frac{1}{\alpha} \right)> 0$, then the function~$\mathsf{L}_{\vect{z}}$ in~\eqref{EqLxy} is separable.
The second part of the proof, consists in showing that if the function~$\mathsf{L}_{\vect{z}}$ is separable, then, for all~$\alpha \in \set{K}_{Q, \vect{z}}$,~$K^{(2)}_{Q, \vect{z}}\left(-\frac{1}{\alpha} \right)> 0$.

The first part is as follows.  
From Lemma~\ref{CorDerivatives}, it holds that for all~$\alpha \in \set{K}_{Q, \vect{z}}$, 
\begin{IEEEeqnarray}{rcl}
\nonumber
& &K^{(2)}_{Q, \vect{z}}\left(-\frac{1}{\alpha} \right)\\
\label{EqChaparral}
& = & \displaystyle\int\left( \mathsf{L}_{\vect{z}}\left( \vect{\theta}\right) -  K^{(1)}_{Q, \vect{z}}\left(-\frac{1}{\alpha} \right)  \right)^2  \mathrm{d}P^{(Q, \alpha)}_{\vect{\Theta}| \vect{Z} = \vect{z}} \left( \vect{\theta} \right) \\
\nonumber
& = & \displaystyle\int_{\set{R}_0\left( \alpha \right)} \left( \mathsf{L}_{\vect{z}}\left( \vect{\theta}\right) -  K^{(1)}_{Q, \vect{z}}\left(-\frac{1}{\alpha} \right)  \right)^2  \mathrm{d}P^{(Q,\alpha)}_{\vect{\Theta}| \vect{Z} = \vect{z}} \left( \vect{\theta} \right) \\
\nonumber
&  &+  \displaystyle\int_{\set{R}_1\left( \alpha \right)} \left( \mathsf{L}_{\vect{z}}\left( \vect{\theta}\right) -  K^{(1)}_{Q, \vect{z}}\left(-\frac{1}{\alpha} \right)  \right)^2  \mathrm{d}P^{(Q,\alpha)}_{\vect{\Theta}| \vect{Z} = \vect{z}} \left( \vect{\theta} \right)\IEEEeqnarraynumspace \\
& & + \displaystyle\int_{\set{R}_2\left( \alpha \right)} \left( \mathsf{L}_{\vect{z}}\left( \vect{\theta}\right) - K^{(1)}_{Q, \vect{z}}\left(-\frac{1}{\alpha} \right) \right)^2  \mathrm{d}P^{(Q,\alpha)}_{\vect{\Theta}| \vect{Z} = \vect{z}} \left( \vect{\theta} \right),\IEEEeqnarraynumspace
\end{IEEEeqnarray}
where the sets~$\set{R}_0(\alpha)$,~$\set{R}_1(\alpha)$, and~$\set{R}_2(\alpha)$ are respectively defined in~\eqref{EqSetR}. 
Hence, 
\begin{IEEEeqnarray}{rcl}
\nonumber
& & K^{(2)}_{Q, \vect{z}}\left(-\frac{1}{\alpha} \right)\\
& = & \displaystyle\int_{\set{R}_1\left( \alpha \right)} \left( \mathsf{L}_{\vect{z}}\left( \vect{\theta}\right) -  K^{(1)}_{Q, \vect{z}}\left(-\frac{1}{\alpha} \right)  \right)^2  \mathrm{d}P^{(Q,\alpha)}_{\vect{\Theta}| \vect{Z} = \vect{z}} \left( \vect{\theta} \right)\\
\label{EqChaparralito}
& & +\displaystyle\int_{\set{R}_2\left( \gamma \right)} \left( \mathsf{L}_{\vect{z}}\left( \vect{\theta}\right) - K^{(1)}_{Q, \vect{z}}\left(-\frac{1}{\alpha} \right) \right)^2  \mathrm{d}P^{(Q,\alpha)}_{\vect{\Theta}| \vect{Z} = \vect{z}} \left( \vect{\theta} \right).\IEEEeqnarraynumspace
\end{IEEEeqnarray}
Under the assumption that for all~$\alpha \in \set{K}_{Q, \vect{z}}$ the function~$K^{(2)}_{Q, \vect{z}}$  in~\eqref{EqK123} satisfies~$K^{(2)}_{Q, \vect{z}}\left(-\frac{1}{\alpha} \right)> 0$, it follows that at least one of the following claims is true: \newline
\textbf{(a)}~$P^{(Q,\alpha)}_{\vect{\Theta}| \vect{Z} = \vect{z}} \left(\set{R}_1(\alpha) \right) > 0$; and \newline
\textbf{(b)}~$P^{(Q,\alpha)}_{\vect{\Theta}| \vect{Z} = \vect{z}} \left(\set{R}_2(\alpha) \right) > 0$.  

Nonetheless, from Lemma~\ref{LemmaR1R2}, it follows that both claims \textbf{(a)} and \textbf{(b)} hold simultaneously.  
Hence, the sets~$\set{R}_1(\alpha)$ and~$\set{R}_2(\alpha)$ are both nonnegligible with respect to~$P^{(Q,\alpha)}_{\vect{\Theta}| \vect{Z} = \vect{z}}~$ and moreover, it holds that for all~$\left( \vect{\nu}_1, \vect{\nu}_2 \right) \in \set{R}_1(\alpha) \times \set{R}_2(\alpha)$, 
\begin{IEEEeqnarray}{rcl}
+\infty > \mathsf{L}_{\vect{z}} \left( \vect{\nu}_1 \right) &> K^{(1)}_{Q, \vect{z}}\left(-\frac{1}{\alpha} \right)  > & \mathsf{L}_{\vect{z}}\left(\vect{\nu}_2\right),
\end{IEEEeqnarray}
where $\mathsf{L}_{\vect{z}} \left( \vect{\nu}_1 \right)  < +\infty$ follows from the fact that $ P^{\left(Q, \lambda\right)}_{\vect{\Theta} | \vect{Z} = \vect{z}} \left( \left\lbrace \vect{\theta} \in \set{M}:  \mathsf{L}_{\vect{z}}\left( \vect{\theta}\right) = +\infty \right\rbrace  \right) =0$ (Lemma~\ref{CorPositive}).
This proves that under the assumption that for all~$\alpha \in \set{K}_{Q, \vect{z}}$,~$K^{(2)}_{Q, \vect{z}}\left(-\frac{1}{\alpha} \right)> 0$,  the function~$\mathsf{L}_{\vect{z}}$ in~\eqref{EqLxy} is separable with respect to~$P^{(Q,\alpha)}_{\vect{\Theta}| \vect{Z} = \vect{z}}$. From Lemma~\ref{CorAC}, it holds that the function~$\mathsf{L}_{\vect{z}}$ is separable with respect to~$Q$.  This completes the first part of the proof. 

The second part of the proof is simpler. Assume that the empirical risk function~$\mathsf{L}_{\vect{z}}$ in~\eqref{EqLxy} is separable with respect to~$P^{(Q,\alpha)}_{\vect{\Theta}| \vect{Z} = \vect{z}}$. That is, for all~$\gamma \in \set{K}_{Q, \vect{z}}$, there exist a positive real~$c_{\gamma} >0$;  and two subsets~$\set{A}(\gamma)$ and~$\set{B}(\gamma)$ of~$\set{M}$ that are nonnegligible with respect to~$P^{\left(Q, \gamma \right)}_{\vect{\Theta} | \vect{Z} = \vect{z}}$ in~\eqref{EqGenpdf} and verify that for all~$(\vect{\nu}_1,\vect{\nu}_2) \in \set{A}(\gamma) \times \set{B}(\gamma)$, 
\begin{IEEEeqnarray}{rcl}
\label{EqTwoNonnegligibleSetsTontin}
+\infty > \mathsf{L}_{\vect{z}} \left( \vect{\nu}_1 \right) &> c_{\gamma} >& \mathsf{L}_{\vect{z}}\left(\vect{\nu}_2\right).
\end{IEEEeqnarray}

From Lemma~\ref{CorDerivatives}, it holds that
\begin{IEEEeqnarray}{rcl}
\nonumber
& &K^{(2)}_{Q, \vect{z}}\left(-\frac{1}{\gamma} \right)\\
& = & \displaystyle\int\left( \mathsf{L}_{\vect{z}}\left( \vect{\theta}\right) -  K^{(1)}_{Q, \vect{z}}\left(-\frac{1}{\gamma} \right)  \right)^2  \mathrm{d}P^{\left(Q, \gamma \right)}_{\vect{\Theta} | \vect{Z} = \vect{z}} \left( \vect{\theta} \right) \\
& = & \displaystyle\int_{\set{A}\left( \gamma \right)} \left( \mathsf{L}_{\vect{z}}\left( \vect{\theta}\right) -  K^{(1)}_{Q, \vect{z}}\left(-\frac{1}{\gamma} \right)  \right)^2  \mathrm{d}P^{\left(Q, \gamma \right)}_{\vect{\Theta} | \vect{Z} = \vect{z}} \left( \vect{\theta} \right) \\
& + &  \displaystyle\int_{\set{B}\left( \gamma \right)} \left( \mathsf{L}_{\vect{z}}\left( \vect{\theta}\right) -  K^{(1)}_{Q, \vect{z}}\left(-\frac{1}{\gamma} \right)  \right)^2  \mathrm{d}P^{\left(Q, \gamma \right)}_{\vect{\Theta} | \vect{Z} = \vect{z}} \left( \vect{\theta} \right)\\
& + & \displaystyle\int_{\set{M}\setminus\left(\set{A}(\gamma) \cup \set{B}(\gamma) \right)} \left( \mathsf{L}_{\vect{z}}\left( \vect{\theta}\right) - K^{(1)}_{Q, \vect{z}}\left(- \,\frac{1}{\gamma} \right) \right)^2  \mathrm{d}P^{\left(Q, \gamma \right)}_{\vect{\Theta} | \vect{Z} = \vect{z}} \left( \vect{\theta} \right) \supersqueezeequ\IEEEeqnarraynumspace\\
\label{EqWildHorses}
& > & 0,
\end{IEEEeqnarray}
where the inequality~\eqref{EqWildHorses} follows from the following facts. 
First, if~$c_{\gamma} < K^{(1)}_{Q, \vect{z}}\left(-\frac{1}{\gamma} \right)$, with~$c_{\gamma}$ in~\eqref{EqTwoNonnegligibleSetsTontin}, then for all~$\vect{\nu} \in  \set{B}(\gamma)$, it holds that $K^{(1)}_{Q, \vect{z}}\left(-\frac{1}{\gamma} \right)> c_{\gamma} > \mathsf{L}_{\vect{z}}\left( \vect{\nu}\right)$, and thus, 
\begin{IEEEeqnarray}{rcl}
\left( \mathsf{L}_{\vect{z}}\left( \vect{\nu}\right) - K^{(1)}_{Q, \vect{z}}\left(-\frac{1}{\gamma} \right) \right)^2  &>& \left( c_{\gamma} - K^{(1)}_{Q, \vect{z}}\left(-\frac{1}{\gamma} \right) \right)^2,\IEEEeqnarraynumspace
\end{IEEEeqnarray}
which implies, 
\begin{IEEEeqnarray}{rcl}
\nonumber
& & \displaystyle\int_{\set{B}\left( \gamma \right)} \left( \mathsf{L}_{\vect{z}}\left( \vect{\theta}\right) -  K^{(1)}_{Q, \vect{z}}\left(-\frac{1}{\gamma} \right)  \right)^2  \mathrm{d}P^{\left(Q, \gamma \right)}_{\vect{\Theta} | \vect{Z} = \vect{z}} \left( \vect{\theta} \right) \\
&  > & \left( c_{\gamma} - K^{(1)}_{Q, \vect{z}}\left(-\frac{1}{\gamma} \right) \right)^2 P^{\left(Q, \gamma \right)}_{\vect{\Theta} | \vect{Z} = \vect{z}} \left( \set{B}\left( \gamma \right)\right)\\
 & >& 0.
\end{IEEEeqnarray}

Second, if~$c_{\gamma} \geqslant K^{(1)}_{Q, \vect{z}}\left(-\frac{1}{\gamma} \right)$ then for all~$\vect{\nu} \in  \set{A}(\gamma)$, it holds that $\mathsf{L}_{\vect{z}}\left( \vect{\nu}\right) > c_{\gamma} \geqslant K^{(1)}_{Q, \vect{z}}\left(-\frac{1}{\gamma} \right)$, and thus, 
\begin{IEEEeqnarray}{rcl}
\left( \mathsf{L}_{\vect{z}}\left( \vect{\nu}\right) - K^{(1)}_{Q, \vect{z}}\left(-\frac{1}{\gamma} \right) \right)^2  &>& \left( c_{\gamma} - K^{(1)}_{Q, \vect{z}}\left(-\frac{1}{\gamma} \right) \right)^2,\IEEEeqnarraynumspace
\end{IEEEeqnarray}
which implies, 
\begin{IEEEeqnarray}{rcl}
\nonumber
& & \displaystyle\int_{\set{A}\left( \gamma \right)} \left( \mathsf{L}_{\vect{z}}\left( \vect{\theta}\right) -  K^{(1)}_{Q, \vect{z}}\left(-\frac{1}{\gamma} \right)  \right)^2  \mathrm{d}P^{\left(Q, \gamma \right)}_{\vect{\Theta} | \vect{Z} = \vect{z}} \left( \vect{\theta} \right) \\
&  > & \left( c_{\gamma} - K^{(1)}_{Q, \vect{z}}\left(-\frac{1}{\gamma} \right) \right)^2 P^{\left(Q, \gamma \right)}_{\vect{\Theta} | \vect{Z} = \vect{z}} \left( \set{A}\left( \gamma \right)\right) \IEEEeqnarraynumspace\\
 & > & 0.
\end{IEEEeqnarray}
Hence, under the assumption that the empirical risk function~$\mathsf{L}_{\vect{z}}$ in~\eqref{EqLxy} is separable, it holds that for all~$\gamma \in \set{K}_{Q, \vect{z}}$, 
$K^{(2)}_{Q, \vect{z}}\left(-\frac{1}{\gamma} \right) > 0$.
This completes the proof.

\section{Proof of Lemma~\ref{LemmaRzLowerBoundB}}\label{AppProofLemmaRzLowerBoundB}

Consider the partition of the set $\set{M}$ formed by the sets $\set{A}_0$, $\set{A}_1$, and $\set{A}_2$ in~\eqref{EqSetsA}.
%
%
From~\eqref{EqNiceK1}, for all~$\lambda \in \set{K}_{Q, \vect{z}}$, with~$\set{K}_{Q, \vect{z}}$ in~\eqref{EqSetKxy}, it holds that, 
\begin{IEEEeqnarray}{rcl}
\nonumber
& & K^{(1)}_{Q, \vect{z}} \left(-\frac{1}{\lambda} \right) \\
& = & \int_{\set{A}_0 } \mathsf{L}_{\vect{z}}\left( \vect{\theta} \right)  \mathrm{d}P^{\left(Q, \lambda\right)}_{\vect{\Theta} | \vect{Z} = \vect{z}}  (\vect{\theta}) 
+ \int_{\set{A}_1} \mathsf{L}_{\vect{z}}\left( \vect{\theta} \right)  \mathrm{d}P^{\left(Q, \lambda\right)}_{\vect{\Theta} | \vect{Z} = \vect{z}}  (\vect{\theta}) \\
& &+ \int_{\set{A}_2} \mathsf{L}_{\vect{z}}\left( \vect{\theta} \right)  \mathrm{d}P^{\left(Q, \lambda\right)}_{\vect{\Theta} | \vect{Z} = \vect{z}}  (\vect{\theta}) \IEEEeqnarraynumspace\\
\label{EqCha1097}
& = &   
\int_{\set{A}_0 } \mathsf{L}_{\vect{z}}\left( \vect{\theta} \right)  \mathrm{d}P^{\left(Q, \lambda\right)}_{\vect{\Theta} | \vect{Z} = \vect{z}}  (\vect{\theta}) 
+ \int_{\set{A}_2} \mathsf{L}_{\vect{z}}\left( \vect{\theta} \right)  \mathrm{d}P^{\left(Q, \lambda\right)}_{\vect{\Theta} | \vect{Z} = \vect{z}}  (\vect{\theta}) \IEEEeqnarraynumspace\\
\label{EqCha1098}
& = &   
\delta_{Q,\vect{z}}^{\star} P^{\left(Q, \lambda\right)}_{\vect{\Theta} | \vect{Z} = \vect{z}}  (\set{L}^{\star}_{Q,\vect{z}}) 
+ \int_{\set{A}_2} \mathsf{L}_{\vect{z}}\left( \vect{\theta} \right)  \mathrm{d}P^{\left(Q, \lambda\right)}_{\vect{\Theta} | \vect{Z} = \vect{z}}  (\vect{\theta})\\
\label{EqCha109A87a}
& \geqslant &   \delta_{Q,\vect{z}}^{\star} P^{\left(Q, \lambda\right)}_{\vect{\Theta} | \vect{Z} = \vect{z}}  (\set{L}^{\star}_{Q,\vect{z}})  
+ \delta_{Q,\vect{z}}^{\star} P^{\left(Q, \lambda\right)}_{\vect{\Theta} | \vect{Z} = \vect{z}}  (\set{A}_{2}) \\
\label{EqCha109A87}
& = & \delta_{Q,\vect{z}}^{\star},
\end{IEEEeqnarray}
where
the equality in~\eqref{EqCha1097} follows by noticing that~$Q\left( \set{A}_1 \right) = 0$, which implies that~$P^{\left(Q, \lambda\right)}_{\vect{\Theta} | \vect{Z} = \vect{z}}  (\set{A}_{1})$ $=$ $0$ (Lemma~\ref{LemmaMutualAC}); 
the equality in~\eqref{EqCha1098} follows from noticing that~$\set{A}_0 = \set{L}^{\star}_{Q,\vect{z}}$, with~$\set{L}^{\star}_{Q,\vect{z}}$ in~\eqref{EqSetLStar}; and 
the equality in~\eqref{EqCha109A87a} follows from~\eqref{EqSetA20987}.
This completes the proof.


 \section{Proof of Theorem~\ref{CorAssymptoticMean}}\label{AppProofCorAssymptoticMean}

From~\eqref{EqCha1098} in the proof of Lemma~\ref{LemmaRzLowerBoundB}, it holds that  
\begin{IEEEeqnarray}{rcl}
\nonumber
& & \lim_{\lambda \rightarrow 0^{+}}  K^{(1)}_{Q, \vect{z}} \left(-\frac{1}{\lambda} \right) \\
\label{EqPeriklesTotis}
& = &   \lim_{\lambda \rightarrow 0^{+}}   
\delta_{Q,\vect{z}}^{\star} P^{\left(Q, \lambda\right)}_{\vect{\Theta} | \vect{Z} = \vect{z}}  (\set{L}^{\star}_{Q,\vect{z}}) 
+ \lim_{\lambda \rightarrow 0^{+}}    \int_{\set{A}_2}\hspace{-1.5ex} \mathsf{L}_{\vect{z}}\left( \vect{\theta} \right)  \mathrm{d}P^{\left(Q, \lambda\right)}_{\vect{\Theta} | \vect{Z} = \vect{z}}  (\vect{\theta}) \supersqueezeequ\IEEEeqnarraynumspace\\
\nonumber
& = &\lim_{\lambda \rightarrow 0^{+}}   
\delta_{Q,\vect{z}}^{\star} P^{\left(Q, \lambda\right)}_{\vect{\Theta} | \vect{Z} = \vect{z}}  (\set{L}^{\star}_{Q,\vect{z}}) \\
\label{EqPeriklesTot}
& &  +  \lim_{\lambda \rightarrow 0^{+}}  \int_{\set{A}_2}  \mathsf{L}_{\vect{z}}\left( \vect{\theta} \right)  \frac{\mathrm{d}P^{\left(Q, \lambda\right)}_{\vect{\Theta} | \vect{Z} = \vect{z}}}{\mathrm{d}Q} \left( \vect{\theta} \right)   \mathrm{d}Q (\vect{\theta}) \\
\nonumber
& = &\lim_{\lambda \rightarrow 0^{+}}   
\delta_{Q,\vect{z}}^{\star} P^{\left(Q, \lambda\right)}_{\vect{\Theta} | \vect{Z} = \vect{z}}  (\set{L}^{\star}_{Q,\vect{z}}) \\
\label{EqPeriklesTos}
& &+    \int_{\set{A}_2}  \mathsf{L}_{\vect{z}}\left( \vect{\theta} \right) \lim_{\lambda \rightarrow 0^{+}} \frac{\mathrm{d}P^{\left(Q, \lambda\right)}_{\vect{\Theta} | \vect{Z} = \vect{z}}}{\mathrm{d}Q} \left( \vect{\theta} \right)   \mathrm{d}Q (\vect{\theta}) \\
\label{EqPeriklesTotis125}
& = & \delta_{Q,\vect{z}}^{\star}\lim_{\lambda \rightarrow 0^{+}}   
 P^{\left(Q, \lambda\right)}_{\vect{\Theta} | \vect{Z} = \vect{z}}  (\set{L}^{\star}_{Q,\vect{z}}) \\
 \label{EqPeriklesTotisa}
 & = & \delta_{Q,\vect{z}}^{\star},
\end{IEEEeqnarray}
where, 
 the equality in~\eqref{EqPeriklesTos} follows from noticing two facts: 
~$(a)$ For all~$\lambda \in \set{K}_{Q, \vect{z}}$, the Randon-Nikodym derivative~$\frac{\mathrm{d}P^{\left(Q, \lambda\right)}_{\vect{\Theta} | \vect{Z} = \vect{z}}}{\mathrm{d}Q}$ is positive and finite (Lemma~\ref{CorPositive}); and 
~$(b)$ For all~$\vect{\theta} \in \set{A}_2$, it holds that~$ \lim_{\lambda \rightarrow 0^{+}}  \frac{\mathrm{d}P^{\left(Q, \lambda\right)}_{\vect{\Theta} | \vect{Z} = \vect{z}}}{\mathrm{d}Q}\left( \vect{\theta} \right) = 0$ (Lemma~\ref{CorAssymptoticZero}).
Hence, the dominated convergence theorem~\cite[Theorem~$1.6.9$]{ash2000probability} holds.
The inequality in~\eqref{EqPeriklesTotis125} follows from Lemma~\ref{TheoSetLStarConcentration}. 
This completes the proof.

\section{Proof of Theorem~\ref{CorDecreasingSet}}\label{AppProofCorDecreasingSet}

From Theorem~\ref{CorDecreasingAverage}, it follows that for all~$(\lambda_1, \lambda_2) \in \set{K}_{Q, \vect{z}}\times\set{K}_{Q, \vect{z}}$ with~$\lambda_1 > \lambda_2$,
\begin{IEEEeqnarray}{rcl}
\nonumber
\displaystyle\int \mathsf{L}_{\vect{z}}\left( \vect{\alpha}\right)   \frac{\mathrm{d}P^{(Q, \lambda_{1})}_{\vect{\Theta}| \vect{Z} = \vect{z}}}{\mathrm{d}Q} \left( \vect{\alpha} \right)  \mathrm{d} Q\left(\vect{\alpha} \right) \supersqueezeequ & \geqslant &
 \displaystyle\int \mathsf{L}_{\vect{z}}\left( \vect{\alpha}\right)   \frac{\mathrm{d}P^{\left(Q, \lambda_2 \right)}_{\vect{\Theta}| \vect{Z} = \vect{z}}}{\mathrm{d}Q} \left( \vect{\alpha} \right)  \mathrm{d} Q\left(\vect{\alpha} \right),\supersqueezeequ\IEEEeqnarraynumspace
\end{IEEEeqnarray}
which implies the following inclusions:
\begin{subequations}\label{EqKeyInclusions}
\begin{IEEEeqnarray}{rcl}
\set{R}_1(\lambda_2) & \subseteq & \set{R}_1(\lambda_1), \mbox{ and }\\
\label{EqKeyFact101}
\set{R}_2(\lambda_1) & \subseteq & \set{R}_2(\lambda_2),
\end{IEEEeqnarray}
\end{subequations}
with the sets~$\set{R}_1(\cdot)$ and~$\set{R}_2(\cdot)$ in~\eqref{EqSetR}.
From~\eqref{EqSetN}, it holds that for all~$i \in \lbrace 1,2 \rbrace$,  
\begin{IEEEeqnarray}{rcl}
\label{EqKeyFact102}
 \set{N}_{Q, \vect{z}}(\lambda_i) = \set{R}_2 (\lambda_i)^{\sf{c}},
\end{IEEEeqnarray}
where the complement is with respect to~$\set{M}$. 
Thus, the inclusion in~\eqref{EqKeyFact101} and the equality in~\eqref{EqKeyFact102} yields,
\begin{IEEEeqnarray}{rcl}
\label{EqKeyFact201}
\set{N}_{Q, \vect{z}} (\lambda_1) \supseteq \set{N}_{Q, \vect{z}} (\lambda_2). 
\end{IEEEeqnarray}
The inclusion~$\set{M} \supseteq \set{N}_{Q, \vect{z}}(\lambda_1)$ follows from~\eqref{EqSetN}.
Alternatively, the inclusion~$\set{N}_{Q, \vect{z}}(\lambda_2) \supseteq \set{N}_{Q, \vect{z}}^{\star}$, follows from Lemma~\ref{LemmaRzLowerBoundB} and from observing that for all~$\vect{\nu} \in  \set{N}_{Q, \vect{z}}^{\star}$, 
\begin{IEEEeqnarray}{rcl}
\mathsf{R}_{\vect{z}}\left(P^{\left(Q, \lambda_2 \right)}_{\vect{\Theta} | \vect{Z} = \vect{z}} \right) 
&  \geqslant  & \delta_{Q,\vect{z}}^{\star}=   \mathsf{L}_{\vect{z}}\left( \vect{\nu}\right),
\end{IEEEeqnarray}
which implies that~$\vect{\nu} \in \set{N}_{Q, \vect{z}}(\lambda_2)$.
This completes the proof of~\eqref{EqDecreasingsets}.

The proof of~\eqref{EqNavir2315} is as follows. From the intermediate value theorem \cite[Theorem $4.23$]{rudin1953bookPrinciples} and the assumption that the empirical risk function~$\mathsf{L}_{\vect{z}}$ in~\eqref{EqLxy} is continuous on~$\set{M}$, it follows that for all~$\lambda \in \set{K}_{Q, \vect{z}}$,  there always exists a model~$\vect{\theta} \in \set{M}$, such that
\begin{IEEEeqnarray}{rcl}
 \mathsf{L}_{\vect{z}}\left( \vect{\theta}\right)  =  \displaystyle\int\mathsf{L}_{\vect{z}}\left( \vect{\alpha}\right)    \mathrm{d} P^{\left(Q, \lambda\right)}_{\vect{\Theta} | \vect{Z} = \vect{z}}\left(\vect{\alpha} \right),
\end{IEEEeqnarray}
which implies that~$\set{R}_0\left( \lambda \right)$ is not empty,  and as a consequence,~$\set{N}_{Q, \vect{z}}\left( \lambda\right) = \set{R}_0\left( \lambda\right) \cup \set{R}_1\left( \lambda\right)$ is not empty. 
Hence,  for all~$\vect{\theta} \in \set{R}_0\left( \lambda_1\right)$ it holds that~$\vect{\theta} \notin \set{N}_{Q, \vect{z}}\left( \lambda_2 \right)$. This proves that the elements of~$\set{R}_0\left( \lambda_1\right)$ are in~$\set{N}_{Q, \vect{z}}\left( \lambda_1\right)$ but not in~$\set{N}_{Q, \vect{z}}\left( \lambda_2\right)$.
This, together with~\eqref{EqKeyFact201}, verifies that 
\begin{IEEEeqnarray}{rcl}
\set{N}_{Q, \vect{z}}\left( \lambda_1 \right)  &\supset& \set{N}_{Q, \vect{z}}\left( \lambda_{2} \right).
\end{IEEEeqnarray}
The strict inclusion~$\set{M} \supset \set{N}_{Q, \vect{z}}(\lambda_1)$ is proved by contradiction. Assume that there exists a~$\lambda \in \set{K}_{Q, \vect{z}}$ such that~$\set{M} = \set{N}_{Q, \vect{z}}(\lambda)$.  Then,~$\set{R}_2\left( \lambda \right) = \emptyset$ and thus,~$P^{\left(Q, \lambda\right)}_{\vect{\Theta} | \vect{Z} = \vect{z}}\left(\set{R}_2\left( \lambda \right)\right) =0$, which together with Lemma~\ref{LemmaR1R2}, implies  that $P^{\left(Q, \lambda\right)}_{\vect{\Theta} | \vect{Z} = \vect{z}}\left(\set{R}_1\left( \lambda \right)\right) =0$ and consequently, 
\begin{equation}
P^{\left(Q, \lambda\right)}_{\vect{\Theta} | \vect{Z} = \vect{z}}\left(\set{R}_0\left( \lambda \right)\right) =1.
\end{equation}
This contradicts the assumption that the function~$\mathsf{L}_{\vect{z}}$ is separable (Definition~\ref{DefSeparableLxy}). Hence,~$\set{M} \supset \set{N}_{Q, \vect{z}}(\lambda_1)$.

Finally, the strict inclusion~$\set{N}_{Q, \vect{z}}(\lambda_2) \supset \set{N}_{Q, \vect{z}}^{\star}$ is proved by contradiction. 
Assume that there exists a~$\lambda \in \set{K}_{Q, \vect{z}}$ such that~$\set{N}_{Q, \vect{z}}^{\star}= \set{N}_{Q, \vect{z}}(\lambda)$. That is, 
\begin{IEEEeqnarray}{lCl}
\nonumber
& &\left\lbrace \vect{\theta} \in \set{M}: \mathsf{L}_{\vect{z}}\left( \vect{\theta}\right)  \leqslant    \delta_{Q,\vect{z}}^{\star}  \right \rbrace \\
& = & \set{N}_{Q, \vect{z}}^{\star} \\
& = & \set{N}_{Q, \vect{z}}(\lambda) \\
& =  & \left\lbrace \vect{\theta} \in \set{M}: \mathsf{L}_{\vect{z}}\left( \vect{\theta}\right)  \leqslant   K^{(1)}_{Q, \vect{z}}\left(-\,\frac{1}{\lambda} \right)\right \rbrace.\IEEEeqnarraynumspace
\end{IEEEeqnarray}

Hence,  three cases might arise: \newline
$(a)$ there exists a~$\lambda\in \set{K}_{Q, \vect{z}}$, such that~$\delta_{Q,\vect{z}}^{\star} <  K^{(1)}_{Q, \vect{z}}\left(-\frac{1}{\lambda} \right)$ and it holds that~$$\left\lbrace \vect{\nu} \in \set{M}: \delta_{Q,\vect{z}}^{\star} < \mathsf{L}_{\vect{z}}\left( \vect{\nu}\right)  \leqslant  K^{(1)}_{Q, \vect{z}}\left(-\frac{1}{\lambda} \right) \right \rbrace = \emptyset;$$
$(b)$ there exists a~$\lambda\in \set{K}_{Q, \vect{z}}$, such that~$\delta_{Q,\vect{z}}^{\star} >  K^{(1)}_{Q, \vect{z}}\left(-\frac{1}{\lambda} \right)$ and it holds that~$$\left\lbrace \vect{\nu} \in \set{M}: K^{(1)}_{Q, \vect{z}}\left(-\frac{1}{\lambda}  \right) < \mathsf{L}_{\vect{z}}\left( \vect{\nu}\right)  \leqslant  \delta_{Q,\vect{z}}^{\star} \right \rbrace = \emptyset;$$ 
or 
$(c)$ there exists a~$\lambda\in \set{K}_{Q, \vect{z}}$, such that~$\delta_{Q,\vect{z}}^{\star} =  K^{(1)}_{Q, \vect{z}}\left(-\frac{1}{\lambda} \right)$. 

The cases~$(a)$ and~$(b)$ are absurd. Hence, the proof is complete only by considering the case~$(c)$. In the case~$(c)$, it holds that,
\begin{IEEEeqnarray}{rcl}
\set{R}_{1}\left( \lambda \right)
& = & \left\lbrace \vect{\nu} \in \set{M}: \mathsf{L}_{\vect{z}}\left( \vect{\nu}\right)  < \delta_{Q,\vect{z}}^{\star} \right \rbrace,
\end{IEEEeqnarray}
and from the definition of~$\delta_{Q,\vect{z}}^{\star}$ in~\eqref{EqDeltaStar}, it holds that
\begin{equation}\label{EqCanapeLit}
P^{\left(Q, \lambda\right)}_{\vect{\Theta} | \vect{Z} = \vect{z}} \left(  \set{R}_1\left( \lambda \right) \right) = 0.
\end{equation}
From Lemma~\ref{LemmaR1R2} and~\eqref{EqCanapeLit}, it follows that,
\begin{equation}\label{EqCanapeLit9876}
P^{\left(Q, \lambda\right)}_{\vect{\Theta} | \vect{Z} = \vect{z}}\left(  \set{R}_2\left( \lambda \right) \right) = 0.
\end{equation}
Finally, by noticing that 
\begin{IEEEeqnarray}{rcl}
\nonumber
1 & = & P^{\left(Q, \lambda\right)}_{\vect{\Theta} | \vect{Z} = \vect{z}}\left(  \set{R}_0\left( \lambda \right) \right)  + P^{\left(Q, \lambda\right)}_{\vect{\Theta} | \vect{Z} = \vect{z}}\left(  \set{R}_1\left( \lambda \right) \right)  \\
& & + P^{\left(Q, \lambda\right)}_{\vect{\Theta} | \vect{Z} = \vect{z}}\left(  \set{R}_2\left( \lambda \right) \right) \\
& = & P^{\left(Q, \lambda\right)}_{\vect{\Theta} | \vect{Z} = \vect{z}}\left(  \set{R}_0\left( \lambda \right) \right),
\end{IEEEeqnarray}
reveals a contradiction to the assumption that the function~$\mathsf{L}_{\vect{z}}$ is separable with respect to~$P^{\left(Q, \lambda\right)}_{\vect{\Theta} | \vect{Z} = \vect{z}}$ (and thus, separable with respect to~$Q$ by Lemma~\ref{CorAC}).
This completes the proof of~\eqref{EqNavir2315}.

\section{Proof of Theorem~\ref{CorDecreasingProbability}}\label{AppProofCorDecreasingProbability}

The proof of~\eqref{EqIncreasingProbabilityA9} is based on the analysis of the derivative of~$P^{\left(Q, \lambda\right)}_{\vect{\Theta} | \vect{Z} = \vect{z}} \left( \set{A} \right)$  with respect to~$\lambda$, for some fixed set~$\set{A} \subseteq \BorSigma{\set{M}}$. 
More specifically, given a~$\gamma \in \set{K}_{Q, \vect{z}}$, it holds that 
\begin{IEEEeqnarray}{rcl}
\label{EqMountainGray}
P^{\left(Q, \gamma \right)}_{\vect{\Theta} | \vect{Z} = \vect{z}} \left( \set{A} \right)  & = & \int_{\set{A}}  \frac{\mathrm{d}P^{\left(Q, \gamma \right)}_{\vect{\Theta} | \vect{Z} = \vect{z}}}{\mathrm{d}Q} \left( \vect{\alpha} \right) \mathrm{d}Q \left( \vect{\alpha} \right),
\end{IEEEeqnarray}
and from the fundamental theorem of calculus \cite[Theorem~$6.21$]{rudin1953bookPrinciples}, it follows that for all~$(\lambda_1, \lambda_2) \in \set{K}_{Q, \vect{z}}\times\set{K}_{Q, \vect{z}}$ with~$\lambda_1 > \lambda_2$,
\begin{IEEEeqnarray}{rcl}
\nonumber
& & P^{\left(Q, \lambda_1 \right)}_{\vect{\Theta}| \vect{Z} = \vect{z}} \left( \set{A}\right) - P^{\left(Q, \lambda_2 \right)}_{\vect{\Theta}| \vect{Z} = \vect{z}} \left( \set{A} \right) \\
& = & \int_{\lambda_{2}}^{\lambda_{1}}   \frac{\mathrm{d}}{\mathrm{d}\gamma}P^{\left(Q, \gamma \right)}_{\vect{\Theta} | \vect{Z} = \vect{z}} \left( \set{A} \right) \mathrm{d}\gamma\\
\label{EqBisson9876}
& = &  \int_{\lambda_{2}}^{\lambda_{1}}  \frac{\mathrm{d}}{\mathrm{d}\gamma}  \int_{\set{A}} \frac{\mathrm{d}P^{\left(Q, \gamma \right)}_{\vect{\Theta} | \vect{Z} = \vect{z}}}{\mathrm{d}Q} \left( \vect{\alpha} \right) \mathrm{d}Q \left( \vect{\alpha} \right) \mathrm{d}\gamma\\
\label{EqBisson9877} 
& = &  \int_{\lambda_{2}}^{\lambda_{1}}  \int_{\set{A}}  \frac{\mathrm{d}}{\mathrm{d}\gamma} \frac{\mathrm{d}P^{\left(Q, \gamma \right)}_{\vect{\Theta} | \vect{Z} = \vect{z}}}{\mathrm{d}Q} \left( \vect{\alpha} \right) \mathrm{d}Q \left( \vect{\alpha} \right) \mathrm{d}\gamma,
\end{IEEEeqnarray}
where 
the equality in~\eqref{EqBisson9876} follows from~\eqref{EqMountainGray}; and
the equality in~\eqref{EqBisson9877} holds  from Lemma~\ref{CorPositive} and the dominated convergence theorem \cite[Theorem~$1.6.9$]{ash2000probability}. 

For all~$\vect{\theta} \in \supp Q$, the following holds,
\begin{IEEEeqnarray}{rcl}
\nonumber
& &\frac{\mathrm{d}}{\mathrm{d} \lambda} \frac{\mathrm{d}P^{\left(Q, \lambda\right)}_{\vect{\Theta} | \vect{Z} = \vect{z}}}{\mathrm{d}Q} \left( \vect{\theta} \right) \\
& = & \frac{\mathrm{d}}{\mathrm{d} \lambda}
\frac{\exp\left( -\frac{\mathsf{L}_{\vect{z}}\left( \vect{\theta}\right)}{\lambda}\right)}{\displaystyle\int\exp\left( -\frac{\mathsf{L}_{\vect{z}}\left( \vect{\nu} \right)}{\lambda} \right)  \mathrm{d} Q\left(\vect{\nu} \right)}\\
\nonumber
& = &\frac{\frac{1}{\lambda^2}\mathsf{L}_{\vect{z}}\left( \vect{\theta}\right) \exp\left( -\frac{\mathsf{L}_{\vect{z}}\left( \vect{\theta}\right)}{\lambda}\right)}{\displaystyle\int\exp\left( -\frac{\mathsf{L}_{\vect{z}}\left( \vect{\nu} \right)}{\lambda} \right)  \mathrm{d} Q\left(\vect{\nu} \right)} \\
& - &  \frac{\frac{1}{\lambda^2} \exp\left( -\frac{\mathsf{L}_{\vect{z}}\left( \vect{\theta}\right)}{\lambda}\right) \displaystyle\int\mathsf{L}_{\vect{z}}\left( \vect{\alpha}\right)  \exp\left( -\frac{\mathsf{L}_{\vect{z}}\left( \vect{\alpha} \right)}{\lambda} \right)  \mathrm{d} Q\left(\vect{\alpha} \right)}{\left( \displaystyle\int\exp\left( -\frac{\mathsf{L}_{\vect{z}}\left( \vect{\nu} \right)}{\lambda} \right)  \mathrm{d} Q\left(\vect{\nu} \right) \right)^2}\\
\nonumber
& = & \frac{1}{\lambda^2}\mathsf{L}_{\vect{z}}\left( \vect{\theta}\right)  \frac{\mathrm{d}P^{\left(Q, \lambda\right)}_{\vect{\Theta} | \vect{Z} = \vect{z}}}{\mathrm{d}Q} \left( \vect{\theta} \right) \\
& -  &\frac{1}{\lambda^2}\frac{\mathrm{d}P^{\left(Q, \lambda\right)}_{\vect{\Theta} | \vect{Z} = \vect{z}}}{\mathrm{d}Q} \left( \vect{\theta} \right) \displaystyle\int\mathsf{L}_{\vect{z}}\left( \vect{\nu}\right)  \frac{\mathrm{d}P^{\left(Q, \lambda\right)}_{\vect{\Theta} | \vect{Z} = \vect{z}}}{\mathrm{d}Q} \left( \vect{\nu} \right) \mathrm{d} Q\left(\vect{\nu} \right)\\
\label{EqDerivativeRND}
& = & \frac{1}{\lambda^2}  \frac{\mathrm{d}P^{\left(Q, \lambda\right)}_{\vect{\Theta} | \vect{Z} = \vect{z}}}{\mathrm{d}Q} \left( \vect{\theta} \right) \left(\mathsf{L}_{\vect{z}}\left( \vect{\theta}\right)  - \displaystyle\int\mathsf{L}_{\vect{z}}\left( \vect{\nu}\right)    \mathrm{d} P^{\left(Q, \lambda\right)}_{\vect{\Theta} | \vect{Z} = \vect{z}}\left(\vect{\nu} \right) \right).\squeezeequ\IEEEeqnarraynumspace
\end{IEEEeqnarray}
Plugging~\eqref{EqDerivativeRND} into~\eqref{EqBisson9877} yields,
\begin{IEEEeqnarray}{lll}
\nonumber
& & P^{\left(Q, \lambda_1 \right)}_{\vect{\Theta}| \vect{Z} = \vect{z}} \left( \set{A}\right) - P^{\left(Q, \lambda_2 \right)}_{\vect{\Theta}| \vect{Z} = \vect{z}} \left( \set{A} \right) \\
\nonumber
& = & \int_{\lambda_{2}}^{\lambda_{1}}  \int_{\set{A}}    \frac{1}{\gamma^2} \frac{\mathrm{d}P^{\left(Q, \gamma \right)}_{\vect{\Theta} | \vect{Z} = \vect{z}}}{\mathrm{d}Q} \left( \vect{\alpha} \right) \\
& & \left(\mathsf{L}_{\vect{z}}\left( \vect{\alpha}\right)  - \displaystyle\int\mathsf{L}_{\vect{z}}\left( \vect{\nu}\right)    \mathrm{d} P^{\left(Q, \gamma \right)}_{\vect{\Theta} | \vect{Z} = \vect{z}}\left(\vect{\nu} \right) \right)  \mathrm{d}Q \left( \vect{\alpha} \right) \mathrm{d}\gamma  \\
\nonumber
& = &\int_{\lambda_{2}}^{\lambda_{1}}  \int_{\set{A}} \frac{1}{\gamma^2}   \left(\mathsf{L}_{\vect{z}}\left( \vect{\alpha}\right)  - \displaystyle\int\mathsf{L}_{\vect{z}}\left( \vect{\nu}\right)    \mathrm{d} P^{\left(Q, \gamma \right)}_{\vect{\Theta} | \vect{Z} = \vect{z}}\left(\vect{\nu} \right) \right)  \mathrm{d} P^{\left(Q, \gamma \right)}_{\vect{\Theta} | \vect{Z} = \vect{z}}\left( \vect{\alpha} \right)  \mathrm{d}\gamma. \supersqueezeequ\\
\label{EqBissonite97}
\end{IEEEeqnarray}
 Note that for all~$\vect{\alpha} \in \set{N}_{Q, \vect{z}}\left( \lambda_2\right)$, it holds that for all~$\gamma \in (\lambda_2, \lambda_1)$,  
 \begin{equation}
\mathsf{L}_{\vect{z}}\left( \vect{\alpha}\right)  - \displaystyle\int\mathsf{L}_{\vect{z}}\left( \vect{\nu}\right)    \mathrm{d} P^{\left(Q, \gamma \right)}_{\vect{\Theta} | \vect{Z} = \vect{z}}\left(\vect{\nu} \right)  \leqslant 0, 
 \end{equation}
and thus, 
\begin{IEEEeqnarray}{rcl}
\nonumber
 \int_{\set{N}_{Q, \vect{z}}\left( \lambda_2 \right)}  \frac{1}{\gamma^2}   \left(\mathsf{L}_{\vect{z}}\left( \vect{\alpha}\right)  - \displaystyle\int\mathsf{L}_{\vect{z}}\left( \vect{\nu}\right)    \mathrm{d} P^{\left(Q, \gamma \right)}_{\vect{\Theta} | \vect{Z} = \vect{z}}\left(\vect{\nu} \right) \right)  \mathrm{d} P^{\left(Q, \gamma \right)}_{\vect{\Theta} | \vect{Z} = \vect{z}}\left( \vect{\alpha} \right) \supersqueezeequ & \leqslant & 0.\supersqueezeequ\\
 \label{EqMountainBlue1}
\end{IEEEeqnarray}
The equalities in~\eqref{EqBissonite97} and~\eqref{EqMountainBlue1}, with~$
\set{A} = \set{N}_{Q, \vect{z}}\left( \lambda \right)$, imply that 
\begin{IEEEeqnarray}{lll}
\label{EqBissonNumber}
P^{\left(Q, \lambda_1 \right)}_{\vect{\Theta}| \vect{Z} = \vect{z}} \left( \set{N}_{Q, \vect{z}}\left( \lambda_2 \right)\right) - P^{\left(Q, \lambda_2 \right)}_{\vect{\Theta}| \vect{Z} = \vect{z}} \left( \set{N}_{Q, \vect{z}}\left( \lambda_2 \right) \right)
& \leqslant & 0.
\end{IEEEeqnarray}

The inequality~$0< P^{\left(Q, \lambda_1 \right)}_{\vect{\Theta}| \vect{Z} = \vect{z}}(\set{N}_{Q, \vect{z}}(\lambda_{2}) )$ in~\eqref{EqIncreasingProbabilityA9} is proved by contradiction. Assume that for some~$\lambda \in \set{K}_{Q, \vect{z}}$
it holds that~$0 =P^{\left(Q, \lambda\right)}_{\vect{\Theta} | \vect{Z} = \vect{z}}(\set{N}_{Q, \vect{z}}(\lambda_{2}))$. Then, $P^{\left(Q, \lambda\right)}_{\vect{\Theta} | \vect{Z} = \vect{z}}(\set{R}_0(\lambda_{2}))$ $+$ $P^{\left(Q, \lambda\right)}_{\vect{\Theta} | \vect{Z} = \vect{z}}(\set{R}_1(\lambda_{2})) = 0$, which implies that~$P^{\left(Q, \lambda\right)}_{\vect{\Theta} | \vect{Z} = \vect{z}}(\set{R}_2(\lambda_{2}))=1$, which is a contradiction. See for instance, 
Lemma~\ref{LemmaGenR1R2}.
This completes the proof of~\eqref{EqIncreasingProbabilityA9}.

The proof of strict inequality in~\eqref{EqIncreasingProbabilityA9} is divided into two parts. 
The first part shows that if for all pairs~$(\lambda_1, \lambda_2)  \in \set{K}_{Q, \vect{z}}\times\set{K}_{Q, \vect{z}}$ with~$\lambda_1 > \lambda_2$,
\begin{IEEEeqnarray}{rcl}
\label{EqPonny}
 P^{\left(Q, \lambda_1 \right)}_{\vect{\Theta}| \vect{Z} = \vect{z}}(\set{N}_{Q, \vect{z}}(\lambda_2) )  < P^{\left(Q, \lambda_2 \right)}_{\vect{\Theta}| \vect{Z} = \vect{z}}(\set{N}_{Q, \vect{z}}(\lambda_2) ),
\end{IEEEeqnarray}
then the function~$\mathsf{L}_{\vect{z}}$ is separable with respect to~$Q$.
The second part of the proof shows that if the function~$\mathsf{L}_{\vect{z}}$ is separable with respect to~$Q$, then, for all pairs~$(\lambda_1, \lambda_2)  \in \set{K}_{Q, \vect{z}}\times\set{K}_{Q, \vect{z}}$ with~$\lambda_1 > \lambda_2$, the inequality in~\eqref{EqPonny} holds. 

The first part is as follows.  In the proof of Theorem~\ref{CorDecreasingSet} it is shown (see~\eqref{EqBissonite97})  that for all pairs~$(\lambda_1, \lambda_2)  \in \set{K}_{Q, \vect{z}}\times\set{K}_{Q, \vect{z}}$ with~$\lambda_1 > \lambda_2$,
\begin{IEEEeqnarray}{lll}
\nonumber
& & P^{\left(Q, \lambda_1 \right)}_{\vect{\Theta}| \vect{Z} = \vect{z}} \left( \set{N}_{Q, \vect{z}} \left( \lambda_2\right) \right) - P^{\left(Q, \lambda_2 \right)}_{\vect{\Theta}| \vect{Z} = \vect{z}} \left(  \set{N}_{Q, \vect{z}} \left( \lambda_2\right) \right) \\
\nonumber
& = & \int_{\lambda_2}^{\lambda_1}  \int_{ \set{N}_{Q, \vect{z}} \left( \lambda_2\right)} \frac{1}{\gamma^2}  \left(\mathsf{L}_{\vect{z}}\left( \vect{\alpha}\right)  - K^{(1)}_{Q, \vect{z}} \left(-\,\frac{1}{\gamma} \right) \right)  \mathrm{d} P^{\left(Q, \gamma \right)}_{\vect{\Theta} | \vect{Z} = \vect{z}}\left( \vect{\alpha} \right)  \mathrm{d}\gamma.\supersqueezeequ\\
\label{EqHivaOa}
\end{IEEEeqnarray}

Assume that for a given pair~$(\lambda_1, \lambda_2)  \in \set{K}_{Q, \vect{z}}\times\set{K}_{Q, \vect{z}}$, with~$\lambda_1 > \lambda_2$, the inequality in~\eqref{EqPonny} holds. Then, from~\eqref{EqHivaOa},
\begin{IEEEeqnarray}{lll}
\nonumber
& & 0 \\
\nonumber
 &>& \int_{\lambda_2}^{\lambda_1}  \int_{ \set{N}_{Q, \vect{z}} \left( \lambda_2\right)} \frac{1}{\gamma^2}  \left(\mathsf{L}_{\vect{z}}\left( \vect{\alpha}\right)  - K^{(1)}_{Q, \vect{z}} \left(-\frac{1}{\gamma} \right)  \right)  \mathrm{d} P^{\left(Q, \gamma \right)}_{\vect{\Theta} | \vect{Z} = \vect{z}}\left( \vect{\alpha} \right)  \mathrm{d}\gamma \squeezeequ\\
 \nonumber
& = &  \int_{\lambda_2}^{\lambda_1}  \int_{\set{R}_1 \left( \lambda_2\right)} \frac{1}{\gamma^2}  \left(\mathsf{L}_{\vect{z}}\left( \vect{\alpha}\right)  - K^{(1)}_{Q, \vect{z}} \left(-\frac{1}{\gamma} \right)\right)  \mathrm{d} P^{\left(Q, \gamma \right)}_{\vect{\Theta} | \vect{Z} = \vect{z}}\left( \vect{\alpha} \right)  \mathrm{d}\gamma, \squeezeequ \\
\label{EqFenuaEnata}
\end{IEEEeqnarray}
where the equality in~\eqref{EqFenuaEnata} follows from noticing that~$\set{R}_0\left( \lambda_2\right)$ and~$\set{R}_1\left( \lambda_2 \right)$ form a partition of~$\set{N}_{Q, \vect{z}}\left( \lambda_2 \right)$, with the sets~$\set{R}_0\left( \lambda_2\right)$,~$\set{R}_1\left( \lambda_2\right)$ and~$\set{N}_{Q, \vect{z}}\left( \lambda_2\right)$ defined in~\eqref{EqSetR0},~\eqref{EqSetR1}, and~\eqref{EqSetN}, respectively.

The inequality in~\eqref{EqFenuaEnata} implies that the set~$\set{R}_1\left( \lambda_2 \right)$ is nonnegligible with respect to ~$P^{\left(Q, \gamma \right)}_{\vect{\Theta} | \vect{Z} = \vect{z}}$, for some~$\gamma \in (\lambda_2, \lambda_1)$.
Hence, from Lemma~\ref{LemmaGenR1R2}, it follows that both sets~$\set{R}_1\left( \lambda_2 \right)$  and~$\set{R}_2\left( \lambda_2 \right)$ are nonnegligible with respect to~$P^{\left(Q, \gamma \right)}_{\vect{\Theta} | \vect{Z} = \vect{z}}$.

From the  arguments above, it has been proved that given a pair~$(\lambda_1, \lambda_2)  \in \set{K}_{Q, \vect{z}}\times\set{K}_{Q, \vect{z}}$ with~$\lambda_1 > \lambda_2$, if 
\begin{IEEEeqnarray}{rcl}
 P^{\left(Q, \lambda_1 \right)}_{\vect{\Theta}| \vect{Z} = \vect{z}}(\set{N}_{Q, \vect{z}}(\lambda_2) )  < P^{\left(Q, \lambda_2 \right)}_{\vect{\Theta}| \vect{Z} = \vect{z}}(\set{N}_{Q, \vect{z}}(\lambda_2) ),
\end{IEEEeqnarray}
then there always exists a positive~$\gamma \in (\lambda_1, \lambda_2)$ such that the sets~$\set{R}_1\left( \lambda_2 \right)$ and~$\set{R}_2\left( \lambda_2 \right)$ are not negligible with respect to~$P^{\left(Q, \gamma \right)}_{\vect{\Theta} | \vect{Z} = \vect{z}}$.
Moreover, such sets~$\set{R}_1\left( \lambda_2 \right)$ and~$\set{R}_2\left( \lambda_2 \right)$  satisfy for all~$(\vect{\nu}_1,\vect{\nu}_2) \in \set{R}_2\left( \lambda \right) \times \set{R}_1\left( \lambda\right)$, 
\begin{IEEEeqnarray}{rcl}
\label{EqTwoNonnegligibleSetsPingPong}
+\infty > \mathsf{L}_{\vect{z}} \left( \vect{\nu}_1 \right) &> K^{(1)}_{Q, \vect{z}}\left( - \frac{1}{\lambda} \right) >& \mathsf{L}_{\vect{z}}\left(\vect{\nu}_2\right),
\end{IEEEeqnarray}
which together with Definition~\ref{CorAC} verify that the function~$\mathsf{L}_{\vect{z}}$ is separable with respect to~$P^{\left(Q, \gamma \right)}_{\vect{\Theta} | \vect{Z} = \vect{z}}$ (and thus, with respect to~$Q$ by Lemma~\ref{CorAC}). This ends the first part of the proof.
 
The second part of the proof is under the assumption that the empirical risk function~$\mathsf{L}_{\vect{z}}$ in~\eqref{EqLxy} is separable with respect to Q (and thus, with respect to~$P^{\left(Q, \gamma \right)}_{\vect{\Theta} | \vect{Z} = \vect{z}}$ by Lemma~\ref{CorAC}). That is, from Definition~\ref{CorAC}, for all~$\gamma \in \set{K}_{Q, \vect{z}}$, there exist a positive real~$c_{\gamma} >0$ and two subsets~$\set{A}(\gamma)$ and~$\set{B}(\gamma)$ of~$\set{M}$ nonnegligible with respect to~$P^{\left(Q, \gamma \right)}_{\vect{\Theta} | \vect{Z} = \vect{z}}$ in~\eqref{EqGenpdf} that verify that for all~$(\vect{\nu}_1,\vect{\nu}_2) \in \set{A}(\gamma) \times \set{B}(\gamma)$, 
\begin{IEEEeqnarray}{rcl}
\label{EqTwoNonnegligibleSets07696}
\mathsf{L}_{\vect{z}} \left( \vect{\nu}_1 \right) &> c_{\gamma} >& \mathsf{L}_{\vect{z}}\left(\vect{\nu}_2\right).
\end{IEEEeqnarray}
In the proof of Theorem~\ref{CorDecreasingSet}, cf.~\eqref{EqBissonite97}, it has been proved that  given a pair $(\alpha_1, \alpha_2)$  $\in$ $\set{K}_{Q, \vect{z}}\times\set{K}_{Q, \vect{z}}$, with~$\alpha_1 >  \gamma > \alpha_2$, it holds that for all subsets~$\set{A}$ of~$\set{M}$, %
\begin{IEEEeqnarray}{lll}
\nonumber
& & P^{\left(Q, \alpha_1 \right)}_{\vect{\Theta} | \vect{Z} = \vect{z}} \left( \set{A}\right) - P^{\left(Q, \alpha_2 \right)}_{\vect{\Theta} | \vect{Z} = \vect{z}} \left( \set{A} \right) \\
\nonumber
& = & \int_{\alpha_{2}}^{\alpha_{1}} \hspace{-1ex} \int_{\set{A}}    \frac{1}{\lambda^2} \frac{\mathrm{d}P^{\left(Q, \lambda\right)}_{\vect{\Theta} | \vect{Z} = \vect{z}}}{\mathrm{d}Q} \left( \vect{\alpha} \right) \left(\mathsf{L}_{\vect{z}}\left( \vect{\alpha}\right)  - K^{(1)}_{Q, \vect{z}} \left(-\frac{1}{\lambda} \right) \right)  \mathrm{d}P \left( \vect{\alpha} \right) \mathrm{d}\lambda \squeezeequ \\
\label{EqIvermectine}
& = &\int_{\alpha_{2}}^{\alpha_{1}}  \hspace{-1ex} \int_{\set{A}} \frac{1}{\lambda^2}   \left(\mathsf{L}_{\vect{z}}\left( \vect{\alpha}\right)  - K^{(1)}_{Q, \vect{z}} \left(-\frac{1}{\lambda} \right) \right)  \mathrm{d} P^{\left(Q, \lambda\right)}_{\vect{\Theta} | \vect{Z} = \vect{z}}\left( \vect{\alpha} \right)  \mathrm{d}\lambda. \squeezeequ \IEEEeqnarraynumspace
\end{IEEEeqnarray}

Hence,  two cases are studied. The first case considers that 
\begin{equation}
\label{EqRainbowRoll}
c_{\gamma}< K^{(1)}_{Q, \vect{z}}\left( - \frac{1}{\gamma} \right),
\end{equation}
 with~$c_{\gamma}$ in~\eqref{EqTwoNonnegligibleSets07696}. The second case considers that \begin{equation}
\label{EqRainbowRoll2345}
c_{\gamma} \geqslant K^{(1)}_{Q, \vect{z}}\left( - \frac{1}{\gamma} \right).
\end{equation}

In the first case, it follows from~\eqref{EqSetN} that
\begin{IEEEeqnarray}{rcl}
\label{EqTilinTilin}
\set{B}\left( \gamma \right) \subset \set{N}_{Q, \vect{z}}\left( \gamma \right),
\end{IEEEeqnarray}
which implies that
\begin{IEEEeqnarray}{rcl}
  \label{EqNuclearIneq0}
    P^{\left(Q, \gamma \right)}_{\vect{\Theta} | \vect{Z} = \vect{z}} \left(  \set{N}_{Q, \vect{z}}\left( \gamma \right) \right) & \geqslant & P^{\left(Q, \gamma \right)}_{\vect{\Theta} | \vect{Z} = \vect{z}} \left( \set{B} \left( \gamma \right) \right) \\
  \label{EqNuclearIneq}
  & > & 0,
\end{IEEEeqnarray}
where, 
the inequality in~\eqref{EqNuclearIneq} follows from the fact that~$ \set{B} \left( \gamma \right)~$ is nonnegligible with respect to~$P^{\left(Q, \gamma \right)}_{\vect{\Theta} | \vect{Z} = \vect{z}}$.
This implies that the set~$\set{N}_{Q, \vect{z}}\left( \gamma \right)$ is not negligible with respect~$P^{\left(Q, \gamma \right)}_{\vect{\Theta} | \vect{Z} = \vect{z}}$. 
Moreover, from~\eqref{EqSetN} and~\eqref{EqTilinTilin}, it follows that for all~$\vect{\alpha} \in \set{N}_{Q, \vect{z}}\left( \gamma \right)$ and for all~$\lambda \in (\gamma, \alpha_1)$,  
\begin{IEEEeqnarray}{rcl}
\label{EqSushiA}
\mathsf{L}_{\vect{z}}\left( \vect{\alpha}\right)  - \displaystyle\int\mathsf{L}_{\vect{z}}\left( \vect{\nu}\right)    \mathrm{d} P^{\left(Q, \lambda\right)}_{\vect{\Theta} | \vect{Z} = \vect{z}}\left(\vect{\nu} \right) &< & \mathsf{L}_{\vect{z}}\left( \vect{\alpha}\right)  - c_{\gamma} \\
\label{EqSushi}
&< &0, 
\end{IEEEeqnarray}
where the inequality in~\eqref{EqSushiA} follows from~\eqref{EqRainbowRoll}; and 
the inequality in~\eqref{EqSushi} follows from~\eqref{EqTwoNonnegligibleSets07696}.
Thus, 
\begin{IEEEeqnarray}{rcl}
\nonumber
 \int_{\gamma}^{\alpha_{1}}  \int_{\set{N}_{Q, \vect{z}}(\gamma)} \frac{1}{\lambda^2}   \left(\mathsf{L}_{\vect{z}}\left( \vect{\alpha}\right)  - K^{(1)}_{Q, \vect{z}} \left(-\frac{1}{\lambda} \right)  \right)  \mathrm{d} P^{\left(Q, \lambda\right)}_{\vect{\Theta} | \vect{Z} = \vect{z}}\left( \vect{\alpha} \right)  \mathrm{d}\lambda \squeezeequ & < & 0,
\end{IEEEeqnarray}
which implies, from~\eqref{EqIvermectine}, that  
\begin{IEEEeqnarray}{rcl}
P^{\left(Q, \alpha_1 \right)}_{\vect{\Theta} | \vect{Z} = \vect{z}} \left( \set{N}_{Q, \vect{z}}\left( \gamma \right)\right) - P^{\left(Q, \gamma \right)}_{\vect{\Theta} | \vect{Z} = \vect{z}} \left( \set{N}_{Q, \vect{z}}\left( \gamma \right) \right) &<& 0. 
\end{IEEEeqnarray}
Assume now that~$c_{\gamma} \geqslant K^{(1)}_{Q, \vect{z}}\left( - \frac{1}{\gamma} \right)$. Hence, the following holds \begin{IEEEeqnarray}{rcl}
\label{EqTilinTilin143}
A\left( \gamma \right)  \subseteq \set{R}_{2}\left( \gamma \right),
\end{IEEEeqnarray}
which implies that
\begin{IEEEeqnarray}{rcl}
  \label{EqNuclearIneq0987}
    P^{\left(Q, \gamma \right)}_{\vect{\Theta} | \vect{Z} = \vect{z}} \left(  \set{R}_{2}\left( \gamma \right) \right) & \geqslant & P^{\left(Q, \gamma \right)}_{\vect{\Theta} | \vect{Z} = \vect{z}} \left( \set{A} \left( \gamma \right) \right) \\
  \label{EqNuclearIneq975}
  & > & 0,
\end{IEEEeqnarray}
where 
the inequality in~\eqref{EqNuclearIneq975} follows from the fact that~$ \set{A} \left( \gamma \right)~$ is nonnegligible with respect to~$P^{\left(Q, \gamma \right)}_{\vect{\Theta} | \vect{Z} = \vect{z}}$.
This implies that the set~$\set{R}_{2}\left( \gamma \right)$ is not negligible with respect~$P^{\left(Q, \gamma \right)}_{\vect{\Theta} | \vect{Z} = \vect{z}}$. From Lemma~\ref{LemmaR1R2}, it follows that both~$\set{R}_{1}\left( \gamma \right)$ and~$\set{R}_{2}\left( \gamma \right)$ are nonnegligible with respect to~$P^{\left(Q, \gamma \right)}_{\vect{\Theta} | \vect{Z} = \vect{z}}$. 
Using this result, the following holds,
\begin{IEEEeqnarray}{rcl}
P^{\left(Q, \gamma \right)}_{\vect{\Theta} | \vect{Z} = \vect{z}} \left( \set{N}_{Q, \vect{z}}\left( \gamma \right) \right) & \geqslant & P^{\left(Q, \gamma \right)}_{\vect{\Theta} | \vect{Z} = \vect{z}} \left( \set{R}_1\left( \gamma \right) \right) \\
& > & 0,
\end{IEEEeqnarray}
which proves  the set~$\set{N}_{Q, \vect{z}}\left( \gamma \right)$ is nonnegligible with respect to~$P^{\left(Q, \gamma \right)}_{\vect{\Theta} | \vect{Z} = \vect{z}}$.

From~\eqref{EqSetN} and Theorem~\ref{CorDecreasingAverage}, it follows that for all~$\vect{\alpha} \in \set{N}_{Q, \vect{z}}\left( \gamma \right)$ and for all~$\lambda \in (\gamma, \alpha_1)$,  
\begin{IEEEeqnarray}{rcl}
\label{EqSpicyTunaRoll}
0 &\geqslant &\mathsf{L}_{\vect{z}}\left( \vect{\alpha}\right)  - \displaystyle\int\mathsf{L}_{\vect{z}}\left( \vect{\nu}\right)    \mathrm{d} P^{\left(Q, \gamma\right)}_{\vect{\Theta} | \vect{Z} = \vect{z}}\left(\vect{\nu} \right)  \\
\label{EqSpicyTunaRollRock}
& > & \mathsf{L}_{\vect{z}}\left( \vect{\alpha}\right)  - \displaystyle\int\mathsf{L}_{\vect{z}}\left( \vect{\nu}\right)    \mathrm{d}  P^{\left(Q, \lambda\right)}_{\vect{\Theta} | \vect{Z} = \vect{z}}\left(\vect{\nu} \right).
\end{IEEEeqnarray}
Thus, 
\begin{IEEEeqnarray}{rcl}
\nonumber
 \int_{\gamma}^{\alpha_1}  \int_{\set{N}_{Q, \vect{z}}(\gamma)} \frac{1}{\lambda^2}   \left(\mathsf{L}_{\vect{z}}\left( \vect{\alpha}\right)  - K^{(1)}_{Q, \vect{z}} \left(-\frac{1}{\lambda} \right) \right)  \mathrm{d} P^{\left(Q, \lambda\right)}_{\vect{\Theta} | \vect{Z} = \vect{z}}\left( \vect{\alpha} \right)  \mathrm{d}\lambda  \squeezeequ & < & 0,
\end{IEEEeqnarray}
which implies, from~\eqref{EqIvermectine}, that  
\begin{IEEEeqnarray}{rcl}
\label{EqCovid19delta}
P^{\left(Q, \alpha_1 \right)}_{\vect{\Theta} | \vect{Z} = \vect{z}} \left( \set{N}_{Q, \vect{z}}\left( \gamma \right)\right) - P^{\left(Q, \gamma \right)}_{\vect{\Theta} | \vect{Z} = \vect{z}} \left( \set{N}_{Q, \vect{z}}\left( \gamma \right) \right) &<& 0. 
\end{IEEEeqnarray}
This completes the proof.

 \section{Proof of Lemma~\ref{CorStrongStrinctDecreasingProbability}}\label{AppProofCorStrongStrinctDecreasingProbability}
The proof is based on the following two observations. First, note  that $\left(\set{N}_{Q, \vect{z}}\left( \lambda_2 \right)\right)^\sfc$ $=$ $\set{R}_2\left( \lambda_{2}\right)$, with the set~$\set{R}_2\left( \cdot \right)$ defined in~\eqref{EqSetR2}.
Second, note that
\begin{IEEEeqnarray}{rcl}
\set{N}_{Q, \vect{z}}\left( \lambda_{1} \right) & = & \set{N}_{Q, \vect{z}}\left( \lambda_{2} \right) \cup \left( \set{N}_{Q, \vect{z}}\left( \lambda_{1}\right) \cap \set{R}_2\left( \lambda_{2}\right)  \right),
\end{IEEEeqnarray}
and the fact that the sets~$ \set{N}_{Q, \vect{z}}\left( \lambda_{2} \right)$ and~$\left( \set{N}_{Q, \vect{z}}\left( \lambda_{1}\right) \cap \set{R}_2\left( \lambda_{2}\right)  \right)$ are disjoint. 
Hence, for all~$i \inCountTwo$,
\begin{IEEEeqnarray}{rcl}
\nonumber
& &P^{(\lambda_i)}_{\vect{\Theta}| \vect{Z} = \vect{z}}(\set{N}_{Q, \vect{z}}(\lambda_1) )  \\
& = & P^{(\lambda_{i})}_{\vect{\Theta}| \vect{Z} = \vect{z}}\bigg( \set{N}_{Q, \vect{z}}\left( \lambda_{2} \right) \cup \left( \set{N}_{Q, \vect{z}}\left( \lambda_{1}\right) \cap \set{R}_2\left( \lambda_{2}\right)  \right) \bigg)\supersqueezeequ\\
\nonumber
& = & P^{(\lambda_{i})}_{\vect{\Theta}| \vect{Z} = \vect{z}}\bigg( \set{N}_{Q, \vect{z}}\left( \lambda_{2} \right)\bigg) \\
\label{EqManuia0824592}
& & + P^{(\lambda_{i})}_{\vect{\Theta}| \vect{Z} = \vect{z}}\bigg( \set{N}_{Q, \vect{z}}\left( \lambda_{1}\right) \cap \set{R}_2\left( \lambda_{2}\right)   \bigg)\\
\label{EqManuia0824593}
& = & P^{(\lambda_{i})}_{\vect{\Theta}| \vect{Z} = \vect{z}}\bigg( \set{N}_{Q, \vect{z}}\left( \lambda_{2} \right)\bigg),  
\end{IEEEeqnarray}
where the equality in~\eqref{EqManuia0824592} follows from Lemma~\ref{LemmaMutualAC} and the equality in~\eqref{EqKeyConditionForStrong}.

Finally, under the assumption that the empirical function~$\mathsf{L}_{\vect{z}}$ in~\eqref{EqLxy}  is separable,  it holds  from Theorem~\ref{CorDecreasingProbability} that
\begin{IEEEeqnarray}{rcl}
\label{EqNsetsEqualA}
P^{\left(Q, \lambda_1 \right)}_{\vect{\Theta}| \vect{Z} = \vect{z}}(\set{N}_{Q, \vect{z}}(\lambda_2) )  < P^{\left(Q, \lambda_2 \right)}_{\vect{\Theta}| \vect{Z} = \vect{z}}(\set{N}_{Q, \vect{z}}(\lambda_2) ).
\end{IEEEeqnarray}
 Plugging~\eqref{EqManuia0824593} into~\eqref{EqNsetsEqualA}, with~$i =1$, yields,
 \begin{IEEEeqnarray}{rcl}
\label{EqNsetsEqualC}
P^{\left(Q, \lambda_1 \right)}_{\vect{\Theta}| \vect{Z} = \vect{z}}(\set{N}_{Q, \vect{z}}(\lambda_1) )  < P^{\left(Q, \lambda_2 \right)}_{\vect{\Theta}| \vect{Z} = \vect{z}}(\set{N}_{Q, \vect{z}}(\lambda_2) ),
\end{IEEEeqnarray}
 and this completes the proof.  
  
\section{Proof of Theorem~\ref{TheoremSensitivityA}}\label{AppProofTheoremSensitivityA}

Consider the following lemma.
\begin{lemma}\label{LemmaSensitivityA}
Given two probability measures~$P_1$ and~$P_2$ over~$\Bormeaspace{\set{M}}$, with~$P_2$ absolutely continuous with respect to~$P_1$, the following holds for all~$\vect{z} \in \left( \set{X} \times \set{Y}\right)^n$,
\begin{IEEEeqnarray}{rcl}
& &  \mathsf{R}_{\vect{z}}\left( P_2 \right) -  \mathsf{R}_{\vect{z}}\left( P_1 \right) \IEEEeqnarraynumspace\\
\nonumber
& \leqslant &  \inf_{t <0} \left( \frac{\KL{P_2}{P_1}+  \log\left( \int \exp\left( t \, \left( \mathsf{L}_{\vect{z}}\left(\vect{\theta}\right)  - \mathsf{R}_{\vect{z}}\left( P_1 \right) \right) \right) \mathrm{d}P_1(\vect{\theta}) \right)  }{t} \right),  \middlesqueezeequ
\end{IEEEeqnarray}
where the function~$\mathsf{L}_{\vect{z}}$ and the functional~$\mathsf{R}_{\vect{z}}$ are defined in~\eqref{EqLxy} and in~\eqref{EqRxy}, respectively.
\end{lemma}
\begin{IEEEproof}
From  \cite[Corollary~$4.15$, Page 100]{boucheron2013book}, it follows that the probability measures~$P_1$ and~$P_2$ in~$\Bormeaspace{\set{M}}$ satisfy the following equality:
\begin{IEEEeqnarray}{rcl}
D\left( P_2\|P_1\right) = \sup_{f} \int f\left(\vect{\theta}\right) \mathrm{d}P_2\left(\vect{\theta}\right) - \log \int \exp\left( f\left(\vect{\theta}\right) \right) \mathrm{d}P_1\left(\vect{\theta}\right),\supersqueezeequ\IEEEeqnarraynumspace
\end{IEEEeqnarray}
where the supremum is over the space of all measurable functions~$f$ with respect to~$\Bormeaspace{\set{M}}$ and~$\Bormeaspace{\reals}$, such that~$\int  \exp \left( f\left(\vect{\theta}\right) \right) \mathrm{d}P_1\left(\vect{\theta}\right) < \infty$.
Hence, for all~$\vect{z} \in \left( \set{X} \times \set{Y}\right)^n$ and for all~$t \in \left( -\infty, 0 \right)$, it follows that the empirical risk function~$\mathsf{L}_{\vect{z}}$ in~\eqref{EqLxy} satisfies that
\begin{IEEEeqnarray}{rcl}
\nonumber
&& D\left( P_2\|P_1\right) \\
\nonumber
& \geqslant & \int t \mathsf{L}_{\vect{z}}\left(\vect{\theta}\right) \mathrm{d}P_2\left(\vect{\theta}\right) - \log \int \exp\left( t \mathsf{L}_{\vect{z}} \left(\vect{\theta}\right) \right) \mathrm{d}P_1\left(\vect{\theta}\right)\\
\nonumber
& \geqslant & \int t \mathsf{L}_{\vect{z}}\left(\vect{\theta}\right) \mathrm{d}P_2\left(\vect{\theta}\right) \\
\nonumber
& & - \log \int \exp\left( t \mathsf{L}_{\vect{z}} \left(\vect{\theta}\right) + t \mathsf{R}_{\vect{z}}\left( P_1 \right) - t\mathsf{R}_{\vect{z}}\left( P_1 \right)\right) \mathrm{d}P_1\left(\vect{\theta}\right)  \\
\nonumber
& = & \int t \mathsf{L}_{\vect{z}}\left(\vect{\theta}\right) \mathrm{d}P_2\left(\vect{\theta}\right) - t\mathsf{R}_{\vect{z}}\left( P_1 \right) \\
\nonumber
& & -  \log \int \exp\left( t \mathsf{L}_{\vect{z}} \left(\vect{\theta}\right) - t \mathsf{R}_{\vect{z}}\left( P_1 \right) \right) \mathrm{d}P_1\left(\vect{\theta}\right)\\
\nonumber
& = & t\mathsf{R}_{\vect{z}}\left( P_2 \right) - t\mathsf{R}_{\vect{z}}\left( P_1 \right)  -  \log \int \exp\left( t \mathsf{L}_{\vect{z}} \left(\vect{\theta}\right) - t \mathsf{R}_{\vect{z}}\left( P_1 \right) \right) \mathrm{d}P_1\left(\vect{\theta}\right),\supersqueezeequ\IEEEeqnarraynumspace
\end{IEEEeqnarray}
which leads to
\begin{IEEEeqnarray}{rcl}
\nonumber
& & \mathsf{R}_{\vect{z}}\left( P_2 \right) - \mathsf{R}_{\vect{z}}\left( P_1 \right) \supersqueezeequ\\
 & \leqslant &\frac{D\left( P_2\|P_1\right)  + \log \int \exp\left( t \left( \mathsf{L}_{\vect{z}} \left(\vect{\theta}\right) -  \mathsf{R}_{\vect{z}}\left( P_1 \right) \right) \right) \mathrm{d}P_1\left(\vect{\theta}\right)}{t}.\supersqueezeequ\IEEEeqnarraynumspace
\end{IEEEeqnarray}
Given that~$t$ can be chosen arbitrarily in~$\left( -\infty, 0 \right)$, it holds that
\begin{IEEEeqnarray}{rcl}
\nonumber
& & \mathsf{R}_{\vect{z}}\left( P_2 \right) - \mathsf{R}_{\vect{z}}\left( P_1 \right)  \supersqueezeequ \\
\nonumber
& \leqslant & \inf_{t \in \left( -\infty, 0 \right) }\frac{D\left( P_2\|P_1\right)  + \log \int \exp\left( t \left( \mathsf{L}_{\vect{z}} \left(\vect{\theta}\right) -  \mathsf{R}_{\vect{z}}\left( P_1 \right) \right) \right) \mathrm{d}P_1\left(\vect{\theta}\right)}{t}, \supersqueezeequ 
\end{IEEEeqnarray}
which completes the proof.
\end{IEEEproof}

From Lemma~\ref{LemmaSensitivityA}, it holds that the probability measure~$P^{\left(Q, \lambda\right)}_{\vect{\Theta} | \vect{Z} = \vect{z}}$ in~\eqref{EqGenpdf}, satisfies for all~$P \in \triangle_{Q}\Bormeaspace{\set{M}}$,
\begin{IEEEeqnarray}{l}
\nonumber 
\mathsf{R}_{\vect{z}}\left( P \right)  -  \mathsf{R}_{\vect{z}}\left( P^{\left(Q, \lambda\right)}_{\vect{\Theta} | \vect{Z} = \vect{z}} \right)   \leqslant     \inf_{t \in \left( -\infty, 0 \right)} \Bigg( \frac{D\left( P\| P^{\left(Q, \lambda\right)}_{\vect{\Theta} | \vect{Z} = \vect{z}} \right)}{t} \\
\label{EqThisIneq53}
+  \frac{ \log\left( \int \exp\left( t \, \left( \mathsf{L}_{\vect{z}}\left(\vect{\theta}\right)  - K^{(1)}_{Q, \vect{z}}\left( - \frac{1}{\lambda} \right)  \right) \right) \mathrm{d}P^{\left(Q, \lambda\right)}_{\vect{\Theta} | \vect{Z} = \vect{z}}(\vect{\theta}) \right)  }{t} \Bigg), \middlesqueezeequ\IEEEeqnarraynumspace
\end{IEEEeqnarray}
where the function~$K^{(1)}_{Q, \vect{z}}$ is defined in~\eqref{EqNiceK1} and satisfies~\eqref{EqRK}.
Moreover, for all~$t \in \left(-\infty, 0 \right)$,
\begin{IEEEeqnarray}{rcl}
\nonumber
& &\log\left( \int \exp\left( t \, \left( \mathsf{L}_{\vect{z}}\left(\vect{\theta}\right)  - K^{(1)}_{Q, \vect{z}}\left( - \frac{1}{\lambda} \right)  \right) \right) \mathrm{d}P^{\left(Q, \lambda\right)}_{\vect{\Theta} | \vect{Z} = \vect{z}}(\vect{\theta}) \right)\\
\label{EqTubuai976}
 & = & \log\left( \int \exp\left( t \, \mathsf{L}_{\vect{z}}\left(\vect{\theta}\right) \right) \mathrm{d}P^{\left(Q, \lambda\right)}_{\vect{\Theta} | \vect{Z} = \vect{z}}(\vect{\theta}) \right)  - t K^{(1)}_{Q, \vect{z}}\left( - \frac{1}{\lambda} \right) \\
\label{EqTubuai977}
& = & J_{\vect{z}, Q, \lambda} (t) - t K^{(1)}_{Q, \vect{z}}\left( - \frac{1}{\lambda} \right)\\
\label{EqTubuai978}
& \leqslant &  \frac{1}{2}t^2  \beta_{Q,\vect{z}}^2 ,
\end{IEEEeqnarray}
where 
the equality in~\eqref{EqTubuai977} follows from~\eqref{EqJflow}; 
the inequality in~\eqref{EqTubuai978} follows from Theorem~\ref{LemmaGaussGen}; and the constant~$\beta_{Q,\vect{z}}$ is defined in~\eqref{EqEatingEarlyToday}.

Plugging~\eqref{EqTubuai978} into~\eqref{EqThisIneq53} yields for all~$t \in \left(-\infty, 0 \right)$,
\begin{IEEEeqnarray}{rcl}
\nonumber
& & \mathsf{R}_{\vect{z}}\left( P \right) - \mathsf{R}_{\vect{z}}\left( P^{\left(Q, \lambda\right)}_{\vect{\Theta} | \vect{Z} = \vect{z}} \right) \\
\label{EqWaitingForDali}
& \leqslant &  \inf_{t \in \left( -\infty, 0 \right)} \frac{D\left( P \| P^{\left(Q, \lambda\right)}_{\vect{\Theta} | \vect{Z} = \vect{z}} \right) + \frac{1}{2}t^2  \beta_{Q,\vect{z}}^2  }{t}. \IEEEeqnarraynumspace
\end{IEEEeqnarray}
Let the~$c \in \reals$ be defined as follows:
\begin{IEEEeqnarray}{cCc}
\label{EqWhatIsThis}
c &\triangleq&  \mathsf{R}_{\vect{z}}\left( P \right) - \mathsf{R}_{\vect{z}}\left( P^{\left(Q, \lambda\right)}_{\vect{\Theta} | \vect{Z} = \vect{z}} \right).
\end{IEEEeqnarray}
Hence, from~\eqref{EqWaitingForDali}, it follows that  for all~$t \in \left(-\infty, 0 \right)$,
\begin{IEEEeqnarray}{rcl}
\label{EqCc}
c \, t -  \frac{1}{2}t^2  \beta_{Q,\vect{z}}^2  \leqslant D\left( P \| P^{\left(Q, \lambda\right)}_{\vect{\Theta} | \vect{Z} = \vect{z}} \right).
\end{IEEEeqnarray}

The rest of the proof consists in finding an explicit expression for the absolute value of~$c$ in \eqref{EqCc}. 
To this aim, consider the function~$\phi: \reals \rightarrow \reals$ such that
\begin{equation}
\label{EqPhi}
\phi(\alpha) =  \frac{1}{2} \alpha^2  \beta_{Q,\vect{z}}^2,
\end{equation}
and note that~$\phi$ is a positive and strictly convex function with~$\phi(0) = 0$.
Let the Legendre-Fenchel transform of~$\phi$ be the function~$\phi^*: \reals \rightarrow \reals$,  and thus for all~$x \in \reals$,
\begin{equation}
\label{EqLFTrans}
\phi^{*}(x) = \max_{t \in \left(-\infty, 0 \right)} x t - \phi(t).
\end{equation}
In particular, note that 
\begin{equation}
\label{EqCrazyStuff23}
\phi^{*}(c) \leqslant D\left( P\| P^{\left(Q, \lambda\right)}_{\vect{\Theta} | \vect{Z} = \vect{z}} \right).
\end{equation}
Note that for all~$x \in \reals$ and for all~$t \in \left(-\infty, 0 \right)$, the function~$\phi^{\star}$ in~\eqref{EqLFTrans} satisfies
\begin{equation}
\label{EqPhiStarfullY}
x \, t -  \frac{1}{2}t^2  \beta_{Q,\vect{z}}^2   \leqslant \phi^{\star}(x) = x \alpha^{\star}(x) - \phi\left( \alpha^{\star}(x) \right),
\end{equation}
where the term~$\alpha^{\star}(x)$ represents the unique solution in~$\alpha$ within the interval~$\left(-\infty,0  \right)$ to
\begin{equation}
\label{EqSimplyDer865}
\frac{\mathrm{d}}{\mathrm{d}\alpha} \left( x \alpha - \phi\left( \alpha \right)  \right) = x - \alpha \beta_{Q,\vect{z}}^2= 0.
\end{equation}
That is, 
\begin{IEEEeqnarray}{rcl}
\label{EqThisDiamond}
\alpha^{\star}(x) & = & \frac{x}{\beta_{Q,\vect{z}}^2}.
\end{IEEEeqnarray}
Plugging~\eqref{EqThisDiamond} into~\eqref{EqPhiStarfullY} yields, 
\begin{IEEEeqnarray}{rcl}
\label{EqHivaOa14}
\phi^{\star}(x) & = &  \frac{x^2}{2\beta_{Q,\vect{z}}^2}.
\end{IEEEeqnarray}
Hence, from~\eqref{EqCrazyStuff23} and~\eqref{EqPhiStarfullY}, given~$c$ in~\eqref{EqWhatIsThis} for all~$t \in \left(-\infty, 0 \right)$,
\begin{IEEEeqnarray}{rcl}
c \, t -  \frac{1}{2}t^2  \beta_{Q,\vect{z}}^2  \leqslant \phi^{\star}(c) \leqslant D\left( P \| P^{\left(Q, \lambda\right)}_{\vect{\Theta} | \vect{Z} = \vect{z}} \right),
\end{IEEEeqnarray}
and thus,  
\begin{equation}
\label{EqThisCrazyIneq865}
\frac{c^2}{2\beta_{Q,\vect{z}}^2}   \leqslant D\left( P \| P^{\left(Q, \lambda\right)}_{\vect{\Theta} | \vect{Z} = \vect{z}} \right).
\end{equation}
This implies that
\begin{equation}
c \leqslant \sqrt{2 \beta_{Q,\vect{z}}^2 D\left( P\| P^{\left(Q, \lambda\right)}_{\vect{\Theta} | \vect{Z} = \vect{z}} \right) }
\end{equation}
and
\begin{equation}
c \geqslant - \sqrt{2 \beta_{Q,\vect{z}}^2 D\left( P\| P^{\left(Q, \lambda\right)}_{\vect{\Theta} | \vect{Z} = \vect{z}} \right) },
\end{equation}
which leads to 
\begin{IEEEeqnarray}{rcl}
\nonumber
& & \abs{ \int \mathsf{L}_{\vect{z}}(\vect{\theta}) \mathrm{d} P (\vect{\theta}) - \int \mathsf{L}_{\vect{z}}(\vect{\theta}) \mathrm{d} P^{\left(Q, \lambda\right)}_{\vect{\Theta} | \vect{Z} = \vect{z}}(\vect{\theta}) } \\ 
& \leqslant &    \sqrt{2 \beta_{Q,\vect{z}}^2 D\left( P\| P^{\left(Q, \lambda\right)}_{\vect{\Theta} | \vect{Z} = \vect{z}} \right) },
\end{IEEEeqnarray}
and completes the proof.

\section{Proof of Theorem~\ref{TheoAminini}}
\label{AppProofOfTheoAminini}

Under the condition that $\lambda \in \set{K}_{Q,P_{\vect{Z}}}$, from Theorem~\ref{TheoremSensitivityEqual} and Definition~\ref{DefGenErr}, it follows that the generalization error $\mathsf{G}_{Q, \lambda}\left(P_{\vect{Z}} \right)$ in~\eqref{EqInsop} satisfies
\begin{IEEEeqnarray}{rCl}
\nonumber
\mathsf{G}_{Q, \lambda}\left(P_{\vect{Z}} \right) & =  &  \lambda \int \Big( \KL{P^{\left(Q, \lambda\right)}_{\vect{\Theta}| \vect{Z} = \vect{\nu}}}{Q} + \KL{P^{\left(Q, \lambda\right)}_{\vect{\Theta}}}{P^{\left(Q, \lambda\right)}_{\vect{\Theta}| \vect{Z} = \vect{\nu}}} \\ 
& &- \KL{P^{\left(Q, \lambda\right)}_{\vect{\Theta}}}{Q} \Big) \mathrm{d} P_{\vect{Z}}(\vect{\nu}) ,\\
\nonumber
& =  & \lambda\Bigg( \int  \KL{P^{\left(Q, \lambda\right)}_{\vect{\Theta}| \vect{Z} = \vect{\nu}} }{P^{\left(Q, \lambda\right)}_{\vect{\Theta}}}  \mathrm{d} P_{\vect{Z}}(\vect{\nu})  \Dsupersqueezeequ\\
\label{EqWorkingInTheGardenOfFlowers}
&   & + \int \KL{P^{\left(Q, \lambda\right)}_{\vect{\Theta}}}{P^{\left(Q, \lambda\right)}_{\vect{\Theta}| \vect{Z} = \vect{\nu}} }  \mathrm{d} P_{\vect{Z}}(\vect{\nu})  \Bigg)  ,\Dsupersqueezeequ \IEEEeqnarraynumspace
\end{IEEEeqnarray}
where the equality in~\eqref{EqWorkingInTheGardenOfFlowers} follows from the fact that
\begin{IEEEeqnarray}{rCl}
\nonumber
& &\int \Big( \KL{P^{\left(Q, \lambda\right)}_{\vect{\Theta}| \vect{Z} = \vect{\nu}}}{Q}  - \KL{P^{\left(Q, \lambda\right)}_{\vect{\Theta}}}{Q} \Big) \mathrm{d} P_{\vect{Z}}(\vect{\nu}) \Dsupersqueezeequ \IEEEeqnarraynumspace\\
& = &\int  \KL{P^{\left(Q, \lambda\right)}_{\vect{\Theta}| \vect{Z} = \vect{\nu}}}{Q} \mathrm{d} P_{\vect{Z}}(\vect{\nu}) -\KL{P^{\left(Q, \lambda\right)}_{\vect{\Theta}}}{Q} 
\Dsupersqueezeequ \IEEEeqnarraynumspace \\
\nonumber
& = &\displaystyle\int \left( \int \log\left(\frac{\mathrm{d} P^{\left(Q, \lambda\right)}_{\vect{\Theta}| \vect{Z} = \vect{\nu}}}{\mathrm{d} Q} (\vect{\theta}) \right) \mathrm{d} P^{\left(Q, \lambda\right)}_{\vect{\Theta}| \vect{Z} = \vect{\nu}}(\vect{\theta}) \right)P_{\vect{Z}}(\vect{\nu}) \Dsupersqueezeequ \\
&  & -\KL{P^{\left(Q, \lambda\right)}_{\vect{\Theta}}}{Q} 
\Dsupersqueezeequ \IEEEeqnarraynumspace \\
\nonumber
& = &\displaystyle\int \left( \int \log\left(\frac{\mathrm{d} P^{\left(Q, \lambda\right)}_{\vect{\Theta}| \vect{Z} = \vect{\nu}}}{\mathrm{d} Q} (\vect{\theta}) \right) \mathrm{d} P^{\left(Q, \lambda\right)}_{\vect{\Theta}| \vect{Z} = \vect{\nu}}(\vect{\theta}) \right) \mathrm{d} P_{\vect{Z}}(\vect{\nu})\Dsupersqueezeequ \IEEEeqnarraynumspace \\
&  & -\int\log\left( \frac{\mathrm{d}P^{\left(Q, \lambda\right)}_{\vect{\Theta}}}{\mathrm{d}Q} (\vect{\theta}) \right)\mathrm{d}P^{\left(Q, \lambda\right)}_{\vect{\Theta}}(\vect{\theta}) 
\Dsupersqueezeequ \IEEEeqnarraynumspace \\
\nonumber
& = &\displaystyle\int \left( \int \log\left(\frac{\mathrm{d} P^{\left(Q, \lambda\right)}_{\vect{\Theta}| \vect{Z} = \vect{\nu}}}{\mathrm{d} Q} (\vect{\theta}) \right) \mathrm{d} P^{\left(Q, \lambda\right)}_{\vect{\Theta}| \vect{Z} = \vect{\nu}}(\vect{\theta}) \right) \mathrm{d}P_{\vect{Z}}(\vect{\nu}) 
\Dsupersqueezeequ \\
\label{EqQuieroDarteUnBeso}
&  & - \int \left( \int\log\left( \frac{\mathrm{d}P^{\left(Q, \lambda\right)}_{\vect{\Theta}}}{\mathrm{d}Q} (\vect{\theta}) \right)\mathrm{d} P^{\left(Q, \lambda\right)}_{\vect{\Theta}| \vect{Z} = \vect{\nu}}(\vect{\theta})  \right) \mathrm{d}P_{\vect{Z}}(\vect{\nu}) 
\Dsupersqueezeequ \IEEEeqnarraynumspace \\
\nonumber
& = &\displaystyle\int \Bigg( \int \Bigg(  \log\left(\frac{\mathrm{d} P^{\left(Q, \lambda\right)}_{\vect{\Theta}| \vect{Z} = \vect{\nu}}}{\mathrm{d} Q} (\vect{\theta}) \right) \Dsupersqueezeequ \\
\label{EqDarteLaFuerza}
& & + \log\left( \frac{\mathrm{d}Q}{\mathrm{d}P^{\left(Q, \lambda\right)}_{\vect{\Theta}}} (\vect{\theta}) \right) \Bigg)
\mathrm{d} P^{\left(Q, \lambda\right)}_{\vect{\Theta}| \vect{Z} = \vect{\nu}}(\vect{\theta}) \Bigg) \mathrm{d}P_{\vect{Z}}(\vect{\nu}) \Dsupersqueezeequ \IEEEeqnarraynumspace \\
& = &\displaystyle\int \Bigg( \int  \log\left(\frac{\mathrm{d} P^{\left(Q, \lambda\right)}_{\vect{\Theta}| \vect{Z} = \vect{\nu}}}{\mathrm{d} P^{\left(Q, \lambda\right)}_{\vect{\Theta}}} (\vect{\theta}) \right) 
\mathrm{d} P^{\left(Q, \lambda\right)}_{\vect{\Theta}| \vect{Z} = \vect{\nu}}(\vect{\theta}) \Bigg) \mathrm{d}P_{\vect{Z}}(\vect{\nu}) \Dsupersqueezeequ \IEEEeqnarraynumspace\\
& = &\displaystyle\int \KL{P^{\left(Q, \lambda\right)}_{\vect{\Theta}| \vect{Z} = \vect{\nu}}}{P^{\left(Q, \lambda\right)}_{\vect{\Theta}}}  \mathrm{d}P_{\vect{Z}}(\vect{\nu}). \Dsupersqueezeequ \IEEEeqnarraynumspace
\end{IEEEeqnarray}
The equality in~\eqref{EqQuieroDarteUnBeso} follows from~\eqref{EqBarPThetaX}; and
the equality in~\eqref{EqDarteLaFuerza} follows from the fact that the measures $Q$ and $P^{\left(Q, \lambda\right)}_{\vect{\Theta}| \vect{Z} = \vect{\nu}}$, with $\vect{\nu}\in\supp P_{\vect{Z}}$, are mutually absolutely continuous (Lemma~\ref{LemmaMutualAC}). 
This completes the proof.
  \end{appendices}

\bibliographystyle{IEEEtran}
\bibliography{references}

\end{document}